%% file: main.tex
\documentclass[b5paper,twoside]{book}
\usepackage[paper=b5paper]{geometry}

\usepackage[utf8]{inputenc}
\usepackage{amssymb}
\usepackage{mathtools}
\usepackage[hidelinks]{hyperref}
\usepackage{tikz-cd}
\usepackage{tikz}
\usepackage[shortlabels]{enumitem}
\usepackage{dsfont}
\usepackage{amsthm}
\usepackage{pdfpages}

\usepackage{personalization/Style}
\usepackage{personalization/FavoritePackages}
\usepackage{personalization/Commands}
\usepackage{subfiles}

\usepackage{imakeidx}

\makeindex[intoc, title=Index of Notation, columnsep=80px, options= -s indexStyle.ist]

\title{Varieties admitting a holomorphic symplectic form: LLV algebras and derived equivalences}
\author{Dion Leijnse}
\date{ }

\begin{document}
\newgeometry{left=2.8cm,bottom=0cm}
\maketitle

\newpage
\thispagestyle{empty}
\restoregeometry
\vspace*{\fill}

\noindent 
This research has been carried out at the Korteweg--de Vries Institute for Mathematics and during a stay at the Hausdorff institute for Mathematics in Bonn and was funded by the European Research Council (ERC), grant 864145, and by the Deutsche Forschungsgemeinschaft (DFG, German Research Foundation) under Germany's Excellence Strategy – EXC-2047/1–390685813.

\restoregeometry

\frontmatter
\tableofcontents

\chapter{Acknowledgements}

\subfile{chapters/dankwoord}

\mainmatter

\chapter{Introduction}
\label{chap: introduction}

\subfile{chapters/introduction}

\chapter{Preliminaries}
\subfile{chapters/chapter2}

\chapter{The LLV Lie algebra of a variety admitting a holomorphic symplectic form}
\label{chap: LLV computation}

\subfile{chapters/chapter4}

\chapter{Computation of the Verbitsky component}
\label{chap: SH computation}

\subfile{chapters/chapter5}

\chapter{Construction of $\psi$}
\label{chap: fatcorwise psi}

\subfile{chapters/chapter6}

\chapter{Construction of $\Psi^{\sH}$}
\label{chap: construction of Psi}

\subfile{chapters/chapter7}

\chapter{Generalizations of the main theorem}
\label{sec: low degree cases}

\subfile{chapters/chapter8}

\chapter{Enriques manifolds}
\label{chap: enriques manifolds}

\subfile{chapters/chapter9}

\backmatter

\chapter{Summary}\label{chap:summary}

\subfile{chapters/summary}

\chapter{Samenvatting}\label{chap:samenvatting}

\subfile{chapters/samenvatting}

\bibliography{PhDReferences.bib}

\printindex

\end{document}

%% file: chapters/dankwoord.tex
During my PhD and my 4 years in Amsterdam, many people have supported me in one way or another, and I would like to use this space to thank those who have supported me throughout this journey and made my time here so enjoyable. I will undoubtedly forget to mention some people here; I am grateful to everyone who supported me, even if you aren't mentioned here.

Als eerste wil ik Lenny bedanken. Dankje dat je me de kans gaf om een promotie bij jou te doen, het was een fantastische ervaring. Er zijn weinigen die zo duidelijk en helder als jij wiskunde begrijpen en het zo overzichtelijk uit kunnen leggen. Het was dan ook een voorrecht om altijd bij je langs te kunnen lopen met mijn wiskundige vragen. Je manier van begeleiden en vermogen om problemen tot de kern te kunnen reduceren waren erg inspirerend en hebben mij veel geholpen.

I am also very grateful to Mingmin for being my copromotor, and to Eric, Hessel, Yagna, Lie, Martijn and Sergey for being in my thesis committee. Thank you Lie for the inspiring discussions we had when I visited you in Strasbourg. Thank you Lie and Martijn for your comments on the text.

It was always a joy to be at work, and this is mainly due to all my great colleagues at the KdvI: Bjarne, Felix, Edward, Wouter, Kees, Jeroen, Emma, Jurre, Georg, Jort, Shin, Vincent, Chemy, Noah, Dhruva, Dhyan, Kaan, Reinier, Yachen, Zekun, Leo, Gabriëlle, Siebe, David, Pjotr, Annika, Thijs, Robin, and I have probably forgotten many people here. I've had lots of fun during all the cookie times, board game nights and bouldering that we've done. Thanks for creating the great atmosphere at the institute. I am also grateful to Bjarne for proofreading the introduction of this thesis and giving valuable feedback, and to Edward and Dhruva for being my paranymphs.

During my PhD, I had the oppurtunity to stay for a few months in London and later also a few months in Bonn. I have had an amazing time in these two cities, and I am very grateful to those who have made this possible.

Daarnaast wil ik ook Joris, Yannick, Arne, Michiel, Etienne, Willem en Thomas bedanken. Ons vele contact, waaronder met variërende deelverzamelingen van deze groep wekelijkse online spellenavonden, spellendagen, wandelweekenden, Dungeons and Dragons sessies en vakanties zorgden altijd voor de nodige afleiding van de PhD. Ook wil ik jullie heel erg bedanken dat jullie altijd voor me klaar stonden en voor de interessante discussies over wiskunde en wetenschap.

Uiteraard kunnen papa en mama hier niet ontbreken. Ik was hier zonder jullie steun nooit gekomen. Ik ben ontzettend dankbaar voor alle kansen en mogelijkheden die jullie me hebben gegeven, en ook voor de liefdevolle opvoeding. Ook wil ik Kenzie, mijn zusje, bedanken voor jouw steun en dat je altijd voor me klaarstaat.

Almost always when doing mathematics or anything else requiring concentration, I am listening to music. By a wide margin, my two most-listened to artists are Nightwish and Two Steps from Hell, and I can listen endlessly to the beautiful songs they have created. I am grateful to these artists for the music they have written.

Als laatste en zeker niet de minste wil ik mijn vriendin Reja bedanken. Je wist me altijd te ondersteunen wanneer dat nodig was, als ik er even doorheen zat of te veel stress had. Je bent fantastisch en weet me altijd op te vrolijken. Ik ben je intens dankbaar, en hoop dat het me tijdens jouw PhD gaat lukken om je net zo goed te steunen als dat je mij bijgestaan hebt.

%% file: chapters/introduction.tex
The derived category of an abelian category was introduced by Grothendieck and Verdier in the 1960s \cite{VerdierDerivedCategories}. It was originally introduced as a natural setting for doing homological algebra, and is very useful for studying derived functors. This has featured prominently in algebraic geometry, an area where homological algebra is used extensively. Here, one often looks at \myindex{$D^b(X)$}, the derived category of coherent sheaves on a variety $X$.

Later, the derived category was studied as an invariant of the variety $X$, for example by Bondal, Mukai and Orlov \cite{MukaiAbelian, OrlovK3, OrlovAbelian, BondalOrlov}. It turns out that many geometric properties of $X$ are captured by its derived category. For example, if two smooth projective varieties $X$ and $Y$ over $\mathbb{C}$ have equivalent derived categories (for us, this will always mean a $\mathbb{C}$-linear additive exact equivalence of triangulated categories $D^b(X) \simeq D^b(Y)$), then they have the same dimension and their canonical bundles have the same order \cite[Proposition~4.1]{Huybrechts}. In fact, Bondal and Orlov proved that if $X$ has either ample or anti-ample canonical bundle and $D^b(X) \simeq D^b(Y)$, then $X$ and $Y$ are isomorphic. Similarly, if $X$ and $Y$ are derived equivalent smooth projective curves over $\mathbb{C}$, then $X$ and $Y$ are isomorphic \cite[Corollary~5.46]{Huybrechts}. 

Despite the above examples, the derived category of a smooth projective variety $X$ does not recover $X$ in general. The earliest example of two non-isomorphic varieties with equivalent derived categories comes from Mukai \cite{MukaiAbelian}, who showed that for every abelian variety $X$, there is an equivalence $D^b(X) \simeq D^b(X^*)$, where $X^*$ denotes the dual abelian variety of~$X$.

Another interesting example in this direction comes from K3 surfaces. For a K3 surface $X$, after identifying $\sH^4(X; \mathbb{Z})$ with $\mathbb{Z}$, the cohomology $\sH^\bullet(X; \mathbb{Z})$ becomes a lattice, where one defines the pairing $\langle -, -\rangle$ by letting $\langle x, y \rangle$ be the degree $4$ part of the cup product $x \cup y$ for $x,y \in \sH^\bullet(X; \mathbb{Z})$. Mukai \cite{MukaiK3} equipped this lattice with the weight two Hodge structure which extends the weight two Hodge structure of $\sH^2(X; \mathbb{Z})$, and where $\sH^0(X; \mathbb{Z})$ and $\sH^4(X; \mathbb{Z})$ are algebraic (i.e. of type $(1,1)$). The resulting structure is called the Mukai lattice of $X$. Orlov \cite{OrlovK3} then showed that two K3 surfaces $X$ and $Y$ are derived equivalent if and only if there is a Hodge isometry between their Mukai lattices, a result which is now known as the derived Torelli theorem.

Despite the fact that two derived equivalent varieties need not be isomorphic, they do still share a lot of structure. An important result in this direction comes from Orlov \cite{OrlovAbelian}, who showed that every equivalence $\Phi \colon D^b(X) \isomto D^b(Y)$ for smooth projective varieties $X$ and $Y$ is given by a so-called Fourier-Mukai transform. If the varieties $X$ and $Y$ are defined over the complex numbers, one obtains from this an induced isomorphism on cohomology\index{$\Phi^{\sH}$}
\[\Phi^{\sH} \colon \sH^\bullet(X; \mathbb{Q}) \to \sH^\bullet(Y; \mathbb{Q}),\]
where $\sH^\bullet(X; \mathbb{Q})$ denotes the singular cohomology of the analytification $X(\mathbb{C})^{\an}$. The isomorphism $\Phi^{\sH}$ does not need to respect the grading. However, it does preserve the Hodge grading, in the sense that for every $n$ we have (after extending scalars to $\mathbb{C}$ and using the Hodge decomposition):
\begin{equation}
    \label{eq: Phi^H respects HS}
    \Phi^{\sH}\left(\bigoplus_{q-p = n} \sH^q(X, \Omega^p_X)\right) = \bigoplus_{q-p = n} \sH^q(Y, \Omega_Y^p).
\end{equation}

Roughly, $\Phi^{\sH}$ respects the columns of the Hodge diamond, but not the rows. It is then natural to ask whether there could also exist an isomorphism $\sH^\bullet(X; \mathbb{Q}) \to \sH^\bullet(Y; \mathbb{Q})$ which preserves the grading, and is an isomorphism of Hodge structures in each degree. The following was conjectured by Orlov:
\begin{conjecture}[Orlov \cite{orlov2005derived}]
    \label{conj: orlov cohomology}
    Let $X$ and $Y$ be two smooth projective varieties over $\mathbb{C}$, and assume that there exists an equivalence of triangulated categories $D^b(X) \simeq D^b(Y)$. Then for every $i$ there is an isomorphism of Hodge structures $\sH^i(X; \mathbb{Q}) \cong \sH^i(Y; \mathbb{Q})$.
\end{conjecture}
In fact, Orlov conjectured something even stronger, namely that derived equivalent varieties have isomorphic motives:
\begin{conjecture}[Orlov \cite{orlov2005derived}]
    Let $X$ and $Y$ be two smooth projective varieties over a field $k$, and assume that there exists an equivalence of triangulated categories $D^b(X) \simeq D^b(Y)$. Then the rational Chow motives of $X$ and $Y$ are isomorphic.
\end{conjecture}

A weaker form of Conjecture \ref{conj: orlov cohomology} just asks for equality of the Hodge numbers $h^{p,q}(X) = h^{p,q}(Y)$, where $h^{p,q}(X) = \dim \sH^q(X, \Omega_X^p)$. The invariance of Hodge numbers under derived equivalences would also follow from conjectures made by Kontsevich in the context of mirror symmetry \cite{KontesvichMS}. Equation \eqref{eq: Phi^H respects HS} implies that for every $n$ the sum $\sum_{q-p = n} h^{p,q}$ is a derived invariant, which is some evidence towards this conjecture. It should be noted that the invariance of Hodge numbers is not true in positive characteristic \cite{DerivedHodgePositiveCharacteristic}.

These conjectures are already known in various cases, which we summarize here:
\begin{enumerate}
    \item If $X$ and $Y$ are smooth projective curves over $\mathbb{C}$, or $X$ has (anti)-ample canonical bundle, and $D^b(X) \simeq D^b(Y)$, then we have seen above that $X \cong Y$, so the conjectures are trivially true.
    \item If $X$ and $Y$ are abelian varieties with duals $X^*$ and $Y^*$ respectively, and $D^b(X) \simeq D^b(Y)$, then there is an isomorphism $X \times X^* \cong Y \times Y^*$ by \cite{OrlovAbelian}. Since $X^*$ and $X$ are isogenous, there is an isomorphism of weight $1$ Hodge structures $\sH^1(X; \mathbb{Q}) \cong \sH^1(X^*; \mathbb{Q})$. The projectivity of $X$ and $Y$ implies that the Hodge structures on their cohomology are polarizable, and hence split as a direct sum of their simple sub-Hodge structures. It follows that there is an isomorphism of weight $1$ Hodge structures $\sH^1(X; \mathbb{Q}) \cong \sH^1(Y; \mathbb{Q})$, and taking the exterior power of this isomorphism gives an isomorphism $\sH^\bullet(X; \mathbb{Q}) \cong \sH^\bullet(Y; \mathbb{Q})$ which respects the grading and is an isomorphism of Hodge structures in each degree.
    \item If $X$ and $Y$ are smooth projective K3 surfaces and $\Phi \colon D^b(X) \isomto D^b(Y)$ is an equivalence, then we saw above that there is a not necessarily graded isomorphism $\Phi^{\sH} \colon \sH^\bullet(X; \mathbb{Q}) \cong \sH^\bullet(Y; \mathbb{Q})$ as above. By applying a Tate twist to the cohomology of $X$ in degrees $0$ and $4$, we can give the cohomology a weight $2$ Hodge structure, and this is then isomorphic to $\mathbb{Q}(-1)^{\oplus 2} \oplus \sH^2(X; \mathbb{Q})$. Hence $\Phi^{\sH}$ gives an isomorphism of weight $2$ Hodge structures
    \[ \mathbb{Q}(-1)^{\oplus 2} \oplus \sH^2(X; \mathbb{Q}) \cong \mathbb{Q}(-1)^{\oplus 2} \oplus \sH^2(Y; \mathbb{Q}).\]
    Since the category of polarizable weight $2$ Hodge structures is semisimple, we can cancel the factors $\mathbb{Q}(-1)^{\oplus 2}$, and this proves Conjecture~\ref{conj: orlov cohomology} for K3 surfaces.

    Note that this example not only works for K3 surfaces, but in fact for all surfaces with $h^1 = 0$, and can even be generalized to all smooth projective surfaces.
    \item The above example was generalized by Taelman \cite{TaelmanDerivedEquivalences} to hyperkähler varieties, who showed that Conjecture \ref{conj: orlov cohomology} holds if $X$ and $Y$ are derived equivalent hyperkähler varieties.
    \item If $X$ and $Y$ are derived equivalent smooth projective varieties over $\mathbb{C}$, then $h^{0,1}(X) = h^{0,1}(Y)$ by \cite[Corollary B]{PopaSchnell}.
    \item The above result also implies that if $X$ and $Y$ are smooth projective varieties over $\mathbb{C}$ of dimension at most $3$ and $D^b(X) \simeq D^b(Y)$, then $X$ and $Y$ have the same Hodge numbers, see \cite[Corollary C]{PopaSchnell}.
\end{enumerate}

Since Conjecture \ref{conj: orlov cohomology} is known for varieties with ample or anti-ample canonical bundle, a logical next step is to try to prove the conjecture for more general varieties with trivial canonical bundle other than abelian varieties or hyperkähler varieties. By the Beauville-Bogomolov decomposition \cite{BeauvilleC1Zero}, any smooth projective complex variety $X$ with trivial canonical bundle admits an étale covering
\begin{equation}
    \label{eq: BB covering}
    A \times \prod X_i \times \prod Y_j \to X,
\end{equation}
where $A$ is an abelian variety, the $X_i$ are hyperkähler varieties and the $Y_j$ are strict Calabi-Yau varieties. One could hope that a derived equivalence between two varieties with trivial canonical bundle can be understood in terms of these simpler factors. Huybrechts and Nieper-Wisskirchen \cite{HuybrechtsNieperWisskirchen} obtained some results in this direction. As a first question, one may wonder whether the derived category can distinguish between the three different kinds of factors that occur. In this direction, they prove the following:

\begin{theorem}[{\cite[Theorem 0.4]{HuybrechtsNieperWisskirchen}}]
    Let $X$ and $Y$ be two derived equivalent smooth projective varieties. If $X$ is either an abelian or a hyperkähler variety, then $Y$ is of the same type.
\end{theorem}

Moreover, when considering more than one factor in the Beauville-Bogomolov decomposition, they obtain the following result:
\begin{theorem}[{\cite[Proposition 4.11]{HuybrechtsNieperWisskirchen}}]
    \label{thm: HNW derived invariance products}
    Let $X_1, \dots, X_m$ and $X'_1, \dots, X'_n$ be hyperkähler varieties, and assume that $D^b(\prod X_i) \simeq D^b(\prod X'_j)$. Then there is a bijection $\sigma \colon \{1, \dots, m\} \to \{1, \dots, n\}$ such that $\dim X_i = \dim X'_{\sigma(i)}$ and $b_2(X_i) = b_2(X'_{\sigma(i)})$ for all $i$.
\end{theorem}

Given these two theorems, there are several interesting paths to explore. On the one hand, since Conjecture \ref{conj: orlov cohomology} is already known for the abelian and hyperkähler factors, one could try to prove it for the strict Calabi-Yau factors. Another possible avenue of research is to combine the results obtained so far, and study what happens if the Beauville-Bogomolov decomposition of $X$ is nontrivial, for example if there is a nontrivial étale covering, or if there are several factors in the Beauville-Bogomolov decomposition of $X$. Note that the Beauville-Bogomolov decomposition is not unique, but that there is not too much variation in the possible étale covers of a given variety with trivial canonical bundle. See Section~\ref{sec: reduction to quotients} for a summary of some results of Beauville in this direction.

This thesis makes contributions in the second direction. Our main result is that Conjecture \ref{conj: orlov cohomology} holds for many varieties which admit a holomorphic symplectic form. If $X$ is a smooth projective variety which admits a holomorphic symplectic form, its Beauville-Bogomolov covering \eqref{eq: BB covering} does not contain any strict Calabi-Yau factors, so we only have to deal with the abelian and hyperkähler factors. The reason that we restrict to varieties with a holomorphic symplectic form is that we can use a theorem from Taelman \cite{TaelmanDerivedEquivalences}, which also underlies his proof of Conjecture~\ref{conj: orlov cohomology} for hyperkähler varieties. This result is about the derived invariance of the so-called LLV Lie algebra that can be associated to a smooth projective variety. The LLV algebra of a smooth projective variety $X$, denoted by $\llv(X; \mathbb{Q})$, is a semisimple Lie algebra over $\mathbb{Q}$, which comes with a representation on $\sH^\bullet(X; \mathbb{Q})$. In Chapter~\ref{chap: LLV computation}, we discuss the LLV algebra in more detail.

The Lie algebra $\llv(X; \mathbb{Q})$ can be computed in many cases. For example, the definition of the Mukai lattice of a K3 surface can be generalized to obtain a Mukai lattice $\widetilde \sH(X; \mathbb{Q})$ of a hyperkähler variety $X$, and there is an isomorphism $\llv(X; \mathbb{Q}) \cong \mathfrak{so}(\widetilde \sH(X; \mathbb{Q}))$, see \cite{LooijengaLunts}, \cite{VerbitskyPhDThesis} and \cite{SoldatenkovHyperkahler}. Taelman \cite{TaelmanDerivedEquivalences} proved that for symplectic varieties, the LLV algebra, together with its representation on cohomology, is a derived invariant:

\begin{theorem}[{\cite[Theorem A]{TaelmanDerivedEquivalences}}]
    \label{thm: Taelman theorem A}
    Let $X$ and $Y$ be two smooth projective varieties over $\mathbb{C}$, both admitting a holomorphic symplectic form, and let $\Phi \colon D^b(X) \simeq D^b(Y)$ be an equivalence of triangulated categories. Then $\Phi$ induces an isomorphism of LLV algebras \index{$\Phi^{\llv}$}$\Phi^{\llv} \colon \llv(X; \mathbb{Q}) \to \llv(Y; \mathbb{Q})$, with the property that $\Phi^{\sH} \colon \sH^\bullet(X; \mathbb{Q}) \to \sH^\bullet(Y; \mathbb{Q})$ is equivariant with respect to the actions of the LLV algebras on the cohomology.
\end{theorem}

The derived invariance of the LLV algebra can be used to prove the derived invariance of other invariants of $X$. For example, we show in Section \ref{subsec: intermezzo, Huybrechts Nieper Wisskirchen} how this result, together with some computations of the LLV algebra and its representation on cohomology, can be used to give an alternative proof of Theorem~\ref{thm: HNW derived invariance products}. 

\subsubsection{Main results}
We now discuss the main results of this thesis, all of which build on Theorem~\ref{thm: Taelman theorem A}. First, we will prove Orlov's conjecture for varieties admitting a holomorphic symplectic form, under some assumptions on their LLV algebra. Take an integer $n \geq 1$. We say that a semisimplie Lie algebra $\mathfrak{g}$ over $\mathbb{Q}$ has a simple factor of type $A_n$ if the base change $\mathfrak{g} \otimes_\mathbb{Q} \mathbb{C}$ has a simple factor of type $A_n$, and similarly for the other simple types $B_n, C_n$, etc. We will require our varieties $X$ to satisfy:

\begin{enumerate}[label=(\textasteriskcentered)]
    \item\label{condition *} The LLV algebra $\llv(X; \mathbb{Q})$ does not have any simple factors of the type $C_2, A_3$ or $D_4$.
\end{enumerate}

\begin{theorem}
    \label{thm: main theorem}
    Let $X$ and $Y$ be smooth projective varieties over $\mathbb{C}$ satisfying condition \ref{condition *}. Assume that $X$ and $Y$ both admit a holomorphic symplectic form, and assume that there exists an equivalence of triangulated categories $D^b(X) \simeq D^b(Y)$. Then there exists an isomorphism $\sH^\bullet(X; \mathbb{Q}) \cong \sH^\bullet(Y; \mathbb{Q})$ respecting the Hodge bigrading (i.e. there exists a degree-preserving isomorphism which, in each degree, is an isomorphism of Hodge structures).
\end{theorem}

The proof of this theorem takes up most of the thesis; we finish the proof in Chapter~\ref{chap: construction of Psi}.

\begin{remark}
    Condition \ref{condition *} suffices to prove Theorem \ref{thm: main theorem}, but in some cases we can slightly weaken this condition. See Chapter \ref{sec: low degree cases} for more details, where we state Theorem~\ref{thm: main theorem} in its most general form in Theorem~\ref{thm: main theorem generalized}.
\end{remark}

Before we continue, let us discuss condition \ref{condition *} a bit. Unfortunately, it is not vacuous. To construct an example of a variety violating this condition, let $A$ be an abelian surface (in particular, $A$ admits a symplectic form) and let $S$ be a K3 surface equipped with a symplectic automorphism $\sigma$ of order $5$ (such a K3 surface exists, see for instance \cite{MukaiSymplecticK3Automorphisms}). Choose an element $a \in A$ of order $5$. Let $G = \langle g \rangle$ be the cyclic group of order $5$ and let $G$ act on $S \times A$ by letting $g (s, x) = (\sigma(s), x + a)$. Then $G$ acts freely on $S \times A$ and preserves a symplectic form, so $X := (S \times A)/G$ admits a holomorphic symplectic form. By \cite[Proposition 1.1]{K3PrimeOrderAutomorphisms}, we have $\dim \sH^2(S; \mathbb{Q})^\sigma = 6$, and we will see (Theorem \ref{thm: summary LLV computation}) that the LLV algebra of $X = (S \times A)/G$ is a product of two factors of type $D_4$. 

However, we are still able to prove Conjecture \ref{conj: orlov cohomology} for varieties of the above type. In fact, by studying the geometry of four-dimensional varieties with holomorphic symplectic forms, we will see that we do not need condition \ref{condition *} in dimension $4$:

\begin{theorem}[{Theorem \ref{thm: main theorem dimension 4}}]
    Let $X$ and $Y$ be two smooth projective varieties of dimension $4$ over $\mathbb{C}$, both admitting a holomorphic symplectic form, and assume that there is a triangulated equivalence $D^b(X) \simeq D^b(Y)$. Then there is an isomorphism $\sH^\bullet(X; \mathbb{Q}) \cong \sH^\bullet(Y; \mathbb{Q})$ which preserves the Hodge bigrading.
\end{theorem}

The above example can be generalized to higher dimensions, for instance by taking the product of $(S \times A)/G$ with another symplectic variety. Therefore there are still varieties admitting a holomorphic symplectic form to which Theorem~\ref{thm: main theorem} and its generalizations do not apply.

We will see in Lemma \ref{lem: * and ** equivalent} that condition \ref{condition *} can be restated in a more explicit manner. Using the Beauville-Bogomolov decomposition, we will see in Section~\ref{sec: reduction to quotients} that the variety $X$ can be written in the form $X \cong (X_0 \times \prod_{i=1}^k X_i)/G$, where $X_0$ is an abelian variety, the $X_i$ are hyperkähler varieties, and $G$ is a finite group acting freely on $X_0 \times \prod_{i=1}^k X_i$. Since $\sH^1(X_0 \times \prod_{i=1}^k X_i; \mathbb{Q}) = \sH^1(X_0; \mathbb{Q})$, we obtain a representation of $G$ on $\sH^1(X_0; \mathbb{Q})$. Furthermore, for every $1 \leq i \leq k$, there is a maximal subgroup $G_i \subseteq G$ such that the action of $G$ on $X_0 \times \prod X_i$ restricts to an action of $G_i$ on $X_i$ (see Section~\ref{sec: reduction to quotients} for more details). Write $\Irr_\mathbb{C}(G)$ for a set of representatives for the isomorphism classes of irreducible $\mathbb{C}$-valued $G$-representations. Then \ref{condition *} is equivalent to:

\begin{enumerate}[label=(\textasteriskcentered \textasteriskcentered)]
    \item\label{condition **} The action of $G$ on $X$ has the following properties:
    \begin{enumerate}
        \item Every self-dual $\rho \in \Irr_\mathbb{C}(G)$ with $(\extp^2 \rho)^G = 0$ that occurs in the representation $G \acts \sH^1(X_0; \mathbb{C})$ does not have multiplicity $3$ or $4$.
        \item Every self-dual $\rho \in \Irr_\mathbb{C}(G)$ with $(\extp^2 \rho)^G = \mathbb{C}$ that occurs in the representation $G \acts \sH^1(X_0; \mathbb{C})$ does not have multiplicity $2$.
        \item Every non-self-dual $\rho \in \Irr_\mathbb{C}(G)$ that occurs in the representation of $G$ on $\sH^1(X_0; \mathbb{C})$ has multiplicity not equal to $2$.
        \item For every $i$ with $1 \leq i \leq k$, the dimension of $\sH^2(X_i;\mathbb{Q})^{G_i}$ is $5$ or at least $7$.
    \end{enumerate}
\end{enumerate}
The dimension of $\sH^2(X_i;\mathbb{Q})^{G_i}$ is at least $3$ (see Lemma \ref{lem: HK H^2 invariants dimension}), so property (d) means that this dimension is not allowed to be $3, 4$ or $6$.

Lastly, in Chapter \ref{chap: enriques manifolds} we will prove Orlov's conjecture for varieties with a torsion canonical bundle whose simply connected covering is a hyperkähler variety. Such varieties are sometimes called \textbf{Enriques manifolds}, as they generalize Enriques surfaces to higher dimensions.

\begin{theorem}[Theorem \ref{thm: orlov conj for enriques}]
    Let $X$ and $Y$ be two Enriques manifolds. Let $\Phi \colon D^b(X) \simeq D^b(Y)$ be an equivalence of triangulated categories. Then for every $i$ there is an isomorphism $\sH^i(X; \mathbb{Q}) \cong \sH^i(Y; \mathbb{Q})$ that preserves the Hodge structure.
\end{theorem}
The proof of this theorem is independent of the theory developed in the earlier chapters of this thesis, and instead builds directly on Theorem \ref{thm: Taelman theorem A}.

Having discussed this, we give a rough outline of the proof of Theorem \ref{thm: main theorem}. First, we will compute the LLV algebras of $X$ and $Y$. Using the Beauville-Bogomolov decomposition, we can write the LLV algebra of $X$ as a product of smaller Lie algebras. We will then explicitly compute these smaller factors by using algebras with involutions. This is done in Chapter~\ref{chap: LLV computation}, with the results summarized in Theorem \ref{thm: summary LLV computation}.

By Theorem \ref{thm: Taelman theorem A}, the equivalence $\Phi \colon D^b(X) \to D^b(Y)$ induces an isomorphism of LLV algebras $\Phi^{\llv} \colon \llv(X; \mathbb{Q}) \to \llv(Y; \mathbb{Q})$. The complexification of the LLV algebra of a symplectic variety $X$ contains two important elements, denoted by $h_X$ and $h_X'$, which measure the cohomological grading and Hodge grading of the cohomology $\sH^\bullet(X; \mathbb{Q})$. In Chapters \ref{chap: fatcorwise psi} and \ref{chap: construction of Psi} we will construct an automorphism $\Psi^{\llv}$ of $\llv(Y; \mathbb{Q})$ and an automorphism $\Psi^{\sH}$ of $\sH^\bullet(Y; \mathbb{Q})$, equivariant with respect to $\Psi^{\llv}$, such that $\Psi^{\sH} \circ \Phi^{\sH}$ preserves the Hodge bigrading:
\begin{equation}
    \label{eq: The diagram}
    \begin{tikzcd}
        \sH^\bullet(X; \mathbb{Q}) \arrow[r, "\Phi^{\sH}"] \arrow[loop below] & \sH^\bullet(Y; \mathbb{Q})  \arrow[loop below] \arrow[r, "\Psi^{\sH}"] & \sH^\bullet(Y; \mathbb{Q}) \arrow[loop below]\\
        \llv(X; \mathbb{Q}) \arrow[r, "\Phi^{\llv}"] & \llv(Y; \mathbb{Q}) \arrow[r, "\Psi^{\llv}"] & \llv(Y; \mathbb{Q}). 
    \end{tikzcd}
\end{equation}

Our strategy will be to construct an algebraic group $\mathcal{G}$ over $\mathbb{Q}$ with an inclusion $\llv(Y; \mathbb{Q}) \subseteq \Lie(\mathcal{G})$ and an element $\psi \in \mathcal{G}(\mathbb{Q})$, such that the map $\Ad(\psi) \circ \Phi^{\llv} \colon \llv(X; \mathbb{Q}) \to \llv(Y; \mathbb{Q})$ sends $h_X$ to $h_Y$ and $h'_X$ to $h'_Y$. This will be done in Chapter \ref{chap: fatcorwise psi}. Ideally, we would then like to integrate the action of $\llv(Y; \mathbb{Q})$ on $\sH^\bullet(Y; \mathbb{Q})$ and use this to obtain a representation of $\mathcal{G}$ on $\sH^\bullet(Y; \mathbb{Q})$. This approach has some problems, and in Chapter \ref{chap: construction of Psi} we will modify the group $\mathcal{G}$ in a way that fixes these problems. The element $\psi \in \mathcal{G}(\mathbb{Q})$ will then induce an automorphism $\Psi^{\sH}$ of $\sH^\bullet(Y; \mathbb{Q})$. The fact that $\Ad(\psi) \circ \Phi^{\llv}$ sends $h_X$ to $h_Y$ and $h'_X$ to $h'_Y$ implies that $\Psi^{\sH} \circ \Phi^{\sH}$ preserves the Hodge bigrading, and this proves Theorem \ref{thm: main theorem}.

%% file: chapters/chapter2.tex
In this chapter we recall some definitions and results that will be used throughout the thesis. Most of the material here is well-known, and we provide proofs where no reference could be found. We also prove a few new results in this chapter which will be used later on.

\section{Lie algebras}

\subsection{Isomorphisms of semisimple Lie algebras}
\label{sec: product of simple lie algs}
Let $\mathfrak{g} = \mathfrak{g}_1 \times \mathfrak{g}_2$ be a product of Lie algebras over a field $F$. Given representations $V_i$ of $\mathfrak{g}_i$ for $i = 1, 2$, we can form the external tensor product $V = V_1 \boxtimes V_2$. This is a representation of $\mathfrak{g}_1 \times \mathfrak{g}_2$, with underlying space $V_1 \otimes V_2$, and action determined by letting $(X, Y) \in \mathfrak{g}_1 \times \mathfrak{g}_2$ act on $v \otimes w \in V_1 \otimes V_2$ by
$$(X, Y) \cdot (v \otimes w) = Xv \otimes w + v \otimes Yw.$$

Recall that for a Lie algebra $\mathfrak{g}$ over a field $F$, a representation $V$ of $\mathfrak{g}$ is called absolutely irreducible if $V \otimes_F \overline{F}$ is an irreducible representation of $\mathfrak{g} \otimes_F \overline{F}$.

The goal of this section is to prove the following theorem: 
\begin{theorem}
    \label{thm: split Lie alg reps}
    Let $\mathfrak{g} = \mathfrak{g}_1 \times \dots \times \mathfrak{g}_n$ and $\mathfrak{h} = \mathfrak{h}_1 \times \dots \times \mathfrak{h}_n$ be products of finite-dimensional simple Lie algebras over a subfield $F$ of $\mathbb{C}$. For each $i$, let $V_i$ be a representation of $\mathfrak{g}_i$, and let $W_i$ be a representation of $\mathfrak{h}_i$. Let $V = V_1 \boxtimes \dots \boxtimes V_n$ and $W = W_1 \boxtimes \dots \boxtimes W_n$, and assume that $V$ and $W$ are absolutely irreducible. Let $\psi \colon \mathfrak{g} \to \mathfrak{h}$ be an isomorphism of Lie algebras, and let $\varphi \colon V \to W$ be a linear isomorphism which is compatible with $\psi$, in the sense that $\varphi(Xv) = \psi(X)\varphi(v)$ for all $X \in \mathfrak{g}$ and $v \in V$. Then there exists a permutation $\sigma \in S_n$ such that: 
    \begin{enumerate}
        \item[(1)] The map $\psi$ is a product of isomorphisms $\psi_i \colon \mathfrak{g}_i \to \mathfrak{h}_{\sigma(i)}$.
        \item[(2)] The map $\varphi$ is a tensor product of isomorphisms $\varphi_i \colon V_i \to W_{\sigma(i)}$, which are compatible with the Lie algebra action under $\psi_i$.
    \end{enumerate} 
\end{theorem}
The first statement of the theorem is known, and follows from \cite[Theorem~1.54]{KnappLieGroups}. I could not find the second statement in the literature, so we prove it below.

We need a small lemma before we prove the above theorem.

\begin{lemma}
    \label{lem: decompose box product}
    Let $\mathfrak{g}_1$ and $\mathfrak{g}_2$ be semisimple Lie algebras over a field $F$ of characteristic~$0$. For $i = 1, 2$, let $V_i$ and $W_i$ be finite-dimensional representations of $\mathfrak{g}_i$, and assume that $V_1$ and $W_1$ are irreducible. Let $V = V_1 \boxtimes V_2$ and $W = W_1 \boxtimes W_2$. If $V \cong W$, then $V_1 \cong W_1$ and $V_2 \cong W_2$.
\end{lemma}
\begin{proof}
    Denote the inclusion $\mathfrak{g}_1 \to \mathfrak{g}_1 \times \mathfrak{g}_2$ by $i_1$. Then precomposition with $i_1$ turns $V$ into a $\mathfrak{g}_1$-representation. Explicitly, this representation is given by $X \cdot (v_1 \otimes v_2) = Xv_1 \otimes v_2$. This implies that $V$, as a $\mathfrak{g}_1$-representation, is isomorphic to $\dim(V_2)$ copies of $V_1$. Similarly, we see that $W$ is isomorphic to $\dim(W_2)$ copies of $W_1$. Hence $V_1^{\dim(V_2)} \cong W_1^{\dim(W_2)}$, and the irreducibility of $V_1$ and $W_1$ then implies that $V_1 \cong W_1$ and $\dim(V_2) = \dim(W_2)$. By an analogous argument we see that $V_2^{\dim(V_1)} \cong W_2^{\dim(W_1)}$ as representations of $\mathfrak{g}_2$. Since $\dim(V_1) = \dim(W_1)$, this implies $V_2 \cong W_2$. 
\end{proof}
\begin{remark}
    The assumption that $V_1$ and $W_1$ are irreducible is crucial. An easy counterexample is given by taking $V_1$ and $W_2$ the trivial representation, and $V_2$ and $W_1$ two copies of the trivial representation.
\end{remark}

We can now prove Theorem \ref{thm: split Lie alg reps}.

\begin{proof}[{Proof of Theorem \ref{thm: split Lie alg reps}}]
    By \cite[Theorem 1.54]{KnappLieGroups}, the minimal ideals of $\mathfrak{g}$ are exactly $\mathfrak{g}_1, \dots, \mathfrak{g}_n$, and the minimal ideals of $\mathfrak{h}$ are exactly $\mathfrak{h}_1, \dots, \mathfrak{h}_n$. Since an isomorphism sends minimal ideals to minimal ideals, we see that there must be a $\sigma \in S_n$ and isomorphisms $\psi_i \colon \mathfrak{g}_i \to \mathfrak{h}_{\sigma(i)}$ such that $\psi$ is the product of the $\psi_i$. This proves the first part of the theorem.

    For the second part, we use Lemma \ref{lem: decompose box product} and Schurs lemma. By reindexing the $\mathfrak{h}_i$, we may assume that $\sigma = \id$. Precomposition with $\psi$ turns $W$ into a $\mathfrak{g}$-representation, where the action is given by $X \cdot w := \psi(X)w$ for $X \in \mathfrak{g}$ and $w \in W$. We denote this representation by $W_\psi$. The map $\varphi \colon V \to W_\psi$ is an isomorphism of $\mathfrak{g}$-representations, since $\varphi(Xv) = \psi(X)\varphi(v)$ for all $X \in \mathfrak{g}$ and $v \in V$. Since $W = W_1 \boxtimes \dots \boxtimes W_n$ and $\psi$ is a direct product of the maps $\psi_i$, we see that $W_\psi = W_{1, \psi_1} \boxtimes \dots \boxtimes W_{n, \psi_n}$, where $W_{i, \psi_i}$ is the representation of $\mathfrak{g}_i$ with underlying vector space $W_i$, and where $X_i \in \mathfrak{g}_i$ acts on $w_i \in W_i$ via $X_i \cdot w_i := \psi_i(X_i)w_i$.

    Therefore we have an isomorphism $V_1 \boxtimes \dots \boxtimes V_n \cong  W_{1, \psi_1} \boxtimes \dots \boxtimes W_{n, \psi_n}$. Since $V$ and $W$ are irreducible, all $V_i$ and products $V_1 \boxtimes \dots \boxtimes V_i$ must also be irreducible, and the same holds for the $W_{i,\psi_i}$ and products $W_{1,\psi_1} \boxtimes \dots \boxtimes W_{i,\psi_i}$. Therefore, we can inductively apply Lemma \ref{lem: decompose box product} to conclude that there are isomorphisms of $\mathfrak{g}_i$-representations $\widetilde \varphi_i \colon V_i \to W_{i,\psi_i}$ for all $i$.

    To finish the proof, it suffices to show that $\varphi = \lambda\widetilde \varphi_1 \otimes \dots \otimes \widetilde \varphi_n$ for some $\lambda \in F^*$. By Schurs lemma, using that $V$ and $W$ are absolutely irreducible, the collection of isomorphisms $V \cong W_\psi$ is isomorphic to $F^*$. Since both $\varphi$ and $\widetilde\varphi_1 \otimes \dots \otimes \widetilde\varphi_n$ are such isomorphisms, we see that they must be a scalar multiple of each other.
\end{proof}

\subsection{The centroid of a simple Lie algebra}
\label{sec: centroid section}
We now recall some facts about the centroid of a simple Lie algebra. Given a Lie algebra $\mathfrak{g}$ over a field $k$, it could be possible that $\mathfrak{g}$ is actually defined over a bigger field. For example, we can consider $\mathfrak{sl}_2(\mathbb{Q}(i))$ as a Lie algebra over $\mathbb{Q}$, but it can also be given the structure of a Lie algebra over $\mathbb{Q}(i)$. If $\mathfrak{g}$ is a simple Lie algebra, the centroid is a field; it turns out to be the biggest field over which $\mathfrak{g}$ is defined. 

We follow the book by Jacobson \cite[Chapter X]{JacobsonLie}, but our definition of the centroid is slightly different from the one given there. The definition given in \cite{JacobsonLie} works for arbitrary (non-associative) algebras, while the definition that we give is specific to Lie algebras.

\begin{definition}
    Let $k$ be a field and let $\mathfrak{g}$ be a Lie algebra over $k$. The \textbf{centroid} of $\mathfrak{g}$, denoted by \myindex{$Z(\mathfrak{g})$}, is the endomorphism algebra of the adjoint representation $\ad \colon \mathfrak{g} \to \mathfrak{gl}(\mathfrak{g})$. Explicitly,
    \[Z(\mathfrak{g}) = \{\alpha \in \End(\mathfrak{g}) : \forall X \in \mathfrak{g}, \ad(X) \circ \alpha = \alpha \circ \ad(X)\}.\]
    \end{definition}
The proof of the following proposition is basically the same as given in \cite[Chapter~X, Lemma~1]{JacobsonLie}, but adapted to our definition of the centroid.
\begin{proposition}
    \label{prop: centroid commutative}
    Let $\mathfrak{g}$ be a Lie algebra over a field $k$, and assume that $\ad \colon \mathfrak{g} \to \mathfrak{gl}(\mathfrak{g})$ is injective. Then the centroid of $\mathfrak{g}$ is a commutative $k$-algebra. 
\end{proposition}

\begin{proof}
    We claim that for all $\alpha \in Z(\mathfrak{g})$ and $X \in \mathfrak{g}$, we have $\ad(X) \circ \alpha = \ad(\alpha(X))$. To see why, note that the equality $\ad(X) \circ \alpha = \alpha \circ \ad(X)$ means that for all $Y \in \mathfrak{g}$ we have
    $$[X, \alpha(Y)] = \alpha([X, Y]).$$
    It follows that also $[\alpha(X), Y] = \alpha([X, Y])$ for all $X, Y \in \mathfrak{g}$, so indeed $\ad(X) \circ \alpha = \ad(\alpha(X))$ for all $X \in \mathfrak{g}$. 

    If we now have $\alpha, \beta \in Z(\mathfrak{g})$, the above claim gives for all $X \in \mathfrak{g}$:
    $$\ad(\alpha\beta(X)) = \alpha \ad(\beta(X)) = \alpha \ad(X) \beta = \ad(\alpha(X))\beta = \ad(\beta\alpha(X)).$$
    The injectivity of $\ad$ then implies that $\alpha \beta(X) = \beta \alpha(X)$ for all $X \in \mathfrak{g}$, so indeed $\alpha\beta = \beta\alpha$.
\end{proof}

\begin{lemma}
    Let $\mathfrak{g}$ be a simple Lie algebra over a field $k$. Then $Z(\mathfrak{g})$ is a field. 
\end{lemma}

\begin{proof}
    Since $\mathfrak{g}$ is simple, the adjoint representation is a faithful irreducible representation. Then Schurs lemma implies that $Z(\mathfrak{g})$ is a division algebra. It is commutative by Proposition \ref{prop: centroid commutative}, so $Z(\mathfrak{g})$ is a field.
\end{proof}

In particular, $\mathfrak{g}$ is a Lie algebra over the field $Z(\mathfrak{g})$. The following result is a variant of \cite[Chapter X, Theorem 3]{JacobsonLie}:

\begin{lemma}
    \label{lem: centroid universal property}
    Let $\mathfrak{g}$ be a finite-dimensional simple Lie algebra over a field $k$. Then $k = Z(\mathfrak{g})$ if and only if $\mathfrak{g} \otimes_k \overline{k}$ is simple.
\end{lemma}

\begin{proof}
    We will first show that $\mathfrak{g} \otimes_{Z(\mathfrak{g})} \overline{Z(\mathfrak{g})}$ is simple. We know that $\mathfrak{g} \otimes_{Z(\mathfrak{g})} \overline{Z(\mathfrak{g})}$ is semisimple, so it is a product of simple Lie algebras. The endomorphisms of the adjoint representation of $\mathfrak{g} \otimes_{Z(\mathfrak{g})} \overline{Z(\mathfrak{g})}$ are given by the base change of the endomorphisms of the adjoint representation of $\mathfrak{g}$, so the endomorphism algebra of the adjoint representation of $\mathfrak{g} \otimes_{Z(\mathfrak{g})} \overline{Z(\mathfrak{g})}$ is $\overline{Z(\mathfrak{g})}$. This is a field, so the adjoint representation of $\mathfrak{g} \otimes_{Z(\mathfrak{g})} \overline{Z(\mathfrak{g})}$ is irreducible, which implies that $\mathfrak{g} \otimes_{Z(\mathfrak{g})} \overline{Z(\mathfrak{g})}$ is a simple Lie algebra.

    Now assume that $\mathfrak{g} \otimes_k \overline{k}$ is simple. For every $\lambda \in k$, multiplication by $\lambda$ is an endomorphism of the adjoint representation, so $k \subseteq Z(\mathfrak{g})$. Since $\mathfrak{g}$ has finite dimension over $k$, the extension $Z(\mathfrak{g}) / k$ is finite, so $\overline{k} = \overline{Z(\mathfrak{g})}$. Therefore, we have a surjective homomorphism of Lie algebras $\mathfrak{g} \otimes_k \overline{k} \to \mathfrak{g} \otimes_{Z(\mathfrak{g})} \overline{Z(\mathfrak{g})}$ given by $X \otimes \alpha \mapsto \alpha X$. If $k \neq Z(\mathfrak{g})$, this homomorphism is not injective for dimension reasons, which implies that $\mathfrak{g} \otimes_k \overline{k}$ is not simple.
\end{proof}

Let $k_0$ be the ground field of $Z(\mathfrak{g})$. If $\mathfrak{g}$ is finite-dimensional over $k_0$, then the condition that $\overline{k} = \overline{Z(\mathfrak{g})}$ is automatically satisfied for any field $k$ over which $\mathfrak{g}$ is defined.

\begin{definition}
    \label{def: type of lie alg}
    Let $k$ be a field and let $\mathfrak{g}$ be a finite-dimensional simple Lie algebra over $k$. If $\mathfrak{g} \otimes_{Z(\mathfrak{g})} \overline{k} \cong \mathfrak{sl}_{n+1}(\overline{k})$, we say that $\mathfrak{g}$ is of \textbf{type $A_{n}$}. We use the analogous terminology for the other simple types, e.g. $B_n, C_n$ and $D_n$. 
\end{definition}

Equivalently, $\mathfrak{g}$ is of type $A_n$ if $\mathfrak{g} \otimes_k \overline{k} \cong \prod_{i=1}^d \mathfrak{sl}_{n+1}(\overline{k})$ for some integer $d$, and we have similar statements for the other simple types.

As an example, consider the Lie algebra $\mathfrak{sl}_2(\mathbb{Q}(i))$ seen as a Lie algebra over $\mathbb{Q}$. This is a simple Lie algebra of type $A_1$, and we have $\mathfrak{sl}_2(\mathbb{Q}(i)) \otimes_\mathbb{Q} \overline{\mathbb{Q}} \cong \mathfrak{sl}_2(\overline{\mathbb{Q}}) \times \mathfrak{sl}_2(\overline{\mathbb{Q}})$.

\section{Algebras with involutions}
\label{sec: algebras with involutions}
In this section, we recall some properties of algebras with involutions. We will see in Chapter \ref{chap: LLV computation} how these can be used to describe the LLV algebra of varieties admitting a holomorphic symplectic form. The exposition here follows the Book of Involutions \cite{BookOfInvolutions}, and we also mention some results from Jacobson \cite[Chapter~X]{JacobsonLie}.

Throughout this section, $k$ will denote a field.

\subsection{Definition and examples}
\begin{definition}
    Let $A$ be an (associative, unital) algebra over $k$. An \textbf{involution} on $A$ is a map $\sigma \colon A \to A$ which is linear, anti-commutative (i.e. $\sigma(ab) = \sigma(b)\sigma(a)$ for $a, b \in A$) and satisfies $\sigma^2 = \id$. We call the pair $(A, \sigma)$ an \textbf{algebra with involution}.
\end{definition}

An involution is not required to be $k$-linear. For example, conjugation is an involution of the $\mathbb{C}$-algebra $\mathbb{C}$.

\begin{example}
    \label{ex: swapping involution}
    Let $B$ be a $k$-algebra. Then $A = B \times B^{\op}$ has an involution $\sigma$ defined by $\sigma(b_1, b_2) = (b_2, b_1)$. The center of $A$ is equal to $Z(B) \times Z(B)$, and $\sigma$ acts on this by swapping the two factors.
\end{example}

The following result can be found in \cite[Chapter 1]{BookOfInvolutions}:
\begin{lemma}
    \label{lem: involution from bilin form}
    Let $V$ be a vector space over $k$, equipped with a nondegenerate bilinear form $b$ that is either symmetric or anti-symmetric. Then for every $a \in \End_k(V)$ there is a unique endomorphism $\sigma_b(a) \in \End_k(V)$ such that for all $v, w \in V$ we have
    $$b(av, w) = b(v, \sigma_b(a)w).$$
    The obtained map \index{$\sigma_b$}$\sigma_b \colon \End_k(V) \to \End_k(V)$ is an involution. \qed
\end{lemma}
We denote the obtained algebra with involution by $(\End_k(V), \sigma_b)$.

\begin{definition}
    Let $(A, \sigma)$ be a $k$-algebra with involution, and let $V$ be an $A$-module. We say that a bilinear form $b \colon V \times V \to k$ is \textbf{$(A, \sigma)$-equivariant} if, for all $a \in A$ and $v, w \in V$ we have
    $$b(av, w) = b(v, \sigma(a)w).$$
\end{definition}

\begin{remark}
    If the involution $\sigma$ is clear from the context, we write $A$-equivariant instead of $(A, \sigma)$-equivariant.
\end{remark}

\begin{lemma}
    \label{lem: restriction of involution of equiv bilin form}
    Let $(A, \sigma)$ be an algebra with involution and let $V$ be an $A$-module with equivariant nondegenerate bilinear form $b$. Then the involution $\sigma_b$ on $\End_k(V)$ restricts to an involution on $\End_A(V)$.
\end{lemma}

We denote the obtained involution on $\End_A(V)$ by $\widetilde \sigma$.

\begin{proof}
    An $f \in \End_k(V)$ is in $\End_A(V)$ if and only if $f \circ a =  a \circ f$ for all $a \in A$. Therefore we want to check that $\sigma_b(f) \circ a = a \circ \sigma_b(f)$ for such $f$. If $f \in \End_A(V)$ and $a \in A$, we see for all $x, y \in V$ that:
    \begin{align*}
        b(a \sigma_b(f)(x), y) &= b(\sigma_b(f)(x), \sigma(a)y)\\
        &= b(x, f(\sigma(a)y)) \\
        &= b(x, \sigma(a)f(y))\\
        &= b(ax, f(y))\\
        &= b(\sigma_b(f)(ax), y),
    \end{align*}
    so we see that indeed $\sigma_b(f) \in \End_A(V)$.
\end{proof}

\begin{definition}[Dual of a module over an algebra with involution]
    \label{def: dual of module over algebra with involution}
    Let $(A, \sigma)$ be an algebra with involution over $k$, and let $V$ be an $A$-module. Then $V^* = \Hom(V, k)$ is a right $A$-module, and we turn this into a left $A$-module via the isomorphism $\sigma \colon A \to A^{\op}$. Explicitly, for $\eta \in V^*$ and $a \in A$, we define $a\eta \in V^*$ by
    $$(a\eta)(v) = \eta(\sigma(a)v)$$
    for all $v \in V$.
\end{definition}

We will often be interested in bilinear forms arising in the following special way. Let $V$ be a vector space over a field $k$. We write $\widetilde V = V \oplus V^*$. Then $\widetilde V$ carries a nondegenerate symmetric bilinear form, given by 
$$b((v, \eta), (v', \eta')) = \eta(v') + \eta'(v),$$
where $v, v' \in V$ and $\eta, \eta' \in V^*$.

\begin{lemma}
    \label{lem: bilin form on V tilde is equivariant}
    Let $(A, \sigma)$ be an algebra with involution over $k$ and let $V$ be an $A$-module. Then the bilinear form $b$ on $\widetilde V = V \oplus V^*$ defined above is $(A,\sigma)$-equivariant.
\end{lemma}
The lemma can be proved by a direct computation, which we leave to the reader. By combining this with Lemma \ref{lem: restriction of involution of equiv bilin form}, we obtain an involution \myindex{$\widetilde \sigma$} on $\End_A(\widetilde V)$.

\subsection{Classification of involutions}
Recall that a $k$-algebra $A$ is central simple if and only if $A \otimes_k \overline{k}$ is of the form $M_n(\overline{k})$ for some integer $n$. If $A$ is a finite-dimensional simple algebra, it is automatically central simple over its center. 

While the algebra $A = B \times B^{\op}$ in Example \ref{ex: swapping involution} is not simple, it does not have any nontrivial $\sigma$-invariant ideals if $B$ is simple. Therefore, from the perspective of algebras with involutions, it is simple in some sense. 

\begin{definition}
    Let $A$ be a semisimple algebra over $k$ and let $\sigma$ be an involution on $A$. If $A$ has no nontrivial two-sided ideals $I \subseteq A$ with $\sigma(I) = I$, we say that the pair $(A, \sigma)$ is a \textbf{simple algebra with involution}.
\end{definition}

\begin{remark}
    It should be noted that a simple algebra with involution need not be simple. The term should be read as: simple algebra-with-involution.
\end{remark}

\begin{definition}
    Let $(A, \sigma)$ be a simple algebra with involution. We say that $\sigma$ is an \textbf{involution of the first kind} if $\sigma|_{Z(A)} = \id_{Z(A)}$. If $\sigma$ does not act as the identity on $Z(A)$, we say that $\sigma$ is an \textbf{involution of the second kind}. 
\end{definition}

Using the classification of semisimple algebras, one easily sees the following:

\begin{lemma}
    Let $(A, \sigma)$ be a simple algebra with involution.
    \begin{enumerate}
        \item If $\sigma$ is an involution of the first kind, then $Z(A)$ is a field.
        \item If $\sigma$ is an involution of the second kind, then $\ell := Z(A)^\sigma$ is a field, and $Z(A)$ is either a degree $2$ field extension of $\ell$ or $Z(A) \cong \ell \times \ell$, with $\sigma|_{Z(A)}$ the swapping involution.
    \end{enumerate}
    If $Z(A)$ is a field, the underlying algebra $A$ is simple, while $A$ is a product of two simple algebras if $Z(A)$ is a product of two fields.\qed
\end{lemma}

For example, if $B$ is a simple $k$-algebra, the swapping involution on $B \times B^{\op}$ from Example \ref{ex: swapping involution} is an involution of the second kind. The involution constructed in Lemma \ref{lem: involution from bilin form} is an involution of the first kind.

\begin{lemma}
    \label{lem: Endomorphism inherits involution type}
    Let $(A, \sigma)$ be a finite-dimensional simple $k$-algebra with involution, and let $V$ be a nonzero finitely generated $A$-module equipped with a nondegenerate symmetric or anti-symmetric $A$-equivariant bilinear form $b$. Then $\sigma$ is an involution of the first kind if and only if the involution $\widetilde \sigma$ on $\End_A(V)$ is an involution of the first kind.
\end{lemma}
\begin{proof}
    We will first show that $Z(A) = Z(\End_A(V))$. If $Z(A)$ is a field, then $A$ is simple, and Wedderburn theory implies that $Z(A) = Z(\End_A(V))$. Otherwise, $Z(A)$ is a product of two fields, and hence $A$ is a product of two simple algebras and the involution $\sigma$ is of the second kind. Then $A$ has two simple modules up to isomorphism, say $V_1$ and $V_2$, and these are dual to each other. We can then write $V \cong V_1^m \oplus V_2^n$ for some integers $m$ and $n$. Since $V$ carries a nondegenerate $A$-equivariant bilinear form, we have an $A$-equivariant isomorphism $V \cong V^*$, which implies that $m = n$. In particular, $m$ and $n$ are both nonzero, so we also have $Z(A) = Z(\End_A(V))$. 
    
    We will show that $\sigma$ and $\widetilde \sigma$ coincide on $Z(A)$, so take a $\lambda \in Z(A)$. Since $b$ is $A$-equivariant, we have $b(\lambda v, w) = b(v, \sigma(\lambda)w)$ for all $v, w \in V$. The involution $\widetilde \sigma$ is the adjoint map for $b$, so we also have $b(\lambda v, w) = b(v, \widetilde \sigma(\lambda)w)$ for all $v, w \in V$. Hence 
    $$b(v, \sigma(\lambda)w) = b(v, \widetilde \sigma(\lambda)w)$$
    for all $v, w \in V$. Since $b$ is nondegenerate, it follows that $\sigma(\lambda) = \widetilde \sigma(\lambda)$. 
\end{proof}

For simple algebras with involutions, we have the following classification theorem, see \cite[Proposition 2.1 and Proposition 2.14]{BookOfInvolutions}:
\begin{theorem}
    \label{thm: involution classification}
    Let $k$ be a field of characteristic unequal to $2$ and let $(A,\sigma)$ be a finite-dimensional simple $k$-algebra with involution.
    \begin{enumerate}
        \item[(I)] If the involution $\sigma$ is of the first kind, then $\ell := Z(A)$ is a field, and either
    \begin{enumerate}
        \item[(1)] there is an integer $n \geq 1$ such that $A \otimes_\ell \overline{\ell} \cong M_n(\overline{\ell})$ and under this isomorphism $\sigma$ is given by the involution $\sigma_b$ associated to some nondegenerate symmetric bilinear form $b$ on $\overline{\ell}^n$, or
        \item[(2)] there is an even integer $n \geq 2$ such that $A \otimes_\ell \overline{\ell} \cong M_n(\overline{\ell})$ and under this isomorphism $\sigma$ is given by the involution $\sigma_b$ associated to some nondegenerate skew-symmetric bilinear form $b$ on $\overline{\ell}^n$.
    \end{enumerate}
    \item[(II)] If the involution $\sigma$ is of the second kind, then $\ell := Z(A)^\sigma$ is a field, and $Z(A) \otimes_{\ell} \overline{\ell} \cong \overline{\ell} \times \overline{\ell}$. Moreover:
    \begin{enumerate}
        \item[(3)] There is an integer $n \geq 1$ such that $A \otimes_{\ell} \overline{\ell} \cong M_n(\overline{\ell}) \times M_n(\overline{\ell})^{\op}$, and under this isomorphism $\sigma$ is given by $\sigma(x, y) = (y, x)$.
    \end{enumerate}
    \end{enumerate}
    
\end{theorem}

Using the classification theorem, we introduce a little bit more terminology:
\begin{definition}
    Let $(A, \sigma)$ be a simple algebra with involution. 
    \begin{enumerate}
        \item If we land in the first case of the above theorem, we say that $\sigma$ is an \textbf{orthogonal involution}.
        \item If we land in the second case of the above theorem, we say that $\sigma$ is a \textbf{symplectic involution}.
        \item If we land in the third case of the above theorem, we say that $\sigma$ is a \textbf{unitary involution}.
    \end{enumerate}
\end{definition}

In particular, an involution is of the second kind if and only if it is unitary, while an involution of the first kind is either orthogonal or symplectic.

\subsection{The group of isometries}
There are various groups that can be associated to an algebra with involution. In this section, we will look at the group of isometries, and consider this as an abstract group. Later, in Section \ref{subsubsec: alg groups of involutions}, we will look at the group of isometries as an algebraic group.
\begin{definition}
    Let $(A, \sigma)$ be an algebra with involution. The \textbf{group of isometries of} $(A, \sigma)$ is defined by\index{$\Iso(A, \sigma)$}
    $$\Iso(A, \sigma) = \{a \in A^\times : a \sigma(a) = 1\},$$
    where $A^\times \subseteq A$ is the group of invertible elements.
\end{definition}

Let $V$ be an $A$-module, and suppose that $V$ has an $A$-equivariant nondegenerate symmetric bilinear form $b$. Then we can look at the $A$-equivariant orthogonal group, defined by\index{$\sO_A(V, b)$}
$$\sO_A(V, b) = \{f \in \sO(V, b) : \forall a \in A, v \in V, f(av) = af(v)\}.$$
Recall from Lemma \ref{lem: restriction of involution of equiv bilin form} that $\End_A(V)$ carries an involution $\widetilde \sigma$.
\begin{lemma}
    \label{lem: SO_A is Iso(A, sigma)}
    There is an equality $\sO_A(V, b) = \Iso(\End_A(V), \widetilde\sigma)$.
\end{lemma}
\begin{proof}
    We have 
    $$\End_A(V)^\times = \{f \in \GL(V) : \forall a \in A, v \in V, f(av) = af(v)\}.$$
    Both $\Iso(\End_A(V), \widetilde \sigma)$ and $\sO_A(V, b)$ are subgroups of $\End_A(V)^\times$; it suffices to show that for $f \in \End_A(V)^\times$ we have $f \in \sO(V, b)$ if and only if $f \circ \widetilde \sigma(f) = 1$.

    An endomorphism $f$ of $V$ is in $\sO(V, b)$ if and only if for all $v, w \in V$ we have $b(f(v), f(w)) = b(v, w)$. But $b(f(v), f(w)) = b(v, (\widetilde\sigma(f)\circ f)(w))$. Therefore, we have 
    $$b(v, w) = b(v, (\widetilde \sigma(f) \circ f)(w))$$
    for all $v, w \in V$. Since $b$ is nondegenerate, this is equivalent to $\widetilde\sigma(f) \circ f = \id_V$, which is in turn equivalent to $f \circ \widetilde \sigma(f) = \id_V$.
\end{proof}

Theorem \ref{thm: involution classification} can be used to compute the group of isometries over an algebraically closed field. The following result can be found in \cite[Section~23.A]{BookOfInvolutions}.
\begin{lemma}
    \label{lem: isometries algebraically closed field}
    Let $k$ be a field of characteristic unequal to $2$, and let $(A, \sigma)$ be a finite-dimensional simple $k$-algebra with involution.
    \begin{enumerate}
        \item[(I)] If the involution $\sigma$ is of the first kind with $\ell := Z(A)$, let $n$ be the integer with $A \otimes_\ell \overline{\ell} \cong M_n(\overline{\ell})$.
    \begin{enumerate}
        \item[(1)] If $\sigma$ is orthogonal, then $\Iso(A \otimes_\ell \overline{\ell}, \sigma)$ is isomorphic to the orthogonal group $\sO(n, \overline{\ell})$.
        \item[(2)] If $\sigma$ is symplectic, then $\Iso(A \otimes_\ell \overline{\ell}, \sigma)$ is isomorphic to the symplectic group $\Sp(n, \overline{\ell})$.
    \end{enumerate}
    \item[(II)] If the involution $\sigma$ is of the second kind, let $\ell := Z(A)^\sigma$ and let $n$ be the integer such that $A \otimes_\ell \overline{\ell} \cong M_n(\overline{\ell}) \times M_n(\overline{\ell})$.
    \begin{enumerate}
        \item[(3)] There is an isomorphism $\Iso(A \otimes_\ell \overline{\ell}, \sigma) \cong \GL_n(\overline{\ell})$.
    \end{enumerate}
    \end{enumerate}
\end{lemma}
\begin{proof}
    This follows directly from Theorem \ref{thm: involution classification}.
\end{proof}

\subsection{Examples coming from representation theory}
\label{sec: alg with involution from representation theory}
Representation theory gives useful examples of algebras with involutions. Indeed, let $k$ be a field and let $G$ be a finite group. Then the group algebra 
$$k[G] = \left\{\sum_{g \in G} a_gg : a_g \in k\right\}$$
carries an involution $\sigma$ defined by $\sigma(g) = g^{-1}$.

We now have two ways of taking the dual of a $G$-representation $V$. On the one hand, we can take the dual as defined in representation theory, where we let $g \in G$ act on $\eta \in V^*$ by requiring that $(g\eta)(v) = \eta(g^{-1}v)$ for all $v \in V$. On the other hand, if we see $V$ as a $k[G]$-module, we can take its dual using the construction of Definition \ref{def: dual of module over algebra with involution}. It is easy to see that these two constructions coincide:
\begin{lemma}
    \label{lem: dual of representation is dual of module}
    Let $V$ be a $G$-representation. Then the dual of $V$ as a $G$-representation coincides with the dual of $V$ as a $k[G]$-module under the equivalence of categories between $G$-representations and $k[G]$-modules.\qed
\end{lemma}

\begin{definition}
    Let $G$ be a finite group and let $k$ be a field. We write \myindex{$\Irr_k(G)$} for a set of representatives of the irreducible $k$-valued $G$-representations.
\end{definition}

If the characteristic of $k$ does not divide $|G|$, the group algebra is semisimple, so we can write $k[G] = A_1 \times \dots \times A_r$ for simple $k$-algebras $A_1, \dots, A_r$. The isomorphism classes of simple $k[G]$-modules are in bijection with $\{A_1, \dots, A_r\}$. Therefore, the isomorphism classes of simple factors of $k[G]$ are in bijection with $\Irr_k(G)$.

\begin{lemma}
    \label{lem: rep self-dual iff involution of first kind}
    Let $G$ be a finite group, let $k$ be a field of characteristic coprime to $|G|$, let $\rho \in \Irr_k(G)$ and let $A$ be the simple factor of $k[G]$ associated to $\rho$. Then the involution $\sigma$ of $k[G]$ restricts to an involution $\sigma|_A$ on $A$ if and only if $\rho$ is self-dual. If $\rho$ is not self-dual, the involution $\sigma$ of $k[G]$ restricts to the swapping involution on $A \times A'$, where $A'$ denotes the simple factor of $k[G]$ associated to $\rho^*$. 
\end{lemma}
\begin{proof}
    This follows directly from Lemma~\ref{lem: dual of representation is dual of module}.
\end{proof}

If $k$ is not algebraically closed, the base change of an irreducible representation $\rho \in \Irr_k(G)$ to $\overline{k}$, denoted by $\rho_{\overline{k}} := \rho \otimes_k \overline{k}$, decomposes further into irreducible representations over $\overline{k}$, say $\rho_{\overline{k}} = \bigoplus_{i=1}^k \rho_i^{n_i}$ for some integers $n_i \in \mathbb{N}$. 

The above lemma implies the following:

\begin{corollary}
    \label{cor: self-dual constituents iff first kind}
    Let $G$ be a finite group, let $k$ be a field of odd characteristic coprime to $|G|$, let $\rho \in \Irr_k(G)$ be a self-dual irreducible representation and let $A$ be the simple factor of $k[G]$ associated to $\rho$. Then $\sigma|_A$ is an involution of the first kind if and only if all irreducible constituents of $\rho_{\overline{k}}$ are self-dual.
\end{corollary}
\begin{proof}
    The base change $A_{\overline{k}} := A \otimes_k \overline{k}$ is a product of simple factors, say $A_{\overline{k}} \cong \prod_i A_i$. Each of these simple factors is isomorphic to a matrix algebra $M_{n_i}(\overline{k})$. It follows from the classification in Theorem \ref{thm: involution classification} that $\sigma|_A$ is of the first kind exactly if it preserves these simple factors. And, as seen in Lemma~\ref{lem: rep self-dual iff involution of first kind}, this happens if and only if all irreducible constituents of $\rho_{\overline{k}}$ are self-dual.
\end{proof}

\begin{lemma}
    \label{lem: self dual can be checked on 1 rep}
    Let $G$ be a finite group, let $\rho \in \Irr_\mathbb{Q}(G)$ and let $\tau$ be an irreducible constituent of $\rho_{\mathbb{C}}$. Then $\tau$ is self-dual if and only if all irreducible constituents of $\rho_{\mathbb{C}}$ are self-dual.
\end{lemma}
\begin{proof}
    Decompose $\rho_{\mathbb{C}} = \bigoplus \rho_i^{n_i}$ into irreducible constituents, so $\tau$ is one of the $\rho_i$. Let $\chi_i$ be the character of $\rho_i$. Then the $\chi_i$ are all conjugate under the action of $\Gal(\overline{\mathbb{Q}}/\mathbb{Q})$, see \cite[Theorem 9.21(c)]{Isaacs}. Now $\rho_i$ is self-dual if and only if $\overline{\chi_i} = \chi_i$, and this property is invariant under the action of $\Gal(\overline{\mathbb{Q}}/\mathbb{Q})$.
\end{proof}

\subsection{The Lie algebra of skew elements}
We will now look at the collection of skew elements of an algebra with involution. This is a Lie algebra, and will play an important role in understanding the LLV algebra of a varieties admitting a holomorphic symplectic form.
\begin{definition}
    Let $(A, \sigma)$ be a $k$-algebra with involution. The \textbf{Lie algebra of skew elements} is the set\index{$\Skew(A, \sigma)$}
    $$\Skew(A, \sigma) := \{a \in A : \sigma(a) = -a\},$$
    equipped with the Lie bracket $[a,b] = ab - ba$. It is not hard to see that $\Skew(A, \sigma)$ is indeed a Lie algebra (over the field $k^\sigma$). We denote its derived Lie algebra by \myindex{$\mathfrak{s}(A, \sigma)$}, so 
    $$\mathfrak{s}(A, \sigma) = [\Skew(A, \sigma), \Skew(A, \sigma)].$$
    In particular, we have a chain of inclusions
    $$\mathfrak{s}(A, \sigma) \subseteq \Skew(A, \sigma) \subseteq A.$$
\end{definition}

For an algebra with involution of the form $(\End_A(V), \widetilde \sigma)$, the Lie algebra of skew elements has an alternative description in terms of the orthogonal or symplectic Lie algebra, similar to Lemma \ref{lem: SO_A is Iso(A, sigma)}:

\begin{lemma}
    \label{lem: Skew = so_A}
    Let $k$ be a field of characteristic not $2$, let $(A, \sigma)$ be an algebra with involution over $k$ and let $V$ be an $A$-module equipped with an $(A,\sigma)$-equivariant nondegenerate symmetric or skew-symmetric bilinear form $b$. Then 
    \begin{equation*}
        \Skew(\End_A(V), \widetilde \sigma) = \begin{cases}\mathfrak{so}(V) \cap \End_A(V) & \text{ if } b \text{ is symmetric}\\
        \mathfrak{sp}(V) \cap \End_A(V) & \text{ if } b \text{ is skew-symmetric}
        \end{cases}
    \end{equation*} 
    as Lie subalgebras of $\mathfrak{gl}(V)$.

    In particular, if $A = k$ and $\sigma = \id_k$, then $\Skew(\End_k(V), \sigma_b)$ is either $\mathfrak{so}(V, b)$ or $\mathfrak{sp}(V, b)$, depending on whether $b$ is symmetric or skew-symmetric.
\end{lemma}
\begin{proof}
    The proof is analogous to that of Lemma \ref{lem: SO_A is Iso(A, sigma)}, where we now observe that an $f \in \End_k(V)$ is in $\mathfrak{so}(V)$ if and only if $b(f(v), w) + b(v, f(w)) = 0$ for all $v, w \in V$, and that this is equivalent to $\widetilde\sigma(f) + f = 0$.
\end{proof}

Let $k$ be a field of characteristic not $2$ and let $(A, \sigma)$ be an algebra with involution over $k$. To understand $\mathfrak{s}(A, \sigma)$, we first compute it after the base change to $\overline k$. The following theorem is a summary of some results that can be found in \cite[Chapter X]{JacobsonLie}:

\begin{theorem}
    \label{Thm: computation of s(A, sigma)}
    Let $k$ be a field of characteristic not equal to $2$ and let $(A, \sigma)$ be a finite-dimensional simple algebra with involution over $k$. Let $n$ be the degree of $A$ over its center, so $n^2 = [A : Z(A)]$.
    \begin{enumerate}
        \item[(I)] If $\sigma$ is an involution of the first kind, then $\ell := Z(A)$ is a field and:
    \begin{enumerate}
        \item[(1)] If $\sigma$ is symplectic, then $n$ is even and $\mathfrak{s}(A \otimes_{\ell} \overline{\ell}, \sigma) \cong \mathfrak{sp}_{n}(\overline \ell)$.
        \item[(2)] If $\sigma$ is orthogonal, then $\mathfrak{s}(A \otimes_{\ell} \overline{\ell}, \sigma) \cong \mathfrak{so}_{n}(\overline \ell)$.
    \end{enumerate}
    \item[(II)] If $\sigma$ is an involution of the second kind, then $\ell := Z(A)^\sigma$ is a field. Moreover:
    \begin{enumerate}
        \item[(3)] The involution $\sigma$ is unitary, and $\mathfrak{s}(A \otimes_{\ell} \overline{\ell}, \sigma) \cong \mathfrak{sl}_{n}(\overline{\ell})$.
    \end{enumerate}
    \end{enumerate}
    In all cases above where $\mathfrak{s}(A \otimes_{\ell} \overline \ell, \sigma)$ is simple, the centroid of $\mathfrak{s}(A, \sigma)$ is equal to $\ell = Z(A)^\sigma$.
\end{theorem}

In the terminology of Definition \ref{def: type of lie alg}, we have the following:

\begin{enumerate}
    \item If $\sigma$ is symplectic, then $\mathfrak{s}(A, \sigma)$ is of type $C_l$ for some $l \in \mathbb{N}$. This is simple if $l \geq 1$.
    \item If $\sigma$ is orthogonal, then $\mathfrak{s}(A, \sigma)$ is either of type $D_l$ or type $B_l$ for some $l \in \mathbb{N}$. If it is of type $D_l$, it is simple for $l \geq 3$, while it is simple for all $l \geq 1$ if it is of type $B_l$.
    \item If $\sigma$ is unitary, then $\mathfrak{s}(A, \sigma)$ is of type $A_{l-1}$ for some $l$. This is simple if $l \geq 2$.
\end{enumerate}

\begin{proof}[{Proof of Theorem \ref{Thm: computation of s(A, sigma)}}]
    To prove the theorem, we first compute $\mathfrak{s}(A \otimes_{\ell} \overline{\ell}, \sigma)$. First assume that $\sigma$ is an involution of the first kind. By Theorem \ref{thm: involution classification}, we have $A \otimes_\ell \overline \ell \cong M_n(\ell) = \End_k(\overline \ell^n)$, and $\sigma = \sigma_b$ for some nondegenerate bilinear form $b$ on $\overline{\ell}^n$. Then Lemma \ref{lem: Skew = so_A} gives that $\mathfrak{s}(A \otimes_{\ell} \overline{\ell}, \sigma)$ is isomorphic to either $\mathfrak{sp}_n(\overline{\ell})$ or $\mathfrak{so}_n(\overline{\ell})$ if $\sigma$ is of the first kind. 
    
    If $\sigma$ is an involution of the second kind, Theorem \ref{thm: involution classification} gives an isomorphism $A \otimes_\ell \overline \ell \cong M_n(\overline{\ell}) \times M_n(\overline{\ell})^{\op}$, with the involution given by swapping the two factors. An explicit computation shows that the embedding $\mathfrak{gl}(\overline{\ell}^n) \to M_n(\overline{\ell}) \times M_n(\overline{\ell})^{\op}$ given by $x \mapsto (x, -x)$ induces an isomorphism of Lie algebras 
    $$\mathfrak{gl}(\overline{\ell}^n) \to \Skew(M_n(\overline{\ell}) \times M_n(\overline{\ell})^{\op}, \sigma),$$ and hence we obtain an isomorphism $\mathfrak{sl}(\overline{\ell}^n) \isomto \mathfrak{s}(M_n(\overline{\ell}) \times M_n(\overline{\ell})^{\op}, \sigma)$.

    It remains to show that $\ell = Z(\mathfrak{s}(A, \sigma))$ if $\mathfrak{s}(A, \sigma) \otimes_{\ell} \overline \ell$ is simple. For any $\lambda \in \ell = Z(A)^\sigma$, the linear map $A \to A$ given by multiplication by $\lambda$ respects the skew elements (if $\sigma(b) = -b$, then also $\sigma(\lambda b) = \sigma(b)\sigma(\lambda) = -b\lambda = -\lambda b$). Then multiplication by $\lambda$ restricts to a map $\mathfrak{s}(A, \sigma) \to \mathfrak{s}(A, \sigma)$, so we have an inclusion $\ell \subseteq Z(\mathfrak{s}(A, \sigma))$. Since $\mathfrak{s}(A, \sigma) \otimes_{\ell} \overline \ell$ is simple, Proposition \ref{lem: centroid universal property} gives an equality $\ell = Z(\mathfrak{s}(A, \sigma))$. 
\end{proof}

A corollary of Theorem \ref{Thm: computation of s(A, sigma)} is that it gives a numerical criterion for detecting whether an involution of the first kind is orthogonal or symplectic. The following result is \cite[Proposition 2.6]{BookOfInvolutions}:
\begin{lemma}
    \label{lem: involution type from skew elements}
    Let $k$ be a field of characteristic~$0$ and let $A$ be a central simple algebra over $k$. Let $n$ be the degree of $A$ over $k = Z(A)$. Assume that $\sigma$ is an involution of the first kind on $A$. Then:
    \begin{itemize}
        \item The dimension of $\Skew(A, \sigma)$ is equal to $\frac{n(n-1)}{2}$ if and only if $\sigma$ is orthogonal.
        \item The dimension of $\Skew(A, \sigma)$ is equal to $\frac{n(n+1)}{2}$ if and only if $\sigma$ is symplectic.
    \end{itemize}
\end{lemma}

\subsection{Reconstructing an algebra with involution from its Lie algebra of skew elements}
An isomorphism $(A, \sigma) \to (B, \tau)$ of algebras with involution induces an isomorphism between the Lie algebras $\mathfrak{s}(A, \sigma) \to \mathfrak{s}(B, \tau)$. Conversely, a result by Jacobson says that, under suitable assumptions on $A$, every such isomorphism comes from an isomorphism of algebras with involution. The following result is \cite[Chapter X, Theorems 11 and 12]{JacobsonLie}:
\begin{theorem}
    \label{thm: isom lifts from skew to algebra}
    Let $k$ be a field of characteristic zero, and let $(A, \sigma)$ and $(B, \tau)$ be two finite-dimensional simple algebras with involution over $k$, and assume that $\sigma$ and $\tau$ restrict to the identity on $k$. Let $\psi \colon \mathfrak{s}(A, \sigma) \to \mathfrak{s}(B, \tau)$ be an isomorphism of Lie algebras over $k$, and assume that we are in one of the following cases:
    \begin{enumerate}
        \item $\sigma$ and $\tau$ are unitary, and $\mathfrak{s}(A, \sigma)$ is of type $A_l$ with $l \geq 2$.
        \item $\sigma$ and $\tau$ are orthogonal, and $\mathfrak{s}(A, \sigma)$ is of type $B_l$ with $l \geq 2$.
        \item $\sigma$ and $\tau$ are symplectic, and $\mathfrak{s}(A, \sigma)$ is of type $C_l$ with $l \geq 2$.
        \item $\sigma$ and $\tau$ are orthogonal, and $\mathfrak{s}(A, \sigma)$ is of type $D_l$ with $l = 3$ or $l \geq 5$.
    \end{enumerate}
    
    Then $\psi$ extends in a unique way to an isomorphism $\widetilde \psi \colon A \to B$ of $k$-algebras, satisfying $\tau \circ \widetilde \psi = \widetilde \psi \circ \sigma$.
\end{theorem}

\begin{proof}
    The proof is already given in the book of Jacobson, except for the cases $C_2$ and $D_3$, and Jacobson does not prove that $\widetilde \psi$ commutes with the involutions.

    Since the proof of the $D_3$ case goes the same as that for $D_l$ with $l \geq 5$, we give a uniform proof for both cases. We omit the proof of the $C_2$ case, which is analogous to the proof of the $C_l$ case with $l \geq 3$ as given in \cite{JacobsonLie}.

    So, assume that $\mathfrak{s}(A, \sigma)$ is of type $D_l$, and let $\psi \colon \mathfrak{s}(A, \sigma) \to \mathfrak{s}(B, \tau)$ be an isomorphism. Let $\ell$ be the center of $A$, and let $\overline \ell$ be an algebraic closure of $\ell$. Then $\ell$ is also the centroid of $\mathfrak{s}(A, \sigma)$ by Theorem \ref{Thm: computation of s(A, sigma)}. Similarly, $Z(B)$ is isomorphic to the centroid of $\mathfrak{s}(B, \tau)$, so $\psi$ induces an isomorphism $\ell \to Z(B)$. Under this isomorphism, the map $\psi \colon \mathfrak{s}(A, \sigma) \to \mathfrak{s}(B,\tau)$ is $\ell$-linear. Throughout the rest of this proof, we will denote a base change from $\ell$ to $\overline{\ell}$ by a subscript $\overline{\ell}$.

    For every integer $n$, we equip $M_n(\overline{\ell})$ with the standard orthogonal involution $\sigma_{\std}$, given by the transpose $X \mapsto X^t$. There exist isomorphisms of algebras with involutions $\varphi_A \colon (A_{\overline \ell}, \sigma) \isomto (M_n(\overline \ell), \sigma_{\std})$ and $\varphi_B \colon (B_{\overline \ell}, \tau) \isomto (M_m(\overline \ell), \sigma_{\std})$ for some integers $n$ and $m$. Since $\mathfrak{s}(A, \sigma)$ and $\mathfrak{s}(B, \tau)$ are isomorphic, we have $m = n$. Then $\varphi_A$ and $\varphi_B$ also induce isomorphisms $\mathfrak{s}(A, \sigma)_{\overline{\ell}} \to \mathfrak{s}(M_n(\overline{\ell}), \sigma_{\std})$ and $\mathfrak{s}(B, \tau)_{\overline{\ell}} \to \mathfrak{s}(M_n(\overline{\ell}), \sigma_{\std})$.

    Let $\chi := \varphi_B \circ \psi \circ \varphi_A^{-1}$, this is an automorphism of $\mathfrak{s}(M_n(\overline{\ell}), \sigma_{\std})$. We now have the solid arrows in the following diagram: 
    \begin{equation*}
        \begin{tikzcd}
            A_{\overline{\ell}} \arrow[r, phantom, "\supseteq"] \arrow[d, dashrightarrow, "\widetilde \psi"] & \mathfrak{s}(A, \sigma)_{\overline{\ell}} \arrow[r, "\varphi_A"', "\sim"] \arrow[d, "\psi", "\sim" {anchor=south, rotate=90}] & \mathfrak{s}(M_n(\overline{\ell}), \sigma_{\std}) \arrow[d, "\chi", "\sim" {anchor=south, rotate=90}] \arrow[r, phantom, "\subseteq"] & M_n(\overline{\ell}) \arrow[d, dashrightarrow, "\widetilde \chi"]\\
            B_{\overline{\ell}} \arrow[r, phantom, "\supseteq"] & \mathfrak{s}(B, \tau)_{\overline{\ell}} \arrow[r, "\varphi_B"', "\sim"] & \mathfrak{s}(M_n(\overline{\ell}), \sigma_{\std}) \arrow[r, phantom, "\subseteq"] & M_n(\overline{\ell}).
        \end{tikzcd}
    \end{equation*}

    The classification of automorphisms of Lie algebras over $\overline \ell$ gives that $\chi$ must be of the form $X \mapsto MXM^{-1}$ for a certain $M \in \sO_{n}(\overline \ell)$ (this is why we exclude the $D_4$ case, where there are more automorphisms due to triality). Then $\chi$ extends to an algebra automorphism $\widetilde \chi$ of $M_{n}(\overline \ell)$, also given by conjugation by $M$.

     Now let $\widetilde \psi = \varphi_B^{-1} \circ \widetilde \chi \circ \varphi_A$, this is an isomorphism $A_{\overline{\ell}} \to B_{\overline{\ell}}$. By \cite[Chapter X, Lemma 4]{JacobsonLie}, the subalgebra of $A$ generated by $\mathfrak{s}(A, \sigma)$ equals all of $A$, and the subalgebra of $B$ generated by $\mathfrak{s}(B, \tau)$ is all of $B$. Since $\widetilde \psi$ extends $\psi$, and $\psi$ sends $\mathfrak{s}(A, \sigma)$ to $\mathfrak{s}(B, \tau)$, it follows that $\widetilde \psi$ restricts to an isomorphism $A \to B$. 

     Since the involution is given by multiplication by $-1$ on $\mathfrak{s}(A, \sigma)$ and $\mathfrak{s}(B, \tau)$, the map $\psi$ respects the involution. Since $A$ and $B$ are generated by $\mathfrak{s}(A, \sigma)$ and $\mathfrak{s}(B, \tau)$ respectively, it follows that $\widetilde \psi$ also respects the involution.
\end{proof}

\begin{remark}
    It seems to me that the cases of type $C_2$ and $D_3$ are excluded in \cite[Chapter X]{JacobsonLie} due to the exceptional isomorphisms $B_2 \cong C_2$ and $A_3 \cong D_3$. Jacobson does not assume that $\sigma$ and $\tau$ are both of the same type, and then such an exceptional isomorphism prevents the proof from working. Since we assume a priori that $\sigma$ and $\tau$ are involutions of the same type, we can also treat these two cases. 
\end{remark}

Nevertheless, a few low-dimensional cases are absent from the above statement, and the theorem is not necessarily true in these cases. This is one reason why we assume the varieties to satisfy condition \ref{condition *} in our main theorem. There are various reasons for the failure of Theorem~\ref{thm: isom lifts from skew to algebra} in the cases that are not covered; we highlight two of them here.

In the $D_4$ case, the proof breaks down due to triality, and this can be used to construct a counterexample. Let $A = B = M_8(\mathbb{C})$, equipped with the orthogonal involution $\sigma$ defined by $X \mapsto X^t$ for $X \in M_8(\mathbb{C})$. If the theorem were true in this case, it would say that every automorphism of $\mathfrak{s}(A, \sigma)$ lifts to an automorphism of $M_8(\mathbb{C})$. The Skolem-Noether theorem implies that $\Aut(M_8(\mathbb{C}), \sigma) \cong \sO_8(\mathbb{C})/Z(\sO_8(\mathbb{C}))$, which has two connected components. This is smaller than $\Aut(\mathfrak{so}_8(\mathbb{C}))$, which has $6$ connected components.

The second counterexample is for the $A_1$ case with $\sigma$ and $\tau$ unitary involutions. We will explicitly write down non-isomorphic algebras with involution $(A, \sigma)$ and $(B, \tau)$, such that $\mathfrak{s}(A, \sigma) \cong \mathfrak{s}(B, \tau)$, where $\mathfrak{s}(A, \sigma)$ is of type $A_1$. 

Let $Q = \mathbb{Q} \oplus \mathbb{Q}i \oplus \mathbb{Q}j \oplus \mathbb{Q}k$ be the quaternion algebra with $i^2 = j^2 = k^2 = -1$. For $\alpha = t + xi + yj + zk \in Q$, let 
$$\overline{\alpha} = t - xi - yj - zk,$$
this defines an involution of $Q$. Let $N$ be a squarefree integer, and let $A = Q \otimes_\mathbb{Q} \mathbb{Q}(\sqrt{N})$. If we write $\tau$ for the nontrivial automorphism of $\mathbb{Q}(\sqrt{N})$, then the map $\sigma \colon A \to A$ defined by $\alpha \otimes z \mapsto \overline{\alpha} \otimes \tau(z)$ is an involution of the second kind on $A$.

Elements $\alpha \in A$ can be represented as $\alpha = a + bi + cj + dk$ with $a,b,c,d \in \mathbb{Q}(\sqrt{N})$. Then $-\sigma(\alpha) = -\tau(a) + \tau(b)i + \tau(c)j + \tau(d)k$, so $\alpha = -\sigma(\alpha)$ exactly if $\alpha$ is of the form
$$\alpha = x\sqrt{N} + ui + vj + wk,$$
where $x, u, v, w \in \mathbb{Q}$. By computing the derived Lie algebra, we see that
$$\mathfrak{s}(A, \sigma) = \{ui + vj + wk : u, v, w \in \mathbb{Q}\}.$$
This Lie algebra is independent of $N$, while the algebras $A$ for different choices of~$N$ are non-isomorphic, since they have different centers. 

The behaviour seen here can be explained by Galois cohomology. Let $(A, \sigma)$ be an algebra with involution, where $\sigma$ is an involution of the second kind. If $\mathfrak{s}(A, \sigma) \otimes_\mathbb{Q} \overline{\mathbb{Q}} \cong \mathfrak{sl}_2(\overline{\mathbb{Q}})$, then $A \otimes_\mathbb{Q} \overline{\mathbb{Q}} \cong M_2(\overline{\mathbb{Q}}) \times M_2(\overline{\mathbb{Q}})^{\op}$, where $\sigma$ becomes the swapping involution $\tau$. The collection of automorphisms of the algebra with involution $(M_2(\overline{\mathbb{Q}}) \times M_2(\overline{\mathbb{Q}})^{\op}, \tau)$ is equal to $\PGL_2(\overline{\mathbb{Q}}) \rtimes \mathbb{Z}/2\mathbb{Z}$, so the twisted forms of the algebra with involution $(A, \sigma)$ are classified by the Galois cohomology set
$$\sH^1(\mathbb{Q}, \PGL_2(\overline{\mathbb{Q}}) \rtimes \mathbb{Z}/2\mathbb{Z}).$$
On the other hand, the twisted forms of $\mathfrak{sl}_2(\mathbb{Q})$ are classified by the cohomology set $\sH^1(\mathbb{Q}, \Aut(\mathfrak{sl}_2(\overline{\mathbb{Q}})))$, where $\Aut(\mathfrak{sl}_2(\overline{\mathbb{Q}})) \cong \PGL_2(\overline{\mathbb{Q}})$. Every automorphism of $(A, \sigma)$ induces an automorphism of $\mathfrak{s}(A, \sigma)$, and this yields a surjective group homomorphism
$$\pi \colon \PGL_2(\overline{\mathbb{Q}}) \rtimes \mathbb{Z}/2\mathbb{Z} \to \PGL_2(\overline{\mathbb{Q}}).$$
The kernel of $\pi$ is isomorphic to $\mathbb{Z}/2\mathbb{Z}$, equipped with the trivial action of $\Gal(\overline{\mathbb{Q}}/\mathbb{Q})$. This yields a long exact sequence on cohomology
\begin{equation*}
    \begin{tikzcd}
        \mathbb{Z}/2\mathbb{Z} \arrow[r] & \PGL_2(\mathbb{Q}) \rtimes \mathbb{Z}/2\mathbb{Z} \arrow[r, "\pi"] \arrow[d, phantom, ""{coordinate, name=W}] & \PGL_2(\mathbb{Q}) \arrow[dll, rounded corners, to path={ -- ([xshift=2ex]\tikztostart.east)
    |- (W) [near end]\tikztonodes
-| ([xshift=-2ex]\tikztotarget.west)
-- (\tikztotarget)}]\\
        \sH^1(\mathbb{Q}, \mathbb{Z}/2\mathbb{Z}) \arrow[r, "\iota"] & \sH^1(\mathbb{Q}, \PGL_2(\overline{\mathbb{Q}}) \rtimes \mathbb{Z}/2\mathbb{Z}) \arrow[r] & \sH^1(\mathbb{Q}, \PGL_2(\overline{\mathbb{Q}})).
    \end{tikzcd}
\end{equation*}
Since $\pi$ is surjective, the kernel of $\iota$ is trivial. Now $\sH^1(\mathbb{Q}, \mathbb{Z}/2\mathbb{Z}) \cong \mathbb{Q}^*/\mathbb{Q}^{*2}$, which is infinite and classifies the quadratic extensions of $\mathbb{Q}$, see \cite{WashingtonGaloisCohomology}. This analysis can be continued to get a more detailed classification of all algebras with involution $(A, \sigma)$ having a fixed twisted form of $\mathfrak{sl}_2(\mathbb{Q})$ as their Lie algebra $\mathfrak{s}(A, \sigma)$.

\section{Algebraic groups}
\label{subsec: algebraic groups}
In this section we will discuss some results about algebraic groups that we will need. We mostly follow the exposition from the Book of Involutions \cite[Chapter VI]{BookOfInvolutions} and Milne \cite{milneALA, milneiAG}. Proposition~\ref{prop: Rmu to RG extension} and Corollary~\ref{cor: extending mu_n-rep to G_m-rep} are new results as far as I know.

Since we are mostly interested in affine algebraic groups, we will use the following definition of an algebraic group:
\begin{definition}
    Let $k$ be a field. An \textbf{algebraic group over} $k$ is an affine group scheme $G$ of finite type over $k$.
\end{definition}

We will denote algebraic groups by an underline. For example, for a field~$k$ and a $k$-vector space~$V$, we let \myindex{$\uGL(V)$} be the algebraic group whose $R$-points for a commutative $k$-algebra~$R$ are given by $\GL(V \otimes_k R)$.

\subsection{Exactness in characteristic $0$}
Let $k$ be a field, and let $f \colon G \to H$ be a morphism of algebraic groups over $k$. Recall that the kernel of $f$ is the algebraic group $\ker(f)$ whose $R$-points for a commutative $k$-algebra $R$ are given by 
$$\ker(f)(R) = \ker(f \colon G(R) \to H(R)).$$
If we write $G = \Spec(A)$ and $H = \Spec(B)$ with $A$ and $B$ Hopf algebras, then $f$ corresponds to a morphism of Hopf algebras $f^* \colon B \to A$. The image of $f$ is the algebraic subgroup of $H$ corresponding to $B/\ker(f^*)$. 

Suppose we now have a sequence of morphisms
\begin{equation*}
    \begin{tikzcd}
        G_1 \arrow[r, "f"] & G_2 \arrow[r, "g"] & G_3.
    \end{tikzcd}
\end{equation*}
We say that this is exact in $G_2$ if $\im(f) = \ker(g)$ as algebraic subgroups of $G_2$.

Recall that for an algebraic group $G$ over a field $k$, its Lie algebra is the tangent space at the identity, and this is a Lie algebra over $k$. We denote this Lie algebra by $\Lie(G)$.

We now recall a criterion to check whether a sequence of algebraic groups is exact. It is probably well-known, but the author could not find a full proof in the literature.

\begin{lemma}
\label{lem: exactness algebraic closure}
    Let $k$ be a field of characteristic zero, and let 
    \begin{equation}
        \label{eq: ses alg grps}
        \begin{tikzcd}
            1 \arrow[r] & G_1 \arrow[r, "\alpha"]& G_2 \arrow[r, "\beta"]& G_3 \arrow[r]& 1
        \end{tikzcd}
    \end{equation}
    be a sequence of algebraic groups over $k$. Then this sequence is exact if and only if the sequence
    \begin{equation}
        \label{eq: ses k bar points}
        \begin{tikzcd}
            1 \arrow[r] & G_1(\overline{k}) \arrow[r] & G_2(\overline{k}) \arrow[r] & G_3(\overline{k}) \arrow[r] & 1
        \end{tikzcd}
    \end{equation}
    is exact.
\end{lemma}

\begin{proof}
    By \cite[Proposition 22.15]{BookOfInvolutions}, taking $\overline{k}$-points is exact, so if \eqref{eq: ses alg grps} is exact, then \eqref{eq: ses k bar points} is also exact.
    
    Now suppose that the sequence of closed points is exact. We will prove that \eqref{eq: ses alg grps} is exact. First observe that $G_2 \to G_3$ is surjective by \cite[Proposition~22.3]{BookOfInvolutions}. Therefore, it suffices to show that $G_1 = \ker(\beta)$. Let $N = \ker(\beta)$, we then have a short exact sequence $1 \to N \to G_2 \to G_3 \to 1$. Since \eqref{eq: ses k bar points} is exact, we have an equality $\im(\alpha)(\overline{k}) = \ker(\beta)(\overline{k}) = N(\overline{k})$, so \cite[Corollary 3.26]{milneiAG} implies that $N = \im(\alpha)$ as subgroups op $G_2$. This gives a morphism $G_1 \to N$, and we obtain a commutative diagram
    \begin{equation*}
        \begin{tikzcd}
            1 \arrow[r] & G_1 \arrow[r, "\alpha"] \arrow[d]& G_2 \arrow[r, "\beta"] \arrow[d]& G_3 \arrow[r] \arrow[d]& 1\\
            1 \arrow[r] & N \arrow[r]& G_2 \arrow[r, "\beta"]& G_3 \arrow[r] & 1,
        \end{tikzcd}
    \end{equation*}
    where the bottom row is exact. We have to prove that the top row is exact. 
    
    By assumption, the map $G_1 \to N$ induces an isomorphism $G_1(\overline{k}) \to N(\overline{k})$. Then \cite[Proposition 22.3]{BookOfInvolutions} gives that the map $G_1 \to N$ is surjective. For injectivity, let $K = \ker(G_1 \to N)$. Then $K(\overline{k}) = \{e\}$ (where $e \in G_1(\overline{k})$ is the unit element), so $\dim(K) = 0$. Since $k$ has characteristic zero, all algebraic groups over $k$ are smooth \cite[Proposition~21.10]{BookOfInvolutions}, so $\Lie(K) = \{0\}$. Therefore, the map $G_1 \to N$ induces an injection $\Lie(G_1) \to \Lie(N)$, so $G_1 \to N$ is injective by \cite[Proposition~22.2]{BookOfInvolutions}. Then \cite[Proposition~22.5]{BookOfInvolutions} implies that $G_1 \to N$ is an isomorphism, so $G_1 = \ker(\beta)$.
\end{proof}

\subsection{Weil restriction}
Given a finite separable extension of fields $F/k$, we will often want to obtain an algebraic group over $k$ from an algebraic group defined over $F$. This can be done by using Weil restriction. 
\begin{definition}
    Let $F/k$ be a finite separable field extension and let $G$ be an algebraic group over $F$. The \textbf{Weil restriction} of $G$, denoted by \index{$\sR_{F/k}$}$\sR_{F/k}G$ is the algebraic group over $k$ representing the functor $A \mapsto G(A \otimes_k F)$ for a commutative $k$-algebra $A$ (this functor is representable by \cite[Lemma 20.6]{BookOfInvolutions}).
\end{definition}

We will need a few properties of algebraic groups and Weil restriction.

\begin{lemma}
    \label{lem: Lie algebra of Weil restriction}
    Let $F/k$ be a finite separable field extension and let $G$ be a smooth algebraic group over $F$ with Lie algebra $\mathfrak{g}$. Then the Lie algebra of $\sR_{F/k}G$ is isomorphic to $\mathfrak{g}$ seen as a Lie algebra over $k$ (instead of a Lie algebra over $F$).
\end{lemma}
\begin{proof}
    See \cite[Proposition 21.4(3)]{BookOfInvolutions}.
\end{proof}

\begin{lemma}[{\cite[Proposition 20.7]{BookOfInvolutions}}]
    \label{lem: Weil restriction F-bar}
    Let $F/k$ be a finite separable field extension of degree $d$, and let $\overline{F}$ be an algebraic closure of $F$. Then $(\sR_{F/k}G)(\overline F) \cong (G(\overline{F}))^d$.
\end{lemma}
\begin{proof}
    This follows from the definition of Weil restriction. We have:
    $$(\sR_{F/k}G)(\overline{F}) = G(F \otimes_k \overline{F}) \cong (G(\overline{F}))^d.$$
\end{proof}

\begin{lemma}
    \label{lem: weil restriction exact}
    Let $F/k$ be a finite extension of fields of characteristic~$0$ and let $1 \to G_1 \to G_2 \to G_3 \to 1$ be a short exact sequence of algebraic groups over $F$. Then the induced sequence
    \begin{equation*}
        \begin{tikzcd}
            1 \arrow[r] &\sR_{F/k}G_1 \arrow[r] &\sR_{F/k}G_2 \arrow[r] &\sR_{F/k}G_3 \arrow[r] &1
        \end{tikzcd}
    \end{equation*}
    is also exact.
\end{lemma}
\begin{proof}
    Since $k$ and $F$ have characteristic~$0$, the result follows from Lemmas \ref{lem: exactness algebraic closure} and \ref{lem: Weil restriction F-bar}.
\end{proof}

\begin{lemma}
    \label{lem: weil restriction product}
    Let $F/k$ be a finite separable field extension and let $G$ and $H$ be algebraic groups over $F$. Then there is a canonical isomorphism $\sR_{F/k}(G \times H) \cong \sR_{F/k} G \times \sR_{F/k}H$.
\end{lemma}
\begin{proof}
    Weil restriction is right adjoint to the functor of extending scalars from $k$ to $F$ \cite[Proposition 20.7]{BookOfInvolutions}, so it preserves limits.
\end{proof}

\subsection{From Lie algebras to algebraic groups}
A lot of the theory of algebraic groups is already captured by their Lie algebras, especially over fields of characteristic~$0$. We will recall here a few results that we will need about this.

\begin{lemma}[{\cite[Proposition II.2.13]{milneALA}}]
    \label{lem: morphism on algebraic groups determined by Lie alg}
    Let $G$ and $H$ be algebraic groups over a field of characteristic~$0$, and let $f, g \colon G \to H$ be two homomorphisms. If $G$ is connected and the induced maps $\Lie(f), \Lie(g) \colon \Lie(G) \to \Lie(H)$ coincide, then $f = g$.\qed
\end{lemma}

Lie's third theorem says that every complex Lie algebra is the Lie algebra of a complex Lie group. There is an analog for semisimple Lie algebras and algebraic groups over fields other than the complex numbers. Before we state this result, we need a few definitions. Recall that an algebraic group $G$ over a field $k$ is called a group of multiplicative type if $G_{k^{\sep}}$ is isomorphic to a product of groups of the form $\mathbb{G}_{m, k^{\sep}}$ and $\mu_{n, k^{\sep}}$.

\begin{definition}
    Let $G$ and $G'$ be connected algebraic groups over a field $k$. A homomorphism $f \colon G' \to G$ is called a \textbf{multiplicative isogeny} if $f$ is surjective and $\ker(f)$ is a finite algebraic group of multiplicative type.
\end{definition}
\begin{definition}
    Let $G$ be a connected semisimple algebraic group over a field $k$. Then $G$ is called \textbf{simply connected} if for every connected semisimple algebraic group $G'$ and multiplicative isogeny $f \colon G' \to G$, the map $f$ is an isomorphism.
\end{definition}

The following is a useful way to prove that an algebraic group is simply connected:
\begin{lemma}
    \label{lem: simply connected iff}
    Let $G$ be a semisimple geometrically connected algebraic group over a subfield $k \subseteq \mathbb{C}$. If $G(\mathbb{C})$ is simply connected in the analytic topology, then $G$ is simply connected as an algebraic group over $k$.
\end{lemma}
\begin{proof}
    Suppose $G$ is not simply connected as an algebraic group. Then there is a nontrivial multiplicative isogeny $G' \to G$. This gives a nontrivial covering $G'(\mathbb{C}) \to G(\mathbb{C})$, so $G(\mathbb{C})$ is also not simply connected, which gives a contradiction.
\end{proof}

\begin{definition}
Let $k$ be a field.
\begin{itemize}
    \item Let $G$ be an algebraic group over $k$. We denote its category of finite-dimensional $k$-representations by $\Rep(G)$.
    \item Let $\mathfrak{g}$ be a Lie algebra over $k$. We denote its category of finite-dimensional $k$-representations by $\Rep(\mathfrak{g})$. 
\end{itemize}  
\end{definition}

With these definitions, we can state the analogue of Lie's third theorem for algebraic groups. The result is mostly proved in a book by Milne \cite{milneALA}.
\begin{theorem}
    \label{thm: integrate Lie algebra action}
    Let $G$ be a semisimple simply connected algebraic group over a field $k$ of characteristic~$0$, and let $\mathfrak{g} = \Lie(G)$. Then taking Lie algebras gives an equivalence of categories
    $$\Rep(G) \isomto \Rep(\mathfrak{g}).$$
\end{theorem}
\begin{proof}
    By \cite[Theorem III.3.12]{milneALA} and \cite[Lemma III.3.8]{milneALA}, there exists a semisimple algebraic group $G'$ with $\Lie(G') = \mathfrak{g}$ and an equivalence $\Rep(G') \cong \Rep(\mathfrak{g})$. Moreover, part (c) of \cite[Theorem III.3.12]{milneALA} gives an isogeny $f \colon G' \to G$. Since $G'$ is semisimple and $k$ has characteristic zero, the kernel of $f$ is of multiplicative type. Since $G$ is simply connected, $f$ must be an isomorphism, and this finishes the proof.
\end{proof}

\subsection{Cohomology of algebraic groups}
Let $G$ be an algebraic group over a field $k$, and let $\ell/k$ be a Galois extension. Then $\Gal(\ell/k)$ acts on $G(\ell)$. We denote the associated Galois cohomology by
$$\sH^i(\ell/k, G) := \sH^i(\Gal(\ell/k), G(\ell)).$$
If $\ell = k^{\sep}$ is the separable closure of $k$, we write $\sH^i(k, G) := \sH^i(k^{\sep}/k, G)$. If $G$ is not abelian, the Galois cohomology is only defined for $i \leq 1$.

We will also need the following result, which generalizes Hilbert 90:
\begin{lemma}
    \label{lem: Hilbert 90}
    Let $F/k$ be a finite separable field extension. Then $$\sH^1(k, \sR_{F/k}\mathbb{G}_{m, F}) = 0.$$
\end{lemma}
\begin{proof}
    Shapiro's Lemma \cite[Lemma~29.6]{BookOfInvolutions} says that
    \[\sH^1(k, \sR_{F/k}\mathbb{G}_{m, F}) \cong \sH^1(F, \mathbb{G}_{m, F}).\]
    The result then follows from Hilbert 90.
\end{proof}

\subsection{Representation theory of some groups of multiplicative type}
We recall some facts about the representation theory of groups of multiplicative type. Using this, we will describe a way of extending representations of algebraic groups with a central factor $\sR_{\ell/k}\mu_{n, \ell}$ for $\ell/k$ a finite seperable field extension, see Corollary \ref{cor: extending mu_n-rep to G_m-rep}.

\begin{definition}
    Let $G$ be an algebraic group over a field $k$. The \textbf{character group} of $G$ is the abelian group\index{$X(G)$}
    $$X(G) := \Hom(G_{k^{\sep}}, \mathbb{G}_{m, k^{\sep}}).$$
\end{definition}
The group $\Gal(k^{\sep} / k)$ acts on $G_{k^{\sep}}$ and $\mathbb{G}_{m, k^{\sep}}$, and this induces an action of $\Gal(k^{\sep} / k)$ on $X(G)$, explicitly given by $\gamma \cdot f = \gamma \circ f \circ \gamma^{-1}$ for $f \in X(G)$ and $\gamma \in \Gal(k^{\sep}/k)$. We denote the set of orbits of the action of 
$\Gal(k^{\sep}/k)$ on $X(G)$ by $X(G)/\Gal(k^{\sep}/k)$.

We will need a few facts about groups of multiplicative type. Their behaviour is governed by their character groups:
\begin{theorem}
    Let $k$ be a field. Then taking the character group induces an equivalence of categories between the category of groups of multiplicative type over $k$ and the category of finitely-generated abelian groups with a continuous action of $\Gal(k^{\sep}/k)$.
\end{theorem}
\begin{proof}
    See \cite[Theorem 12.23]{milneiAG}.
\end{proof}

\begin{lemma}
    \label{lem: Q-reps of multiplicative groups}
    Let $k$ be a field and let $G$ be a group of multiplicative type over~$k$. Let $\Irr_k(G)$ be a set of representatives for the isomorphism classes of finite-dimensional irreducible $G$-representations over $k$. Then the category of finite-dimensional $G$-representations is semisimple and there is a bijection $\Irr_k(G) \to X(G)/\Gal(k^{\sep}/k)$, given by sending $V \in \Irr_k(G)$ to the collection of characters occurring in $V_{k^{\sep}}$.
\end{lemma}
\begin{proof}
    See \cite[Theorem 12.30]{milneiAG}.
\end{proof}

If $\ell/k$ is a separable field extension and $n$ a positive integer, then $\sR_{\ell/k}\mathbb{G}_{m,\ell}$ and $\sR_{\ell/k}\mu_{n, \ell}$ are groups of multiplicative type.  

\begin{lemma}
    \label{lem: characters mu_n to G_m surjective}
    Let $n \geq 1$ be an integer and let $\ell/k$ be a finite separable field extension. The inclusion $\sR_{\ell/k}\mu_{n, \ell} \to \sR_{\ell/k}\mathbb{G}_{m, \ell}$ induces a $\Gal(k^{\sep}/k)$-equivariant surjection $X(\sR_{\ell/k}\mathbb{G}_{m, \ell}) \to X(\sR_{\ell/k}\mu_{n, \ell})$.
\end{lemma}
\begin{proof}
    Since $\sR_{\ell/k}\mu_{n, \ell}$ and $\sR_{\ell/k}\mathbb{G}_{m, \ell}$ are algebraic groups over $k$, the inclusion $(\sR_{\ell/k}\mu_{n, \ell})_{k^{\sep}} \to (\sR_{\ell/k}\mathbb{G}_{m, F})_{k^{\sep}}$ is Galois equivariant, so the functoriality of $\Hom$ implies that the map $X(\sR_{\ell/k}\mathbb{G}_{m, \ell}) \to X(\sR_{\ell/k}\mu_{n, \ell})$ is equivariant with respect to the action of $\Gal(k^{\sep}/k)$.

    To see that this map is surjective, observe that we have a commutative diagram
    \begin{equation}
        \label{eq: Weil restriction Qbar diagram}
        \begin{tikzcd}
            (\sR_{\ell/k}\mu_{n, \ell})_{k^{\sep}} \arrow[r] \arrow[d, "\sim" {anchor=north, rotate=90}] & (\sR_{\ell/k}\mathbb{G}_{m, \ell})_{k^{\sep}} \arrow[d, "\sim" {anchor=north, rotate=90}]\\
            \mu_{n, k^{\sep}}^d \arrow[r, "\iota^d"] & \mathbb{G}_{m, k^{\sep}}^d,
        \end{tikzcd}
    \end{equation}
    where $\iota \colon \mu_{n, k^{\sep}} \to \mathbb{G}_{m, k^{\sep}}$ is the standard inclusion and $d = [\ell : k]$. There are identifications $X(\mu_{n, k^{\sep}}) \cong \mathbb{Z}/n\mathbb{Z}$ and $X(\mathbb{G}_{m, k^{\sep}}) \cong \mathbb{Z}$, where $\iota$ induces the standard projection $\mathbb{Z} \to \mathbb{Z}/n\mathbb{Z}$. Under these identifications we see that the map $\sR_{\ell/k} \mu_{n, \ell} \to \sR_{\ell/k}\mathbb{G}_{m, \ell}$ induces the surjection $\mathbb{Z}^d \to (\mathbb{Z}/n\mathbb{Z})^d$ on character groups.
\end{proof}

\begin{lemma}
    \label{lem: Galois action on X of Weil restriction}
    Let $\ell/k$ be a finite Galois extension, and let $G$ be either $\mathbb{G}_{m, \ell}$ or $\mu_{n, \ell}$. Let $S$ be the set of $k$-linear field embeddings $\ell \to \overline{k}$. Then $X(\sR_{\ell/k} G) \cong \bigoplus_{\sigma \in S} X(G)$, and the Galois group $\Gal(\overline{k}/k)$ acts on $X(\sR_{\ell/k} G)$ by permuting the factors; explicitly a $\gamma \in \Gal(\overline{k}/k)$ sends $(x_\sigma)_\sigma \in \bigoplus_{\sigma \in S} X(G)$ to $(x_{\gamma^{-1}\sigma})_\sigma$.
\end{lemma}
There is a bijection $S \cong \Gal(\ell/k)$, and the action of $\Gal(\overline{k}/k)$ on $S$ is given by taking the surjection $\Gal(\overline{k}/k) \to \Gal(\ell/k)$ and composing this with the action of $\Gal(\ell/k)$ on $S$ given by composition.

If we take $G = \mathbb{G}_{m, \ell}$ in the theorem, the induced representation of $\Gal(\ell/k)$ on $\mathbb{Z}^S$ is just the regular representation.

\begin{proof}
    Write $G = \Spec(A)$, where $A$ is either the $\ell$-algebra $\ell[t, t^{-1}]$ or $\ell[t]/(t^n)$. The proof of \cite[Lemma 20.6]{BookOfInvolutions} provides the following description of the Weil restriction: for a $k$-linear inclusion $\sigma \colon \ell \to \overline{k}$, let $A_\sigma = A \otimes_{\ell,\sigma} \overline{k}$, then \[\sR_{\ell/k}G \times_{\Spec k} \Spec \overline{k} \cong \Spec\Bigl(\bigotimes_{\sigma \in S} A_\sigma\Bigl).\]
    The Galois group $\Gal(\overline{k}/k)$ acts on $\bigotimes_{\sigma \in S} A_\sigma$ by letting $\gamma \in \Gal(\overline{k}/k)$ act as the tensor product of the maps $A_\sigma \to A_{\gamma \sigma}$ given by sending $a \otimes x$ to $a \otimes \gamma(x)$ for $a \in A$ and $x \in \overline{k}$.

    We have
    $$X(\sR_{\ell/k} G) = \Hom((\sR_{\ell/k} G)_{\overline{k}}, \mathbb{G}_{m, \overline{k}}) = \prod_{\sigma \in S} \Hom_{\text{Hopf}, \overline{k}}(\overline{k}[t, t^{-1}],  A \otimes_{\ell, \sigma} \overline{k}),$$
    where $\Hom_{\text{Hopf}, \overline{k}}$ denotes the collection of $\overline{k}$-linear Hopf algebra morphisms. Then a $\gamma \in \Gal(\overline{k}/k)$ acts on this by sending a map $f \colon \overline{k}[t] \to A_{\sigma}$ to the map $f^\gamma \colon \overline{k}[t] \to A_{\gamma \sigma}$ defined by $f^\gamma(p) = \gamma(f(\gamma^{-1}(p)))$ for $p \in \overline{k}[t, t^{-1}]$.

    We have isomorphisms 
    \begin{equation*}
        \mathbb{Z} \cong \Hom_{\text{Hopf}, \overline{k}}(\overline{k}[t, t^{-1}], \overline{k}[t, t^{-1}]),
    \end{equation*}
    given by $m \mapsto (t \mapsto t^m)$ and
    \begin{equation*}
        \mathbb{Z}/n\mathbb{Z} \cong \Hom_{\text{Hopf}, \overline{k}}(\overline{k}[t, t^{-1}], \overline{k}[t]/(t^n)),
    \end{equation*}
    given by $m \mapsto (t \mapsto t^m)$. All $\gamma \in \Gal(\overline{k}/k)$ send $t$ to itself, so in both cases, we see that $f^\gamma = f$ for $f \in \Hom_{\text{Hopf}, \overline{k}}(\overline{k}[t, t^{-1}], A_\sigma)$. It follows that the Galois group acts on $X(\sR_{\ell/k} G) = \bigoplus_\sigma X(G)$ by permuting the factors.
\end{proof}

Let $W$ be a vector space of dimension $n$ over a field $k$. There is an inclusion $\uSL(W) \to \uGL(W)$ of algebraic groups over $k$. Given a representation $V$ of $\uSL(W)$, one may wonder whether it can be extended to a representation of $\uGL(W)$. One way to achieve this is to use that $\uGL(W)$ is isomorphic to the quotient $(\mathbb{G}_{m,k} \times \uSL(W))/\mu_{n, k}$, so we have to find a representation of $\mathbb{G}_{m,k}$ on $V$ that commutes with the $\uSL(W)$-representation such that the subgroup $\mu_{n,k} \subseteq \mathbb{G}_{m,k} \times \uSL(V)$ acts trivially on $W$. Similar extension problems can be phrased in the same way, for example for the inclusions $\uSp(V) \subseteq \uGSp(V)$ for $V$ a symplectic space or $\uSpin(V) \subseteq \uGSpin(V)$ for $V$ a quadratic space. We will generalize this in Corollary~\ref{cor: extending mu_n-rep to G_m-rep}.

For an algebraic group $G$ over a field $k$ and a representation $\rho \colon G \to \uGL(V)$ on a $k$-vector space $V$, we write \myindex{$C_{\uGL(V)}(G)$} for the centralizer of $\rho(G)$ in $\uGL(V)$, this is an algebraic subgroup of $\uGL(V)$.

\begin{lemma}
    \label{lem: centralizer criterion}
    Let $G$ be an algebraic group of multiplicative type over a field $k$, and let $\rho \colon G \to \uGL(V)$ be a finite-dimensional $k$-representation. Decompose $V_{k^{\sep}}$ into isotypic components $V_{k^{\sep}} = \bigoplus_{\alpha \in X(G)} V_\alpha$. Then
    $$C_{\uGL(V)}(G)(k^{\sep}) = \{g \in \uGL(V)(k^{\sep}) : \forall \alpha \in X(G),\text{ } gV_\alpha \subseteq V_\alpha\}.$$
\end{lemma}
\begin{proof}
    Take a $g \in C_{\uGL(V)}(G)(k^{\sep})$, then $g \rho(\lambda) = \rho(\lambda)g$ for all $\lambda \in G(k^{\sep})$. Then for $\alpha \in X(G)$ and $v \in V_\alpha$, we see that
    $$\rho(\lambda)gv = g\rho(\lambda)v = g\alpha(\lambda) v = \alpha(\lambda) gv,$$
    so we see that $gv \in V_\alpha$, which implies that $gV_\alpha \subseteq V_\alpha$.
    
    For the other inclusion, suppose that we have a $g \in \uGL(V)(k^{\sep})$ with the property that $gV_\alpha \subseteq V_\alpha$ for all $\alpha \in X(G)$. We need to show that $\rho(\lambda)g = g\rho(\lambda)$ for all $\lambda \in G(k^{\sep})$. To show this, it suffices to show that $\rho(\lambda)gv = g\rho(\lambda)v$ for all $v \in V_{k^{\sep}}$. Since $V_{k^{\sep}} = \bigoplus_\alpha V_\alpha$, it suffices to show the equality $\rho(\lambda)gv = g\rho(\lambda)v$ for all $\alpha \in X(G)$ and $v \in V_\alpha$. Since $gV_\alpha \subseteq V_\alpha$ by assumption, we have for $v \in V_\alpha$: 
    $$\rho(\lambda) gv = \alpha(\lambda)gv = g\alpha(\lambda)v = g\rho(\lambda)v.$$
    Hence $g \in C_{\uGL(V)}(G)(k^{\sep})$.
\end{proof}

\begin{proposition}
    \label{prop: Rmu to RG extension}
    Let $\ell/k$ be a finite Galois extension of fields, and let $n \geq 1$ be an integer. Let $\rho \colon \sR_{\ell/k}\mu_{n, \ell} \to \uGL(V)$ be a finite-dimensional $k$-representation. Then $\rho$ extends to a representation $\widetilde \rho \colon \sR_{\ell/k} \mathbb{G}_{m, \ell} \to \uGL(V)$, with the property that
    $$C_{\uGL(V)}(\sR_{\ell/k}\mu_{n, \ell})(k^{\sep}) = C_{\uGL(V)}(\sR_{\ell/k}\mathbb{G}_{m, \ell})(k^{\sep}).$$
\end{proposition}

\begin{proof}
    To ease notation, we will abbreviate $\sR_{\ell/k}\mu_{n, \ell}$ to $\sR\mu_n$ and $\sR_{\ell/k}\mathbb{G}_{m, \ell}$ to $\sR\mathbb{G}_m$. By Lemma \ref{lem: Q-reps of multiplicative groups}, we can decompose $V_{k^{\sep}} = \bigoplus_{T} V_T^{n_T}$ into irreducible representations of $\sR\mu_{n}$, where the direct sum runs over $T$ in $X(\sR\mu_{n})/\Gal(k^{\sep}/k)$.
    
    By Lemma \ref{lem: characters mu_n to G_m surjective}, there is a surjective map $X(\sR \mathbb{G}_{m}) \to X(\sR\mu_{n})$. By looking at the explicit description of this Galois action given in Lemma \ref{lem: Galois action on X of Weil restriction}, we see that for every Galois orbit $T \subseteq X(\sR\mu_{n})$, there exists a Galois orbit $T' \subseteq X(\sR\mathbb{G}_{m})$ mapping bijectively onto $T$. For every $T$, we choose one such lift $T'$ of $T$ and fix it throughout the rest of the proof.

    So, for every $T$, we see that $V_T$ is isomorphic to the restriction of some $V_{T'}$, and this gives us a representation of $\sR\mathbb{G}_{m}$ on $V$ extending the $\sR\mu_{n}$-representation. For every $\alpha \in X(\sR \mu_n)$, let $\alpha' \in X(\sR \mathbb{G}_m)$ be its lift. Then $V_\alpha = V_{\alpha'}$ as subspaces of $V_{k^{\sep}}$ for all $\alpha \in X(\sR \mu_n)$. Then Lemma \ref{lem: centralizer criterion} gives that the centralizers of $\sR\mu_n$ and $\sR \mathbb{G}_m$ in $\uGL(V)$ are equal.    
\end{proof}

\begin{corollary}
    \label{cor: extending mu_n-rep to G_m-rep}
    Let $k$ be a field of characteristic~$0$ and let $\Gamma$ be an algebraic group over $k$. Suppose we have an algebraic subgroup $\sR_{\ell/k}\mu_{n, \ell} \subseteq Z(\Gamma)$ for some finite Galois extension $\ell/k$ and integer $n \geq 1$. Let $\rho \colon \Gamma \to \uGL(V)$ be a finite-dimensional $k$-representation of $\Gamma$. Then there exists a representation of $\sR_{\ell/k}\mathbb{G}_{m, \ell} \times \Gamma$ on $V$ such that the algebraic subgroup $\sR_{\ell/k}\mu_{n, \ell} \to \sR_{\ell/k}\mathbb{G}_{m, \ell} \times \Gamma$ embedded via $\lambda \mapsto (\lambda^{-1}, \lambda)$ acts trivially on $V$, and with the property that
    \[C_{\uGL(V)}(\sR_{\ell/k}\mathbb{G}_{m, \ell} \times \Gamma) = C_{\uGL(V)}(\Gamma).\]
\end{corollary}

The corollary implies that in the diagram below, with $\iota$ induced by the inclusion of the second factor, one can always find a dashed arrow making the diagram commute:
\begin{equation*}
    \begin{tikzcd}
        \Gamma \arrow[r, "\rho"]\arrow[d, "\iota"] & \uGL(V).\\
        (\sR_{\ell/k} \mathbb{G}_{m, \ell} \times \Gamma)/\sR_{\ell/k}\mu_{n,\ell} \arrow[ur, dashrightarrow]
    \end{tikzcd}
\end{equation*}

\begin{proof}[{Proof of Corollary \ref{cor: extending mu_n-rep to G_m-rep}}]
    As in the proof of Proposition \ref{prop: Rmu to RG extension}, we will abbreviate the Weil restrictions of $\mu_{n, \ell}$ and $\mathbb{G}_{m, \ell}$ along $\ell/k$ to $\sR\mu_n$ and $\sR\mathbb{G}_m$. By Proposition \ref{prop: Rmu to RG extension}, the restriction of $\rho$ to $\sR\mu_n$ lifts to a representation $\sigma \colon \sR\mathbb{G}_m \to \uGL(V)$, such that the images of $\sR\mu_n$ and $\sR\mathbb{G}_m$ have the same centralizer in $\uGL(V)(k^{\sep})$. Since $k$ has characteristic~$0$, using \cite[Corollary 3.26]{milneiAG}, we see that $C_{\uGL(V)}(\sR \mu_n) = C_{\uGL(V)}(\sR \mathbb{G}_m)$ as algebraic subgroups of $\uGL(V)$.

    Since $\sR\mu_n \subseteq Z(\Gamma)$, we have an inclusion $\Gamma \subseteq C_{\uGL(V)}(\sR\mu_n)$, and hence also $\Gamma \subseteq C_{\uGL(V)}(\sR\mathbb{G}_m)$. We now define a map $\widetilde \rho \colon \sR \mathbb{G}_m \times \Gamma \to \uGL(V)$, which for every commutative $k$-algebra $R$ sends a pair $(\lambda, g) \in (\sR \mathbb{G}_m)(R) \times \Gamma(R)$ to
    $$\widetilde \rho(\lambda, g) = \sigma(\lambda)\rho(g).$$
    The fact that $\Gamma \subseteq C_{\uGL(V)}(\sR\mathbb{G}_m)$ implies that $\widetilde \rho$ is a group homomorphism. By construction, the subgroup $\sR \mu_n \to \sR\mathbb{G}_m \times \Gamma$ embedded via $\lambda \mapsto (\lambda^{-1}, \lambda)$ is in the kernel of $\widetilde \rho$.

    It remains to show that the centralizers of $\sR\mathbb{G}_m \times \Gamma$ and $\Gamma$ in $\uGL(V)$ coincide. By using that $\Gamma \subseteq C_{\uGL(V)}(\sR\mathbb{G}_m)$, we see:
    \begin{equation*}
        C_{\uGL(V)}(\sR\mathbb{G}_m \times \Gamma) = C_{\uGL(V)}(\sR\mathbb{G}_m) \cap C_{\uGL(V)}(\Gamma) = C_{\uGL(V)}(\Gamma).
    \end{equation*}
\end{proof}

\section{Equivariant Hodge structures}
In this section, we look at Hodge structures equipped with an action of a finite group. I could not find a reference for the result of Proposition~\ref{prop: category of polarizable G-equiv HS is semisimple}, so a proof is given here.

\begin{definition}
    Let $G$ be a finite group. A \textbf{polarizable $G$-equivariant Hodge structure} is a polarizable rational Hodge stucture $V$ which is also a $G$-representation, such that for every $g \in G$, the map $g \colon V \to V$ is a morphism of Hodge structures. 

    If $V$ and $W$ are two $G$-equivariant Hodge structures, we say that a linear map $f \colon V \to W$ is a \textbf{morphism of $G$-equivariant Hodge structures} if $f$ is both a morphism of Hodge structures and a morphism of $G$-representations.
    
    These objects and morphism form the category of \textbf{polarizable $G$-equivariant Hodge structures}.
\end{definition}

\begin{proposition}
    \label{prop: category of polarizable G-equiv HS is semisimple}
    Let $G$ be a finite group. Then:
    \begin{enumerate}
        \item The category of polarizable $G$-equivariant Hodge structures is semisimple.
        \item  Let $V$ be a polarizable $G$-equivariant Hodge structure of weight $k$. Then there is a $G$-equivariant isomorphism of Hodge structures $V \cong V^*(-k)$.
    \end{enumerate}
\end{proposition}
\begin{proof}
    Let $V$ be a polarizable $G$-equivariant Hodge structure, and for every $k \in \mathbb{Z}$, let $V_k \subseteq V$ be the subspace of elements of weight $k$. Then $V \cong \bigoplus_{k \in \mathbb{Z}} V_k$, and also $gV_k \subseteq V_k$ for every $g \in G$ and $k \in \mathbb{Z}$. Hence it suffices to show that every $k$ the category of polarizable $G$-equivariant Hodge structures of weight $k$ is semisimple.
    
    So, let $V$ be a polarizable $G$-equivariant Hodge structure of weight $k$ and let $Q$ be a polarization on $V$. Using Maschke's theorem, we obtain a $G$-invariant bilinear form $Q'$ defined by
    $$Q'(v, w) = \frac{1}{|G|}\sum_{g \in G}Q(gv, gw).$$
    It is easy to see that $Q'$ is again a polarization and that it is also $G$-equivariant.

    If we now take a $G$-equivariant sub-Hodge structure $W \subseteq V$, then its orthogonal complement $W^\perp$ with respect to $Q'$ is also a $G$-equivariant sub-Hodge structure of $V$ and satisfies $V = W \oplus W^\perp$, which proves the first statement.

    For the second statement, observe that $Q'$ is a $G$-equivariant perfect pairing $V \otimes V \to \mathbb{Q}(-k)$, and hence induces a $G$-equivariant isomorphism $V \cong V^*(-k)$. 
\end{proof}

\section{The Beauville-Bogomolov decomposition}
\label{sec: reduction to quotients}
In this section we discuss the Beauville-Bogomolov decomposition, which can be used to reduce the proof of Theorem \ref{thm: main theorem} to the case where $X$ and $Y$ have the specific form $X = (X_0 \times \prod_{i=1}^k X_i)/G$ and $Y = (Y_0 \times \prod_{j=1}^l \times Y_j)/H$, where the $X_i$ and $Y_j$ for $i,j \geq 1$ are hyperkähler varieties, $X_0$ and $Y_0$ are even-dimensional abelian varieties, and $G$ and $H$ are finite groups that act freely and preserve at least one holomorphic symplectic form. After that, we will mention a few constraints on the possible group actions on varieties of the form $X_0 \times \dots \times X_k$ where $X_0$ is an abelian variety and the $X_i$ for $i \geq 1$ are hyperkähler varieties. All results in this section are due to Beauville \cite{BeauvilleKatata}.

\begin{theorem}
    \label{thm: symplectic variety is quotient}
    Let $X$ be a smooth connected projective variety over $\mathbb{C}$ admitting a holomorphic symplectic form. Then there exist hyperkähler varieties $X_1, \dots, X_k$, an even-dimensional abelian variety $X_0$ and a finite group $G$ acting freely on $X_0 \times \dots \times X_k$, such that $X \cong (\prod_{i=0}^k X_i)/G$.
\end{theorem}
\begin{proof}
    Since $X$ is symplectic, it has trivial canonical bundle. Hence it admits an étale covering $A' \times \prod X'_i \times \prod Y_j \to X$ where $A'$ is an abelian variety, the $X_i'$ are hyperkähler varieties and the $Y_j$ are strict Calabi-Yau varieties \cite{BeauvilleKatata}. Since the pullback of a symplectic form along an étale morphism is a symplectic form, the étale covering is also symplectic. Hence we do not have any strict Calabi-Yau factors and $A'$ has even dimension. The étale covering $A' \times \prod X_i' \to X$ can be refined to a Galois covering $X_0 \times \prod_{i=1}^k X_i\to X$ with Galois group $G$ (see the proof of \cite[Proposition~3]{BeauvilleKatata}), so that $X \cong (X_0 \times \prod_{i=1}^k X_i)/G$. The action of $G$ on $X$ is free because the map $X \to X/G$ is étale.
\end{proof}

The Beauville-Bogomolov decomposition $X_0 \times X_1 \times \dots \times X_k \to X$ is not unique. However, it turns out that there is not too much wiggle room. The following two results follow from results by Beauville \cite{BeauvilleKatata}:

\begin{proposition}[{\cite[Proposition 3]{BeauvilleKatata}}]
    \label{prop: Beauville minimal split covering}
    Let $X$ be a smooth connected projective variety over $\mathbb{C}$, admitting a holomorphic symplectic form. Then there exists a Galois covering $\pi \colon \prod_{i=0}^k X_i \to X$, with the property that any other étale covering $\prod_{i=0}^k X'_i \to X$ factors via $\pi$. Moreover, the action of $G$ on $X_0 \times X_1 \times \dots \times X_k$ does not have any automorphism of the form 
    \[(x_0, x_1, \dots, x_k) \mapsto (x_0 + a, x_1, \dots, x_k)\]
    for any $a \in X_0$. \qed
\end{proposition}

\begin{proposition}[{\cite[Theorem]{BeauvilleKatata}}]
    Let $X$ be a smooth connected projective variety over $\mathbb{C}$ admitting a holomorphic symplectic form. Let $X \cong (\prod_{i=0}^k X_i)/G$ and $X \cong (\prod_{i=0}^{k'} X'_i)/G'$ be two ways of writing $X$ as in Theorem \ref{thm: symplectic variety is quotient}. Then $k = k'$, and there is a bijection $\sigma \colon \{1, \dots, k\} \cong \{1, \dots, k'\}$ such that $X_i \cong X'_{\sigma(i)}$ for all $1 \leq i \leq k$. \qed
\end{proposition}

The following theorem by Beauville describes the possible actions of $G$ on the product $X_0 \times \prod_{i=1}^k X_i$. Note that we do not require automorphisms of abelian varieties to respect the group structure.

\begin{theorem}[{\cite[Section 3]{BeauvilleKatata}}]
    \label{thm: HKA automorphisms}
    Let $X_1, \dots, X_k$ be hyperkähler varieties and let $X_0$ be an abelian variety. Let $X = X_0 \times X_1 \times \dots \times X_k$ and let $g$ be an automorphism of $X$. Then there exists a permutation $\sigma \in S_k$, isomorphisms $g_i \colon X_{\sigma^{-1}(i)} \to X_i$ and an automorphism $g_0$ of $X_0$ such that 
    $$g(x_0, x_1, \dots, x_k) = (g_0(x_0), g_1(x_{\sigma^{-1}(1)}), \dots, g_k(x_{\sigma^{-1}(k)}))$$ 
    for all elements $(x_0, x_1, \dots, x_k) \in X_0 \times X_1 \times \dots \times X_k$. \qed
\end{theorem}

\begin{remark}
    There is a lot more that can be said about free symplectic group actions on varieties of the form $X_0 \times X_1 \times \dots \times X_k$ with $X_0$ an even-dimensional abelian variety and the $X_i$ for $i \geq 1$ hyperkähler varieties. See Section~\ref{sec: geometry of symplectic varieties} for a more detailed study of such group actions.
\end{remark}

We end this chapter by introducing some notation that will be used throughout the thesis. Consider a group $G$ acting on $X_0 \times X_1 \times \dots \times X_k$, where $X_0$ is an abelian variety and the $X_i$ for $i \geq 1$ are hyperkähler varieties. For $g \in G$, let $\sigma_g$ be the permutation of $\{1, \dots, k\}$ obtained from the above theorem. We then get a group homomorphism $\sigma \colon G \to S_k$ by mapping $g$ to $\sigma_g$. Let \index{$I$}$I = \{1, \dots, k\}$, then $G$ acts on $I$ via the homomorphism $G \to S_k$. For $i \in I$, we write $[i] \subseteq I$ for the orbit of $i$ under the action of $G$, and define \index{$G_i$}$G_i = \{g \in G : \sigma_g(i) = i\}$. In particular, note that $G_i$ acts on the variety $X_i$.

%% file: chapters/chapter4.tex
In this chapter, we will recall the definition of the LLV algebra of a projective variety, due to Looijenga, Lunts and Verbitsky \cite{LooijengaLunts, VerbitskyLLV}. We then extensively study the theory of LLV algebras. Some of the results presented here are due to Looijenga and Lunts or Verbitsky, but this chapter also contains many new results. This study of LLV algebras culminates in the computation of the LLV algebra of a smooth projective variety admitting a holomorphic symplectic form. Throughout this chapter, $K$ will be a field of characteristic~$0$. All vector spaces and algebras over $K$ will be assumed to be finite-dimensional.

\section{LLV algebras of (Jordan-)Lefschetz modules}
\label{sec: LLV properties}
In this section we recall the definition of Lefschetz modules $(\mathfrak{a}, M)$ and their LLV algebras, as introduced by Looijenga, Lunts and Verbitsky \cite{LooijengaLunts, VerbitskyLLV}. We then study some of their properties.

\begin{definition}
    Let $M$ be a graded vector space of finite dimension over a field $K$. We say that a linear endomorphism $e \colon M \to M$ of degree $2$ has the \textbf{hard Lefschetz property} if for every $k \in \mathbb{N}$ the $k$-fold composite $e^k \colon M_{-k} \to M_k$ is an isomorphism.
\end{definition}

Let $h \colon M \to M$ be the linear endomorphism which acts as multiplication by $k$ on $M_k$. Given a degree $2$ operator $e \colon M \to M$ with the hard Lefschetz property, the Jacobson-Morozov theorem \cite[VIII.11.2, Proposition~2]{BourbakiLie7and8} says that there is a unique endomorphism $f \colon M \to M$ of degree $-2$ giving an isomorphism $\langle e, h, f \rangle \cong \mathfrak{sl}_2(K)$, where $e, h$ and $f$ map to the standard generators $\begin{psmallmatrix} 0 & 1\\ 0 & 0\ \end{psmallmatrix}, \begin{psmallmatrix} 1 & 0\\0 & -1\end{psmallmatrix}$ and $\begin{psmallmatrix} 0 & 0\\ 1 & 0 \end{psmallmatrix}$ of $\mathfrak{sl}_2(K)$.

In general, we can consider a $K$-vector space $\mathfrak{a}$ and a map $e \colon \mathfrak{a} \to \End(M)_2$ whose image consists of mutually commuting linear endomorphisms of degree $2$. We will assume that at least one element of the image of $e$ has the hard Lefschetz property. Since being invertible is a Zariski-open condition, a Zariski-open collection of elements of $\im(e)$ has the hard Lefschetz property. For $a \in \mathfrak{a}$, we will denote the associated endomorphism of $M$ by $e_a$, and denote the `dual' operator by $f_a$ (if it exists).

\begin{definition}
    Let $\mathfrak{a}$ be a $K$-vector space, $M$ a graded $K$-vector space and $e \colon \mathfrak{a} \to \End(M)_2$ a linear map whose image consists of mutually commuting endomorphisms of $M$. The \textbf{LLV algebra} of $(\mathfrak{a}, M)$, denoted by \myindex{$\llv(\mathfrak{a}, M)$}, is the Lie subalgebra of $\mathfrak{gl}(M)$ generated by all Lie subalgebras $\mathfrak{sl}_2(K) \cong \langle e_a, h, f_a \rangle \subseteq \mathfrak{gl}(M)$ arising from $a \in \mathfrak{a}$ for which $e_a$ has the hard Lefschetz property.
\end{definition}

\begin{definition}
    Let $\mathfrak{a}$ be a $K$-vector space, $M$ a graded $K$-vector space and $e \colon \mathfrak{a} \to \End(M)_2$ a linear map whose image consists of mutually commuting endomorphisms of $M$. We say that $(\mathfrak{a}, M)$ is a \textbf{Lefschetz module} if the Lie algebra $\llv(\mathfrak{a}, M)$ is semisimple.
\end{definition}

For a graded vector space $M$ and $n \in \mathbb{Z}$, we will write \myindex{$M[n]$} for $M$ with the grading shifted by $n$, so $M[n]_k = M_{n+k}$. The motivating example of LLV algebras comes from the cohomology of compact Kähler manifolds. 
\begin{example}
    \label{ex: LLV of var}
    Let $X$ be a compact connected Kähler manifold of complex dimension $n$, let $M = \sH^\bullet(X; \mathbb{Q})[n]$ and let $\mathfrak{a} = \sH^2(X; \mathbb{Q})$. Then $\mathfrak{a}$ acts on $M$ by the cup product, and this indeed gives a collection of mutually commuting degree $2$ endomorphisms of $M$. We denote $\llv(\mathfrak{a}, M)$ by $\llv(X; \mathbb{Q})$. This Lie algebra is nonzero, since any Kähler class has the hard Lefschetz property. Moreover, $\llv(X; \mathbb{Q})$ is semisimple by \cite[(1.9)]{LooijengaLunts}, so $(\mathfrak{a}, M)$ is a Lefschetz module.

    We will also be interested in \myindex{$\llv(X; K)$} for other fields $K$ of characteristic~$0$. Since $\mathbb{Q} \subseteq K$, there is still a Kähler class in $\sH^2(X; K)$, so $\llv(X; K)$ is still a nonzero semisimple Lie algebra over $K$. There is a natural map $\llv(X; \mathbb{Q}) \otimes_\mathbb{Q} K \to \llv(X; K)$, and it is not hard to see that this map is an isomorphism. Therefore it suffices to know $\llv(X; \mathbb{Q})$.
\end{example}

If $M$ is a graded $K$-vector space, then $\End_K(M)$ obtains a grading, the degree $n$ part is equal to:
\[ \End_K(M)_n = \bigoplus_i \Hom_K(M_i, M_{i + n}).\]
Then $\llv(\mathfrak{a}, M)$ inherits a grading via the inclusion $\llv(\mathfrak{a}, M) \subseteq \End_K(M)$.

The LLV algebras of the varieties that interest us in this thesis behave quite nicely, they are all so-called Jordan-Lefschetz algebras. They have been studied quite extensively by Looijenga--Lunts. For the proof that the properties in the following definition are indeed equivalent, see \cite[Proposition 2.1]{LooijengaLunts}.
\begin{definition}
    Let $(\mathfrak{a}, M)$ be a Lefschetz module. We say that $(\mathfrak{a}, M)$ is a \textbf{Jordan-Lefschetz pair} if the following equivalent properties hold:
    \begin{enumerate}
        \item $\llv(\mathfrak{a}, M)$ only lives in degrees $-2, 0$ and $2$.
        \item For any $a, b \in \mathfrak{a}$ that both have the hard Lefschetz property, the operators $f_a$ and $f_b$ commute.
    \end{enumerate}
\end{definition}

Looijenga and Lunts have classified all possible Jordan-Lefschetz pairs over an algebraically closed field of characteristic zero \cite[Corollary 2.6]{LooijengaLunts}.

A first property of Jordan-Lefschetz pairs is that injectivity of morphisms from their LLV algebras can be checked by only looking in degree $2$.

\begin{lemma}
    \label{lem: check inj in deg 2}
    Let $(\mathfrak{a}, M)$ be a Jordan-Lefschetz pair and let $\mathfrak{g} = \llv(\mathfrak{a}, M)$. Suppose that $\mathfrak{h}$ is a graded Lie algebra and that $\varphi \colon \mathfrak{g} \to \mathfrak{h}$ is a degree-preserving homomorphism. If the restriction $\varphi_2 \colon \mathfrak{g}_2 \to \mathfrak{h}_2$ to degree $2$ parts is injective, then $\varphi$ is injective.
\end{lemma}
\begin{proof}
    Since $\mathfrak{g}$ is semisimple, we can write $\mathfrak{g} = \prod_i \mathfrak{g}^i$ as a product of simple Lie algebras. We start by observing that every $\mathfrak{g}^i$ is graded. Indeed, every $\mathfrak{g}^i$ is an $\mathfrak{g}$-representation and $h$ is a semisimple element of $\mathfrak{g}$, so $\mathfrak{g}^i$ is the sum of its eigenspaces for the action of $h$. Hence $\mathfrak{g}_k = \bigoplus_i \mathfrak{g}^i_k$ for $k = -2, 0, 2$.
    
    Write $\mathfrak{g}'$ for the product of those $\mathfrak{g}^i$ which have a nonzero component in degree $2$. By \cite[Proposition 2.1]{TaelmanDerivedEquivalences}, any factor $\mathfrak{g}^i$ that does not contain an element of degree $2$ also does not contain any elements of degree $-2$, so the only way in which $\varphi$ can fail to be injective is if there is a simple factor living purely in degree $0$. Hence, we are done if we show that $\mathfrak{g} = \mathfrak{g}'$. 

    By assumption, $\mathfrak{g}'$ contains all the $e_a$ with $a \in \mathfrak{a}$. Since it is an ideal of $\mathfrak{g}$, it also contains $h = [e_a, f_a]$. But then also all $f_a$ lie in $\mathfrak{g}'$ by \cite[Proposition~2.1]{TaelmanDerivedEquivalences}. Hence all generators of $\mathfrak{g}$ lie in $\mathfrak{g}'$, so $\mathfrak{g} = \mathfrak{g}'$.
\end{proof}

\section{Tensor products and the LLV algebra}
In this section we discuss the LLV algebra of a tensor product of Lefschetz modules.

Given two graded vector spaces $M_1$ and $M_2$ with vector spaces $\mathfrak{a}_1$ and $\mathfrak{a}_2$ acting on $M_1$ and $M_2$ respectively (with mutually commuting actions by degree $2$ endomorphisms), we obtain an action of $\mathfrak{a}_1 \oplus \mathfrak{a}_2$ on the tensor product $M_1 \otimes M_2$ given by
$$(a_1, a_2) \cdot (m_1 \otimes m_2) = a_1m_1 \otimes m_2 + m_1 \otimes a_2 m_2.$$
This is also an action by mutually commuting degree $2$ endomorphisms. 

We equip $\mathfrak{gl}(M_1) \times \mathfrak{gl}(M_2)$ with a grading by setting $$(\mathfrak{gl}(M_1) \times \mathfrak{gl}(M_2))_k = \mathfrak{gl}(M_1)_k \times \mathfrak{gl}(M_2)_k,$$
and then we obtain a degree-preserving inclusion $\mathfrak{gl}(M_1) \times \mathfrak{gl}(M_2) \to \mathfrak{gl}(M_1 \otimes M_2)$, where we let $(\alpha, \beta) \in \mathfrak{gl}(M_1) \times \mathfrak{gl}(M_2)$ act on $M_1 \otimes M_2$ by sending $m_1 \otimes m_2$ to $\alpha(m_1) \otimes m_2 + m_1 \otimes \beta(m_2)$. 
In the discussion above Lemma 1.2 in the paper by Looijenga and Lunts \cite{LooijengaLunts}, it is observed that the LLV algebra of such a product is just the product of the LLV algebras:
\begin{lemma}
    \label{lem: tensor prod LLV}
    In the above situation, if $M_1$ and $M_2$ are nonzero and the actions of $\mathfrak{a}_1$ and $\mathfrak{a}_2$ on $M_1$ and $M_2$ respectively have the hard Lefschetz property, then the action of $\mathfrak{a}_1 \oplus \mathfrak{a}_2$ on $M_1 \otimes M_2$ also has the hard Lefschetz property. Moreover, the natural degree-preserving inclusion $\mathfrak{gl}(M_1) \times \mathfrak{gl}(M_2) \to \mathfrak{gl}(M_1 \otimes M_2)$ induces an isomorphism $\llv(\mathfrak{a}_1, M_1) \times \llv(\mathfrak{a}_2, M_2) \cong \llv(\mathfrak{a}_1 \oplus \mathfrak{a}_2, M_1 \otimes M_2)$ of graded Lie algebras. \qed
\end{lemma}

It is not difficult to see that the tensor product of Jordan-Lefschetz pairs is again a Jordan-Lefschetz pair:

\begin{corollary}
    \label{cor: product of JL is JL}
    Let $(\mathfrak{a}_1, M_1)$ and $(\mathfrak{a}_2, M_2)$ be two Jordan-Lefschetz pairs. Then $(\mathfrak{a}_1 \oplus \mathfrak{a}_2, M_1 \otimes M_2)$ is also a Jordan-Lefschetz pair.
\end{corollary}
\begin{proof}
    Let $\llv_i = \llv(\mathfrak{a}_i, M_i)$ for $i = 1, 2$. Then $\llv_1$ and $\llv_2$ are both semisimple, so $\llv_1 \times \llv_2$ is also semisimple. Furthermore, if we denote by $\llv_{i,k}$ the degree $k$ part of $\llv_i$, then $(\llv_1 \times \llv_2)_k = \llv_{1, k} \times \llv_{2, k}$, so it immediately follows that the pair $(\mathfrak{a}_1 \oplus \mathfrak{a}_2, M_1 \otimes M_2)$ is Jordan-Lefschetz.
\end{proof}

\section{LLV algebras of graded commutative algebras}
We will now turn our attention to Lefschetz pairs coming from a graded commutative (unital, associative) finite-dimensional $K$-algebra $R$, where we can define the so-called Verbitsky component. We will see later that, under some assumptions, the LLV algebra can be computed from just the Verbitsky component (Proposition \ref{prop: gR isom to gSR}). Recall that a graded $K$-algebra $R$ is called graded commutative if, for all integers $n, m$ and elements $a \in R_n$ and $b \in R_m$ we have $ab = (-1)^{mn}ba$.

\begin{definition}
    Let $R$ be a graded $K$-algebra. If there is an integer $n$ such that $R$ is supported in degrees $0$ up to $2n$ and $R_{2n} \neq 0$, we say that $R$ has \textbf{depth} $n$.
\end{definition}

One could also wish to define the depth of a ring $R$ supported in degrees $0$ up to $2n$ to be $2n$, but we stick with this convention as has been set by Looijenga-Lunts.

\begin{definition}
    Let $R$ be a graded commutative $K$-algebra of depth $n$. We say that $R$ has the \textbf{hard Lefschetz property} if at least one $a \in R_2$ has the hard Lefschetz property for the action on $R[n]$.
\end{definition}

\begin{example}
    \label{ex: kahler manifold cohom HL}
    Let $X$ be a compact connected Kähler manifold of dimension $n$ over $\mathbb{C}$. Then $\sH^\bullet(X; K)$ is a graded commutative $K$-algebra of depth $n$. Moreover, the hard Lefschetz theorem says that the Kähler class in $\sH^2(X; \mathbb{Q}) \subseteq \sH^2(X; K)$ has the hard Lefschetz property. 
\end{example}

\begin{definition}
    Let $R$ be a graded commutative $K$-algebra of depth $n$ with the hard Lefschetz property. Let $\mathfrak{a} = R_2$ and $M = R[n]$, where $a \in \mathfrak{a}$ acts on $R[n]$ by multiplication. We denote the resulting Lie algebra $\llv(\mathfrak{a}, M)$ by \myindex{$\llv(R)$}, and call this the \textbf{LLV algebra of $R$}.
\end{definition}
Recall that we require that the action of $\mathfrak{a}$ on $M$ is by mutually commuting endomorphisms. This is indeed the case in the above, since $R$ is graded commutative.

\begin{example}
    Continuing Example \ref{ex: kahler manifold cohom HL}, let $X$ be a compact connected Kähler manifold of dimension $n$ over $\mathbb{C}$ and let $R = \sH^\bullet(X; K)$. Then the LLV algebra $\llv(R)$ equals the Lie algebra $\llv(X; K)$ from Example \ref{ex: LLV of var}.
\end{example}

\begin{lemma}
    \label{lem: R_2 to g(R)_2 injective}
    Let $R$ be a graded commutative $K$-algebra of depth $n$ with the hard Lefschetz property. Then there is an injective map $R_2 \to \llv(R)_2$ sending $r \in R_2$ to the operator $e_r \colon R \to R$ of multiplication by $r$.
\end{lemma}
\begin{proof}
    Let $U \subseteq R_2$ be the set of $r \in R_2$ with the hard Lefschetz property. Then $U$ is Zariski-open in $R_2$, so it is dense since $K$ is infinite (recall that we assume $K$ to be of characteristic~$0$). Hence the map $e \colon U \to \llv(R)_2$ which sends $r$ to $e_r$ extends to all of $R_2$ by linearity. If $r, s \in R_2$ and $r \neq s$, then $e_r \neq e_s$, so the map $e$ is injective.
\end{proof}

One could wonder whether the natural map $R_2 \to \llv(R)_2$ also has to be surjective. This is not the case, it is possible to construct examples where this map is not surjective.  

\begin{lemma}
    \label{lem: f_a of subring}
    Let $R$ be a graded commutative $K$-algebra with the hard Lefschetz property and let $S \subseteq R$ be a subalgebra of $R$. Assume that $R$ and $S$ both have depth $n$, and assume there is an $a \in S_2 \subseteq R_2$ with the hard Lefschetz property for both $R[n]$ and $S[n]$. Write $e_a \colon R \to R$ for the operator of multiplication by $a$, and $f_a \colon R \to R$ for its dual operator (so $(e_a, h, f_a)$ is an $\mathfrak{sl}_2$-triple). Similarly, write $(e'_a, h, f'_a)$ for the $\mathfrak{sl}_2$-triple associated to $e'_a \colon S \to S$. Then $f'_a = f_a|_S$.
\end{lemma}
\begin{proof}
    This can be seen using the structure of irreducible $\mathfrak{sl}_2$-representations. The representation $S[n]$ of $\mathfrak{sl}_2 = \langle e_a', h, f_a'\rangle$ splits as a direct sum $\bigoplus_i V_i$ of irreducible representations $V_i$. Each such representation $V_i$ can be understood as having a piece of dimension $1$ in degrees $-k, -k + 2, \dots, k - 2, k$, where $e_a' \colon V_i \to V_i$ raises the degree of each element by $2$ and $f_a'$ lowers the degree by $2$. Then $[e_a', f_a'] = h$. Since $e_a|_S = e_a'$ and $V_i$ is generated by repeatedly applying $e_a'$ to an element of lowest degree, we see that the $V_i \subseteq R$ are also irreducible representations of $\mathfrak{sl}_2 = \langle e_a, h, f_a \rangle$. Hence, we see that $f_a|_{V_i} = f_a'|_{V_i}$ for each $i$, which implies that $f'_a = f_a|_S$.
\end{proof}

\begin{corollary}
    \label{cor: subring HL iff sl_2-invariant}
    Let $R$ be a graded commutative $K$-algebra with the hard Lefschetz property and let $S \subseteq R$ be a subalgebra of $R$. Assume that $R$ and $S$ both have depth $n$, and assume that there is an $a \in S_2 \subseteq R_2$ with the hard Lefschetz property for $R[n]$. Then the following are equivalent:
    \begin{enumerate}
        \item The action of $a$ on $S[n]$ has the hard Lefschetz property.
        \item The subspace $S[n] \subseteq R[n]$ is an invariant subspace for $\langle e_a, h, f_a \rangle$.
    \end{enumerate}
\end{corollary}
\begin{proof}
    If the action of $a$ on $S[n]$ has the hard Lefschetz property, then Lemma~\ref{lem: f_a of subring} implies that the $\mathfrak{sl}_2$-triple of the action of $a$ on $R[n]$ gives the $\mathfrak{sl}_2$-triple for the action of $a$ on $S[n]$, so $S[n]$ is indeed an invariant subspace.

    Now suppose that $S[n]$ is an invariant subspace for $\langle e_a, h, f_a \rangle$. Then $S[n]$ is a representation of $\mathfrak{sl}_2(K)$. The symmetry of $\mathfrak{sl}_2$-representations implies that $\dim S_{n-k} = \dim S_{n+k}$ for all $0 < k \leq n$. Since $S \subseteq R$ and $R$ has the hard Lefschetz property, we know that $a^k \colon R_{n-k} \to R_{n+k}$ is an isomorphism for all $k$, so it induces injections $S_{n-k} \to S_{n+k}$, which must be isomorphisms by the equality of dimensions. Hence the action of $a$ on $S$ has the hard Lefschetz property. 
\end{proof}

\begin{definition}
    Let $R$ be a graded commutative $K$-algebra with the hard Lefschetz property, and assume that $R_0 = K$. The \textbf{Verbitsky component} of $R$ is the subalgebra \myindex{$\SH(R)$} of $R$ generated by $R_2$.
\end{definition}

\begin{lemma}
    Let $R$ be a graded commutative $K$-algebra with $R_0 = K$ and depth $n$, and assume that $R$ has the hard Lefschetz property. Then $\SH(R)$ also has depth $n$.
\end{lemma}
\begin{proof}
    Take an $a \in R_2$ with the hard Lefschetz property. Then $R_{2n} = K \cdot a^n$, so $a^n \neq 0$. But $a^n \in \SH(R)$, so $\SH(R)$ also has depth $n$.
\end{proof}

By applying Corollary \ref{cor: subring HL iff sl_2-invariant} to $\SH(R)$, we obtain the following result. The fact that the first statement implies the second is also proven by Verbitsky \cite[Proposition~8.4]{VerbitskyPhDThesis}.

\begin{lemma}
    \label{lem: SH subrep iff HL}
    Let $R$ be a graded commutative $K$-algebra with $R_0 = K$ and depth $n$, and assume that $R$ has the hard Lefschetz property. The following are equivalent:
    \begin{enumerate}
        \item The $K$-algebra $\SH(R)$ has the hard Lefschetz property.
        \item The subspace $\SH(R)[n] \subseteq R[n]$ is $\llv(R)$-invariant.
    \end{enumerate}
\end{lemma}
\begin{proof}
    If $\SH(R)$ has the hard Lefschetz property, then Corollary \ref{cor: subring HL iff sl_2-invariant} implies that it is preserved by the action of all $\mathfrak{sl}_2$-triples. Therefore it is closed under the action of all generators of $\llv(R)$.

    Conversely, if $\SH(R)[n]$ is $\llv(R)$-invariant, it is in particular invariant under the action of all $\mathfrak{sl}_2$-triples associated to elements with the hard Lefschetz property on $R$. Hence Corollary \ref{cor: subring HL iff sl_2-invariant} implies that $\SH(R)$ has the hard Lefschetz property.
\end{proof}

The following notions were introduced by Looijenga and Lunts and by Verbitsky \cite{LooijengaLunts, VerbitskyPhDThesis}.
\begin{definition}
    Let $R$ be a graded commutative $K$-algebra. We say that $R$ is \textbf{Lefschetz} (of depth $n$) if $R$ has the following properties:
    \begin{enumerate}
        \item $R_0 = K$.
        \item $R$ has depth $n$.
        \item $R$ has the hard Lefschetz property.
        \item The LLV algebra $\llv(R)$ is a semisimple Lie algebra over $K$.
    \end{enumerate}
    If $(R_2, R[n])$ is in addition a Jordan-Lefschetz pair, we say that the algebra $R$ is \textbf{Jordan-Lefschetz} (of depth $n$).
\end{definition}

\begin{example}
    Let $X$ be a compact connected Kähler manifold of dimension $n$ over $\mathbb{C}$. Then $\sH^\bullet(X; K)$ is Lefschetz of depth $n$, this follows from Example~\ref{ex: LLV of var}.
\end{example}
Moreover, it turns out that $\sH^\bullet(X; \mathbb{Q})$ is Jordan-Lefschetz if $X$ is a hyperkähler or abelian variety (this was proven by both Looijenga-Lunts and Verbitsky for hyperkähler varieties \cite{LooijengaLunts, VerbitskyPhDThesis}, and by Looijenga-Lunts for abelian varieties \cite[Proposition~3.2]{LooijengaLunts}).

\begin{lemma}
    \label{lem: SH hard Lefschetz then irreducible}
    Let $R$ be a graded commutative Lefschetz algebra of depth $n$ and assume that $\SH(R)$ also has the hard Lefschetz property. Then $\SH(R)[n]$ is an irreducible representation of $\llv(R)$. 
\end{lemma}
\begin{proof}
    The submodule $\SH(R)[n] \subseteq R[n]$ has the hard Lefschetz property, so it is a representation of $\llv(R)$ by Lemma \ref{lem: SH subrep iff HL}.

    Next, observe that any subrepresentation $V \subseteq R[n]$ is graded, in the sense that if $v \in V$ can be written as $v = \sum_m v_m$ with $v_m \in R[n]_m$, we have $v_m \in V$ for each $m$. This follows from the fact that the grading operator $h \in \llv(R)$ is semisimple, so its action on $V$ is diagonalizable. 

    Since $R_0 = K$, there must be a unique irreducible subrepresentation of $R[n]$ containing $R_0$. Since $R$ is Lefschetz, the Lie algebra $\llv(R)$ semisimple, so the representation $\SH(R)[n]$ is completely reducible. If we have a decomposition $\SH(R)[n] = U \oplus V$, then without loss of generality $U_{-n} \neq 0$ and $V_{-n} = 0$. But since $\SH(R)$ is generated by $1$ and $R_2$, we see that $\SH(R)[n] \subseteq U$, so $\SH(R)[n]$ is irreducible.
\end{proof}

\begin{lemma}
    \label{lem: SH(R) absolutely irreducible}
    Let $R$ be a graded commutative Lefschetz algebra of depth $n$ and assume that $\SH(R)[n]$ is an irreducible representation of $\llv(R)$. Then $\SH(R)[n]$ is absolutely irreducible.
\end{lemma}
\begin{proof}
    Since $R_0 = K$, we also have $\SH(R)_{0} = K$. Let $A$ be the ring of endomorphisms of the $\llv(R)$-representation $\SH(R)[n]$, then Schurs lemma implies that $A$ is a division algebra over $K$. Now take $a \in A$ and $x \in \SH(R)_0$, and let $h \in \llv(R)$ be the grading operator. Then
    $$h(ax) = ahx = -nax,$$
    so also $ax \in \SH(R)_0$. It follows that $\SH(R)_0 = K$ is an $A$-module. Since $A$ is a division algebra over $K$, this is only possible if $A = K$. This implies that the endomorphism algebra of $\SH(R)[n] \otimes_K \overline{K}$ is $\overline{K}$, so Schurs lemma implies that $\SH(R)[n] \otimes_K \overline{K}$ is also irreducible.
\end{proof}

Let $R$ be a graded commutative Jordan-Lefschetz algebra and assume that $\SH(R)$ has the hard Lefschetz property. Since $\SH(R)_2 = R_2$, any element of $R_2$ that has the hard Lefschetz property for $R$ also has it for $\SH(R)$, so by Lemma \ref{lem: f_a of subring} we have a restriction morphism $\llv(R) \to \llv(\SH(R))$. The existence of such a morphism was also shown by Verbitsky \cite[Corollary 8.2]{VerbitskyPhDThesis}. In the next proposition, we will see that under suitable hypotheses this map is an isomorphism, and we will later use this to compute some LLV algebras.

\begin{proposition}
    \label{prop: gR isom to gSR}
    Let $R$ be a graded commutative Jordan-Lefschetz algebra. Assume that the natural map $R_2 \to \llv(R)_2$ from Lemma \ref{lem: R_2 to g(R)_2 injective} is bijective and that $\SH(R)$ has the hard Lefschetz property. Then $\SH(R)$ is also Jordan-Lefschetz, and the natural map $\llv(R) \to \llv(\SH(R))$ is an isomorphism.
\end{proposition}
\begin{proof}
    Since $\SH(R)_2 = R_2$, the natural map $\llv(R) \to \llv(\SH(R))$ is surjective. Because $R$ is Jordan-Lefschetz, its LLV algebra only lives only in degrees $-2, 0$ and $2$, so the same holds for $\llv(\SH(R))$, which implies that $\SH(R)$ is Jordan-Lefschetz as well.

    Since $\SH(R)$ has the hard Lefschetz property, Lemma \ref{lem: R_2 to g(R)_2 injective} gives an injection $R_2 \to \llv(\SH(R))_2$. The assumption that $\llv(R)_2 = R_2$ then implies that the map $\llv(R)_2 \to \llv(\SH(R))_2$ is injective. Then Lemma \ref{lem: check inj in deg 2} implies that the map $\llv(R) \to \llv(\SH(R))$ is injective. Since it is also surjective, we see that it must be an isomorphism.
\end{proof}

We end this section with the following result, which will be useful in combination with Lemma \ref{lem: tensor prod LLV} to decompose the LLV algebra of a tensor product into smaller factors.
\begin{lemma}
    \label{lem: factors of HL are HL}
    Let $R$ and $S$ be two graded commutative $K$-algebras of depth $m$ and $n$ respectively. Let $(r, s) \in R_2 \oplus S_2$ be an element whose action on $R \otimes S$ has the hard Lefschetz property. Then the actions of $r$ on $R$ and $s$ on $S$ both have the hard Lefschetz property.
\end{lemma}
\begin{proof}
    Consider the subalgebra \[R(s) := R \otimes K[s]/(s^{n+1}) \subseteq R \otimes S.\]
    The operator $e_{(r,s)} \colon R \otimes S \to R \otimes S$ restricts to a map $e_{(r, s)} \colon R(s) \to R(s)$. Since $e_{(r,s)}$ has the hard Lefschetz property on $R \otimes S$, the map 
    $$e^i_{(r,s)} \colon R(s)_{m+n-i} \to R(s)_{m+n+i}$$ 
    must be injective for all $i$, so $\dim R_{m-i} \leq \dim R_{m+i}$ for all $i$. By switching the roles of $R$ and $S$, we see that also $\dim S_{n-i} \leq \dim S_{n+i}$ for all $i$. Since $R \otimes S$ has the hard Lefschetz property, this is only possible if $\dim R_{m-i} = \dim R_{m + i}$ and $\dim S_{n-i} = \dim S_{n+i}$ for all $i$. 
    
    The fact that the map $e^i_{(r, s)} \colon R(s)_{m+n-i} \to R(s)_{m+n+i}$ is injective for all $i$ forces the map $e_r^i \colon R_{m - i} \to R_{m + i}$ to also be injective for all $i$. The equality of dimensions then implies that $r$ has the hard Lefschetz property for its action on $R$. Similarly, the action of $s$ on $S$ has the hard Lefschetz property.
\end{proof}

\begin{corollary}
    \label{prop: factors of tensor prod HL}
    Let $R$ and $S$ be two graded commutative $K$-algebras of depth $m$ and $n$ respectively, with $R_0 = K$ and $S_0 = K$. Assume that $R_1 = 0$, and that $R \otimes S$ has the hard Lefschetz property. Then $R$ and $S$ both have the hard Lefschetz property.
\end{corollary}
\begin{proof}
    Since $R_1 = 0$ and $R_0 = S_0 = K$, we have $(R \otimes S)_2 = R_2 \oplus S_2$. Therefore there is an element $(r, s) \in R_2 \oplus S_2 = (R \otimes S)_2$ with the hard Lefschetz property, and the result follows from Lemma~\ref{lem: factors of HL are HL}.
\end{proof}

\section{LLV algebras of type $A_1$}
Let $K$ be a field, and let $R$ be a graded commutative Jordan-Lefschetz algebra over $K$. Then $\llv(R)$ is semisimple, so it is a direct product of simple Lie algebras over $K$. Let $\mathfrak{g}$ be such a simple factor, and assume that $\mathfrak{g}$ is of type $A_1$. Recall from Definition~\ref{def: type of lie alg} that this means that there is a finite field extension $F/K$ such that $\mathfrak{g}$ is defined over $F$, and $\mathfrak{g} \otimes_F \overline{F} \cong \mathfrak{sl}_2(\overline{F})$. In this section we will show that $\mathfrak{g}$ is then actually isomorphic to $\mathfrak{sl}_2(F)$, so it cannot be a nontrivial twisted form of $\mathfrak{sl}_2(F)$.

\begin{lemma}
    \label{lem: simple factor A1 of llv is sl2(F)}
    Let $K$ be a field of characteristic $0$, let $R$ be a graded commutative Jordan-Lefschetz algebra over $K$, and let $\mathfrak{g}$ be a simple factor of $\llv(R)$ of type $A_1$. Then $\mathfrak{g} \cong \mathfrak{sl}_2(F)$, where $F$ is the centroid of $\mathfrak{g}$.
\end{lemma}
\begin{proof}
    As argued in the proof of Lemma~\ref{lem: check inj in deg 2}, the simple factor $\mathfrak{g}$ of $\llv(R)$ is a graded Lie algebra, with grading induced by the grading on $\llv(R)$. Moreover, $\mathfrak{g}$ can not live purely in degree $0$. In particular, there is a nonzero element $e \in \mathfrak{g}_2$.
    
    Then $e$ is a nilpotent element of $\mathfrak{g}$, since it acts nilpotently on $R$. Since $F$ is the centroid of $\mathfrak{g}$, we can see $\mathfrak{g}$ as a Lie algebra over $F$. By applying the Jacobson-Morozov theorem over $F$ to the element $e \in \mathfrak{g}$, we obtain an injective homomorphism $j \colon \mathfrak{sl}_2(F) \to \mathfrak{g}$, see \cite[VIII.11.2, Proposition~2]{BourbakiLie7and8}. Since $\mathfrak{g}$ is of type $A_1$ and has centroid $F$, it has dimension $3$ over $F$, so $j$ is an isomorphism $\mathfrak{sl}_2(F) \cong \mathfrak{g}$.
\end{proof}

\section{Group actions on algebras with the hard Lefschetz property}
\label{sec: group actions on Lefschetz algebras}
In this section we consider graded commutative Lefschetz algebras $R$ equipped with an action of a group $G$ and study the LLV algebra of $R^G$. We will put a few assumptions on the action of $G$ on $R$:

\begin{setup}
    \label{stp: R^G in R Lefschetz}
    Let $R$ be a graded commutative Lefschetz algebra of depth $n$, and let $G$ be a group acting on $R$ by algebra automorphisms preserving the grading. Assume the following:
    \begin{itemize}
        \item The subalgebra $R^G \subseteq R$ is also Lefschetz of depth $n$.
        \item There is an $a \in R^G_2$ that has the hard Lefschetz property for both $R$ and~$R^G$.
    \end{itemize}
\end{setup}

The above conditions on the action of $G$ on $R$ are motivated by the cohomology of compact Kähler manifolds:
\begin{example}
    Let $X$ be a compact Kähler manifold, and let $G$ be a finite group acting freely on $X$. Then we will see in Lemma \ref{lem: X/G setup 4.20} that the action of $G$ on $\sH^\bullet(X; \mathbb{Q})$ satisfies the conditions of Setup \ref{stp: R^G in R Lefschetz}.
\end{example}

\begin{definition}
    \label{def: g_G(R)}
    Let $R$ and $G$ be as in Setup \ref{stp: R^G in R Lefschetz}. We denote by \myindex{$\llv^{\pre}_G(R)$} the Lie subalgebra of $\llv(R)$ generated by the elements $e_a$ and $f_a$ where $a$ ranges over the elements of $R_2^G$ that have the hard Lefschetz property for both $R$ and $R^G$.
\end{definition}

We will see that there is a natural surjection $\llv^{\pre}_G(R) \to \llv(R^G)$, and that this is an isomorphism under some mild hypotheses. The action of $G$ on $R$ induces an action of $G$ on $\mathfrak{gl}(R)$ by conjugation, and this restricts to an action of $G$ on $\llv(R)$.

\begin{lemma}
    \label{lem: g_G(R) in g(R)^G}
    Let $R$ and $G$ be as in Setup \ref{stp: R^G in R Lefschetz}. Then there is an inclusion $\llv^{\pre}_G(R) \subseteq \llv(R)^G$.
\end{lemma}
\begin{proof}
    It suffices to show the $G$-invariance of all $e_a$ and $f_a$ for $a$ running over the elements of $R_2^G$ having the hard Lefschetz property. For the $e_a$ the $G$-invariance is clear. 

    To see that $f_a$ is also invariant, we use the fact that $(\ad e_a)^2 \colon \mathfrak{gl}(R)_{-2} \to \mathfrak{gl}(R)_2$ is an isomorphism sending $f_a$ to $-2e_a$, see the proof of \cite[Proposition~2.1]{TaelmanDerivedEquivalences}. Since $e_a$ is $G$-invariant, $(\ad e_a)^2$ restricts to an isomorphism $\mathfrak{gl}(R)_{-2}^G \to \mathfrak{gl}(R)_2^G$. Hence $f_a$ is also $G$-invariant.
\end{proof}

Let $R$ and $G$ be as in Setup \ref{stp: R^G in R Lefschetz}. Then Lemma \ref{lem: f_a of subring} implies that restriction gives a homomorphism of Lie algebras $\llv^{\pre}_G(R) \to \llv(R^G)$. Since all generators of $\llv(R^G)$ are in the image of this map, it is a surjective homomorphism. 

\begin{lemma}
    \label{lem: inv of JL is JL}
    Let $R$ and $G$ be as in Setup \ref{stp: R^G in R Lefschetz}. If the algebra $R$ is Jordan-Lefschetz, then $R^G$ is also Jordan-Lefschetz.
\end{lemma}
\begin{proof}
    Since $R$ is Jordan-Lefschetz, $\llv(R)$ only lives in degrees $-2, 0$ and $2$, so the same holds for its subalgebra $\llv^{\pre}_G(R)$. Since we have a degree-preserving surjection $\llv^{\pre}_G(R) \to \llv(R^G)$, the Lie algebra $\llv(R^G)$ also has elements only in degrees $-2, 0$ and $2$, so $R^G$ is Jordan-Lefschetz.
\end{proof}

\begin{lemma}
    \label{lem: natural map g_G to g}
    Let $R$ and $G$ be as in Setup \ref{stp: R^G in R Lefschetz}, and assume that $R$ is Jordan-Lefschetz. If the map $R_2 \to \llv(R)_2$ is bijective and $\llv^{\pre}_G(R)$ is semisimple, the natural map $\llv^{\pre}_G(R) \to \llv(R^G)$ is an isomorphism.
\end{lemma}
\begin{proof}
    Let $U$ be the Zariski-open subset of $R^G_2$ consisting of the elements with the hard Lefschetz property for $R^G$, and let $V \subseteq R^G_2$ be the Zariski-open subset of elements with the hard Lefschetz property for both $R$ and $R^G$. Then the linear spans of $U$ and $V$ in $R_2^G$ are both equal to $R_2^G$, because $K$ is infinite. This implies that $R_2^G \subseteq \llv^{\pre}_G(R)_2$ and $R^G_2 \subseteq \llv(R^G)_2$.

    Now Lemma \ref{lem: g_G(R) in g(R)^G} implies that $\llv^{\pre}_G(R)_2 \subseteq \llv(R)^G_2 = R_2^G$, so we have an equality $\llv^{\pre}_G(R)_2 = R_2^G$. Hence the restriction map $\llv^{\pre}_G(R) \to \llv(R^G)$ is injective in degree $2$. Then Lemma \ref{lem: check inj in deg 2} implies that the map $\llv^{\pre}_G(R) \to \llv(R^G)$ is injective, where we use the assumption that $\llv^{\pre}_G(R)$ is semisimple. Since this map is also surjective, it is an isomorphism.
\end{proof}

\begin{lemma}
    \label{lem: tensor factor invariants HL}
    Let $R$ and $S$ be graded commutative algebras of depth $n$ and $m$, both with the hard Lefschetz property. Let $G$ be a group acting on $R$ and $S$ by algebra automorphisms preserving the grading. Assume that $(R_1 \otimes S_1)^G = 0$, that $R_0 = K$ and $S_0 = K$, and that the action of $G$ on $R \otimes S$ satisfies the conditions of Setup~\ref{stp: R^G in R Lefschetz}. Then $R^G$ and $S^G$ also have the hard Lefschetz property, and the actions $G \acts R$ and $G \acts S$ also satisfy the conditions in Setup \ref{stp: R^G in R Lefschetz}.
\end{lemma}
\begin{proof}
    First observe that the assumption $(R_1 \otimes S_1)^G = 0$ implies that $(R \otimes S_2)^G \cong R_2^G \oplus S_2^G$. Since the action of $G$ on $R \otimes S$ satisfies the conditions of Setup~\ref{stp: R^G in R Lefschetz}, there is an element $(r, s) \in R_2^G \oplus S_2^G$ whose action on $R \otimes S$ has the hard Lefschetz property for both $R \otimes S$ and $(R \otimes S)^G$.

    Lemma~\ref{lem: factors of HL are HL} implies that the actions of $r$ on $R$ and $s$ on $S$ have the hard Lefschetz property. Let $h_R$ be the grading operator of $R$ and $h_S$ the grading operator of $S$. We then have $\mathfrak{sl}_2$-triples $(e_r, h_R, f_r)$ and $(e_s, h_S, f_s)$ acting on $R$ and $S$.

    Since $r \in R_2^G$, the map $e_r \colon R \to R$ is $G$-invariant. Since $h_R$ is also $G$-invariant, it follows that $f_r$ is also $G$-invariant. In particular, $R^G \subseteq R$ is an invariant subspace under the action of $\langle e_r, h_R, f_r \rangle \cong \mathfrak{sl}_2(K)$, so the action of $r$ on $R^G$ has the hard Lefschetz property by Corollary~\ref{cor: subring HL iff sl_2-invariant}. Analogously, the action of $G$ on $S$ also satisfies the conditions of Setup~\ref{stp: R^G in R Lefschetz}.
\end{proof}
 
\section{LLV algebras of compact Kähler manifolds}
We will now study the LLV algebras of compact connected Kähler manifolds. Recall from Example \ref{ex: LLV of var} that for a compact connected Kähler manifold $X$, its LLV algebra $\llv(X; K)$ is defined as $\llv(\sH^2(X;K), \sH^\bullet(X; K))$. We will write \myindex{$\SH(X; K)$} for the Verbitsky component $\SH(\sH^\bullet(X; K))$, which is the subalgebra of $\sH^\bullet(X;K)$ generated by $\sH^2(X;K)$. We will now discuss the LLV algebras of hyperkähler manifolds and abelian varieties, these are two fundamental examples for the rest of this thesis. By using the decomposition from Theorem~\ref{thm: symplectic variety is quotient}, we will reduce the computation of the LLV algebra of a holomorphic symplectic variety to these two cases.

\begin{definition}
    Let $A$ be an abelian variety, and let $V = \sH^1(A; K)$. Let $\widetilde V = V \oplus V^*$, equipped with the bilinear form $b$ given by
    $$b((v, \eta), (v', \eta')) = \eta(v') + \eta'(v).$$
\end{definition}

We write $\sH^{\ev}(A; K) = \bigoplus_i \sH^{2i}(A; K)$ for the even part of the cohomology of $A$. In the next result, we use the spinor representation of an orthogonal Lie algebra. See \cite[Chapter 20]{Fulton-Harris} for more details on spinor representations.

\begin{proposition}
    \label{prop: LLV of AV}
    Let $A$ be an abelian variety over $\mathbb{C}$, and let $V = \sH^1(A; K)$. Then there is an isomorphism
    $$\llv(A; K) \cong \mathfrak{so}(\widetilde V, b).$$
    Moreover, the Verbitsky component of $A$ is equal to $\sH^{\ev}(A; K)$, and is isomorphic to the even spinor representation of $\mathfrak{so}(\widetilde V, b)$.
\end{proposition}
\begin{proof}
    See \cite[Proposition 3.3]{LooijengaLunts}, where this is proven with real coefficients. Their proof also works with an arbitrary field of characteristic $0$ as coefficients.
\end{proof}

\begin{definition}
    \label{def: Mukai lattice}
    Let $Y$ be a hyperkähler variety. The \textbf{Mukai lattice} of $Y$ is defined by\index{$\widetilde\sH(Y; K)$} 
    $$\widetilde\sH(Y; K) = K\alpha \oplus \sH^2(Y; K) \oplus K \beta.$$
    We equip $\widetilde \sH(Y; K)$ with a bilinear form $b$ by taking the orthogonal sum of the Beauville-Bogomolov-Fujiki form on $\sH^2(Y; K)$ and the hyperbolic plane $K\alpha \oplus K\beta$ (where $b(\alpha, \beta) = -1$). The grading on the Mukai lattice is given by giving $\alpha$ degree $-2$, giving $\beta$ degree $2$ and by giving $\sH^2(Y; K)$ degree $0$.
\end{definition}

\begin{proposition}
    \label{prop: LLV of HK}
    Let $Y$ be a hyperkähler variety. Then there is an isomorphism 
    \[\llv(Y; K) \cong \mathfrak{so}(\widetilde \sH(Y; K)).\]
\end{proposition}
\begin{proof}
    See {\cite[Proposition 4.5]{LooijengaLunts}} and {\cite{VerbitskyPhDThesis}} for a proof with $K = \mathbb{R}$. For the proof over $\mathbb{Q}$, see {\cite[Proposition 2.9]{SoldatenkovHyperkahler}} or {\cite[Theorem 2.7]{GreenHKLLV}}. By base change, the result is then true over any field of characteristic $0$.
\end{proof}

\begin{lemma}
    \label{lem: llv of product}
    Let $X$ and $Y$ be compact connected Kähler manifolds, and assume that $\sH^1(X; K) = 0$. Then the natural map $ \llv(X; K) \times \llv(Y; K) \to \llv(X \times Y; K)$ is an isomorphism.
\end{lemma}
\begin{proof}
    The Künneth formula implies that $\sH^\bullet(X \times Y; K) \cong \sH^\bullet(X; K) \otimes \sH^\bullet(Y; K)$, and the assumptions $\sH^0(X; K) = \sH^0(Y;K) = K$ and $\sH^1(X; K) = 0$ give $\sH^2(X \times Y; K) \cong \sH^2(X; K) \oplus \sH^2(Y; K)$. The result then follows from Lemma \ref{lem: tensor prod LLV}
\end{proof}

The assumption that $\sH^1(X; K) = 0$ is necessary. For example, consider two abelian varieties $A$ and $B$ of dimension at least $1$. Then Proposition \ref{prop: LLV of AV} implies that $\llv(A \times B; K)$ is simple, so it cannot be isomorphic to $\llv(A; K) \times \llv(B; K)$.

Let $X$ be a compact connected Kähler manifold, and let a group $G$ act on $X$ from the left. We then obtain a right action of $G$ on $\sH^\bullet(X; K)$ by letting $g \in G$ send $\omega \in \sH^\bullet(X; K)$ to $g^*\omega$. 

\begin{lemma}
    \label{lem: X/G setup 4.20}
    Let $X$ be a compact connected Kähler manifold, and let $G$ be a finite group acting freely on $X$. Then the action of $G$ on $\sH^\bullet(X; K)$ satisfies the conditions of Setup \ref{stp: R^G in R Lefschetz}.  
\end{lemma}
\begin{proof}
    We first prove this for $K = \mathbb{Q}$ by using Kähler classes. By base change, the result is then also true over an arbitrary field $K$ of characteristic $0$.

    Let $\pi \colon X \to X/G$ be the quotient map, then $\pi^* \colon \sH^\bullet(X/G; \mathbb{Q}) \to \sH^\bullet(X; \mathbb{Q})$ is injective and identifies $\sH^\bullet(X/G;\mathbb{Q})$ with $\sH^\bullet(X;\mathbb{Q})^G$. Since $X/G$ is also a Kähler manifold of the same dimension as $X$, the algebra $\sH^\bullet(X/G;\mathbb{Q})$ is Lefschetz of depth $\dim X$. 
    
    Therefore it suffices to show that there is an $\omega \in \sH^2(X;\mathbb{Q})^G$ that has the hard Lefschetz property for both $\sH^\bullet(X;\mathbb{Q})$ and $\sH^\bullet(X;\mathbb{Q})^G$. Let $\omega$ be a Kähler class on $X/G$. Then $\omega$ has the hard Lefschetz property for $\sH^\bullet(X/G; \mathbb{Q})$, and the pullback $\pi^*\omega$ is a Kähler class on $X$. Hence $\pi^*\omega$ has the hard Lefschetz property for both $\sH^\bullet(X;\mathbb{Q})$ and $\sH^\bullet(X;\mathbb{Q})^G$.
\end{proof}

Let $X$ be a compact connected Kähler manifold, and let $G$ be a finite group acting freely on $X$. By taking $R = \sH^\bullet(X; K)$ in Definition \ref{def: g_G(R)}, we get the Lie algebra $\llv^{\pre}_G(\sH^\bullet(X; K))$, which is the Lie subalgebra of $\llv(X; K)$ generated by the $\mathfrak{sl}_2$-triples coming from elements of $\sH^2(X; K)^G$ having the hard Lefschetz property for both $X$ and $X/G$.

By applying the results of Section \ref{sec: group actions on Lefschetz algebras}, we conclude the following:
\begin{corollary}
    Let $X$ be a compact connected Kähler manifold, and let $G$ be a finite group acting freely on $X$. Then there is a natural surjection 
    \[\llv^{\pre}_G(\sH^\bullet(X;K)) \to \llv(\sH^\bullet(X/G;K)).\] \qed
\end{corollary}

\begin{corollary}
    \label{cor: quot of JL is JL}
    Let $X$ be a compact connected Kähler manifold, and let $G$ be a finite group acting freely on $X$. If $\sH^\bullet(X; K)$ is Jordan-Lefschetz then $\sH^\bullet(X/G;K)$ is Jordan-Lefschetz as well. Moreover, if $\llv^{\pre}_G(\sH^\bullet(X;K))$ is semi-simple and $\llv(X; K)_2 = \sH^2(X;K)$, then the map \[\llv^{\pre}_G(\sH^\bullet(X; K)) \to \llv(\sH^\bullet(X/G; K))\]
    is an isomorphism.
\end{corollary}
\begin{proof}
    By Lemma \ref{lem: X/G setup 4.20} the action of $G$ on $\sH^\bullet(X; K)$ satisfies the conditions of Setup \ref{stp: R^G in R Lefschetz}, and the result then follows from Lemmas \ref{lem: inv of JL is JL} and \ref{lem: natural map g_G to g}.
\end{proof}

\section{LLV algebras of varieties admitting a holomorphic symplectic form}
We will now restrict ourselves to smooth projective varieties over $\mathbb{C}$ which can be written as a quotient of a product of abelian and hyperkähler varieties. By Theorem~\ref{thm: symplectic variety is quotient}, this includes all varieties which admit a holomorphic symplectic form. Therefore, we fix the following notation for the rest of this chapter:

\begin{setup}
    \label{stp: symplectic var}
    Let $X = (\prod_{i = 0}^k X_i)/G$, where $G$ is a finite group acting freely, $X_0$ is an abelian variety and $X_1, \dots, X_k$ are hyperkähler varieties. We denote by \index{$X_{\HK}$}$X_{\HK} = X_1 \times \dots \times X_k$ the \textbf{hyperkähler part} of $X$. By Theorem \ref{thm: HKA automorphisms}, we have actions of $G$ on $X_0$ and $X_{\HK}$.
\end{setup}

\begin{remark}
    Note that the varieties of the above form are more general than varieties admitting a holomorphic symplectic form, since we do not require the abelian variety $X_0$ to have a symplectic form, and also do not require the action of $G$ to preserve a symplectic form. It is only in Chapters \ref{chap: fatcorwise psi} and \ref{chap: construction of Psi} that we will actually require our varieties to be symplectic.
\end{remark}

\begin{proposition}
    \label{prop: HKA is JL}
    Let $X = (\prod_i X_i)/G$ be as in Setup \ref{stp: symplectic var}. Then $\sH^\bullet(X; K)$ is Jordan-Lefschetz.
\end{proposition}
\begin{proof}
     First note that $\sH^\bullet(X_0; K)$ is Jordan-Lefschetz by \cite[Proposition 3.3]{LooijengaLunts}, and all the $\sH^\bullet(X_i; K)$ are Jordan-Lefschetz by \cite[Proposition~4.5]{LooijengaLunts} (they only prove these statements for $K = \mathbb{R}$, but this directly implies the statements for any field $K$ of characteristic $0$). Therefore $\sH^\bullet(X_0 \times \dots \times X_k; K)$ is Jordan-Lefschetz by Lemma \ref{lem: llv of product} and Corollary \ref{cor: product of JL is JL}. Then Corollary~\ref{cor: quot of JL is JL} gives that $\sH^\bullet(X; K)$ is also Jordan-Lefschetz.
\end{proof}

\begin{lemma}
    \label{lem: prod X_i g_2 = H^2}
    Let $X_0$ be an abelian variety and let $X_1, \dots, X_k$ be hyperkähler varieties. Then $\llv(\prod_{i=0}^k X_i; K)_2 = \sH^2(\prod_{i=0}^k X_i; K)$.
\end{lemma}
\begin{proof}
    We know that $\llv(\prod X_i; K) = \prod_i \llv(X_i; K)$ from Lemma \ref{lem: llv of product}, and we have $\sH^2(\prod X_i; K) = \bigoplus_i \sH^2(X_i; K)$. From Propositions \ref{prop: LLV of AV} and \ref{prop: LLV of HK} we see that for every $i$ we have an equality $\llv(X_i; K)_2 = \sH^2(X_i; K)$, and this implies the result.
\end{proof}

We will now show that $\SH(X; K)$ has the hard Lefschetz property for $X$ as in Setup \ref{stp: symplectic var}. This result is already known for hyperkähler varieties and abelian varieties, and our proof will be by reducing to these cases. The following proof is based on the ideas from the proof of \cite[Theorem 2.15]{GreenHKLLV}.
\begin{proposition}
    \label{prop: SH(X) is HL}
    Let $X = (\prod_i X_i)/G$ be as in Setup \ref{stp: symplectic var}. Then $\SH(X; K)$ has the hard Lefschetz property.
\end{proposition}

Then Lemma \ref{lem: SH hard Lefschetz then irreducible} gives that $\SH(X; K)$ is irreducible, and hence absolutely irreducible by Lemma \ref{lem: SH(R) absolutely irreducible}. This implies the following:
\begin{corollary}
    \label{cor: SH(X) is absolutely irreducible}
    Let $X = (\prod_i X_i)/G$ be as in Setup \ref{stp: symplectic var}. Then $\SH(X; K)$ is absolutely irreducible as an $\llv(X; K)$-representation.\qed
\end{corollary}

\begin{proof}[{Proof of Proposition \ref{prop: SH(X) is HL}}]
    Throughout the proof, we will identify $\sH^\bullet(X;K)$ with the subring 
    $$\sH^\bullet\left(\prod_{i=0}^k X_i; K\right)^G \subseteq \sH^\bullet\left(\prod_{i=0}^k X_i;K\right).$$
    Recall that $\llv^{\pre}_G(\prod X_i; K)$ is the Lie subalgebra of $\llv(\prod X_i; K)$ generated by the $\mathfrak{sl}_2$-triples coming from elements of $\sH^2(X;K)$ having the hard Lefschetz property for both $\sH^\bullet(X; K)$ and $\sH^\bullet(\prod X_i; K)$. We will show that $\SH(X; K)$ is invariant under the action of $\llv^{\pre}_G(\prod X_i)$. It then follows that $\SH(X; K)$ is also invariant under the action of $\llv(X;K)$, and then Lemma \ref{lem: SH subrep iff HL} implies that $\SH(X; K)$ indeed has the hard Lefschetz property.
    
    By Corollary \ref{cor: quot of JL is JL}, $\llv^{\pre}_G(\prod X_i; K)$ only has elements in degrees $-2, 0$ and $2$. Furthermore, by applying Lemma \ref{lem: prod X_i g_2 = H^2} to Corollary \ref{cor: quot of JL is JL}, we obtain an isomorphism $\llv^{\pre}_G(\prod X_i; K)_{2} \cong \sH^2(X; K)$. This acts on $\sH^\bullet(\prod X_i; K)$ by multiplication, so it is clear that $\SH(X; K)$ is preserved under the action of $\llv^{\pre}_G(\prod X_i; K)_{2}$.

    Next, we will show that $\SH(X; K)$ is preserved under the action of the degree $0$ part $\llv^{\pre}_G(\prod X_i ; K)_{0}$. We know from Corollary~\ref{lem: llv of product} that $\llv(\prod X_i; K) \cong \prod_i \llv(X_i; K)$. For $i = 1, \dots, k$, the Lie algebra $\llv(X_i; K)_0$ is isomorphic to $\mathfrak{so}(\sH^2(X_i; K)) \times Kh_i$, where $\mathfrak{so}(\sH^2(X_i;K))$ acts on $\sH^\bullet(X_i;K)$ by derivations \cite[Proposition 4.5]{LooijengaLunts} and $h_i$ is the grading operator on $\sH^\bullet(X_i;K)$. An easy verification shows that $h_i$ is also a derivation (this fact was observed by Verbitsky \cite[Chapter 8]{VerbitskyPhDThesis}). For $X_0$, the situation is similar. It follows from Looijenga-Lunts \cite[Proposition 3.2]{LooijengaLunts} that for every element $u \in \llv(X_0; K)_0$, the endomorphism of $\sH^\bullet(X_0;K)$ is a sum $u_0 + \lambda \cdot \id$ with $\lambda \in K$ and $u_0$ a derivation.

    Now take a $u \in \llv^{\pre}_G(\prod X_i; K)_{0}$. It follows from the above that we can write $u = \sum_{i=0}^k D_i + \lambda$, where each $D_i$ is a derivation on $\sH^\bullet(X_i;K)$ and $\lambda$ acts as multiplication by $\lambda$ on $\sH^\bullet(X; K)$. It is not hard to see that $\SH(X; K)$ is preserved by derivations, and it is also preserved by scalar multiples of the identity. Therefore $u$ preserves $\SH(X; K)$.

    Finally, let us take an element $f_a \in \llv^{\pre}_G(\prod X_i; K)_{-2}$ where $a \in \sH^2(X;K)$ and take $x = x_1 \cdots x_k \in \SH(X; K)_{2k}$ (the argument now closely follows \cite[Theorem 2.15]{GreenHKLLV}). We will show by induction on $k$ that $f_a(x) \in \SH(X; K)_{2k-2}$. It is clear for $k = 0$, since then $f_a(x) = 0$, so let us take $k > 0$. We then have $x = e_{x_1}(x_2 \cdots x_k)$, so:
    $$f_a(x) = f_a(e_{x_1}(x_2 \cdots x_k)) = [e_{x_1}, f_a](x_2 \cdots x_k) - e_{x_1}(f_a(x_2 \cdots x_k)).$$
    We have already seen that $[e_{x_1}, f_a](x_2 \cdots x_k) \in \SH(X; K)$, and $f_a(x_2 \cdots x_k) \in \SH(X; K)$ by the induction hypothesis, so $f_a(x) \in \SH(X; K)$. 
\end{proof}

\section{Decomposing the LLV algebra of $X$}
\label{sec: LLV is HK times Ab}
In this section we split the LLV algebra of a variety $X = (\prod X_i)/G$ as in Setup~\ref{stp: symplectic var} into a direct product of smaller factors, and in Sections \ref{sec: HK LLV} and \ref{sec: abelian LLV} we will explicitly compute these factors. Recall that $K$ denotes a field of characteristic $0$.

\begin{lemma}
    \label{lem: HHK and HAb are JL}
    Let $X$ be as in Setup \ref{stp: symplectic var}. Then the actions of $G$ on $\sH^\bullet(X_0; K)$ and $\sH^\bullet(X_{\HK}; K)$ satisfy the conditions of Setup \ref{stp: R^G in R Lefschetz}, and the graded commutative rings $\sH^\bullet(X_0; K)^G$ and $\sH^\bullet(X_{\HK}; K)^G$ are Jordan-Lefschetz.
\end{lemma}
\begin{proof}
    The action of $G$ on $\sH^\bullet(\prod_{i=0}^k X_i; K)$ satisfies the conditions of Setup~\ref{stp: R^G in R Lefschetz} by Lemma \ref{lem: X/G setup 4.20}, and $\sH^\bullet(\prod_{i=0}^k X_i; K)$ is Jordan-Lefschetz by Proposition~\ref{prop: HKA is JL}. Lemma~\ref{lem: tensor factor invariants HL} implies that the actions of $G$ on $\sH^\bullet(X_0; K)$ and $\sH^\bullet(X_{\HK}; K)$ satisfy the conditions of Setup~\ref{stp: R^G in R Lefschetz}. Then $\sH^\bullet(X_0; K)^G$ and $\sH^\bullet(X_{\HK}; K)^G$ are Jordan-Lefschetz by Lemma~\ref{lem: inv of JL is JL}.
\end{proof}

In particular, the rings $\sH^\bullet(X_0; K)^G$ and $\sH^\bullet(X_{\HK}; K)^G$ have the hard Lefschetz property, so their LLV algebras are defined.

\begin{definition}
\label{def: ab and hk origin}
    Let \index{$\llv_{\Ab}(X; K)$}$\llv_{\Ab}(X; K) = \llv(\sH^\bullet(X_0; K)^G)$, we call this the \textbf{abelian part} of the LLV algebra of $X$. Similarly, let \index{$\llv_{\HK}(X; K)$}$\llv_{\HK}(X; K) = \llv(\sH^\bullet(X_{\HK}; K)^G)$, we call this the \textbf{hyperkähler part} of $\llv(X; K)$. We call the simple factors of $\llv_{\Ab}(X; K)$ the factors of \textbf{abelian origin}, and the simple factors of $\llv_{\HK}(X; K)$ the factors of \textbf{hyperkähler origin}.
\end{definition}

A priori, the Lie algebras $\llv_{\Ab}(X; K)$ and $\llv_{\HK}(X; K)$ depend on the chosen presentation of $X$ as a quotient $X = (\prod X_i)/G$. However, we will see in Lemma~\ref{lem: llvAb and llvHK independent} below that $\llv_{\Ab}(X; K)$ and $\llv_{\HK}(X; K)$ do not depend on this presentation. Therefore the Definition \ref{def: ab and hk origin} is independent of the chosen quotient presentation of $X$. 

Our aim in this section is to prove that $\llv(X; K) \cong \llv_{\Ab}(X; K) \times \llv_{\HK}(X; K)$. The inclusion
\begin{equation}
    \label{eq: inclusion HAb otimes HHK in H}
    \sH^\bullet(X_0; F)^G \otimes \sH^\bullet(X_{\HK}; K)^G \subseteq \sH^\bullet(X; K)^G
\end{equation}
is a strict inclusion in general, so we cannot directly apply Lemma \ref{lem: tensor prod LLV}. However, by looking at Verbitsky components, we can still compute $\llv(X; K)$. Therefore, we introduce the following notation: 
\begin{definition}
    Let $X$ be as in Setup \ref{stp: symplectic var}. We write \myindex{$\SH_{\Ab}(X; K)$} for the Verbitsky component of the ring $\sH^\bullet(X_0; K)^G$, and \myindex{$\SH_{\HK}(X; K)$} for the Verbitsky component of $\sH^\bullet(X_{\HK}; K)^G$.
\end{definition}

\begin{proposition}
    \label{prop: LLV is HK times Ab}
    Let $X$ be as in Setup \ref{stp: symplectic var}. Then the inclusion \eqref{eq: inclusion HAb otimes HHK in H} induces an isomorphism
    $$\SH_{\Ab}(X; K) \otimes \SH_{\HK}(X; K) \cong \SH(X; K).$$
\end{proposition}
\begin{proof}
    The inclusion \eqref{eq: inclusion HAb otimes HHK in H} restricts to an injective ring homomorphism
    $$\iota \colon \SH_{\Ab}(X; K) \otimes \SH_{\HK}(X; K) \to \sH^\bullet(X; K).$$
    Since $\sH^1(X_{\HK}; K) = 0$, we have $\sH^2(X; K) \cong \sH^2(X_0; K)^G \oplus \sH^2(X_{\HK}; K)^G$, so~$\iota$ is an isomorphism in degree $2$. Since the Verbitsky component is exactly the subring generated by the elements in degree $2$, it follows that the image of $\iota$ is equal to $\SH(X; K)$.
\end{proof}

\begin{corollary}
    \label{cor: LLV HK and Ab}
    Let $X$ be as in Setup \ref{stp: symplectic var}. Then there is an isomorphism 
    $$\llv(X; K) \cong \llv_{\Ab}(X; K) \times \llv_{\HK}(X; K).$$
\end{corollary}
\begin{proof}
    By Proposition \ref{prop: SH(X) is HL}, $\SH(X; K)$ has the hard Lefschetz property. By taking $G$-invariants in Lemma \ref{lem: prod X_i g_2 = H^2}, we see that $\llv(X; K)_2 = \sH^2(X; K)$. Hence we can apply Proposition~\ref{prop: gR isom to gSR} to conclude that $\llv(X; K) \cong \llv(\SH(X; K))$. Corollary~\ref{prop: factors of tensor prod HL} implies that both $\SH_{\Ab}(X; K)$ and $\SH_{\HK}(X; K)$ have the hard Lefschetz property. Then Lemma~\ref{lem: tensor prod LLV} and Proposition~\ref{prop: LLV is HK times Ab} imply that 
    \[\llv(X; K) \cong \llv(\SH_{\Ab}(X; K)) \times \llv(\SH_{\HK}(X; K)).\]
    Then Proposition \ref{prop: gR isom to gSR} gives isomorphisms $\llv(\SH_{\Ab}(X; K)) \cong \llv(\sH^\bullet(X_0; K)^G)$ and $\llv(\SH_{\HK}(X; K)) \cong \llv(\sH^\bullet(X_{\HK}; K)^G)$, which finishes the proof.
\end{proof}

\begin{lemma}
    \label{lem: llvAb and llvHK independent}
    Let $X$ be a smooth connected projective variety over $\mathbb{C}$ admitting a holomorphic symplectic form. Suppose we have two ways of writing $X$ as in Setup \ref{stp: symplectic var}, say $X = (\prod X_i)/G$ and $X = (\prod X_i')/G'$. Write $\llv_{\Ab}(X'; K)$ and $\llv_{\HK}(X'; K)$ for the LLV algebras obtained from the presentation $X = (\prod X_i')/G'$. Then the two decompositions
    \begin{equation*}
        \llv(X; K) = \llv_{\Ab}(X; K) \times \llv_{\HK}(X; K)
    \end{equation*}
    and
    \begin{equation*}
        \llv(X; K) = \llv_{\Ab}(X'; K) \times \llv_{\HK}(X'; K)
    \end{equation*}
    coincide, i.e. we have equalities $\llv_{\Ab}(X; K) = \llv_{\Ab}(X'; K)$ and $\llv_{\HK}(X; K) = \llv_{\HK}(X'; K)$.
\end{lemma}
\begin{proof}
    In \cite[Proposition 3]{BeauvilleKatata}, Beauville proves that $X$ admits a unique (up to a non-unique automorphism) minimal Galois covering $\widetilde X \to X$, with $\widetilde X = \widetilde X_0 \times \widetilde X_1 \times \dots \times \widetilde X_l$ a product of an abelian variety $\widetilde X_0$ and hyperkähler varieties $\widetilde X_i$ for $i \geq 1$. Then $\widetilde X$ has the property that every Galois covering $\overline X \to X$, with $\overline X = \prod_j \overline X_j$ a product of an abelian variety and hyperkähler varieties, factors through $\widetilde X \to X$:
    \begin{equation*}
        \begin{tikzcd}
            \overline X \arrow[r] \arrow[dr] & \widetilde X \arrow[d] \\
            & X.
        \end{tikzcd}
    \end{equation*}
    
    Furthermore, for any Galois covering $\overline X \to X$, the map $\overline X \to \widetilde X$ is a Galois covering, say with Galois group $N$, and Beauville argues that every element $g \in N$ acts on $\overline X$ by an action of the form $g = (s, \id) \in \Aut(\overline X_0) \times \Aut(\overline X_{\HK})$, where $s \in \Aut(\overline X_0)$ is a translation. But then we see that $N$ acts trivially on $\sH^\bullet(\overline{X}; K)$. It follows that $\overline X \to \widetilde X$ induces an isomorphism on cohomology, which identifies $\sH^\bullet(\widetilde X_0; K)$ with $\sH^\bullet(\overline X_0; K)$ and $\sH^\bullet(\widetilde X_{\HK}; K)$ with $\sH^\bullet(\overline X_{\HK}; K)$. 

    Let $\SH_{\Ab}(\overline{X}; K) = \SH(\sH^\bullet(\overline X_0; K)^{\overline G})$, and define $\SH_{\HK}(\overline X; K)$ etc. similarly. Then the above implies that the two tensor product decompositions 
    $$\SH(X; K) \cong \SH_{\Ab}(\overline{X}; K) \otimes \SH_{\HK}(\overline{X}; K)$$
    and
    \begin{equation}
        \label{eq: SH widetilde X}
        \SH(X; K) \cong \SH_{\Ab}(\widetilde{X}; K) \otimes \SH_{\HK}(\widetilde{X}; K)
    \end{equation}
    coincide. It then follows that also the product decompositions
    $$\llv(X; K) = \llv_{\Ab}(\overline X; K) \times \llv_{\HK}(\overline X; K)$$
    and
    \begin{equation}
    \label{eq: llv widetilde X}
        \llv(X; K) = \llv_{\Ab}(\widetilde  X; K) \times \llv_{\HK}(\widetilde X; K)
    \end{equation}
    also coincide. 

    However, since composing the map $\widetilde X \to X$ with an automorphism of $\widetilde X$ also gives a minimal Galois covering, we must also check that modifying $\widetilde X$ by an automorphism does not change the decomposition. For this, note that any automorphism of $\widetilde X$ preserves the decomposition $\widetilde X = \widetilde X_0 \times \widetilde X_{\HK}$ by \cite[Section~3]{BeauvilleKatata}. Hence also the tensor product decomposition~\eqref{eq: SH widetilde X} stays invariant under all automorphisms of $\widetilde X$. The same must then hold for the product decomposition~\eqref{eq: llv widetilde X}.

    To finish the proof of the lemma, observe that the Galois coverings $\prod X_i \to X$ and $\prod X'_i \to X$ both factor throuh the minimal Galois presentation $\widetilde X \to X$, and hence the induced product decompositions of $\llv(X; K)$ must coincide.
\end{proof}

\section{The hyperkähler part of the LLV algebra of a symplectic variety}
\label{sec: HK LLV}
We continue with the notation from Setup \ref{stp: symplectic var}. By Corollary \ref{cor: LLV HK and Ab}, the LLV algebra of $X$ is the product of the hyperkähler part and the abelian part. In this section, we explicitly compute the hyperkähler part.

Recall that $X_{\HK} = X_1 \times \dots X_k$. Let $I = \{1, \dots, k\}$. Recall from Section~\ref{sec: reduction to quotients} that the action of $G$ on $X_{\HK}$ induces an action of $G$ on $I$. 
\begin{definition}
    Let $X$ be as in Setup \ref{stp: symplectic var}, and let $[i] \in I/G$ be an orbit. We write \index{$X_{[i]}$}$X_{[i]} = \prod_{j \in [i]} X_j$.
\end{definition}

Then $X_{\HK} = \prod_{[i] \in I/G} X_{[i]}$, and the action of $G$ on $X_{\HK}$ preserves the factors of this product. Similar to Lemma \ref{lem: HHK and HAb are JL}, by repeatedly applying Lemma \ref{lem: tensor factor invariants HL} and using Lemma~\ref{lem: inv of JL is JL} we see:
\begin{lemma}
    \label{lem: H(X_[i]) stp 4.20}
    Let $X$ be as in Setup \ref{stp: symplectic var}. Then for all $[i] \in I/G$, the action of $G$ on $\sH^\bullet(X_{[i]}; K)$ satisfies the conditions of Setup \ref{stp: R^G in R Lefschetz}, and $\sH^\bullet(X_{[i]}; K)^G$ is Jordan-Lefschetz. \qed
\end{lemma}

\begin{definition}
    Let $X$ be as in Setup \ref{stp: symplectic var}. Take $[i] \in I/G$. We then define\index{$\llv_G(X_{[i]}; K)$}\index{$\SH_G(X_{[i]}; K)$}:
    \begin{align*}
        \llv_G(X_{[i]}; K) &:= \llv(\sH^\bullet(X_{[i]}; K)^G),\\
        \SH_G(X_{[i]}; K) &:= \SH(\sH^\bullet(X_{[i]}; K)^G).
    \end{align*}
\end{definition}

\begin{proposition}
    \label{prop: HK SH and LLV decomposition}
    Let $X$ be as in Setup \ref{stp: symplectic var}. Then there is an isomorphism
    $$\SH_{\HK}(X; K) \cong \bigotimes_{[i] \in I/G} \SH_G(X_{[i]}; K),$$
    and this induces an isomorphism
    $$\llv_{\HK}(X; K) \cong \prod_{[i] \in I/G} \llv_G(X_{[i]}; K).$$
\end{proposition}
\begin{proof}
    Since the action of $G$ on $X_{\HK}$ is the products of actions of $G$ on the $X_{[i]}$, we have an isomorphism
    $$\sH^2(X_{\HK}; K) \cong \bigoplus_{[i] \in I/G} \sH^2(X_{[i]}; K)^G.$$
    This implies that the inclusion homomorphism
    $$\iota \colon \bigotimes_{[i] \in I/G} \SH_G(X_{[i]}; K) \to \sH^\bullet(X_{\HK}; K)^G$$
    is an isomorphism in degree $2$, and hence induces an isomorphism onto $\SH_{\HK}(X; K)$, which proves the first part of the proposition.

    The proof of the isomorphism of LLV algebras is then analogous to Corollary~\ref{cor: LLV HK and Ab}.
\end{proof}

Recall that we say that all the $\llv_G(X_{[i]}; K)$ are simple factors of $\llv(X; K)$ of hyperkähler origin.

We will now compute the factors $\llv_G(X_{[i]}; K)$ explicitly. For this, we start with some preliminary results.

\begin{lemma}
    \label{lem: generators of Lie algebra roots}
    Let $\mathfrak{g}$ be a semisimple Lie algebra over an algebraically closed field. Let $\mathfrak{h}$ be a Cartan subalgebra of $\mathfrak{g}$, with associated root system $R$. Let $R' \subseteq R$ be a subset of roots with the property that every element of $R$ is a finite sum of elements of $R'$. Then the Lie algebra $\mathfrak{g}$ is generated by the subset $\oplus_{\alpha \in R'} \mathfrak{g}_\alpha$.
\end{lemma}
\begin{proof}
    Write $\mathfrak{g}'$ for the Lie subalgebra of $\mathfrak{g}$ generated by $\oplus_{\alpha \in R'} \mathfrak{g}_\alpha$. We know that $\mathfrak{g} = \mathfrak{h} \oplus \bigoplus_{\alpha \in R} \mathfrak{g}_\alpha$. For any two roots $\alpha$ and $\beta$ of $\mathfrak{g}$ with $\alpha + \beta$ also a root, the inclusion $[\mathfrak{g}_\alpha, \mathfrak{g}_\beta] \subseteq \mathfrak{g}_{\alpha + \beta}$ is an equality, see \cite[VIII.2.4, Lemma 4]{BourbakiLie7and8}. Since $R'$ additively generates $R$ and all root spaces are one-dimensional, we see that every root space $\mathfrak{g}_\alpha$ for $\alpha \in R$ is contained in $\mathfrak{g}'$.

    It remains to show that $\mathfrak{h} \subseteq \mathfrak{g}'$. The Cartan subalgebra is spanned by the coroots $t_\alpha$ for $\alpha \in R$, and for every root $\alpha$ and nonzero elements $x_\alpha \in \mathfrak{g}_\alpha$ and $x_{-\alpha} \in \mathfrak{g}_{-\alpha}$, the element $[x_\alpha, x_{-\alpha}]$ is a nonzero multiple of $t_\alpha$ by \cite[Proposition~8.3]{HumphreysRepresentationTheory}. Hence $\mathfrak{h} \subseteq \mathfrak{g}'$.
\end{proof}

The following lemma is proven by Verbitsky \cite[Theorem 9.1]{VerbitskyPhDThesis}, but we also give a self-contained proof here. 
\begin{lemma}
    \label{lem: so(U+V) generators}
    Let $K$ be a field of characteristic zero, let $U = K\alpha \oplus K\beta$ be the hyperbolic plane with quadratic form $\begin{psmallmatrix} 0 & 1 \\ 1 & 0 \end{psmallmatrix}$ with respect to the basis $\{\alpha, \beta\}$ and let $(V,b)$ be a non-degenerate quadratic space over $K$ with $\dim(V) \geq 2$. Let $U \oplus V$ be the graded quadratic space where $V$ sits in degree $0$ and where $\alpha$ and $\beta$ have degrees $-2$ and $2$ respectively. Then the graded Lie algebra $\mathfrak{so}(U \oplus V)$ is generated by its elements in degrees $-2$ and $2$.
\end{lemma}
\begin{proof}
    It suffices to show this after passing to an algebraic closure of $K$, where we may put the quadratic form on $U \oplus V$ into standard form given by either the block matrix
    \begin{equation*}
        \begin{pmatrix}
            0 & I_n \\
            I_n& 0
        \end{pmatrix} \text{ or }
        \begin{pmatrix}
            0 & I_n & 0\\
            I_n & 0 & 0\\
            0 & 0 & 1
        \end{pmatrix},
    \end{equation*}
    where $\dim (U \oplus V)$ is either $2n$ or $2n+1$ depending on the parity of $\dim V$. In either case, we can identify $U$ with the subspace $\langle e_1, e_{n+1}\rangle$. With respect to the standard root system (see \cite[Section 18]{Fulton-Harris}), the space $\mathfrak{so}(U \oplus V)_2 \oplus \mathfrak{so}(U \oplus V)_{-2}$ is the direct sum over the root spaces corresponding to the following lists of roots:
    \begin{itemize}
        \item ($\dim V$ even): the roots $\pm(L_1 - L_i)$ and $\pm(L_1 + L_i)$ for $i \neq 1$.
        \item ($\dim V$ odd): the roots $\pm(L_1 - L_i)$ and $\pm(L_1 + L_i)$ for $i \neq 1$, and the roots $\pm L_1$. 
    \end{itemize}
    In both the case where $\dim V$ is even and where it is odd, it is clear that every root of the root system of $\mathfrak{so}(U \oplus V)$ can be written as a finite sum of the roots in the above lists, so we are done by Lemma \ref{lem: generators of Lie algebra roots}.
\end{proof}

We will also need the following geometric ingredients:
\begin{lemma}
    \label{lem: HK H^2 invariants dimension}
    Let $H$ be a finite group of symplectic automorphisms of a hyperkähler variety $Y$. Then $\sH^2(Y; K)^H$ has dimension at least $3$.
\end{lemma}
\begin{proof}
    See for example \cite[Lemma 3.5(1)]{MongardiK32}. The dimension of $\sH^2(Y; K)^H$ is independent of the field $K$ (as long as $K$ is of characteristic $0$). Therefore, this lemma can be proven by providing~$3$ linearly independent elements of $\sH^2(Y; \mathbb{C})^H$. We can take a symplectic form~$\omega$, its conjugate~$\overline{\omega}$, and the first Chern class of an $H$-equivariant ample line bundle on~$Y$.
\end{proof}

\begin{lemma}
    \label{lem: H^2 of X_[i]}
    Let $X$ and $G$ be as in Setup \ref{stp: symplectic var}, and take $i \in I$. For each $j \in [i]$, choose a $g_j \in G$ such that, under the induced action of $G$ on $[i]$, the element $g_j$ sends $j$ to $i$. Denote the induced isomorphism $X_j \to X_i$ by $g_j$. Then the map
    \begin{equation} 
        \label{eq: isom H^2 of X_[i]}
        \sH^2(X_i; K)^{G_i} \to \left( \bigoplus_{j \in [i]} \sH^2(X_j; K)\right)^G
    \end{equation}
    given by $\omega \mapsto (g_j^*\omega)_{j \in [i]}$ is an isomorphism, and does not depend on the chosen $g_j$.
\end{lemma}
\begin{proof}
    Recall that $G_i \subseteq G$ is the subgroup of elements that stabilize $i$ under the action on $I$. First observe that if $g_j' \in G$ is another element which sends $j$ to $i$, then $g_j'g_j^{-1} \in G_i$. So if $\omega \in \sH^2(X_i; K)^{G_i}$, then $g_j^*\omega = (g_j')^* \omega$. It follows that the map $\omega \mapsto (g_j^*\omega)_j$ is independent of the chosen $g_j$.

    First, we need to verify that for each $\omega \in \sH^2(X_i; K)^{G_i}$ the element $(g_j^*\omega)_j \in \bigoplus_{j \in [i]} \sH^2(X_j; K)$ is actually $G$-invariant, so take a $g \in G$. Take $j \in [i]$, and define $j' \in [i]$ to be the index such that $g$ sends $j'$ to $j$. Then the component $g_j^*\omega$ gets mapped by $g$ to $g^*g_j^* \omega \in \sH^2(X_{j'}; K)$. But the above argument then implies that $g^*g_j^*\omega = g^*_{j'}\omega$. Since this holds for all $j$, we see that $(g_j^*\omega_j)_j$ is indeed $G$-invariant.

    It is clear that the map \eqref{eq: isom H^2 of X_[i]} is injective. For surjectivity, take an element $(\omega_j)_j \in ( \bigoplus_{j \in [i]} \sH^2(X_j; K))^G$. We then have $\omega_j = g^*_j \omega_i$ for all $j \in [i]$, so we see that the map \eqref{eq: isom H^2 of X_[i]} is also surjective.
\end{proof}

Recall from Definition \ref{def: Mukai lattice} that the Mukai lattice $\widetilde \sH(Y; K)$ of a hyperkähler variety $Y$ is the direct sum of $\sH^2(Y; K)$ and a hyperbolic plane $U = K\alpha \oplus K\beta$. Then $\widetilde \sH(Y; K)$ is a non-degenerate graded quadratic space, with $\sH^2(Y; K)$ in degree $0$. For a group $H$ of automorphisms of $Y$, we define an action of $H$ on $\widetilde \sH(Y; K)$ by taking the direct sum of the induced action of $H$ on $\sH^2(Y; K)$ and the trivial action of $H$ on the hyperbolic plane $U$.

In Proposition \ref{prop: LLV of HK} we saw that $\llv(Y;K) \cong \mathfrak{so}(\widetilde \sH(Y; K))$ for a hyperkähler variety $Y$. We will now use this to compute the hyperkähler factors of the LLV algebra of a symplectic variety. 

\begin{proposition}
    \label{prop: LLV of HK factor}
    Let $X$ be as in Setup \ref{stp: symplectic var}. Take an $i \in I$. Let $[i] \subseteq I$ be the orbit of $i$ under the induced action of $G$ on $I$, and let $G_i$ be the stabilizer of $i$. Then there is an isomorphism $\llv_G(X_{[i]}; K) \cong \mathfrak{so}(\widetilde \sH(X_i; K)^{G_i})$, where $G_i$ acts on $\widetilde \sH(X_i; K)$ as described above.
\end{proposition}

\begin{proof}
    Recall that $\llv_G(X_{[i]}; K)$ is defined as the LLV algebra of $\sH^\bullet(X_{[i]}; K)^G$. By Lemma~\ref{lem: H(X_[i]) stp 4.20}, the action of $G$ on $\sH^\bullet(X_{[i]}; K)$ satisfies the conditions of Setup \ref{stp: R^G in R Lefschetz}. Then Lemma~\ref{lem: natural map g_G to g} gives that the surjection $\llv^{\pre}_G(\sH^\bullet(X_{[i]}; K)) \to \llv_G(X_{[i]}, K)$ is an isomorphism if we show that $\llv^{\pre}_G(\sH^\bullet(X_{[i]}; K))_2 = \sH^2(X_{[i]}; K)^G$ and that $\llv^{\pre}_G(\sH^\bullet(X_{[i]}; K))$ is a semisimple Lie algebra. Therefore we will now compute $\llv^{\pre}_G(\sH^\bullet(X_{[i]}; K))$ and observe that it satisfies these two properties.
    
    Lemma \ref{lem: tensor prod LLV} and Proposition \ref{prop: LLV of HK} imply that
    \begin{equation}
        \label{eq: hk factor LLV to prod so}
        \llv\left(\sH^\bullet(X_{[i]}; K)\right) \cong \llv\left(\bigotimes_{j \in [i]} \sH^\bullet(X_j; K)\right) \cong \prod_{j \in [i]} \mathfrak{so}(\widetilde \sH(X_j; K)).
    \end{equation}

    By Lemma \ref{lem: H^2 of X_[i]}, there is an isomorphism
    \[\sH^2(X_i; K)^{G_i} \isomto \Bigl( \bigoplus_{j \in [i]} \sH^2(X_j; K)\Bigl)^G, \text{ } \omega \mapsto (g_j^* \omega)_j.\]
    By extending this map with the diagonal inclusion $U \to \bigoplus_{j \in [i]} U$ of hyperbolic planes, we get a diagonal inclusion $\widetilde \sH(X_i; K)^{G_i} \subseteq \bigoplus_{j \in [i]} \widetilde \sH(X_j; K)$. From this, we obtain an inclusion $\mathfrak{so}(\widetilde \sH(X_i; K)^{G_i}) \subseteq \prod_{j \in [i]}\mathfrak{so}(\widetilde \sH(X_j; K))$. Under the isomorphism of Equation \eqref{eq: hk factor LLV to prod so}, in degrees $2$ and $-2$ the image of this inclusion equals the degree $2$ and $-2$ parts of $\llv^{\pre}_G(\sH^\bullet(X_{[i]}; K))$.

    Lemma \ref{lem: HK H^2 invariants dimension} gives that $\dim \sH^2(X_i; K)^{G_i} \geq 3$, so Lemma \ref{lem: so(U+V) generators} implies that $\mathfrak{so}(\widetilde \sH(X_i; K)^{G_i})$ is generated by its elements in degrees $2$ and $-2$. By definition, $\llv_G(\sH^\bullet(X_{[i]}; K))$ is also generated by its elements in degrees $2$ and $-2$. Therefore, these two Lie algebras must be equal. Note that $\mathfrak{so}(\widetilde \sH(X_i; K)^{G_i})$ is semisimple and that
    $$ \mathfrak{so}(\widetilde \sH(X_i; K)^{G_i})_2 = \sH^2(X_i; K)^{G_i} = \Bigl(\bigoplus_{j \in [i]} \sH^2(X_j; K)\Bigl)^G,$$
    so the map $\llv^{\pre}_G(\sH^\bullet(X_{[i]}; K)) \to \llv_G(X_{[i]}; K)$ is indeed an isomorphism.
\end{proof}

\section{The abelian part of the LLV algebra of a symplectic variety}
\label{sec: abelian LLV}
Let $X$ be as in Setup \ref{stp: symplectic var}, and recall that $K$ is a field of characteristic~$0$. In this section, we will compute $\llv_{\Ab}(X; K)$ in terms of algebras with involution. 

Recall that for a finite group $G$, we write $\Irr_K(G)$ for a set of representatives of the isomorphism classes of irreducible $K$-valued $G$-representations. For $\rho \in \Irr_K(G)$, we write $\rho^*$ for its dual.
\begin{definition}
    Let $G$ be a finite group. We let $\sim$ be the equivalence relation on $\Irr_K(G)$ defined by $\rho \sim \rho'$ if and only if either $\rho \cong \rho'$ or $\rho^* \cong \rho'$.
\end{definition}
In particular, every equivalence class has either $1$ or $2$ elements, depending on whether $\rho$ is self-dual or not. For a $G$-representation $V$ defined over $K$ and an equivalence class $[\rho] \in \Irr_K(G)/\sim$\index{$\Irr_K(G)/\sim$}, we define the $[\rho]$\textbf{-isotypic component} \myindex{$V_{[\rho]}$} to be the direct sum 
\begin{equation*}
    V_{[\rho]} := \bigoplus_{\sigma \in [\rho]} V_\sigma = \begin{cases}
    V_{\rho} & \text{ if } \rho \cong \rho^*\\
    V_{\rho} \oplus V_{\rho^*} & \text{ if } \rho \ncong \rho^*
\end{cases}
\end{equation*}
where $V_\sigma$ denotes the $\sigma$-isotypic component of $V$ for $\sigma \in \Irr_K(G)$.

We now let $V = \sH^1(X_0; K)$, then $V$ is a $G$-representation. Since $X_0$ is an abelian variety, we have an isomorphism of graded rings $\sH^\bullet(X_0; K) \cong \extp^\bullet V$. The decomposition $V \cong \bigoplus_{[\rho]} V_{[\rho]}$ gives an isomorphism $\extp^\bullet V \cong \bigotimes_{[\rho]} \extp^\bullet V_{[\rho]}$. Since $(\extp^\bullet V)^G$ has the hard Lefschetz property, Lemma~\ref{lem: tensor factor invariants HL} implies that for every $[\rho] \in \Irr_K(G)/\sim$ the $K$-algebra $(\extp^\bullet V_{[\rho]})^G$ has the hard Lefschetz property.
\begin{definition}
    \label{def: SH and LLV of V_{[rho]}}
    For $X = (\prod_i X_i)/G$ as in Setup \ref{stp: symplectic var}, let $V = \sH^1(X_0; K)$ and take $[\rho] \in \Irr_K(G)/\sim$. We introduce the following notations\index{$\SH_G(V_{[\rho]}; K)$}\index{$\llv_G(V_{[\rho]}; K)$}:
    \begin{align*}
        \SH_G(V_{[\rho]}; K) &:= \SH\left(\left(\extp^\bullet V_{[\rho]}\right)^G\right),\\
        \llv_G(V_{[\rho]}; K) &:= \llv\left(\left(\extp^\bullet V_{[\rho]}\right)^G\right).
    \end{align*}
\end{definition}

\begin{proposition}
    \label{prop: Ab SH and LLV decomposition}
    Let $X$ be as in Setup \ref{stp: symplectic var}. Then there is an isomorphism
    $$\SH_{\Ab}(X; K) \cong \bigotimes_{[\rho] \in \Irr_K(G)/\sim} \SH_G(V_{[\rho]}; K),$$
    and this induces an isomorphism
    $$\llv_{\Ab}(X; K) \cong \prod_{[\rho] \in \Irr_K(G)/\sim} \llv_G(V_{[\rho]}; K).$$
\end{proposition}
\begin{proof}
    The proof is analogous to the proof of Proposition \ref{prop: HK SH and LLV decomposition}. The decomposition $V = \bigoplus_{\rho \in \Irr_K(G)} V_\rho$ into isotypic components gives a decomposition
    $$\extp^2 V \cong \bigoplus_{\rho \in \Irr_K(G)} \extp^2 V_\rho \oplus \bigoplus_{\rho, \rho' \in \Irr_K(G), \rho \ncong \rho'} V_\rho \otimes V_{\rho'}.$$
    Note that $(V_\rho \otimes V_{\rho'})^G$ is nonzero if and only if $\rho$ and $\rho'$ are dual to each other. Hence we see that
    $$\sH^2(X_0; K)^G = \bigoplus_{[\rho] \in \Irr_K(G)/\sim}\left( \extp^2 V_{[\rho]} \right)^G.$$
    Furthermore, the decomposition $V = \bigoplus_{[\rho]} V_{[\rho]}$ induces a tensor product decomposition 
    $$\extp^\bullet V \cong \bigotimes_{[\rho] \in \Irr_K(G)/\sim} \extp^\bullet V_{[\rho]},$$
    and it follows that the natural map 
    $$\bigotimes_{[\rho] \in \Irr_K(G)/\sim} \SH_G(V_{[\rho]}; K) \to \SH_{\Ab}(X; K)$$
    is an isomorphism. The proof that this tensor product decomposition of the Verbitsky component induces a product decomposition of the LLV algebra is analogous to the proof of Corollary \ref{cor: LLV HK and Ab}.
\end{proof}

\begin{remark}
    The decomposition in this proposition simplifies a bit if $K = \mathbb{Q}$. In this case all representations are self-dual, so $[\rho] = \{\rho\}$ for all $\rho \in \Irr_\mathbb{Q}(G)$.
\end{remark}

We will now explicitly compute the factors $\llv_G(V_{[\rho]}; K)$ that occur in the decomposition of Proposition \ref{prop: Ab SH and LLV decomposition}. For this, we will work with the following data:
\begin{setup}
    \label{stp: isotypic hard Lefschetz algebra}
    Let $G$ be a finite group acting on a finite-dimensional vector space $V$ over $K$. Assume the following:
    \begin{itemize}
        \item There is an element of $(\extp^2 V)^G$ which has the hard Lefschetz property for both $(\bigwedge^\bullet V)^G$ and $\bigwedge^\bullet V$ (i.e. the action of $G$ on $\extp^\bullet V$ satisfies the assumptions of Setup \ref{stp: R^G in R Lefschetz}).
        \item There is a $[\rho] \in \Irr_K(G)/\sim$ with the property that $V \cong \rho^n$ or $V \cong \rho^n \oplus (\rho^*)^n$ for some integer $n$, depending on whether $\rho$ is self-dual or not.
    \end{itemize}
    We obtain the following objects from this:
    \begin{itemize}
        \item The factor $A$ of $K[G]$ associated to $[\rho]$, this is an algebra with involution, see Section \ref{sec: alg with involution from representation theory}. Note that $A$ is either simple or a product of two simple factors, depending on whether $\rho$ is self-dual or not. In either case, $A$ is simple as an algebra with involution.
        \item The $A$-module $\widetilde V = V \oplus V^*$ (where $A$ acts on $V^*$ as described in Definition~\ref{def: dual of module over algebra with involution}), equipped with the bilinear form $b$ described above Lemma~\ref{lem: bilin form on V tilde is equivariant}. This bilinear form $b$ is $A$-equivariant.
        \item The algebra \index{$\widetilde A$}$\widetilde A = \End_A(\widetilde V)$, equipped with the involution $\widetilde \sigma$ corresponding to the bilinear form~$b$.
    \end{itemize}
\end{setup}
If we are in this setting, the LLV algebra of $(\bigwedge^\bullet V)^G$ is defined, and we denote this LLV algebra by $\llv_G(V; K)$. Note that this is consistent with the notation $\llv_G(V_{[\rho]}; K)$ introduced in Definition \ref{def: SH and LLV of V_{[rho]}}.

We already know that $\llv(\extp^\bullet V) \cong \mathfrak{so}(\widetilde V)$, and this is the starting point for the computation of $\llv_G(V; K)$. If we have $G$ and $V$ as in Setup \ref{stp: isotypic hard Lefschetz algebra}, then the action of $G$ on $\extp^\bullet V$ satisfies the conditions of Setup \ref{stp: R^G in R Lefschetz}. Therefore we can also consider the LLV algebra $\llv_G^{\pre}(\extp^\bullet V) \subseteq \llv(\extp^\bullet V)$, and we have a surjection 
$$\pi \colon \llv_G^{\pre}\left(\extp^\bullet V\right) \surjto \llv_G(V; K).$$

\begin{proposition}
    \label{prop: inclusion llv_Gpre in s(widetilde A)}
    Let $V$ and $G$ be as in Setup \ref{stp: isotypic hard Lefschetz algebra}. Then the isomorphism $\llv(\extp^\bullet V) \cong \mathfrak{so}(\widetilde V)$ restricts to an injective homomorphism of Lie algebras $$j \colon \llv_G^{\pre}\left(\extp^\bullet V\right) \injto \mathfrak{s}(\widetilde A, \widetilde \sigma),$$
    where $\mathfrak{s}(\widetilde A, \widetilde \sigma) = \Skew(\widetilde A, \widetilde \sigma)'$.
\end{proposition}
It follows that we have the following diagram:
\begin{equation}
\label{eq: j pi diagram}
    \begin{tikzcd}
        \llv_G^{\pre}(\extp^\bullet V) \arrow[r, hookrightarrow, "j"] \arrow[d, twoheadrightarrow, "\pi"] & \mathfrak{s}(\widetilde A, \widetilde \sigma)\\
        \llv_G(V; K).
    \end{tikzcd}
\end{equation}
We will show that $\pi$ is always an isomorphism, and that, with one exception, the map $j$ is also surjective.

\begin{proof}[{Proof of Proposition \ref{prop: inclusion llv_Gpre in s(widetilde A)}}]
    Let $\mathfrak{so}_G^{\pre}(\widetilde V) \subseteq \mathfrak{so}(\widetilde V)$ be the Lie subalgebra generated by the $G$-invariant elements in degrees $2$ and $-2$. Then the isomorphism $\llv(\extp^\bullet V) \cong \mathfrak{so}(\widetilde V)$ restricts to an isomorphism $\llv_G^{\pre}(\extp^\bullet V) \cong \mathfrak{so}_G^{\pre}(\widetilde V)$. All the generators of $\mathfrak{so}^{\pre}_G(\widetilde V)$ are $G$-invariant, so we have an inclusion $\mathfrak{so}_G^{\pre}(\widetilde V) \subseteq \mathfrak{so}(\widetilde V)^G$. Lemma \ref{lem: Skew = so_A} gives an equality $\Skew(\widetilde A, \widetilde \sigma) = \mathfrak{so}(\widetilde V)^G$. It follows from the definition that $\mathfrak{so}_G^{\pre}(\widetilde V)$ is equal to its own derived Lie algebra, so the inclusion $\mathfrak{so}_G^{\pre}(\widetilde V) \subseteq \Skew(\widetilde A, \widetilde \sigma)$ implies that we have an inclusion $\mathfrak{so}_G^{\pre}(\widetilde V) \subseteq \mathfrak{s}(\widetilde A, \widetilde \sigma)$.
\end{proof}

Recall that if we have a field $F$ and a central simple algebra $A$ over $F$, then $\dim_F(A) = d^2$ for some integer $d$, and we call $d$ the \textbf{degree} of $A$. The degree is invariant under base change, i.e. we have \index{$\deg(A)$}$\deg(A) = \deg(A \otimes_F L)$ for every field extension $L/F$.

If $(A, \sigma)$ is a simple algebra with involution which comes from an equivalence class $[\rho]$ where $\rho$ is not self-dual, we can write $A \cong A_1 \times A_2$ for two simple algebras $A_1$ and $A_2$. Since the induced involution on $A$ is the swapping involution, there is an isomorphism $A_1 \cong A_2^{\op}$. In particular, we have $\deg(A_1) = \deg(A_2)$. We then define the degree of $(A, \sigma)$ in the following way (see \cite[Section 2.B]{BookOfInvolutions}):
\begin{definition}
    Let $(A, \sigma)$ be a simple algebra with involution, and assume that $A = A_1 \times A_2$ is a product of two simple algebras, with $\sigma$ the swapping involution. We define the \textbf{degree} of $(A, \sigma)$ to be \index{$\deg(A, \sigma)$}$\deg(A, \sigma) = \deg(A_1)$.
\end{definition}
We will generally remove the $\sigma$ from the notation, so we just write $\deg(A)$ for $\deg(A, \sigma)$ if this causes no confusion. Note that the index is stable under base change in the following sense: if $(A, \sigma)$ is an algebra with involution of the second kind, and $F := Z(A)$ is a field, then for every field extension $L$ of $F^\sigma$, we have 
\[\deg(A) = \deg(A \otimes_{F^\sigma} L),\]
where the index on the right hand side is understood to be $\deg(A \otimes_{F^\sigma} L, \sigma \otimes \id)$ if $A \otimes_{F^\sigma} L$ is a product of two simple algebras.

We will need the following two technical results, which imply that $\llv_G^{\pre}(\extp^\bullet V)$ is always semisimple, and that the map $j$ in diagram \eqref{eq: j pi diagram} is surjective unless $\deg(\widetilde A) = 4$ and $\widetilde \sigma$ is orthogonal.

\begin{lemma}
    \label{lem: type of s(A, sigma)}
    Let $V, G$ and $[\rho]$ be as in Setup \ref{stp: isotypic hard Lefschetz algebra}, and let $\widetilde A = \End_A(\widetilde V)$, equipped with the involution $\widetilde \sigma$. Then $\deg(\widetilde A)$ is even, and $\mathfrak{s}(\widetilde A, \widetilde \sigma)$ is a simple Lie algebra of the following type, where $l = \deg(\widetilde A)/2$:
    \begin{enumerate}
        \item If $\widetilde \sigma$ is symplectic, then $\mathfrak{s}(\widetilde A, \widetilde \sigma)$ is of type $C_l$.
        \item If $\widetilde \sigma$ is orthogonal, then $\mathfrak{s}(\widetilde A, \widetilde \sigma)$ is of type $D_l$.
        \item If $\widetilde \sigma$ is unitary, then $\mathfrak{s}(\widetilde A, \widetilde \sigma)$ is of type $A_{2l-1}$.
    \end{enumerate}
\end{lemma}
\begin{proof}
    We first prove that $\deg(\widetilde A)$ is even. First suppose that $A$ is simple as a $K$-algebra. Then there are an integer $r$ and a division algebra $D$ such that $A \cong M_r(D)$. Moreover, $A$ has, up to isomorphism, a unique simple module $V_0$. Then $\widetilde V \cong V_0^{2n}$, and $\widetilde A \cong M_{2n}(D^{\op})$, which implies that $\deg(A)$ is even.

    If $A \cong A_1 \times A_2$ is a product of two simple algebras, then $A$ has, up to isomorphism, two simple $A$-modules $V_1$ and $V_2$, where each $V_i$ is an $A_i$-module. Then $\widetilde V \cong V_1^{2n} \oplus V_2^{2n}$ (this follows from the assumptions in Setup~\ref{stp: isotypic hard Lefschetz algebra}). Let $D$ be the division algebra such that there exists an integer $r$ with $A_1 \cong M_r(D)$. Then $\widetilde A \cong M_{2n}(D) \times M_{2n}(D^{\op})$, from which it follows that $\deg(\widetilde A, \widetilde \sigma)$ is even. 

    Now that we have shown that $\deg(\widetilde A)$ is even, we will compute $\mathfrak{s}(\widetilde A, \widetilde \sigma)$. It suffices to compute this Lie algebra after base change to an algebraic closure, so let $F = Z(A)$ and let $\mathfrak{g} = \mathfrak{s}(\widetilde A, \widetilde \sigma) \otimes_{F^\sigma} \overline{F}$. Note that $F^\sigma$ is the centroid of $\mathfrak{s}(\widetilde A, \widetilde \sigma)$.

    Let $d = \deg(\widetilde A)$. Observe that
    \begin{equation*}
        \widetilde A \otimes_{F^\sigma} \overline{F} \cong \begin{cases}
            M_d(\overline{F}) & \text{if } \widetilde \sigma \text{ is of the first kind} \\
            M_d(\overline{F}) \times M_d(\overline{F}) & \text{if } \widetilde \sigma \text{ is of the second kind.}
        \end{cases}
    \end{equation*}
    
    It now follows from Theorem \ref{Thm: computation of s(A, sigma)} that the Lie algebra $\mathfrak{g}$ is of type $C_l, D_l$ or $A_{2l - 1}$ with $l = d/2$, depending on whether $\widetilde \sigma$ is symplectic, orthogonal or unitary.
\end{proof}

\begin{lemma}
    \label{lem: derived skew generated in degrees +-2}
    Let $V, G$ and $[\rho]$ be as in Setup \ref{stp: isotypic hard Lefschetz algebra}, let $\widetilde A = \End_A(\widetilde V)$ and let $F = Z(A)$. Then:
    \begin{enumerate}
        \item If we are not in the case that $(\widetilde A, \widetilde \sigma)$ is of orthogonal type and $\deg(\widetilde A) = 4$, then $\mathfrak{s}(\widetilde A, \widetilde \sigma)$ is generated by its elements in degrees $2$ and $-2$.
        \item If $\deg(\widetilde A) = 4$ and the involution $\widetilde \sigma$ is orthogonal, then the Lie subalgebra of $\mathfrak{s}(\widetilde A, \widetilde \sigma)$ generated by the elements in degree $2$ and $-2$ is isomorphic to $\mathfrak{sl}_2(F)$.
    \end{enumerate}
\end{lemma}

Note that in the case where $\deg(\widetilde A) = 4$ and $\widetilde \sigma$ is of orthogonal type, $\mathfrak{s}(\widetilde A, \widetilde \sigma)$ is a Lie algebra of type $D_2$. Over an algebraically closed field we have $\mathfrak{so}_4 \cong \mathfrak{sl}_2 \times \mathfrak{sl}_2$, and only one of these copies of $\mathfrak{sl}_2$ is visible in the LLV algebra, because the other $\mathfrak{sl}_2$ factor lives purely in degree zero.

\begin{proof}[{Proof of Lemma \ref{lem: derived skew generated in degrees +-2}}]
    A set of elements of $\mathfrak{s}(\widetilde A, \widetilde \sigma)$ generates $\mathfrak{s}(\widetilde A, \widetilde \sigma)$ if and only if it does so over an algebraic closure of $F^\sigma$, so it suffices to show that $\mathfrak{g} := \mathfrak{s}(\widetilde A, \widetilde \sigma) \otimes_{F^\sigma} \overline{F}$ is generated by its elements in degrees $2$ and $-2$.

    We will now describe a Cartan subalgebra of $\mathfrak{g}$ and give the associated root system for each type that occurs, and use Lemma \ref{lem: generators of Lie algebra roots} to show that $\mathfrak{g}$ is indeed generated by its elements in degrees $-2$ and $2$. For all of the following cases, recall that $V$ lives in degree $1$ and $V^*$ lives in degree $-1$, so that $\End(\widetilde V)$ is graded, with $\End(\widetilde V)_{2} = \Hom(V^*, V)$ and $\End(\widetilde V)_{-2} = \Hom(V, V^*)$. Throughout the rest of this proof, we let $\overline{A} = A \otimes_{F^\sigma} \overline{F}$, and let $l = \deg(\widetilde A)/2$.

    We will first assume that $\widetilde \sigma$ is of the first kind, then there is a unique (up to isomorphism) simple $\overline{A}$-module $V_0$ and $V \otimes_F \overline{F} \cong V_0^l$. Furthermore, the dual $V_0^*$ also is a simple $\overline{A}$-module (via the involution), so there is an $\overline{A}$-equivariant isomorphism $\varphi \colon V_0 \isomto V_0^*$. This induces an isomorphism $\widetilde V \otimes_F \overline{F} \cong V_0^{2l}$. Choose a basis of $V_0$. By taking $2l$ copies of this basis we obtain a basis of $\widetilde V \otimes_F \overline{F}$. If we let $r = \dim_{\overline{F}}(V_0)$, we can write $\{e_1, \dots, e_{2rl}\}$ for the obtained basis of $\widetilde V  \otimes_ F \overline{F}$. Let $h_i = \sum_{j = 1}^r e_{ir + j}$, so $h_i$ is the identity on the $i$-th copy of $V_0$ in $\widetilde V \otimes_F \overline{F}$ and zero on the other $2l-1$ copies of $V_0$.

    The bilinear form on $\widetilde V$ is given by
    $$b((v, \eta), (v', \eta')) = \eta(v') + \eta'(v)$$
    for $v, v' \in V$ and $\eta, \eta' \in V^*$. Under the isomorphism $\widetilde V \otimes_F \overline{F} \cong V_0^{2l}$, this will correspond to the bilinear form $b'$ on $V_0^{2l}$ given by
    $$b'((v_1, \dots, v_{l}, w_1, \dots, w_{l}), (v_1', \dots, v_{l}', w_1', \dots, w_{l}')) = \sum_{i = 1}^{l}(\varphi(w_i)(v_i') + \varphi(w_i')(v_i)$$
    for $(v_1, \dots, w_{l})$ and $(v'_1, \dots, w_{l}')$ in $V_0^{2l}$, where $\varphi \colon V_0 \isomto V_0^*$ is the isomorphism introduced above.

    Now let $h'_i = h_i - h_{l+i} \in \End_{\overline{A}}(V_0^{2l})$. A small computation shows that $h'_i \in \mathfrak{g}$ for every $i = 1, \dots, l$, and in fact $\mathfrak{h} := \langle h'_1, \dots, h'_l\rangle$ is a Cartan subalgebra of $\mathfrak{g}$. Let $L'_i \in \mathfrak{h}^*$ be the dual operator of $h'_i$, so $L_i'(h_j') = \delta_{ij}$.

    \textbf{Type $D_l$}: assume that $\mathfrak{g}$ is of type $D_l$. Then $\mathfrak{g}$ has the root system $R = \{\pm L'_i \pm L'_j\}_{i \neq j}$. The set of roots associated to the subspaces of degree $2$ and $-2$ is then equal to $\{L'_i + L'_j\}_{i \neq j} \cup \{-L'_i - L'_j\}_{i \neq j}$ (see \cite[Chapter 18]{Fulton-Harris}). If $l \geq 3$, this collection of roots additively generates $R$, so Lemma~\ref{lem: generators of Lie algebra roots} implies that $\mathfrak{g}$ is indeed generated by its elements in degrees $2$ and $-2$.

    The case $l = 2$ corresponds to the case $\deg(\widetilde A) = 4$. In this case we see that the set of roots generated by the elements in degrees $2$ and $-2$ is just $\{L'_1 + L'_2, -L'_1 -L'_2\}$. This is isomorphic to the root system $A_1$, so the Lie subalgebra of of $\mathfrak{g}$ generated by the elements in degrees $2$ and $-2$ is isomorphic to $\mathfrak{sl}_2(\overline{F})$. Then Lemma \ref{lem: simple factor A1 of llv is sl2(F)} implies that the Lie subalgebra of $\mathfrak{s}(\widetilde A, \widetilde \sigma)$ generated by the elements in degree $2$ and $-2$ is isomorphic to $\mathfrak{sl}_2(F)$.
    
    If $l = 1$, then $\mathfrak{s}(\End_A(\widetilde V), \widetilde \sigma) = 0$.

    \textbf{Type $C_l$}: now let $\mathfrak{g}$ be of type $C_l$. The root system is $R = \{\pm L'_i \pm L'_j\}_{1 \leq i, j \leq l}$, and the set of roots of the elements of degree $2$ or $-2$ is given by $\{L'_i + L'_j\} \cup \{ -L'_i - L'_j\}_{1 \leq i,j \leq n}$, see \cite[Chapter 16]{Fulton-Harris}. We see that this set of roots generates all of $R$, so we are done.

    \textbf{Type $A_{2l-1}$}: now let $\mathfrak{g}$ be of type $A_{2l-1}$. Since $\overline{A} \cong M_r(\overline{F}) \times M_r(\overline{F})$ for some integer $r$, there are, up to isomorphism, two simple $\overline{A}$-modules $V_0$ and $V_1$, and they are dual to each other. We can write $V \otimes_{F^\sigma} \overline{F} \cong V_0^l \oplus V_1^l$, and then $V^* \otimes_{F^\sigma} \overline{F} \cong (V_0^*)^l \oplus (V_1^*)^l$. Since $V_0$ and $V_1$ are not isomorphic, there is a well-defined projection
    $$\End_{\overline{A}}(V_0^l \oplus V_1^l \oplus (V_0^*)^l \oplus (V_0^*)^l) \to \End_{\overline{A}}(V_0^l \oplus (V_1^*)^l),$$
    and under this projection $\mathfrak{g}$ maps isomorphically onto $\mathfrak{sl}_{\overline{A}}(V_0^l \oplus (V_1^*)^l)$ (the $\overline{A}$-equivariant trace $0$ endomorphisms of $V_0^l \oplus (V_1^*)^l$), and this is isomorphic to $\mathfrak{sl}_{2l}(\overline{F})$. The degree $2$ part maps to $\Hom_{\overline{A}}((V_1^*)^l, V_0^l)$ and the degree $-2$ part maps to $\Hom_{\overline{A}}(V_0^l, (V_1^*)^l)$. If we take the standard Cartan subalgebra 
    $$\mathfrak{h} = \left\{a_1E_{1,1} + \dots + a_{2l}E_{2l, 2l} : \sum a_i = 0\right\}$$
    of $\mathfrak{sl}_{2l}(\overline{F})$ with dual $\mathfrak{h}^* = \langle L_1, \dots, L_{2l}\rangle / (\sum L_i = 0)$, then $R = \{L_i - L_j\}_{i \neq j}$, see \cite[Chapter 15]{Fulton-Harris}. We see that the elements of degrees $-2$ and $2$ correspond to the subset $\{L_i - L_j\}_{|i - j| \geq l}$ of $R$. A quick computation shows that this subset additively generates $R$, which is what we had to show.
\end{proof}

\begin{theorem}
    \label{thm: LLV computation}
    Let $G$ and $V$ be as in Setup \ref{stp: isotypic hard Lefschetz algebra}. Then $\llv_G(V; K) \cong \mathfrak{s}(\widetilde A, \widetilde \sigma)$ unless $\deg(\widetilde A) = 4$ and $\widetilde \sigma$ is orthogonal, in which case $\llv_G(V; K)\cong \mathfrak{sl}_2(F)$, where $F = Z(\widetilde A)$.
\end{theorem}

\begin{proof}
    We will use the maps $j$ and $\pi$ in Diagram \eqref{eq: j pi diagram}. It follows from Lemma~\ref{lem: derived skew generated in degrees +-2} that the image of $j$ is $\mathfrak{s}(\widetilde A, \widetilde \sigma)$, unless $\deg(\widetilde A) = 4$ and $\widetilde \sigma$ is orthogonal. In both cases, it follows from this computation and the injectivity of $j$ that $\llv_G^{\pre}(\extp^\bullet V)$ is semisimple. Then Lemma \ref{lem: natural map g_G to g} implies that $\pi$ is an isomorphism, and this finishes the proof. 
\end{proof}

In order to use these computations, we need to know the type of the involution $\widetilde \sigma$ on $\widetilde A$. For simplicity, we will only do this for $K = \mathbb{Q}$. We first prove a small lemma about the Lie algebra of skew elements of a product of algebras with involution.
Given two algebras with involution $(A, \sigma)$ and $(B, \tau)$, the product $A \times B$ carries a natural involution $\sigma \times \tau$, given by $(\sigma \times \tau)(a, b) = (\sigma(a), \tau(b))$.
\begin{lemma}
    \label{lem: Skew of product}
    Let $(A, \sigma)$ and $(B, \tau)$ be algebras with involution. Then
    $$\Skew(A \times B, \sigma \times \tau) \cong \Skew(A, \sigma) \times \Skew(B, \tau).$$
\end{lemma}
\begin{proof}
    For $a \in A$ and $b \in B$, the condition that $(\sigma \times \tau)(a, b) = -(a, b)$ is equivalent to $\sigma(a) = -a$ and $\tau(b) = -b$.
\end{proof}

Let $G$ be a finite group and let $\rho \in \Irr_\mathbb{Q}(G)$. Choose an irreducible constituent $\tau$ of $\rho_\mathbb{C}$. Then we have already seen that $\tau$ is self-dual if and only if all irreducible constituents of $\rho_\mathbb{C}$ are self-dual (Lemma \ref{lem: self dual can be checked on 1 rep}). If $\tau$ is self-dual, then character theory implies that $(\tau \otimes \tau)^G = \mathbb{C}$. Since $\extp^2 \tau \subseteq \tau \otimes \tau$ is a subrepresentation, one may ask whether the trivial subrepresentation of $\tau \otimes \tau$ already occurs in $\extp^2 \tau$. 

\begin{lemma}
    \label{lem: exterior power independent of i}
    Let $G$ be a finite group, let $\rho \in \Irr_\mathbb{Q}(G)$ and let $\tau$ be an irreducible constituent of $\rho_\mathbb{C}$. Assume that $\tau$ is self-dual. Then $(\extp^2 \tau)^G = 0$ if and only if all irreducible constituents $\rho_i$ of $\rho_\mathbb{C}$ satisfy $(\extp^2 \rho_i)^G = 0$.
\end{lemma}
\begin{proof}
    Let $\mathbf{1}$ be the constant character of $G$ with value $1$. Let $\chi_i$ be the character of $\rho_i$, and let $\wedge^2\chi_i$ be the character of $\extp^2 \rho_i$. For $g \in G$, we have
    $$\wedge^2\chi_i(g) = \frac{1}{2}(\chi_i(g)^2 -\chi_i(g^2)),$$
   see \cite[Page 51]{Isaacs}. There is a Galois extension $F/\mathbb{Q}$ over which the $\rho_i$ are already defined, and $\Gal(F/\mathbb{Q})$ acts transitively on the collection of irreducible constituents of $\rho_F$ \cite[Theorem 9.21]{Isaacs}. Given two irreducible constituents $\rho_i$ and $\rho_j$ of $\rho_{\mathbb{C}}$, let $\sigma \in \Gal(F/\mathbb{Q})$ be an automorphism mapping $\rho_i$ to $\rho_j$. The above formula for $\wedge^2\chi_i$ then implies that $\sigma \circ \wedge^2\chi_i = \wedge^2\chi_j$. Write $\langle \wedge^2\chi_i, \mathbf{1} \rangle$ for the inner product of the characters $\wedge^2\chi_i$ and $\mathbf{1}$. Since the inner product of two characters is always an integer, we see that:
    $$\langle \wedge^2 \chi_j, \mathbf{1}\rangle = \langle \sigma \circ \wedge^2 \chi_i, \mathbf{1} \rangle = \sigma(\langle \wedge^2\chi_i, \mathbf{1}\rangle) = \langle \wedge^2\chi_i, \mathbf{1}\rangle,$$
    and this implies that the dimension of $(\extp^2 \rho_i)^G$ does not depend on the chosen irreducible constituent $\rho_i$ of $\rho_\mathbb{C}$.
\end{proof}

Let $\rho \in \Irr_\mathbb{Q}(G)$ be an irreducible representation of a finite group $G$. Then one can decompose the base change $\rho \otimes_\mathbb{Q} \mathbb{C} \cong \bigoplus_{i=1}^k \rho_i^{r_i}$ into irreducible representations, where the $\rho_i$ are in $\Irr_\mathbb{C}(G)$. It turns out that the multiplicities $r_i$ are independent of $i$, so there is an $r \in \mathbb{N}$ such that $\rho \otimes_\mathbb{Q} \mathbb{C} \cong \bigoplus_{i=1}^k \rho_i^r$, see \cite[Theorem~9.21(a)]{Isaacs}.

\begin{lemma}
    \label{lem: widetilde sigma type computation}
    Let $G$ be a finite group and let $\rho \in \Irr_\mathbb{Q}(G)$. Decompose $\rho \otimes_\mathbb{Q} \mathbb{C} = \bigoplus_{i=1}^k \rho_i^r$ as a direct sum of irreducible representations $\rho_i \in \Irr_\mathbb{C}(G)$. Let $V = \rho^n$ for an integer $n \geq 1$. Then the involution $\widetilde \sigma$ on $\widetilde A = \End_G(\widetilde V)$ is of the following type:
    \begin{itemize}
        \item If the $\rho_i$ are not all self-dual, then $\widetilde \sigma$ is an involution of the second kind.
        \item If the $\rho_i$ are all self-dual and $(\bigwedge^2 \rho_i)^G = \mathbb{C}$ for all $i$, then $\widetilde \sigma$ is a symplectic involution.
        \item If the $\rho_i$ are all self-dual and $(\bigwedge^2 \rho_i)^G = 0$ for all $i$, then $\widetilde \sigma$ is an orthogonal involution.
    \end{itemize}
\end{lemma}

Lemma \ref{lem: exterior power independent of i} implies that the conditions in this lemma need to be checked for just a single irreducible constituent of $\rho$.

\begin{proof}
    Let $A$ be the simple factor of $\mathbb{Q}[G]$ associated to $\rho$. By Corollary~\ref{cor: self-dual constituents iff first kind}, the involution on $A$ is of the second kind if and only if all irreducible constituents of $\rho_i$ are not self-dual. Lemma~\ref{lem: Endomorphism inherits involution type} gives that $\sigma$ and $\widetilde \sigma$ are involutions of the same kind. Hence we see that $\widetilde \sigma$ is an involution of the second kind if and only if the $\rho_i$ are not self-dual.

    So from now on, we will assume that all $\rho_i$ are self-dual and that $\widetilde \sigma$ is an involution of the first kind. Using Lemma \ref{lem: involution type from skew elements}, we can determine the type of $\widetilde \sigma$ by looking at the dimension of $\Skew(\widetilde A, \widetilde \sigma)$. This can be done after passing to an algebraic closure, so we will compute the dimension of $\Skew(\widetilde A \otimes_\mathbb{Q} \mathbb{C}, \widetilde \sigma)$. 

    The decomposition $\rho \otimes_\mathbb{Q} \mathbb{C} = \bigoplus_{i = 1}^k \rho_i^r$ induces a decomposition 
    \begin{align*}
        \widetilde A \otimes_\mathbb{Q} \mathbb{C} &= \End_G(\rho^n \oplus (\rho^{n})^*)\otimes_\mathbb{Q} \mathbb{C}\\
            &\cong \End_G((\rho^n \oplus (\rho^n)^*) \otimes_\mathbb{Q} \mathbb{C})\\
            &\cong \prod_{i=1}^k \End_G( \rho_i^{rn} \oplus (\rho_i^*)^{rn})
    \end{align*}   
    Let $\widetilde A_i = \End_G(\rho_i^{rn} \oplus (\rho_i^*)^{rn})$. Since $\widetilde \sigma$ is an involution of the first kind, it restricts to involutions $\widetilde \sigma_i$ on all the simple factors $\widetilde A_i$. Then by Lemma \ref{lem: Skew of product}, 
    $$\Skew(\widetilde A, \widetilde \sigma) \otimes_\mathbb{Q} \mathbb{C} \cong \Skew\left(\prod_{i=1}^k \widetilde A_i, \prod_{i=1}^k \widetilde \sigma_i\right) \cong \prod_{i=1}^k \Skew(\widetilde A_i, \widetilde \sigma_i).$$
    We will now compute the dimension of $\Skew(\widetilde A_i, \widetilde \sigma)$, depending on whether the dimension of $(\bigwedge^2 \rho_i)^G$ is $0$ or $1$. 

    Observe that $\Skew(\widetilde A_i, \widetilde \sigma_i) = \mathfrak{so}(\widetilde{\rho_i^{rn}})^G$ by Lemma \ref{lem: Skew = so_A} (where $\widetilde{\rho_i^{rn}} = \rho_i^{rn} \oplus (\rho_i^{rn})^*$). By \cite[Formula (20.4)]{Fulton-Harris}, there is an isomorphism $\mathfrak{so}(\widetilde{\rho_i^{rn}}) \cong \bigwedge^2 \widetilde{\rho_i^{rn}}$, which is easily seen to be $G$-equivariant. Therefore there is an isomorphism $\mathfrak{so}(\widetilde{\rho_i^{rn}})^G \cong (\bigwedge^2 \widetilde{\rho_i^{rn}})^G$. In particular, we have an equality
    $$\dim \Skew(\widetilde A_i, \widetilde \sigma_i) = \dim \left(\extp^2 \widetilde{\rho_i^{rn}}\right)^G.$$
    Since $\rho_i$ is self-dual, we have $\dim (\rho_i \otimes \rho_i^*)^G = 1$, so the dimension of the right hand side is 
    \begin{equation*}
        \dim \left(\extp^2 \widetilde{\rho_i^{rn}}\right)^G = \begin{cases}
            \binom{2rn}{2} &\text{ if } (\bigwedge^2 \rho_i)^G = 0\\
            \binom{2rn}{2} + 2rn &\text{ if } (\bigwedge^2 \rho_i)^G = \mathbb{C}.
        \end{cases}
    \end{equation*}
    Then Lemma \ref{lem: involution type from skew elements} gives that $\widetilde \sigma_i$ is orthogonal if $(\bigwedge^2 \rho_i)^G = 0$ and symplectic if $(\bigwedge^2 \rho_i)^G = \mathbb{C}$.
\end{proof}

Now that we have computed the LLV algebra more explicitly, we can also prove the equivalence of the conditions \ref{condition *} and \ref{condition **} introduced in Chapter \ref{chap: introduction}.
\begin{lemma}
    \label{lem: * and ** equivalent}
    Let $X = (X_0 \times \prod_{i=1}^k X_i)/G$, where $X_0$ is an abelian variety, the $X_i$ for $i = 1, \dots, k$ are hyperkähler varieties and $G$ is a finite group whose action on $\prod_{i=0}^k X_i$ preserves a symplectic form. Then $X$ satisfies condition \ref{condition *} if and only if it satisfies condition \ref{condition **}.
\end{lemma}
\begin{proof}
    The equivalence of these two conditions follows from studying in which ways different simple factors arise. Write $\llv(X; \mathbb{Q}) = \mathfrak{h}_1 \times \dots \times \mathfrak{h}_m$ as a product of simple Lie algebras. 
    
    Let $\mathfrak{h}$ be one of the $\mathfrak{h}_\alpha$ for $1 \leq \alpha \leq m$, and first assume that $\mathfrak{h} \cong \mathfrak{so}(\widetilde \sH(X_i;\mathbb{Q})^{G_i})$ for some $i$. This is of type $B_n$ or $D_n$. Since $B_2 \cong C_2$ and $D_3 \cong A_3$, we see that condition \ref{condition *} says that $\mathfrak{h}$ is of type $B_n$ with $n \geq 3$ or of type $D_n$ with $n \geq 5$. This is equivalent to the requirement that $\dim \sH^2(X_i; \mathbb{Q})^{G_i}$ is $5$ or at least $7$, which is exactly what is demanded by condition \ref{condition **}.

    Now suppose we are in the other case. Then there is a $\rho \in \Irr_\mathbb{Q} (G)$ with associated simple factor $A$ of $\mathbb{Q}[G]$, and for $V$ the $\rho$-isotypic component of $\sH^1(X_0; \mathbb{Q})$ we have an isomorphism $\mathfrak{h} \cong \mathfrak{s}(\End_A(\widetilde V), \widetilde \sigma)$ (unless $\mathfrak{h}$ is of type $A_1$). Then for every irreducible constituent $\rho_i$ of $\rho_\mathbb{C}$, the multiplicity of $\rho_i$ in $\sH^1(X_0; \mathbb{C})$ is equal to the integer $l = \deg(\widetilde A)/2$ from Lemma~\ref{lem: type of s(A, sigma)}. Hence $\mathfrak{h}$ is of type $D_l, C_l$ or $A_{2l-1}$, depending on whether $\widetilde \sigma$ is orthogonal, symplectic or unitary. The equivalence of conditions \ref{condition *} and \ref{condition **} now follows from Lemma~\ref{lem: widetilde sigma type computation}.
\end{proof}

\section{Summary}
Over the course of Sections \ref{sec: LLV is HK times Ab}, \ref{sec: HK LLV} and \ref{sec: abelian LLV}, we computed the LLV algebra of an arbitrary holomorphic symplectic variety. We will now summarize the results obtained so far. 

Let $X = (\prod_{i=0}^k X_i)/G$ be as in Setup \ref{stp: symplectic var}. Let $I = \{1, \dots, k\}$, so that $X = (X_0 \times \prod_{i \in I} X_i)/G$. The action of $G$ on $X$ induces an action of $G$ on $I$, and for $i \in I$ we denote its $G$-orbit by $[i]$. Let $K$ be a field of characteristic~$0$. For $[i] \in I/G$, recall that $X_{[i]} = \prod_{j \in [i]} X_j$, and that $\SH_G(X_{[i]}; K)$ denotes the Verbitsky component of the graded ring $\sH^\bullet(X_{[i]}; K)^G$. We write $\llv_G(X_{[i]}; K)$ for the LLV algebra of $\sH^\bullet(X_{[i]}; K)^G$. We have seen in Proposition \ref{prop: LLV of HK} that $\llv_G(X_{[i]}; K) \cong \mathfrak{so}(\widetilde \sH(X_i; K)^{G_i})$, where $G_i$ is the stabilizer of $i$ for the action of $G$ on $I$, and $\widetilde \sH(X_i; K)$ is the Mukai lattice of $X_i$.

Furthermore, the action of $G$ on $\prod_{i=0}^k X_i$ turns $\sH^1(X_0; K) = \sH^1(\prod_i X_i; K)$ into a $G$-representation, which we denote by $V$. For $[\rho] \in \Irr_K(G)/\sim$, we have the $[\rho]$-isotypic component $V_{[\rho]} \subseteq V$. Associated to such a $[\rho]$-isotypic component, recall that $\SH_G(V_{[\rho]}; K)$ is the Verbitsky component of the graded ring $(\extp^\bullet V_{[\rho]})^G$ and $\llv_G(V_{[\rho]}; K)$ denotes the LLV algebra of this ring. By Theorem \ref{thm: LLV computation}, except for some small cases, $\llv_G(V_{[\rho]}; K)$ is isomorphic to $\mathfrak{s}(\End_A(\widetilde{V_{[\rho]}}), \widetilde \sigma) := \Skew(\End_A(\widetilde{V_{[\rho]}}), \widetilde \sigma)'$, where $\widetilde {V_{[\rho]}} = V_{[\rho]} \oplus V_{[\rho]}^*$, where $A$ is the simple factor (or product of two simple factors) of $K[G]$ associated to $[\rho]$, and where $\widetilde \sigma$ is the involution on $\End_A(\widetilde{V_{[\rho]}})$ associated to the bilinear form from Lemma~\ref{lem: bilin form on V tilde is equivariant}.

\begin{theorem}
    \label{thm: summary LLV computation}
    Let $X = (\prod_i X_i)/G$ be as in Setup \ref{stp: symplectic var} and let $K$ be a field of characteristic zero. Then there is an isomorphism
    \begin{equation}
        \label{eq: SH decomposition}
        \SH(X; K) \cong \bigotimes_{[i] \in I/G} \SH_G(X_{[i]}; K) \otimes \bigotimes_{[\rho] \in \Irr_K(G)/\sim} \SH_G(V_{[\rho]}; K),
    \end{equation}
    which induces an isomorphism
    \begin{equation}
        \label{eq: llv decomposition}
        \llv(X; K) \cong \prod_{[i] \in I/G}\llv_G(X_{[i]}; K) \times \prod_{[\rho] \in \Irr_K(G)/\sim} \llv_G(V_{[\rho]}; K),
    \end{equation}
    and this is a decomposition of $\llv(X; K)$ into simple factors. For every such simple factor, the associated factor in the factorization of the Verbitsky component~\eqref{eq: SH decomposition} is an irreducible representation.
    
    Furthermore, the simple factors of the LLV algebra can be computed as follows:
    \begin{enumerate}
        \item $\llv_G(X_{[i]}; K) \cong \mathfrak{so}(\widetilde \sH(X_i; K)^{G_i})$ for $[i] \in I/G$.
        \item $\llv_G(V_{[\rho]}; K) \cong \mathfrak{s}(\End_A(\widetilde{V_{[\rho]}}), \widetilde \sigma)$ for $[\rho] \in \Irr_K(G)/\sim$, unless $\widetilde \sigma$ is orthogonal and $\deg(\End_A(\widetilde{V_{[\rho]}})) = 4$.
        \item $\llv_G(V_{[\rho]}; K) \cong \mathfrak{sl}_2(F)$, where $F = Z(A)$, if $\widetilde \sigma$ is an orthogonal involution and $\deg(\widetilde A) = 4$.
    \end{enumerate}
\end{theorem}
\begin{proof}
    The product decomposition of $\llv(X; K)$ follows from combining Corollary \ref{cor: LLV HK and Ab},  Proposition \ref{prop: HK SH and LLV decomposition} and Proposition \ref{prop: Ab SH and LLV decomposition}. Similarly, the tensor product factorization of $\SH(X; K)$ follows from Propositions~\ref{prop: LLV is HK times Ab}, \ref{prop: HK SH and LLV decomposition} and \ref{prop: Ab SH and LLV decomposition}. Note that $\SH(X; K)$ is absolutely irreducible as an $\llv(X; K)$-representation by Corollary~\ref{cor: SH(X) is absolutely irreducible}, so also all the individual factors of $\SH(X; K)$ are irreducible.

    The computation of the individual factors of the LLV algebra is Proposition~\ref{prop: LLV of HK factor} for the factors of hyperkähler origin, and Theorem \ref{thm: LLV computation} for the factors of abelian origin. All these factors are indeed simple, so Equation \eqref{eq: llv decomposition} indeed gives a decomposition of $\llv(X; K)$ into simple factors.
\end{proof}

%% file: chapters/chapter5.tex
Let $X = (\prod_i X_i)/G$ be as in Setup \ref{stp: symplectic var}, and assume that $X$ satisfies condition~\ref{condition *}. Let $K$ be a field of characteristic~$0$. In the previous chapter, we computed the LLV algebra $\llv(X; K)$, and also decomposed the Verbitsky component $\SH(X; K)$ as a tensor product of smaller factors, see Theorem \ref{thm: summary LLV computation}. In this chapter, we will study these factors in more detail, and describe them as representations of the associated simple factors of the LLV algebra.

We saw in Theorem \ref{thm: summary LLV computation} that the factors of $\SH(X; \mathbb{Q})$ are of two types, either of the form $\SH_G(X_{[i]}; \mathbb{Q})$ for an orbit $[i] \in I/G$, or $\SH_G(V_{[\rho]}; \mathbb{Q})$ for an isotypic component $V_{[\rho]}$ of the $G$-representation $\sH^1(X_0; \mathbb{Q})$. We will study the factors of type $\SH_G(X_{[i]}; \mathbb{Q})$ in Section \ref{subsec: SH of HK factor}, and the factors of type $\SH_G(V_{[\rho]}; \mathbb{Q})$ in Section~\ref{subsec: SH of abelian factors}.

\section{The Verbitsky component of the factors of hyperkähler origin}
\label{subsec: SH of HK factor}
Let $X = (\prod_i X_i)/G$ be as in Setup \ref{stp: symplectic var}, satisfying condition \ref{condition *}. We start by analyzing the Verbitsky components $\SH_G(X_{[i]}; \mathbb{Q})$ associated to the hyperkähler factors. Recall that $\llv_G(X_{[i]}; \mathbb{Q}) \cong \mathfrak{so}(\widetilde \sH(X_i; \mathbb{Q})^{G_i})$, where $\widetilde \sH(X_i; \mathbb{Q}) = \mathbb{Q}\alpha \oplus \sH^2(X_i; \mathbb{Q}) \oplus \mathbb{Q}\beta$ is the Mukai lattice of the hyperkähler variety $X_i$. The computation of the Verbitsky component $\SH_G(X_{[i]}; \mathbb{Q})$ associated to $\llv_G(X_{[i]}; \mathbb{Q})$ generalizes the computation of \cite[Proposition 3.5]{TaelmanDerivedEquivalences}.

For any $D \in \mathbb{N}$, the space $\Sym^D\widetilde \sH(X_i; \mathbb{Q})^{G_i}$ is naturally a representation of $\mathfrak{so}(\widetilde \sH(X_{i}; \mathbb{Q})^{G_i})$, so it also becomes a representation of $\llv_G(X_{[i]}; \mathbb{Q})$ via the isomorphism $\llv_G(X_{[i]}; \mathbb{Q}) \cong \mathfrak{so}(\widetilde \sH(X_{[i]}; \mathbb{Q})^{G_i})$.

\begin{theorem}
    \label{thm: computation of SH of HK factor}
    Let $X = (\prod_{i=0}^k X_i)/G$ be as in Setup \ref{stp: symplectic var}. For every $i = 1, \dots, k$, define $d_i$ by $2d_i = \dim(X_i)$. Then for every $i = 1, \dots, k$ there is a unique injective homomorphism of $\llv_G(X_{[i]}; \mathbb{Q})$-modules 
    $$\SH_G(X_{[i]}; \mathbb{Q}) \to \Sym^{d_i \cdot \#[i]} \widetilde \sH(X_i; \mathbb{Q})^{G_i},$$
    which sends $1$ to $\alpha^{d_i \cdot \#[i]}/(d_i!)^{\#[i]}$.
\end{theorem}
\begin{proof}
     By definition, $\SH_G(X_{[i]}; \mathbb{Q})$ is the Verbitsky component of the algebra $(\bigotimes_{j \in [i]} \sH^\bullet(X_j; \mathbb{Q}))^G$, so we have an inclusion
     \begin{equation}
        \label{eq: SH [i] to otimes SH j}
         \SH_G(X_{[i]}; \mathbb{Q}) \injto \bigotimes_{j \in [i]} \SH(X_j; \mathbb{Q}).
     \end{equation}
    
    For every $j \in [i]$, we have an inclusion of $\llv(X_j; \mathbb{Q})$-modules $\SH(X_j; \mathbb{Q}) \to \Sym^{d_i} \widetilde \sH(X_j; \mathbb{Q})$ sending $1$ to $\alpha^{d_i}/d_i!$ by \cite[Proposition 3.5]{TaelmanDerivedEquivalences}, so we obtain an inclusion
    \begin{equation}
        \label{eq: tensor SH HK}
        \bigotimes_{j \in [i]} \SH(X_j; \mathbb{Q}) \to \bigotimes_{j \in [i]} \Sym^{d_i} \widetilde \sH(X_j; \mathbb{Q}),
    \end{equation}
    which is equivariant with respect to the action of $\prod_{j \in [i]} \mathfrak{so}(\widetilde \sH(X_j; \mathbb{Q}))$. 
    
    For every $j \in [i]$, choose a $g_j \in G$ such that $g_j$ maps $X_i$ to $X_j$, this induces an isomorphism $g_j^* \colon \sH^2(X_j; \mathbb{Q}) \to \sH^2(X_i; \mathbb{Q})$. Extending by the identity on $\mathbb{Q}\alpha \oplus \mathbb{Q}\beta$ gives an isomorphism $\widetilde \sH(X_j; \mathbb{Q}) \to \widetilde \sH(X_i; \mathbb{Q})$. These isomorphisms give us a (non-canonical) isomorphism 
    \begin{equation}
        \label{eq: Sym homogeneization}
        \bigotimes_{j \in [i]} \Sym^{d_i} \widetilde \sH(X_j; \mathbb{Q}) \isomto \bigotimes_{j \in [i]} \Sym^{d_i} \widetilde \sH(X_i; \mathbb{Q}).
    \end{equation}
    Finally, there is a natural map
    \begin{equation}
        \label{eq: tensor sym quotient}
        \bigotimes_{j \in [i]} \Sym^{d_i} \widetilde \sH(X_i; \mathbb{Q}) \to \Sym^{d_i \cdot \#[i]} \widetilde \sH(X_i; \mathbb{Q}).
    \end{equation}
    By composing the maps \eqref{eq: SH [i] to otimes SH j}, \eqref{eq: tensor SH HK}, \eqref{eq: Sym homogeneization} and \eqref{eq: tensor sym quotient}, we obtain an $\llv_G(X_{[i]}; \mathbb{Q})$-equivariant map
    \begin{equation}
        \label{eq: SH [i] to Sym di[i]}
        \SH_G(X_{[i]}; \mathbb{Q}) \to \Sym^{d_i \cdot \#[i]} \widetilde \sH(X_i; \mathbb{Q}).
    \end{equation}
    By using the description of the maps $\SH(X_j; \mathbb{Q}) \to \Sym^{d_j} \widetilde \sH(X_j; \mathbb{Q})$ from \cite[Proposition 3.5]{TaelmanDerivedEquivalences}, we see that the map \eqref{eq: SH [i] to Sym di[i]} sends $1$ to $\alpha^{d_i\cdot \#[i]}/(d_i!)^{\#[i]}$, and that its image is contained in $\Sym^{d_i \cdot \#[i]} \widetilde \sH(X_i; \mathbb{Q})^{G_i}$. It remains to show that this map is injective. The map is nonzero, since $1 \mapsto \alpha^{d_i\cdot \#[i]}/(d_i!)^{\#[i]}$. Since $\SH_G(X_{[i]}; \mathbb{Q})$ is an absolutely irreducible $\llv_G(X_{[i]}; \mathbb{Q})$-representation by Corollary \ref{cor: SH(X) is absolutely irreducible}, Schurs lemma implies that this map must be injective.

    Lastly, the absolute irreducibility of $\SH_G(X_{[i]}; \mathbb{Q})$ and Schurs lemma imply that the the morphism $\SH_G(X_{[i]}; \mathbb{Q}) \to \Sym^{d_i\cdot\#[i]}\widetilde \sH(X_i; \mathbb{Q})^{G_i}$ is uniquely determined up to scalars, so it is uniquely determined by the extra assumption that $1$ gets mapped to $\alpha^{d_i \cdot \#[i]}/(d_i!)^{\#[i]}$.
\end{proof}

Write $b$ for the bilinear form on $\widetilde \sH(X_i; \mathbb{Q})^{G_i}$. Let $D = d_i\cdot\#[i]$ and consider the contraction operator
$$\Delta \colon \Sym^D \widetilde \sH(X_i; \mathbb{Q})^{G_i} \to \Sym^{D - 2} \widetilde \sH(X_i; \mathbb{Q})^{G_i},$$
which sends $x_1 \cdots x_D$ to $\sum_{k < l}b(x_k, x_l) x_1 \cdots \hat x_k \cdots \hat x_l \cdots x_D$.

\begin{lemma}
    \label{lem: SH is S_{[d]}}
    There is a short exact sequence of $\llv_G(X_{[i]}; \mathbb{Q})$-representations
    \begin{equation*}
    \begin{tikzcd}
        0 \arrow[r] &\SH_G(X_{[i]}; \mathbb{Q}) \arrow[r] &\Sym^{D} \widetilde \sH(X_i)^{G_i} \arrow[r, "\Delta"] &\Sym^{D - 2} \widetilde \sH(X_i)^{G_i} \arrow[r] & 0.
    \end{tikzcd}
    \end{equation*}
\end{lemma}
\begin{proof}
    The proof is exactly the same as \cite[Lemma 3.7]{TaelmanDerivedEquivalences}.
\end{proof}

\section{The Verbitsky component of the factors of abelian origin}
\label{subsec: SH of abelian factors}
We will now compute the Verbitsky component $\SH_G(V_{[\rho]}; \mathbb{C})$ associated to simple factors of $\llv(X; \mathbb{C})$ of the type $\llv_G(V_{[\rho]}; \mathbb{C})$. We will do this by computing the highest weight of $\SH_G(V_{[\rho]}; \mathbb{C})$ (this is the reason why we work over $\mathbb{C}$ instead of $\mathbb{Q}$ in this section). We will only need this computation in the cases where the Lie algebra $\llv_G(V_{[\rho]}; \mathbb{C})$ is of type $A_n$ or $D_n$ for some $n > 1$, but the computation for type $C_n$ is similar and can also be performed. We start with some technical preliminaries.

\subsection{Representation theory of Lie subalgebras}
We start by setting up some notation. Suppose we have an inclusion $\mathfrak{g}_1 \to \mathfrak{g}_2$ of semisimple Lie algebras over an algebraically closed field $k$ of characteristic~$0$. Let $\mathfrak{h}_2$ be a Cartan subalgebra of $\mathfrak{g}_2$, and let $\mathfrak{h}_1 = \mathfrak{g}_1 \cap \mathfrak{h}_2$. Assume that $\mathfrak{h}_1$ is a Cartan subalgebra of $\mathfrak{g}_1$. The inclusion $\mathfrak{h}_1 \to \mathfrak{h}_2$ induces a surjection $f \colon \mathfrak{h}_2^* \to \mathfrak{h}_1^*$, given by $f(\beta) = \beta|_{\mathfrak{h}_1}$ for $\beta \in \mathfrak{h}_2^*$. Let $R_1$ and $R_2$ be the root systems of $\mathfrak{g}_1$ and $\mathfrak{g}_2$ respectively, and choose systems of positive roots $R_1 = R_1^+ \cup R_1^-$ and $R_2 = R_2^+ \cup R_2^-$. Assume that these choices are compatible, in the sense that $f^{-1}(R_1^+) \cap R_2 \subseteq R_2^+$ (explicitly, this means that for every $\beta \in R_2$ with $f(\beta) \in R_1^+$, we have $\beta \in R_2^+$).

From $R_1$ and $R_2$ we obtain root space decompositions $\mathfrak{g}_1 = \mathfrak{h}_1 \oplus \bigoplus_{\alpha \in R_1}\mathfrak{g}_{1, \alpha}$ and $\mathfrak{g}_2 = \mathfrak{h}_2 \oplus \bigoplus_{\beta \in R_2} \mathfrak{g}_{2, \beta}$. We furthermore consider the positive subspaces $\mathfrak{g}_1^+ := \oplus_{\alpha \in R_1^+} \mathfrak{g}_{1, \alpha}$ and $\mathfrak{g}_2^+ := \oplus_{\beta \in R_2^+} \mathfrak{g}_{2, \beta}$.

\begin{lemma}
    \label{lem: positive subspaces inclusion}
    Under the above assumptions, we have an inclusion $\mathfrak{g}_1^+ \subseteq \mathfrak{g}_2^+$.
\end{lemma}
\begin{proof}
    Take an $\alpha \in R_1^+$ and an $X \in \mathfrak{g}_{1, \alpha}$. We are done if we show that $X \in \mathfrak{g}_2^+$. As an element of $\mathfrak{g}_2$, decompose $X = h + \sum_{\beta \in R_2} X_\beta$ with $h \in \mathfrak{h}_2$ and each $X_\beta \in \mathfrak{g}_{2, \beta}$.

    For any $H \in \mathfrak{h}_1$, we now have:
    $$\alpha(H)X = [H,X] = [H, h] + \sum_\beta [H,X_\beta] = \sum_\beta \beta(H)X_\beta,$$
    where $[H, h] = 0$ since $\mathfrak{h}_1 \subseteq \mathfrak{h}_2$ and $\mathfrak{h}_2$ is abelian. This implies that $\alpha(H) = \beta(H)$ for all $\beta$ with $X_\beta \neq 0$, so $h = 0$ and $\alpha = \beta|_{\mathfrak{h}_1} = f(\beta)$ for all $\beta \in R_2$ with $X_\beta \neq 0$. In particular, if $X_\beta \neq 0$, then $\beta \in f^{-1}(R_1^+)$. The assumption $f^{-1}(R_1^+) \cap R_2 \subseteq R_2^+$ then gives that $\beta \in R_2^+$ for all $\beta \in R_2$ with $X_\beta \neq 0$, and this implies that $X \in \mathfrak{g}_2^+$. 
\end{proof}

Now, let $V$ be an irreducible representation of $\mathfrak{g}_2$ with highest weight vector $v$ of weight $\alpha \in R_2$. By restriction, $V$ becomes a representation of $\mathfrak{g}_1$, and we can consider the subrepresentation $\mathfrak{g}_1v \subseteq V$ obtained by repeatedly applying elements from $\mathfrak{g}_1$ to $v$.  

\begin{lemma}
    \label{lem: weight of rep of sub lie alg}
    Under the above assumptions, the $\mathfrak{g}_1$-representation $\mathfrak{g}_1v \subseteq V$ is irreducible and has $v$ as a highest weight vector with weight $f(\beta)$. In particular, $\mathfrak{g}_1v$ is the irreducible representation of $\mathfrak{g}_1$ corresponding to the weight $f(\beta)$.
\end{lemma}
\begin{proof}
    To show that $v$ is a highest weight vector with weight $f(\beta)$, we must show that $v$ has weight $f(\beta)$ and that $\mathfrak{g}_1^+v = 0$.

    To see that $v$ has weight $f(\beta)$, observe that for any $H \in \mathfrak{h}_1$, we have 
    \[Hv = \beta(H)v = \beta|_{\mathfrak{h}_1}(H)v = f(\beta)(H)v,\] 
    because $v$ has weight $\beta$ in the $\mathfrak{g}_2$-representation~$V$.

    We already know that $\mathfrak{g}_2^+v = 0$, since $v$ is a highest weight vector for the $\mathfrak{g}_2$-representation $V$. The inclusion $\mathfrak{g}_1^+ \subseteq \mathfrak{g}_2^+$ from Lemma \ref{lem: positive subspaces inclusion} then implies that also $\mathfrak{g}_1^+v = 0$. 

    The representation $\mathfrak{g}_1v$ is irreducible by \cite[Proposition 14.13(ii)]{Fulton-Harris}.
\end{proof}

\subsection{Computing the Verbitsky component}
Let $X$ be as in Setup \ref{stp: symplectic var}, and assume that $X$ satisfies condition~\ref{condition *}. Now suppose we have a simple factor $\mathfrak{g}$ of $\llv_{\Ab}(X; \mathbb{C})$ which is not of type $A_1$, so there is a $[\rho] \in \Irr_\mathbb{C}(G)$ such that $\mathfrak{g}$ is of the form $\mathfrak{s}(\End_A(\widetilde{V_{[\rho]}}), \widetilde \sigma)$ with $V = \sH^1(X_0; \mathbb{C})$. We wish to understand the associated representation $\SH_G(V_{[\rho]}; \mathbb{C})$ of $\mathfrak{g}$.

Recall that Proposition \ref{prop: Ab SH and LLV decomposition} gives a decomposition 
$$\llv_{\Ab}(X; \mathbb{C}) \cong \prod_{[\rho] \in \Irr_\mathbb{C}(G)/\sim} \llv_G(V_{[\rho]}; \mathbb{C})$$
and a tensor product factorization
$$\SH_{\Ab}(X; \mathbb{C}) \cong \bigotimes_{[\rho] \in \Irr_\mathbb{C}(G)/\sim} \SH_G(V_{[\rho]}; \mathbb{C}).$$

Pick a class $[\rho] \in \Irr_\mathbb{C}(G)/\sim$, then its associated Verbitsky component $\SH_G(V_{[\rho]}; \mathbb{C})$ is defined as $\SH((\extp^\bullet V_{[\rho]})^G)$, and this is a subalgebra of $\extp^\bullet V_{[\rho]}$. Furthermore, $\llv_G(V_{[\rho]}; \mathbb{C})$ is defined as $\llv((\extp^\bullet V_{[\rho]})^G)$. We already know that $\llv(\extp^\bullet V_{[\rho]}) = \mathfrak{so}(\widetilde {V_{[\rho]}})$, and the Verbitsky component of $\extp^\bullet V_{[\rho]}$, seen as a representation of $\llv(\extp^\bullet V_{[\rho]})$, is the even spinor representation \cite[Section 3]{LooijengaLunts}. We will now use Lemma \ref{lem: weight of rep of sub lie alg} to compute the weight of a highest weight vector of $\SH_G(V_{[\rho]}; \mathbb{C})$.

\begin{proposition}
    \label{prop: weight of Dn SH}
    Take a self-dual $\rho \in \Irr_\mathbb{C}(G)$ with $\llv_G(V_{[\rho]}; \mathbb{C})$ of type $D_n$. Let $\omega'$ be the fundamental root corresponding to the even spinor representation of $\llv_G(V_{[\rho]}; \mathbb{C})$, and let $k = \dim \rho$. Then $\SH_G(V_{[\rho]}; \mathbb{C})$ is the representation of highest weight $k\omega'$.
\end{proposition}
Recall that $[\rho] = \{\rho\}$ since $\rho$ is self-dual.
\begin{proof}
    Theorem \ref{thm: LLV computation} implies that $\llv_G(V_{[\rho]}; \mathbb{C}) = \mathfrak{so}(\widetilde{V_{[\rho]}})^G$. Since $\rho$ is self-dual, we can write $V_{[\rho]} = V_\rho \cong \rho^n$. 

   Now choose a basis of $V_{[\rho]}$, where the first $k$ basis elements generate the first copy of $\rho$, the second $k$ basis elements generate the second copy of $\rho$ and so forth. We extend this to a basis of $\widetilde{V_{[\rho]}}$ by letting the first $kn$ basis elements be the basis we constructed for $V_{[\rho]}$, and letting the second $kn$ basis elements be its dual basis. For $1 \leq i,j \leq 2kn$, let $E_{i, j}$ be the matrix with a $1$ in position $(i, j)$ and zeroes elsewhere. For $1 \leq i \leq kn$, let $H_i = E_{i,i} - E_{kn+i, kn+i}$. Then $\mathfrak{h} = \langle H_1, \dots, H_{kn} \rangle$ is a Cartan subalgebra of $\mathfrak{so}(\widetilde{V_{[\rho]}})$. Write $L_1, \dots, L_{kn}$ for the dual basis of $\mathfrak{h}^*$. Note that for every $1 \leq j \leq n$, the sum $H'_{j} := H_{k(j-1) + 1} + \dots + H_{k(j-1) + k}$ is $G$-equivariant, and hence gives an element of $\mathfrak{h}^G$, which is a Cartan subalgebra of $\llv_G(V_{[\rho]}; \mathbb{C})$. Let $\{L_1', \dots, L_n'\}$ be the induced dual basis of $(\mathfrak{h}^G)^*$. 

    The inclusion $\mathfrak{h}^G \subseteq \mathfrak{h}$ induces a restriction map $\mathfrak{h}^* \to (\mathfrak{h}^G)^*$, which sends the first $k$ basis elements $L_1, \dots, L_k$ to $L_1'$, the second $k$ to $L_2'$ and so forth. If we write $\omega = \frac{1}{2}(L_1 + \dots + L_{kn})$ for the highest weight of the even spinor representation of $\mathfrak{so}(\widetilde V)$, and let $\omega' = \frac{1}{2}(L_1' + \dots + L_n')$, then we see that the restriction map $\mathfrak{h}^* \to (\mathfrak{h}^G)^*$ sends $\omega$ to $k\omega'$.

    The set of roots of $\mathfrak{so}(\widetilde{V_{[\rho]}})$ is equal to $R_2 = \{\pm L_i \pm L_j\}_{i \neq j}$. We choose $R_2^+ = \{L_i \pm L_j\}_{i < j}$ as the collection of positive roots. The space of invariants $\mathfrak{so}(\widetilde{V_{[\rho]}})^G$ has root system $R_1 = \{\pm L_i' \pm L_j'\}_{i \neq j}$, and we choose the set $R_1^+ = \{L_i'\pm L_j'\}_{i < j}$ as positive roots. Then the restriction map $f \colon \mathfrak{h}^* \to (\mathfrak{h}^G)^*$ satisfies $f^{-1}(R_1^+) \cap R_2 \subseteq R_2^+$. 
    
    The highest weight space of the even spinor representation $\extp^{\ev} V_{[\rho]}$ is $\extp^{kn} V_{[\rho]}$, so we can choose a nonzero vector $v \in \extp^{kn} V_{[\rho]}$ as a highest weight vector. We now apply Lemma \ref{lem: weight of rep of sub lie alg} to this situation, and conclude that the subspace $\mathfrak{so}(\widetilde{V_{[\rho]}})^Gv \subseteq \extp^\bullet V_{[\rho]}$ is an irreducible subrepresentation with highest weight $\omega|_{(\mathfrak{h}^G)^*} = k\omega'$. 

    To finish the proof, we need to see that $\mathfrak{so}(\widetilde{V_{[\rho]}})^G\cdot v = \SH_G(V_{[\rho]}; \mathbb{C})$. To see this, observe that both these spaces are irreducible representations of $\mathfrak{so}(\widetilde{V_{[\rho]}})^G$ containing $\extp^{kn}V_{[\rho]}$, so they must be equal.
\end{proof}

We will now look at the case where $\llv_G(V_{[\rho]}; \mathbb{C})$ comes from a unitary involution and is of type $A_{2l-1}$ for some integer~$l$. This result is not needed for the proof of Theorem \ref{thm: main theorem}, but only for a generalization of Theorem \ref{thm: main theorem} in Chapter \ref{sec: low degree cases}. In this case, $\llv_G(V_{[\rho]}; \mathbb{C}) \cong \mathfrak{sl}_{2l}(\mathbb{C})$. To state the next proposition, we need to recall the construction of the weight lattice of this Lie algebra, see \cite[Chapter 15]{Fulton-Harris} for more details. We take the Cartan subalgebra $\mathfrak{h}$ consisting of the diagonal matrices of trace zero in $\mathfrak{sl}_{2l}(\mathbb{C})$. Let $L_i \colon \mathfrak{sl}_{2l}(\mathbb{C}) \to \mathbb{C}$ be defined by sending $E_{i,i}$ to $1$ and all other elementary matrices to $0$. Then $\mathfrak{h}^* = \langle L_1, \dots, L_{2l} \rangle / \langle L_1 + \dots + L_{2l}\rangle$. The root lattice is $R = \{L_i - L_j\}_{i \neq j}$, with standard choice of positive roots $\{L_i - L_j\}_{i < j}$. The fundamental weights are then $\omega_1, \dots, \omega_{2l-1}$, where $\omega_i = L_1 + \dots + L_i$, and these can be used to describe the irreducible representations.

\begin{proposition}
    \label{prop: weight of An SH}
    Let $[\rho] \in \Irr_\mathbb{C}(G)/\sim$ with $\rho$ not self-dual and $\llv_G(V_{[\rho]}; \mathbb{C})$ of type $A_{2l-1}$. Let $\omega_1, \dots, \omega_{2l-1}$ be the fundamental weights of $\mathfrak{sl}_{2l}(\mathbb{C})$ as described above. Then $\SH_G(V_{[\rho]}; \mathbb{C})$ is the representation of highest weight $k \omega_{l}$, where $k = \dim \rho$.
\end{proposition}
Recall that the Dynkin diagram of $A_{2l-1}$ consists of $2l-1$ points on a line. It has one nontrivial symmetry, namely by flipping the whole diagram. The weight $\omega_{l}$ sits exactly in the middle of the Dynkin diagram, and is therefore invariant under this symmetry.
\begin{proof}
    Since $\rho$ is not self-dual, the algebra $\End(\widetilde{V_{[\rho]}})^G$ has an involution of the second kind. Because $[\rho] = \{\rho, \rho^*\}$ consists of two elements, we have $V_{[\rho]} = V_\rho \oplus V_{\rho^*}$, and $\widetilde {V_{[\rho]}} = V_{\rho} \oplus V_{\rho^*} \oplus V_{\rho}^* \oplus V_{\rho^*}^*$. By Schurs lemma, 
    $$\End(\widetilde{V_{[\rho]}})^G \cong \End(V_{\rho} \oplus V_{\rho^*}^*)^G \times \End(V_{\rho^*} \oplus V_\rho^*)^G \cong M_{2l}(\mathbb{C}) \times M_{2l}(\mathbb{C}),$$
    where $l$ is the natural number with $V_{\rho} \cong \rho^l$. Choose a basis $\{e_1, \dots, e_k\}$ of $\rho$. By taking $l$ copies of this basis and $l$ copies of its dual basis, we obtain a basis $\{e_1, \dots, e_{2kl}\}$ of $V_{[\rho]}$. We can extend this with its dual basis to obtain a basis $\{e_1, \dots, e_{4kl}\}$ of $\widetilde{V_{[\rho]}}$. We can represent elements of $\mathfrak{sl}(V_\rho \oplus V_{\rho^*}^*)^G$ by block matrices $\begin{psmallmatrix} A & B \\ C & D\end{psmallmatrix}$ coming from the direct sum decomposition, so $A \colon V_\rho \to V_\rho$ and $D \colon V_{\rho^*}^* \to V_{\rho^*}^*$ for example. Using the direct sum decomposition $\widetilde {V_{[\rho]}} = V_{\rho} \oplus V_{\rho^*} \oplus V_{\rho}^* \oplus V_{\rho^*}^*$, elements of $\End(\widetilde{V_{[\rho]}})$ can be represented as four by four block matrices, and we have an inclusion $j \colon \mathfrak{sl}(V_\rho \oplus V_{\rho^*}^*)^G \to \mathfrak{so}(\widetilde{V_{[\rho]}}) \subseteq \End(\widetilde{V_{[\rho]}})$ given by
    \begin{equation*}
        \begin{pmatrix}
            A & B \\ C & D
        \end{pmatrix}
        \mapsto
        \begin{pmatrix}
            A & 0 & 0 & B \\
            0 & -D^t & -B^t & 0\\
            0 & -C^t & -A^t & 0\\
            C & 0 & 0 & D
        \end{pmatrix}\lowdot
    \end{equation*}
    The image of this inclusion is exactly $\llv_G(V_{[\rho]}; \mathbb{C})$, so we have an isomorphism $\llv_G(V_{[\rho]}; \mathbb{C}) \cong \mathfrak{sl}_{2l}(\mathbb{C})$, which is of type $A_{2l-1}$. For $1 \leq i \leq 2kl$, let $H_i = E_{i,i} - E_{2kl+i, 2kl+i} \in \mathfrak{so}(\widetilde{V_{[\rho]}})$. Then a Cartan subalgebra of $\mathfrak{so}(\widetilde{V_{[\rho]}})$ is given by $\mathfrak{h}_2 = \langle H_1, \dots, H_{2kl}\rangle$. In $\mathfrak{sl}(V_\rho \oplus V_{\rho^*}^*)$, the standard Cartan subalgebra is given by $$\left\{a_1 E_{1,1} + \dots + a_{2kl}E_{2kl} : \sum a_i = 0\right\}.$$
    We define $\mathfrak{h}_1$ to be the $G$-invariants of this, so if we let 
    $$E_{i, i}' = E_{k(i-1)+1, k(i-1)+1} + \dots + E_{k(i-1)+k, k(i-1)+k}$$
    for $1 \leq i \leq l$, then
    $$\mathfrak{h}_1 = \left\{a_1 E'_{1,1} + \dots + a_{2l}E'_{2l, 2l}: \sum a_i = 0\right\}.$$
    We see that $j(\mathfrak{h}_1) \subseteq \mathfrak{h}_2$. We write $L_i \in \mathfrak{h}_2^*$ for the element with $L_i(H_j) = \delta_{ij}$. The root system of $\mathfrak{so}(\widetilde{V_{[\rho]}})$ is $R_2 = \{\pm L_i \pm L_j\}_{i \neq j}$, and we take the positive roots $R_2^+ = \{L_i \pm L_j\}_{i < j}$. To get a root system for $\mathfrak{sl}(V_\rho \oplus V_{\rho^*}^*)^G$, we observe that the dual to $\mathfrak{h}_1$ is the quotient
    $$\mathfrak{h}_1^* = \langle L_1', \dots, L_{2l}'\rangle / \langle L_1' + \dots + L'_{2l} = 0\rangle,$$
    where $L_i'$ is defined by $L_i'(E'_{jj}) = \delta_{ij}$. Then the root system for $\mathfrak{sl}(V_\rho \oplus V_{\rho^*}^*)^G$ is given by $R_1 = \{L_i' - L_j'\}_{i \neq j}$.

    Let $\varphi \colon \{1, \dots, 2kl\} \to \{1, \dots, 2l\}$ be the function sending the first $k$ numbers to~$1$, the second $k$ numbers to $2$ etc., so $\varphi(i) = \left \lceil \frac{i}{k} \right \rceil$.
    The restriction map $f \colon \mathfrak{h}_2^* \to \mathfrak{h}_1^*$ is then given by 
    \begin{equation*}
        f(L_i) = \begin{cases}
            L'_{\varphi(i)} \text{ if } \varphi(i) \leq n\\
            -L'_{\varphi(i)} \text{ if } \varphi(i) > n.
        \end{cases}    
    \end{equation*}
    Choose a general collection of real numbers $c_1 > c_2 > \dots > c_l > c_{2l} > c_{2l-1} > \dots > c_{l+1}$ summing to zero. This gives a functional $\lambda \colon \mathfrak{h}_1^* \to \mathbb{R}$ defined by $\lambda(\sum a_i L'_i) = \sum a_ic_i$, which induces a division of $R_1$ into positive and negative roots. The associated set of positive roots is 
    $$R_1^+ = \{L'_i - L'_j : i < j \text{ and } i \leq l\} \cup \{L'_j - L'_i : l < i < j\}.$$ 
    This satisfies $f^{-1}(R_1^+) \cap R_2 \subseteq R_2^+$. Just as in the proof of Proposition \ref{prop: weight of Dn SH}, we see that any nonzero vector $v \in \extp^{2kl}V_{[\rho]}$ is a highest weight vector of the representation $\mathfrak{so}(\widetilde{V_{[\rho]}}) \acts \extp^\bullet V_{[\rho]}$, it has weight $\omega = \frac{1}{2}\sum_{i = 1}^{2kl} L_i$. We will now apply Lemma \ref{lem: weight of rep of sub lie alg} to this situation. We see that the subrepresentation $\mathfrak{sl}(V_\rho \oplus V_{\rho^*}^*)^Gv$ of $\extp^\bullet V_{[\rho]}$ generated by $v$ and $\mathfrak{sl}(V_\rho \oplus V_{\rho^*}^*)^G$ is irreducible of weight $f(\omega)$. This subrepresentation is equal to $\SH_G(V_{[\rho]}; \mathbb{C})$, because they are both irreducible representations of $\mathfrak{sl}(V_\rho \oplus V_{\rho^*}^*)^G$ and have a non-empty intersection.

    We have
    $$f(\omega) = \frac{1}{2}\sum_{i=1}^{2kl} f(L_i) = \frac{k}{2}\sum_{i=1}^l L'_i - \frac{k}{2}\sum_{i=l+1}^{2l}L'_i.$$
    By using the relation $L'_1 + \dots + L'_{2l} = 0$ to eliminate $L'_{2l}$ in the above, we see that this simplifies to
    $$f(\omega) = \frac{k}{2} \sum_{i=1}^l L'_i - \frac{k}{2}\sum_{i= l+1}^{2l-1}L'_i + \frac{k}{2}\sum_{i=1}^{2l-1} L'_i = k\sum_{i=1}^l L'_i.$$

    Our chosen ordering of the roots in $R_1$ is not the standard one often found in the literature (e.g. as in \cite[Chapter 15]{Fulton-Harris}). The standard weights are given by $\omega_1, \dots, \omega_{2l-1}$, where $\omega_i = L'_1 + \dots + L'_i$. The Weyl group of $\mathfrak{sl}(V_\rho \oplus V_{\rho^*}^*)^G$ is the symmetric group $S_{2l}$, and this acts on $\mathfrak{h}_1^*$ by just permuting the elements of the set $\{L'_1, \dots, L'_{2l}\}$. The permutation that will turn our choice of positive roots into the standard one is given by leaving $1, \dots, l$ fixed, and reversing the sequence $l+1, \dots, 2l$. We see that this brings our weight $f(\omega) = k(L'_1 + \dots + L'_l)$ to the weight $k\omega_l$, which finishes the proof.
\end{proof}

\section{Hodge theory of the Verbitsky component}
The existence of a symplectic form on a variety $X$ implies that $h^{2d,0}(X) \neq 0$, where $2d = \dim(X)$. We will need a similar version of this result also for the Verbitsky components associated to the simple factors of $\llv(X; \mathbb{Q})$.

\begin{proposition}
    \label{prop: HS of SH_i is wide}
    Let $X$ be a smooth projective variety over $\mathbb{C}$ admitting a holomorphic symplectic form. Write $\llv(X; \mathbb{Q}) = \prod_\alpha \mathfrak{g}_\alpha$ and $\SH(X; \mathbb{Q}) = \bigotimes_\alpha \SH_\alpha(X)$ for the decompositions from Theorem \ref{thm: summary LLV computation}. For $j \in \mathbb{N}$, let $\SH_\alpha^j(X)$ be the degree $j$ part of $\SH_\alpha(X)$. Let $d_\alpha \in \mathbb{N}$ be maximal with $\SH_\alpha^{2d_\alpha}(X) \neq 0$. Then the Hodge structure on $\SH_\alpha(X)$ satisfies $\SH_\alpha(X)^{d_\alpha, 0} \neq 0$.
\end{proposition}
\begin{proof}
    Let $\sigma \in \sH^2(X; \mathbb{C})$ be a symplectic form on $X$. Then $\sigma^d$ is a nonzero element of $\sH^{2d}(X; \mathbb{C})$. Since $\sigma^d$ is of Hodge type $(2d, 0)$ and is in $\SH(X)$, we see that $\SH(X)^{2d,0} \neq 0$.

    Since $\SH_\alpha(X)$ is generated by elements in degree $2$, which are of Hodge type either $(2,0), (1,1)$ or $(0, 2)$, it follows that $\SH_\alpha(X)^{p,0} = 0$ if $p > d_\alpha$. Therefore $\SH(X)^{2d,0} = \bigotimes_\alpha \SH_\alpha(X)^{d_\alpha, 0}$, and this implies that all $\SH_\alpha(X)^{d_\alpha,0}$ are nonzero.
\end{proof}

We also have the following result, whose second part generalizes \cite[Theorem~4.4]{TaelmanDerivedEquivalences}:
\begin{theorem}
    \label{thm: LLV alg derived invariant}
    Let $X$ and $Y$ be smooth projective varieties over $\mathbb{C}$ carrying a holomorphic symplectic form, and let $\Phi \colon D^b(X) \to D^b(Y)$ be a triangulated equivalence. Then there is an isomorphism of Lie algebras $\Phi^{\llv} \colon \llv(X; \mathbb{Q}) \to \llv(Y; \mathbb{Q})$ so that $\Phi^{\sH} \colon \sH^\bullet(X;\mathbb{Q}) \to \sH^\bullet(Y; \mathbb{Q})$ is equivariant with respect to $\Phi^{\llv}$. Furthermore, the map $\Phi^{\sH}$ restricts to an isomorphism of Verbitsky components \index{$\Phi^{\SH}$}$\Phi^{\SH} \colon \SH(X; \mathbb{Q}) \to \SH(Y; \mathbb{Q})$.
\end{theorem}
\begin{proof}
    The existence of $\Phi^{\llv}$ and the proof that $\Phi^{\sH}$ is equivariant with respect to $\Phi^{\llv}$ is \cite[Theorem A]{TaelmanDerivedEquivalences}. To see that $\Phi^{\sH}$ restricts to an isomorphism of the Verbitsky components, we can use the same proof as \cite[Theorem 4.4]{TaelmanDerivedEquivalences}. Recall that $\SH(X; \mathbb{Q})$ is an irreducible $\llv(X; \mathbb{Q})$-representation by Lemma~\ref{lem: SH(R) absolutely irreducible}. Let $d = \dim(X)$, and let $\omega$ be a symplectic form on $X$. Then multiplication with $\omega^{d/2}$ induces an isomorphism $\mathcal{O}_X \cong \Omega_X^d$, so $\sH^0(X, \Omega_X^d) \cong \mathbb{C}$. We have $\omega^{d/2} \in \SH(X; \mathbb{Q})$, so the Verbitsky component is the unique irreducible subrepresentation of $\sH^\bullet(X;\mathbb{Q})$ that contains the classes of Hodge type $(d,0)$ and $(0,d)$. Since $\Phi^{\sH}$ respects the columns of the Hodge diamond and is equivariant with respect to $\Phi^{\llv}$, we see that it must send $\SH(X; \mathbb{Q})$ to $\SH(Y; \mathbb{Q})$.
\end{proof}

\section{Intermezzo: alternative proof of a result by Huybrechts and Nieper-Wisskirchen}
\label{subsec: intermezzo, Huybrechts Nieper Wisskirchen}
In this section, we use our computations of the LLV algebra and the derived invariance of the LLV algebra to give an alternative proof of (a slight generalization of) Theorem \ref{thm: HNW derived invariance products}, which was originally proven by Huybrechts and Nieper-Wisskirchen. It should be noted that both their and our argument in the end boil down to the derived invariance of Hochschild (co)homology. The argument by Huybrechts and Nieper-Wisskirchen uses this rather directly, while in our proof this fact is used by Taelman in the proof of Theorem \ref{thm: Taelman theorem A}.

\begin{theorem}
    Let $X_0$ and $Y_0$ be abelian varieties, let $X_1, \dots, X_m$ and $Y_1, \dots, Y_n$ be hyperkähler varieties, and assume that $D^b(\prod_{i=0}^m X_i) \simeq D^b(\prod_{j=0}^n Y_j)$. Then  $\dim X_0 = \dim Y_0$, and there is a bijection $\sigma \colon \{1, \dots, m\} \to \{1, \dots, n\}$ such that $\dim X_i = \dim Y_{\sigma(i)}$ and $b_2(X_i) = b_2(Y_{\sigma(i)})$ for all $i$.
\end{theorem}
\begin{remark}
    The following proof becomes quite a bit shorter if we do not have the abelian factors $X_0$ and $Y_0$.
\end{remark}

\begin{proof}
    Write $X = \prod_{i=0}^m X_i$ and $Y = \prod_{j=0}^n Y_j$. By the Künneth formula, we have $\dim \sH^1(X; \mathbb{Q}) = 2 \dim X_0$, so \cite[Corollary B]{PopaSchnell} implies that $X_0$ and $Y_0$ have the same dimension. If they both have even dimension, the varieties $X$ and $Y$ admit a holomorphic symplectic form. If not, we can take an elliptic curve $E$ and replace $X$ by $X \times E$ and $Y$ by $Y \times E$. Once we have obtained an isomorphism $\sH^\bullet(X; \mathbb{Q}) \otimes \sH^\bullet(E; \mathbb{Q}) \cong \sH^\bullet(Y; \mathbb{Q}) \otimes \sH^\bullet(E; \mathbb{Q})$ which preservers the Hodge bigrading, we can then construct an isomorphism $\sH^\bullet(X; \mathbb{Q}) \cong \sH^\bullet(Y; \mathbb{Q})$ which preserves the Hodge bigrading. Therefore we may assume that $X$ and $Y$ are symplectic.

    By Theorem \ref{thm: summary LLV computation} the LLV algebra of $X$ decomposes as
    \[\llv(X; \mathbb{Q}) \cong \mathfrak{so}(\widetilde \sH^1(X_0; \mathbb{Q})) \times \prod_{i=1}^m \mathfrak{so}(\widetilde \sH(X_i;\mathbb{Q})),\]
    where $\widetilde \sH^1(X_0; \mathbb{Q}) = \sH^1(X_0; \mathbb{Q}) \oplus \sH^1(X_0; \mathbb{Q})^*$. Similarly, we have
    \[\llv(Y; \mathbb{Q}) \cong \mathfrak{so}(\widetilde \sH^1(Y_0; \mathbb{Q})) \times \prod_{i=1}^n \mathfrak{so}(\widetilde \sH(Y_i;\mathbb{Q})).\]
    Theorem \ref{thm: Taelman theorem A} gives us an isomorphism $\Phi^{\llv} \colon \llv(X; \mathbb{Q}) \to \llv(Y; \mathbb{Q})$, such that the isomorphism $\Phi^{\sH} \colon \sH^\bullet(X; \mathbb{Q}) \to \sH^\bullet(Y; \mathbb{Q})$ is equivariant with respect to $\Phi^{\llv}$. By Theorem \ref{thm: split Lie alg reps}, the isomorphism $\Phi^{\llv}$ is a product of isomorphisms between the simple factors.

    We first claim that $\Phi^{\llv}$ maps the simple factor $\mathfrak{so}(\widetilde \sH^1(X_0; \mathbb{Q}))$ of $\llv(X; \mathbb{Q})$ to the simple factor $\mathfrak{so}(\widetilde \sH^1(Y_0; \mathbb{Q}))$ of $\llv(Y; \mathbb{Q})$. To see why, we look at the odd cohomology. Write 
    \[\sH^\bullet(X_0; \mathbb{Q}) = \sH^{\ev}(X_0; \mathbb{Q}) \oplus \sH^{\odd}(X_0; \mathbb{Q})\] for the decomposition of $\sH^\bullet(X_0; \mathbb{Q})$ into the even and odd parts. Then $\sH^{\ev}(X_0; \mathbb{Q})$ and $\sH^{\odd}(X_0; \mathbb{Q})$ are both irreducible representations of $\llv(X_0; \mathbb{Q})$, they are the even and odd spinor representation. Since $\sH^1(X_i; \mathbb{Q}) = 0$ for all $i \geq 1$ and by using the symmetries of the Hodge diamond coming from the symplectic form, we see that the $\llv(X; \mathbb{Q})$-subrepresentation
    \[ \sH^{\odd}(X_0; \mathbb{Q}) \otimes \bigotimes_{i=1}^m \SH(X_i; \mathbb{Q}) \subseteq \sH^\bullet(X; \mathbb{Q})\]
    is the only irreducible subrepresentation of $\sH^\bullet(X; \mathbb{Q})$ with elements that are in the $(d-1)$-th column of the Hodge diamond of $X$, where $d$ is the integer such that $2d = \dim X$. The analogous statement holds for $Y$. Since $\Phi^{\sH}$ respects the columns of the Hodge diamonds, it restricts to a $\Phi^{\llv}$-equivariant isomorphism 
    \[ \sH^{\odd}(X_0; \mathbb{Q}) \otimes \bigotimes_{i=1}^m \SH(X_i; \mathbb{Q}) \isomto \sH^{\odd}(Y_0; \mathbb{Q}) \otimes \bigotimes_{j=1}^n \SH(Y_j; \mathbb{Q}).\]
    Since the even and odd spinor representations are not isomorphic, the simple factor $\mathfrak{so}(\widetilde \sH^1(X_0; \mathbb{Q}))$ is the only simple factor of $\llv(X; \mathbb{Q})$ whose representation on $\SH(X; \mathbb{Q})$ differs from the representation on $\sH^{\odd}(X_0; \mathbb{Q}) \otimes \bigotimes_{i=1}^m \SH(X_i; \mathbb{Q})$. Since the same holds for $\mathfrak{so}(\widetilde \sH^1(Y_0; \mathbb{Q}))$, we conclude that $\Phi^{\llv}$ indeed maps $\mathfrak{so}(\widetilde \sH^1(X_0; \mathbb{Q}))$ to $\mathfrak{so}(\widetilde \sH^1(Y_0; \mathbb{Q}))$.

    By returning to the decomposition of $\Phi^{\llv}$ as a product of isomorphisms between the simple factors, we see that there must be a bijection $\sigma \colon \{1, \dots, m\} \to \{1, \dots, n\}$ and isomorphisms 
    \[ \mathfrak{so}(\widetilde \sH(X_i; \mathbb{Q})) \isomto \mathfrak{so}(\widetilde \sH(Y_{\sigma(i)}; \mathbb{Q})).\]
    This immediately implies that $b_2(X_i) = b_2(Y_{\sigma(i)})$ for all $i$. It remains to show that $\dim X_i = \dim Y_{\sigma(i)}$ for all $i$. To see this, we use Theorem \ref{thm: split Lie alg reps} to conclude that the isomorphism $\Phi^{\SH} \colon \SH(X; \mathbb{Q}) \to \SH(Y; \mathbb{Q})$ is the tensor product of the isomorphism $\SH(X_0; \mathbb{Q}) \to \SH(Y_0; \mathbb{Q})$ and the isomorphisms $\SH(X_i; \mathbb{Q}) \to \SH(Y_{\sigma(i)}; \mathbb{Q})$. But the dimension of $X_i$ is also the width of the Hodge structure on $\SH(X_i; \mathbb{Q})$, see Proposition \ref{prop: HS of SH_i is wide}. Since the isomorphism $\SH(X_i; \mathbb{Q}) \to \SH(Y_{\sigma(i)}; \mathbb{Q})$ preserves the columns of the Hodge diamond, we conclude that $\dim X_i = \dim Y_{\sigma(i)}$ for all $i$, and this finishes the proof.
\end{proof}

%% file: chapters/chapter6.tex
Let $X$ and $Y$ be two smooth projective varieties over $\mathbb{C}$, admitting a holomorphic symplectic form and satisfying condition \ref{condition *} from the introduction, and let $\Phi \colon D^b(X) \to D^b(Y)$ be an equivalence of triangulated categories. By \cite[Theorem~A]{TaelmanDerivedEquivalences}, there is an induced isomorphism of LLV algebras $\Phi^{\llv} \colon \llv(X; \mathbb{Q}) \to \llv(Y; \mathbb{Q})$, with the property that the isomorphism $\Phi^{\sH} \colon \sH^\bullet(X; \mathbb{Q}) \to \sH^\bullet(Y; \mathbb{Q})$ is equivariant with respect to $\Phi^{\llv}$. Recall that there is an element\index{$h_X$} $h_X \in \llv(X; \mathbb{Q})$ which acts on $\sH^i(X; \mathbb{Q})$ as multiplication by $\dim(X) - i$, and there is an element\index{$h'_X$} $h'_X \in \llv(X; \mathbb{C})$ which acts on $\sH^{p,q}(X)$ as multiplication by by $q-p$. Together, the actions of $h_X$ and $h_X'$ determine the Hodge bigrading on $\sH^\bullet(X; \mathbb{Q})$ (we discuss $h'_X$ in more detail in Section~\ref{sec: Hodge operator}). The goal of this chapter is to prove the following theorem:

\begin{theorem}
    \label{thm: construction of psi}
    Let $X$ and $Y$ be smooth projective varieties over $\mathbb{C}$, admitting a holomorphic symplectic form and satisfying condition \ref{condition *}, and let $\Phi \colon D^b(X) \to D^b(Y)$ be an equivalence of triangulated categories. Then there exist:
    \begin{enumerate}
        \item A reductive algebraic group \myindex{$\mathcal{H}$} over $\mathbb{Q}$,
        \item an inclusion $\llv(Y; \mathbb{Q}) \subseteq \Lie(\mathcal{H})$,
        \item and an element $\psi \in \mathcal{H}(\mathbb{Q})$,
    \end{enumerate}
    satisfying the following properties:
    \begin{enumerate}
        \item The automorphism $\Ad(\psi) \in \Aut(\Lie(\mathcal{H}))$ sends $\llv(Y;\mathbb{Q})$ to itself.
        \item The composition $\Ad(\psi) \circ \Phi^{\llv} \colon \llv(X; \mathbb{Q}) \to \llv(Y; \mathbb{Q})$ sends $h_X$ to $h_Y$ and $h'_X$ to $h'_Y$.
    \end{enumerate}
\end{theorem}

The obtained automorphism $\Ad(\psi)$ of $\llv(Y; \mathbb{Q})$ gives us the map $\Psi^{\llv}$ in diagram~\eqref{eq: The diagram} in the introduction.

We will prove this theorem by constructing $\mathcal{H}$ and $\psi$ `factorwise'. In chapters \ref{chap: LLV computation} and \ref{chap: SH computation}, we computed the LLV algebras and Verbitsky components of $X$ and $Y$, and obtained an explicit description of the simple factors of $\llv(X; \mathbb{Q})$ and $\llv(Y; \mathbb{Q})$. Write $\llv(X; \mathbb{Q}) = \prod_\alpha \mathfrak{g}_\alpha$ and $\llv(Y; \mathbb{Q}) = \prod_\beta \mathfrak{h}_\beta$ for the decompositions into simple factors obtained in Theorem \ref{thm: summary LLV computation}. By Theorem \ref{thm: split Lie alg reps}, after reordering the $\mathfrak{h}_\beta$, the isomorphism $\Phi^{\llv}$ is a product of isomorphisms $\Phi^{\llv}_\alpha \colon \mathfrak{g}_\alpha \to \mathfrak{h}_{\alpha}$. We will then construct for every $\alpha$ a reductive algebraic group $\mathcal{H}_\alpha$ over $\mathbb{Q}$ with $\mathfrak{h}_\alpha \subseteq \Lie(\mathcal{H}_\alpha)$ and an element $\psi_\alpha \in \mathcal{H}_\alpha(\mathbb{Q})$ such that $\Ad(\psi_\alpha) \circ \Phi^{\llv}_\alpha$ gives an isomorphism $\mathfrak{g}_\alpha \to \mathfrak{h}_\alpha$ which preserves the Hodge bigrading. The theorem then follows by taking $\mathcal{H} = \prod \mathcal{H}_\alpha$ and $\psi = (\psi_\alpha)_\alpha$. 

The Lie algebras $\Lie(\mathcal{H}_\alpha)$ do not differ too much from the $\mathfrak{h}_\alpha$: for every $\alpha$, the Lie algebra $\Lie(\mathcal{H}_\alpha)$ will either be equal to $\mathfrak{h}_\alpha$, or it will be the direct product of $\mathfrak{h}_\alpha$ and an abelian Lie algebra. In particular, $\Lie(\mathcal{H}_\alpha)$ is a reductive Lie algebra with semisimple part isomorphic to $\mathfrak{h}_\alpha$.

\section{Reduction to simple factors}
Recall from Theorem \ref{thm: summary LLV computation} that for a variety $X = (\prod_{i=0}^k X_i)/G$ as in Setup~\ref{stp: symplectic var}, which additionally satisfies condition \ref{condition *}, the simple factors $\mathfrak{g}$ of the LLV algebra can be of two possible forms. Either there is a hyperkähler factor $X_i$ such that $\mathfrak{g} \cong \mathfrak{so}(\widetilde \sH(X_i; \mathbb{Q})^{G_i})$, in which case we say that $\mathfrak{g}$ is of hyperkähler origin. Otherwise, there is an isotypic component $V_\rho \subseteq \sH^1(X_0; \mathbb{Q})$ with associated simple algebra $A \subseteq \mathbb{Q}[G]$ such that $\mathfrak{g} \cong \mathfrak{s}(\End_A(\widetilde V_{\rho}), \widetilde \sigma)$. In this case we say that $\mathfrak{g}$ is of abelian origin. In Proposition~\ref{prop: matching of LLV factors} below, we will show that if the isomorphism $\Phi^{\llv}$ matches up two simple factors $\mathfrak{g}$ and $\mathfrak{h}$ of $\llv(X; \mathbb{Q})$ and $\llv(Y; \mathbb{Q})$ respectively, they must both be of hyperkähler origin or both of abelian origin. For this, we need the following result about the signature of the Mukai lattice: 

\begin{lemma}
    \label{lem: H tilde invariants signature}
    Let $Z$ be a hyperkähler variety and let $G$ be a finite group of symplectic automorphisms of $Z$. Then $\widetilde \sH(Z; \mathbb{R})^G$ has signature $(4, k)$ for some $k \in \mathbb{N}$.
\end{lemma}
\begin{proof}
    The quadratic space $\sH^2(Z; \mathbb{R})^G$ has signature $(3, l)$ for some $l \in \mathbb{N}$ by \cite[Lemma 24]{MongardiK3n}. The result then follows from the fact that $\widetilde \sH(Z; \mathbb{R})^G$ is the direct sum of a hyperbolic plane and $\sH^2(Z; \mathbb{R})^G$.
\end{proof}

We will also need the following result. It follows from \cite[Proposition~C.3.14]{ConradReductiveGroupSchemes}, where it is proven in much more generality (including the case $n=8$). However, it is not difficult to give a proof using Theorem \ref{thm: isom lifts from skew to algebra}.
\begin{lemma}
    \label{lem: so determines b}
    Let $(V, b)$ and $(W, b')$ be two quadratic spaces over a field $k$ of characteristic zero, and let $\varphi \colon \mathfrak{so}(V, b) \to \mathfrak{so}(W, b')$ be an isomorphism of Lie algebras. Let $n = \dim_k(V)$. If $n \geq 5$ and $n \neq 8$, then there is a scalar $\lambda \in k^*$ such $(V, b) \cong (W, \lambda b')$ as quadratic spaces.
\end{lemma}
\begin{proof}
    Let $A = \End_k(V)$ and $B = \End_k(W)$, equipped with the adjoint involutions $\sigma$ and $\tau$ corresponding to $b$ and $b'$ respectively. By Theorem \ref{thm: isom lifts from skew to algebra}, the map $\varphi$ lifts to an involution-preserving isomorphism $\widetilde \varphi \colon A \isomto B$ (here we use the assumptions on $\dim_k(V)$). The Skolem-Noether theorem implies that the isomorphism $\widetilde \varphi$ is of the form $\widetilde \varphi = \Ad(f)$ for some isomorphism $f \colon V \isomto W$. The result now follows from the fact that, up to multiplication by a scalar, the bilinear form can be recovered from its adjoint involution \cite[Page~1]{BookOfInvolutions}.
\end{proof}

\begin{proposition}
    \label{prop: matching of LLV factors}
    Let $X$ and $Y$ be smooth projective varieties over $\mathbb{C}$, admitting a holomorphic symplectic form and satisfying condition \ref{condition *}. Let $\mathfrak{g}$ and $\mathfrak{h}$ be simple factors of $\llv(X; \mathbb{Q})$ and $\llv(Y; \mathbb{Q})$ respectively. Assume that $\mathfrak{g} \cong \mathfrak{h}$ and that $\mathfrak{g}$ is not of type $A_1$. 
    \begin{enumerate}
        \item If $\mathfrak{g}$ is of hyperkähler origin, then so is $\mathfrak{h}$. 
        \item If $\mathfrak{g} \cong \mathfrak{s}(\End_A(\widetilde V), \widetilde \sigma)$ for an isotypic component $V$ of $\sH^1(X_0; \mathbb{Q})$, then $\mathfrak{h} \cong \mathfrak{s}(\End_B(\widetilde W), \widetilde \tau)$ for an isotypic component $W$ of $\sH^1(Y_0; \mathbb{Q})$, and the involutions $\widetilde \sigma$ and $\widetilde \tau$ are of the same type.
    \end{enumerate} 
\end{proposition}
\begin{proof}
    This follows from listing all possible isomorphism types of the simple factors of the LLV algebra.

    If a simple factor of $\llv(X; \mathbb{Q})$ is of hyperkähler origin, it is isomorphic to $\mathfrak{so}(\widetilde \sH(X_j; \mathbb{Q})^{G_j})$ for some hyperkähler variety $X_j$. This Lie algebra is either of type $D_n$ or $B_n$ for some~$n$. We know from Lemma~\ref{lem: HK H^2 invariants dimension} that $\dim \widetilde \sH(X_j; \mathbb{Q})^{G_j} \geq 5$.

    If $\mathfrak{g}$ is a factor of abelian origin, then we have seen that it can be of type $A_n, C_n$ or $D_n$, depending on the type of the involution. By condition \ref{condition *} and the assumption that $\mathfrak{g}$ is not of type $A_1$, it cannot be of type $A_1, C_2, A_3$ or $D_4$ (and hence also not of type $D_3$, because $A_3 \cong D_3$). We have also seen that type $D_2$ does not occur, since in this case only one factor of type $A_1$ ends up in the LLV algebra (see Theorem \ref{thm: LLV computation}). Similarly, since $C_1 \cong A_1$, a simple factor of type $C_n$ must have $n \geq 3$. Since this exhausts all exceptional isomorphisms between the types $A_n, C_n$ and $D_n$, we see that if $\mathfrak{g}$ and $\mathfrak{h}$ are two simple factors of abelian origin, they must come from the same type of involution.

    Up to symmetry, the only way in which an isomorphism $\mathfrak{g} \cong \mathfrak{h}$ can occur with $\mathfrak{g}$ and $\mathfrak{h}$ of different origin is if $\mathfrak{g}$ is of hyperkähler origin and $\mathfrak{h}$ is of abelian origin, and both are of type $D_n$ for some $n \geq 5$. We will show that such an isomorphism cannot occur over $\mathbb{R}$.
    
    We know that $\mathfrak{g}_{\mathbb{R}} \cong \mathfrak{so}(\widetilde \sH(X_j; \mathbb{R})^{G_j})$ for some hyperkähler factor $X_j$ of $X$. By Lemma \ref{lem: H tilde invariants signature}, the quadratic form on $\widetilde \sH(X_j; \mathbb{R})^{G_j}$ has signature $(4, k)$ for some $k$. Then $k \geq 6$ by our assumption that $\mathfrak{g}$ is of type $D_n$ with $n \geq 5$. Since $\mathfrak{h}$ is of abelian origin, its base change to $\mathbb{R}$ has the form $\mathfrak{h}_{ \mathbb{R}} = \mathfrak{s}(\End_B(\widetilde W_\mathbb{R}), \widetilde \tau) = \Skew(\End_B(\widetilde W_\mathbb{R}), \widetilde \tau)$ for a semisimple $\mathbb{R}$-algebra $B$ with involution (the second equality follows from the fact that $\Skew(\End_B(\widetilde W_\mathbb{R}), \widetilde \tau)$ is already simple, since $\widetilde \tau$ is orthogonal). Since $\mathfrak{h}_{\mathbb{R}} \cong \mathfrak{g}_{\mathbb{R}}$ and $\mathfrak{g}_{\mathbb{R}}$ is simple, Theorem~\ref{Thm: computation of s(A, sigma)} gives that $Z(B) = Z(\mathfrak{g}_{\mathbb{R}}) = \mathbb{R}$, so $B$ is isomorphic to either $M_l(\mathbb{R})$ or $M_l(\mathbb{H})$ for some $l \in \mathbb{N}$. Then $\End_B(\widetilde W_\mathbb{R})$ is isomorphic to $M_{2n}(\mathbb{R})$ or $M_n(\mathbb{H})$.

    Let $\sigma_b$ be the adjoint involution on $\End_\mathbb{R}(\widetilde \sH(X_j; \mathbb{R})^{G_j})$ associated to the Mukai pairing on $\widetilde \sH(X_j; \mathbb{R})^{G_j}$. Then $\mathfrak{g}_{\mathbb{R}} \cong \Skew(\End_\mathbb{R}(\widetilde \sH(X_j; \mathbb{R})^{G_j}), \sigma_b)$. By Theorem~\ref{thm: isom lifts from skew to algebra}, the isomorphism $\mathfrak{g}_{\mathbb{R}} \cong \mathfrak{h}_{\mathbb{R}}$ lifts to an isomorphism of algebras with involution $(\End_\mathbb{R}(\widetilde \sH(X_j; \mathbb{R})^{G_j}), \sigma_b) \cong (\End_B(\widetilde W_\mathbb{R}), \widetilde \tau)$. It follows that $\End_B(\widetilde W_\mathbb{R})$ must be isomorphic to $M_{2n}(\mathbb{R})$. 

    The bilinear form on $\widetilde W_\mathbb{R}$ is of signature $(nl, nl)$, so the induced involution $\widetilde \sigma$ on $\End_B(\widetilde W_\mathbb{R}) \cong M_{2n}(\mathbb{R})$ corresponds to a bilinear form of signature $(n, n)$ on $\mathbb{R}^{2n}$. So if $\mathfrak{g} \cong \mathfrak{h}$, we would have an isomorphism of real Lie algebras $\mathfrak{so}(4, k) \cong \mathfrak{so}(n, n)$, where $k \geq 6$. This is not possible by Lemma~\ref{lem: so determines b}, which finishes the proof.
\end{proof}

\begin{remark}
    There is an additional piece of information that we could use in the above proof. Namely, associated to the decomposition $\llv(X; \mathbb{Q}) = \prod_\alpha \mathfrak{g}_\alpha$ into simple factors, there is a factorization $\SH(X; \mathbb{Q}) \cong \bigotimes_\alpha \SH_\alpha(X; \mathbb{Q})$, and each $\SH_\alpha(X; \mathbb{Q})$ is an absolutely irreducible representation of $\mathfrak{g}_\alpha$. If we also have a simple factor $\mathfrak{h}_\alpha$ of $\llv(Y; \mathbb{Q})$ and an isomorphism $\Phi^{\llv}_\alpha \colon \mathfrak{g}_\alpha \to \mathfrak{h}_\alpha$, then there is an induced $\Phi^{\llv}_\alpha$-equivariant isomorphism $\Phi^{\SH}_\alpha \colon \SH_\alpha(X; \mathbb{Q}) \to \SH_\alpha(Y; \mathbb{Q})$ (see Theorem~\ref{thm: split Lie alg reps}). One could hope that we could use this to weaken condition \ref{condition *}.

    For example, suppose we have a simple factor $\mathfrak{g}_\alpha$ of $\llv(X; \mathbb{Q})$ of type $C_2$ of abelian origin, and a simple factor $\mathfrak{h}_\alpha$ of $\llv(Y; \mathbb{Q})$ of type $B_2$ of hyperkähler origin, then it could be possible that $\Phi^{\llv}_\alpha$ sends $\mathfrak{g}_\alpha$ to $\mathfrak{h}_\alpha$, because there is an exceptional isomorphism $B_2 \cong C_2$. One could hope that the structure of $\SH_\alpha(X; \mathbb{Q})$ and $\SH_\alpha(Y; \mathbb{Q})$ prevents the existence of the associated isomorphism $\Phi^{\SH}_\alpha$. Unfortunatly, one can do a computation similar to Proposition~\ref{prop: weight of Dn SH} to show that, after complexification, the representations $\SH_\alpha(X; \mathbb{Q})$ and $\SH_\alpha(Y; \mathbb{Q})$ both have their highest weight a multiple of the same fundamental weight under the isomorphism $B_2 \cong C_2$ (with the precise coefficients depending on the situation).

    However, there is one case where this trick does work. Suppose $\mathfrak{g}_\alpha$ is a simple factor of $\llv(X; \mathbb{Q})$ of type $A_3$. If $\mathfrak{g}_\alpha$ is of abelian origin, there are two ways in which it can be constructed, either from an orthogonal involution or from a unitary involution, due to the isomorphism $A_3 \cong D_3$. One can use Propositions~\ref{prop: weight of Dn SH} and~\ref{prop: weight of An SH} to distinguish between these two cases, see Lemma~\ref{lem: Phi llv orthogonal vs unitary} for more details.    
\end{remark}

\section{The Hodge operator}
\label{sec: Hodge operator}
Recall that for a graded vector space $V$ over a field $k$, we have a grading operator $h$ which, for every integer $i$, acts on the graded piece $V_i$ as multiplication by $i$. Similarly, if $V$ is equipped with a Hodge structure, we have the Hodge operator~$h'$:
\begin{definition}
    Let $V$ be a $\mathbb{Q}$-Hodge structure. The Hodge operator $h'$ is the linear endomorphism $h'$ of $V_\mathbb{C}$ which for every pair of integers $(p, q)$ acts as multiplication by $(q-p)$ on $V^{p,q}$.
\end{definition}

The eigenspaces of $h'$ are the columns of the Hodge diamond of $V$. Let $V$ and $W$ be two vector spaces over $\mathbb{Q}$, both equipped with a pure Hodge structure, and with Hodge operators $h'_V$ and $h'_W$ respectively. Then a linear map $f \colon V \to W$ respects the Hodge grading if and only if $f \circ h_V' = h_W' \circ f$, and $f$ respects the grading if and only if $f \circ h_V = h_W \circ f$. Then $f$ is a morphism of Hodge structures if and only if $f$ commutes with both $h$ and $h'$. To stress this, we say that such an $f$ respects the \textbf{Hodge bigrading}.

\begin{lemma}[{\cite[Corollary 2.12]{TaelmanDerivedEquivalences}}]
    Let $X$ be a smooth projective variety over $\mathbb{C}$ admitting a holomorphic symplectic form and let $h_X' \in \End(\sH^\bullet(X; \mathbb{C}))$ be the Hodge operator. Then $h_X' \in \llv(X; \mathbb{C})$. \qed
\end{lemma}

Under the product decomposition $\llv(X; \mathbb{Q}) = \prod \mathfrak{g}_\alpha$ from Theorem \ref{thm: summary LLV computation}, the Hodge operator of $X$ decomposes as $h_X' = (h_{X, \alpha}')_\alpha$, where $h'_{X,\alpha} \in \mathfrak{g}_{\alpha, \mathbb{C}} := \mathfrak{g}_{\alpha} \otimes_\mathbb{Q} \mathbb{C}$ for all~$\alpha$. 
\begin{lemma}
    \label{lem: equality of Hodge operators}
    Let $X$ be a smooth projective variety over $\mathbb{C}$ admitting a holomorphic symplectic form, let $\mathfrak{g}_\alpha$ be a simple factor of $\llv(X, \mathbb{Q})$, and let $h'_{X,\alpha} \in \mathfrak{g}_{\alpha, \mathbb{C}}$ be the $\alpha$-th component of the Hodge operator.
    \begin{enumerate}
        \item If $\mathfrak{g}_\alpha \cong \mathfrak{so}(\widetilde \sH(X_j; \mathbb{Q})^{G_j})$ is of hyperkähler origin, then $h_{X,\alpha}'$ as an endomorphism of $\widetilde \sH(X_j; \mathbb{Q})^{G_j}$ acts as the Hodge operator of $\widetilde \sH(X_j; \mathbb{Q})^{G_j}$.
        \item If $V$ is an isotypic component of $\sH^1(X_0; \mathbb{Q})$ and $\mathfrak{g}_\alpha \cong \mathfrak{s}(\End_A(\widetilde V), \widetilde \sigma) \subseteq \mathfrak{gl}(\widetilde V)$, then $h_{X,\alpha}'$ acts on $\widetilde V$ as the Hodge operator of $\widetilde V$.
    \end{enumerate}
\end{lemma}
\begin{proof}
    For the first statement, recall that the isomorphism $\mathfrak{g}_\alpha \cong \mathfrak{so}(\widetilde \sH(X_j; \mathbb{Q})^{G_j})$ is the restriction of the isomorphism $\llv(X_{[j]}; \mathbb{Q}) \cong \prod_{k \in [j]} \mathfrak{so}(\widetilde \sH(X_k; \mathbb{Q}))$ to the Lie subalgebras generated by $G$-invariant elements in degrees $-2$ and $2$ (see Proposition \ref{prop: LLV of HK factor}). Therefore it suffices to show that for every hyperkähler variety $X_k$, under the isomorphism $\llv(X_k; \mathbb{Q}) \cong \mathfrak{so}(\widetilde \sH(X_k; \mathbb{Q}))$, the Hodge operator of $X_k$ acts on $\widetilde \sH(X_k; \mathbb{Q})$ as the Hodge operator.

    Let $h_k' \in \mathfrak{gl}(\widetilde \sH(X_k; \mathbb{Q}))$ be the Hodge operator of $\widetilde \sH(X_k; \mathbb{Q})$. Since the Beauville-Bogomolov-Fujiki form on $\widetilde \sH(X_k; \mathbb{Q})$ respects the Hodge structure, we see that $h'_k \in \mathfrak{so}(\widetilde \sH(X_k; \mathbb{Q}))$. There is a commutative triangle
    \begin{equation*}
        \begin{tikzcd}
            \mathfrak{so}(\widetilde \sH(X_k; \mathbb{Q})) \arrow[rr, "\sim"] \arrow[rd] & & \mathfrak{llv}(X_k; \mathbb{Q}) \arrow[ld]\\
            & \mathfrak{gl}( \SH(X_k; \mathbb{Q})),
        \end{tikzcd}
    \end{equation*}
    where the action of $\mathfrak{so}(\widetilde \sH(X_k; \mathbb{Q}))$ on $\SH(X_k; \mathbb{Q})$ comes from the injective map $\SH(X_k; \mathbb{Q}) \to \Sym^d \widetilde \sH(X_k; \mathbb{Q})$ from \cite[Proposition 3.5]{TaelmanDerivedEquivalences}, where $d$ is defined by $2d = \dim X_k$. The Verbitsky component $\SH(X_k; \mathbb{Q})$ is a faithful representation of $\llv(X_k; \mathbb{Q})$, so the two diagonal arrows are injective. The image of $h'_k$ in $\mathfrak{gl}( \SH(X_k; \mathbb{Q}))$ acts as the Hodge operator by \cite[Lemma~4.8]{TaelmanDerivedEquivalences}. The Hodge operator of $\llv(X_k; \mathbb{Q})$ is by definition the Hodge operator of $\sH^\bullet(X_k; \mathbb{Q})$, so it also acts on $\SH(X_k; \mathbb{Q})$ as the Hodge operator. Therefore, the isomorphism $\mathfrak{so}(\widetilde \sH(X_k; \mathbb{Q})) \to \llv(X_k; \mathbb{Q})$ sends $h'_k$ to the Hodge operator of $X_k$, and this finishes the proof of the first statement of the lemma. 

    For the second statement, note that the isomorphism $\mathfrak{g}_{\alpha, \mathbb{C}} \cong \mathfrak{s}(\End_A(\widetilde V), \widetilde \sigma)_{\mathbb{C}}$ is obtained as the restriction of the isomorphism $\llv(\extp^\bullet V) \cong \mathfrak{so}(\widetilde V)$ of \cite[Chapter 3]{LooijengaLunts}, so it suffices to show that this isomorphism sends $h'_{\bigwedge^\bullet V}$ to $h'_{\widetilde V}$, where $h'_{\bigwedge^\bullet V}$ to $h'_{\widetilde V}$ are the Hodge operators on $\extp^\bullet V$ and $\widetilde V$ respectively. We have $h'_{\widetilde V} \in \mathfrak{so}(\widetilde V)_0 \cong \mathfrak{gl}(V)$, and \cite[Proposition 3.2]{LooijengaLunts} says that under the induced map $\mathfrak{gl}(V) \to \llv(\bigwedge^\bullet V)$ any $x \in \mathfrak{gl}(V)$ acts on $\bigwedge^\bullet V$ as $\widetilde x - \frac{1}{2} \Tr(x) \id$, where $\widetilde x$ is the unique derivation of $\bigwedge^\bullet V$ extending $x$. The isomorphism $\mathfrak{so}(\widetilde V)_0 \cong \mathfrak{gl}(V)$ sends $h'_{\widetilde V}$ to $h'_V$. The symmetry of the Hodge structure implies that $\Tr(h'_V) = 0$, and a small computation shows that $h'_{\bigwedge^\bullet V}$ is a derivation of $\bigwedge^\bullet V$. Therefore the isomorphism $\mathfrak{so}(\widetilde V) \cong \llv(\bigwedge^\bullet V)$ sends $h'_{\widetilde V}$ to $h'_{\bigwedge^\bullet V}$.
\end{proof}

\section{Proof of Theorem \ref{thm: construction of psi}}
In this section, we will prove Theorem \ref{thm: construction of psi}. We will prove it factorwise, as explained in the introduction to this chapter. In Section \ref{sec: factorwise psi abelian}, we will prove it for factors of abelian origin, except for those of type $A_1$. We prove it for factors of type $A_1$ in Section \ref{sec: factorwise psi A1}, and lastly we prove it for simple factors of hyperkähler origin in Section \ref{sec: factorwise psi HK}. We then summarize these results in Section \ref{sec: factorwise psi summary}.

We continue with the decompositions $\llv(X; \mathbb{Q}) = \prod_\alpha \mathfrak{g}_\alpha$ and $\llv(Y; \mathbb{Q}) = \prod_\beta \mathfrak{h}_\beta$. The Verbitsky component $\SH(X; \mathbb{Q})$ is isomorphic to the tensor product $\bigotimes_\alpha \SH_\alpha(X)$ for certain irreducible representations $\SH_\alpha(X)$ of $\mathfrak{g}_\alpha$, and similarly $\SH(Y; \mathbb{Q}) \cong \bigotimes_\beta \SH_\beta(Y)$ for irreducible representations $\SH_\beta(Y)$ of $\mathfrak{h}_\beta$. The Verbitsky components $\SH(X; \mathbb{Q})$ and $\SH(Y; \mathbb{Q})$ are absolutely irreducible by Corollary~\ref{cor: SH(X) is absolutely irreducible}, so using Theorem~\ref{thm: split Lie alg reps} we obtain the following from a derived equivalence $\Phi \colon D^b(X) \to D^b(Y)$:
\begin{itemize}
    \item After reordering the $\mathfrak{h}_\beta$, we have for every $\alpha$ an isomorphism \index{$\Phi^{\llv}_\alpha$}$\Phi^{\llv}_\alpha \colon \mathfrak{g}_\alpha \to \mathfrak{h}_\alpha$ with the property that $\Phi_\alpha^{\llv}(h'_{X,\alpha}) = h'_{Y,\alpha}$, and $\Phi^{\llv} \colon \llv(X; \mathbb{Q}) \to \llv(Y; \mathbb{Q})$ is the product of the isomorphisms $\Phi^{\llv}_\alpha$.
    \item For every $\alpha$, there is an associated isomorphism\index{$\Phi^{\SH}_\alpha$} $\Phi^{\SH}_\alpha \colon \SH_\alpha(X) \to \SH_\alpha(Y)$ that is equivariant with respect to $\Phi^{\llv}_\alpha$. The isomorphism $\Phi^{\SH}$ is the tensor product of the maps $\Phi^{\SH}_\alpha$.
    \item Proposition \ref{prop: matching of LLV factors} gives that, for every $\alpha$, either $\mathfrak{g}_\alpha$ and $\mathfrak{h}_\alpha$ are both of hyperkähler origin or both of abelian origin. Moreover, if they are both of abelian origin and not of type $A_1$, they both come from an algebra with the same type of involution.
\end{itemize}

\subsection{Factors of abelian origin}
\label{sec: factorwise psi abelian}
Let $\mathfrak{g}_\alpha$ be a simple factor of $\llv(X; \mathbb{Q})$, let $\mathfrak{h}_\alpha$ be the corresponding simple factor of $\llv(Y; \mathbb{Q})$, and let $\Phi_\alpha^{\llv} \colon \mathfrak{g}_\alpha \to \mathfrak{h}_\alpha$ be the isomorphism induced by the derived equivalence $\Phi$. Assume that $\mathfrak{g}_\alpha$ and $\mathfrak{h}_\alpha$ are both of abelian origin, so $\mathfrak{g}_\alpha \cong \mathfrak{s}(\End_A(\widetilde V), \widetilde \sigma)$ and $\mathfrak{h}_\alpha \cong \mathfrak{s}(\End_B(\widetilde W), \widetilde \tau)$ for some algebras $A$ and $B$, an $A$-module $V$ and a $B$-module $W$. We will show that the isomorphism $\Phi_\alpha^{\llv} \colon \mathfrak{g}_\alpha \to \mathfrak{h}_\alpha$ is of the form $\Ad(f)$ for an isomorphism $f \colon \widetilde V \to \widetilde W$, and then construct a $\psi_\alpha \colon \widetilde W \to \widetilde W$ such that $\psi_\alpha \circ f$ respects the Hodge bigrading.

We start by proving some algebraic results, and in Corollary \ref{cor: psi_i for ab var} we will apply these results to the maps $\Phi_\alpha^{\llv}$. Our algebraic results will require the following data:
\begin{setup}
    \label{stp: algebraic data}
    Consider the following data:
    \begin{itemize}
        \item Simple $\mathbb{Q}$-algebras $A$ and $B$ with involutions $\sigma$ and $\tau$ respectively.
        \item An $A$-module $V$ and a $B$-module $W$.
        \item An isomorphism of Lie algebras $\phi \colon \mathfrak{s}(\widetilde A, \widetilde \sigma) \to \mathfrak{s}(\widetilde B, \widetilde \tau)$, where $\widetilde A = \End_A(\widetilde V)$ and $\widetilde B = \End_B(\widetilde W)$, and where $\widetilde \sigma$ and $\widetilde \tau$ are the involutions associated to the bilinear forms on $\widetilde V$ and $\widetilde W$ respectively.
    \end{itemize}
    Assume that this data satisfies the following:
    \begin{itemize}
        \item The modules $V$ and $W$ have the same dimensions as vector spaces over $\mathbb{Q}$.
        \item The Lie algebra $\mathfrak{s}(\widetilde A, \widetilde \sigma)$ is simple of one of the following three types:
        \begin{itemize}
            \item Type $A_n$ with $n \geq 2$.
            \item Type $C_n$ with $n \geq 3$.
            \item Type $D_n$ with $n \geq 5$.
        \end{itemize}
    \end{itemize}
\end{setup}
Note that the restriction of the possible simple types of $\mathfrak{s}(\widetilde A, \widetilde \sigma)$ corresponds to condition \ref{condition *} and requiring that $\mathfrak{s}(\widetilde A, \widetilde \sigma)$ is not of type $A_1$.

\begin{theorem}
    \label{thm: involution morphism from LLV isom}
    Assume that we are in the setting of Setup \ref{stp: algebraic data}. Then there exists an isomorphism $\beta \colon A \to B$ respecting the involutions.
\end{theorem}

\begin{proof}
    The assumptions on $\mathfrak{s}(\widetilde A, \widetilde \sigma)$ imply that Theorem \ref{thm: isom lifts from skew to algebra} can be used to lift the isomorphism $\phi$ to an isomorphism of algebras $\widetilde \phi \colon \widetilde A \to \widetilde B$, preserving the involutions. Let $p \colon Z(\widetilde A) \to Z(\widetilde B)$ be the induced isomorphism on the centers.
    
    Write $A = M_d(D)$ and $B = M_e(E)$ for certain division algebras $D$ and $E$ and integers $d$ and $e$. Let $V_0$ be the unique (up to isomorphism) simple $A$-module and write $\widetilde V \cong V_0^{2m}$. Similarly, let $\widetilde W \cong W_0^{2n}$ for the unique (up to isomorphism) simple $B$-module $W_0$. Then $\widetilde A \cong M_{2m}(D^{\op})$ and $\widetilde B \cong M_{2n}(E^{\op})$. By Wedderburns theorem, the isomorphism $\widetilde A \to \widetilde B$ induces an isomorphism $D \cong E$, and we have $m = n$. We want this isomorphism $D \cong E$ to be linear with respect to $p \colon Z(A) \to Z(B)$.
    
    This isomorphism $D \cong E$ can be made explicit by taking a simple $\widetilde A$-module $V_1$ and a simple $\widetilde B$-module $W_1$. Then $W_1$ becomes a simple $\widetilde A$-module via the action $a \cdot w := \widetilde \phi(a)w$ for $a \in \widetilde A$ and $w \in W_1$. Therefore, there is an isomorphism $\varphi \colon V_1 \to W_1$ with the property that $\varphi(av) = \widetilde\phi(a)\varphi(v)$ for $v \in V_1$ and $a \in \widetilde A$. Conjugation by $\varphi$ gives an isomorphism $D = \End_{\widetilde A}(V_1) \to \End_{\widetilde B}(W_1) = E$. Observe that $\varphi$ is linear over $p$, i.e. we have $\varphi(\lambda v) = p(\lambda)\varphi(v)$ for $\lambda \in Z(A)$ and $v \in V_1$. It follows that the induced map $\End_{\widetilde A}(V_1) \to \End_{\widetilde B}(W_1)$ is also linear over $p$, and therefore we obtain an isomorphism $D \to E$ that is linear over $p$.

    We know that $V_0 \cong D^d$ and $W_0 \cong E^e$. Let $\delta = \dim_\mathbb{Q}(D) = \dim_\mathbb{Q}(E)$, then we have $\dim_\mathbb{Q} V = \delta md$ and $\dim_\mathbb{Q} W = \delta ne$. Since $\dim_\mathbb{Q} V = \dim_\mathbb{Q} W$ by assumption, we have $\delta md = \delta ne$. The fact that $m = n$ then implies $d = e$, and therefore $A \cong M_d(D) \cong M_e(E) \cong B$, and this isomorphism $\beta' \colon A \to B$ is linear over $p$. 

    Next, we want to modify this isomorphism $\beta'$ to an isomorphism $\beta$ that preserves the involutions. The isomorphism $\widetilde \phi \colon \widetilde A \to \widetilde B$ preserves the involution. Lemma \ref{lem: Endomorphism inherits involution type} shows that the involutions on $A$ and $\widetilde A$ are of the same kind, and the same holds for $B$ and $\widetilde B$. Together, this gives that the involutions on $A$ and $B$ both are of the same kind. If the involutions on $A$ and $B$ are both of the first kind, then one can modify the isomorphism $\beta'$ to one that respects the involutions by \cite[Proposition 2.7]{BookOfInvolutions}. If the involutions on $A$ and $B$ are both of the second kind, then Lemma \ref{lem: Endomorphism inherits involution type} shows that they induce the same automorphism of the center, after identifying the centers via $p$. Therefore, we can modify $\beta'$ to an isomorphism that respects the involution by \cite[{Proposition 2.18}]{BookOfInvolutions}.
\end{proof}

For the next result, we will need the following version of the Skolem--Noether theorem. It can easily be deduced from the more absolute version where $k = k'$, as stated in \cite[\href{https://stacks.math.columbia.edu/tag/074Q}{Theorem 074Q}]{stacks}.
\begin{theorem}
    \label{thm: relative Skolem-Noether}
    Let $\alpha \colon k \to k'$ be an isomorphism of fields, let $A$ be a simple $k$-algebra and let $B$ be a finite-dimensional central simple $k'$-algebra. Denote the embeddings by $\iota \colon k \to A$ and $\iota' \colon k' \to B$. Suppose we have two algebra homomorphisms $f,g \colon A \to B$ with the property that $f \circ \iota = \iota' \circ \alpha$ and $g \circ \iota = \iota' \circ \alpha$. Then there is an invertible element $b \in B$ with the property that $f(a) = bg(a)b^{-1}$ for all $a \in A$. \qed
\end{theorem}

Given an isomorphism $f \colon \widetilde V \to \widetilde W$, conjugation by $f$ induces an isomorphism $\Ad(f) \colon \widetilde A \to \widetilde B$.
\begin{theorem}
    \label{thm: similitude from LLV isom}
    Assume that we are in the setting of Setup \ref{stp: algebraic data} and let $\beta \colon A \to B$ be the isomorphism from Theorem~\ref{thm: involution morphism from LLV isom}. Then there exists an isomorphism $f \colon \widetilde V \to \widetilde W$ with the following properties:
    \begin{enumerate}
        \item $f$ is equivariant with respect to the isomorphism $\beta \colon A \to B$.
        \item The restriction $\Ad(f)|_{\mathfrak{s}(\widetilde A, \widetilde \sigma)} \colon \mathfrak{s}(\widetilde A, \widetilde \sigma) \to \mathfrak{s}(\widetilde B, \widetilde \tau)$ is equal to $\phi$.
        \item The map $\Ad(f) \colon \widetilde A \to \widetilde B$ respects the involution, i.e. $\Ad(f) \circ \widetilde \sigma = \widetilde \tau \circ \Ad(f)$.
    \end{enumerate}
\end{theorem}

\begin{proof}
    Let $V_0$ be the (up to isomorphism) unique simple $A$-module, and let $W_0$ be the (up to isomorphism) unique simple $B$-module. Then $W_0$ also becomes a simple $A$-module under $\beta$, so there exists an isomorphism $f_0 \colon V_0 \to W_0$ with the property that $f_0(av) = \beta(a)f_0(v)$ for all $a \in A$ and $v \in V_0$. Let $m$ and $n$ be the integers such that $V \cong V_0^m$ and $W \cong W_0^n$. We saw in the proof of Theorem \ref{thm: involution morphism from LLV isom} that $m = n$, so by taking the direct sum of $n$ copies of $f_0$ we obtain an equivariant isomorphism $f_1 \colon V \to W$. Since $\beta$ respects the involutions, the induced dual map $f_1^* \colon W^* \to V^*$ is equivariant with respect to $\beta^{-1}$. Taking the direct sum of $f_1$ and $(f_1^*)^{-1}$, we obtain an equivariant isometry $\widetilde f_1 \colon \widetilde V \to \widetilde W$.

    As in the proof of Theorem \ref{thm: involution morphism from LLV isom}, let $\widetilde \phi \colon \widetilde A \to \widetilde B$ be the map extending $\phi$. The isomorphism $\Ad(\widetilde f_1) \colon \widetilde A \to \widetilde B$ is linear with respect to $p \colon Z(A) \to Z(B)$, so we can apply Theorem \ref{thm: relative Skolem-Noether} to the maps $\Ad(\widetilde f_1)$ and $\widetilde \phi$ to get an invertible $f_2 \in \widetilde B = \End_B(\widetilde W)$ with $\widetilde \phi(X) = f_2 \widetilde f_1 \circ X \circ \widetilde f_1^{-1} f_2^{-1}$. If we let $f = f_2 \circ \widetilde f_1$, then we have found the right isomorphism $\widetilde V \to \widetilde W$. It is clear from the construction that $f$ is equivariant with respect to $\beta$.

    We know that $\Ad(f) = \widetilde \phi$ as isomorphisms $\widetilde A \to \widetilde B$, so $\Ad(f)$ indeed extends~$\phi$. We know from Theorem~\ref{thm: isom lifts from skew to algebra} that $\widetilde \phi \circ \widetilde \sigma = \widetilde \tau \circ \widetilde \phi$, so we also have $\Ad(f) \circ \widetilde \sigma = \widetilde \tau \circ \Ad(f)$.
\end{proof}

\begin{remark}
    \label{rem: widetilde V HS}
    So far, we defined $\widetilde V$ to be $V \oplus V^*$, where $V$ sits in degree $1$ and $V^*$ in degree $-1$. From now on, $V$ will be equipped with a weight $1$ Hodge structure, and then $V^*$ becomes a weight $-1$ Hodge structure. To obtain a weight $1$ Hodge structure again, we let $\widetilde V = V \oplus V^*(-1)$ from now on. However, we still let the grading operator $h_{\widetilde V}$ act on $V$ as multiplication by $1$, and on $V^*(-1)$ as multiplication by $-1$.
\end{remark}

Now suppose that $V$ and $W$ are equipped with an $A$-equivariant (respectively $B$-equivariant) Hodge structure of weight $1$. This induces a weight $1$ Hodge structure on $\widetilde V = V \oplus V^*(-1)$. Let $h'_{\widetilde V} \in \mathfrak{s}(\widetilde A, \widetilde \sigma)_{\mathbb{C}}$ be the Hodge operator on $\widetilde V$, and similarly let $h'_{\widetilde W} \in \mathfrak{s}(\widetilde B, \widetilde \tau)_\mathbb{C}$ be the Hodge operator on $\widetilde W$. The Hodge operator $h'_{\widetilde V}$ then gives $\widetilde V \otimes_\mathbb{Q} \mathbb{C}$ the Hodge grading, where $V^{0,1} \oplus (V^{1, 0})^*$ has Hodge degree $1$ and $V^{1, 0} \oplus (V^{0,1})^*$ has Hodge degree $-1$.

\begin{lemma}
    \label{lem: f : widetilde V to widetilde W respects HS}
    Assume that we are in the setting of Setup \ref{stp: algebraic data}. Furthermore, assume that $V$ and $W$ have an equivariant Hodge structure of weight $1$, and assume that $\phi \colon \mathfrak{s}(\widetilde A, \widetilde \sigma) \to \mathfrak{s}(\widetilde B, \widetilde \tau)$ sends $h'_{\widetilde V}$ to $h'_{\widetilde W}$. Then the isomorphism $f \colon \widetilde V \to \widetilde W$ from Theorem~\ref{thm: similitude from LLV isom} respects the Hodge grading.
\end{lemma}
\begin{proof}
    By construction, we have $\Ad(f) = \phi$. Hence the assumption $\phi(h'_{\widetilde V}) = h'_{\widetilde W}$ implies that $f \circ h'_{\widetilde V} = h'_{\widetilde W} \circ f$, and this is equivalent to $f$ respecting the Hodge grading.
\end{proof}

We will now modify $f$ by a similitude to obtain a map which preserves the Hodge bigrading (so both the normal and the Hodge grading). The group of similitudes of an algebra with involution $(A, \sigma)$ is defined as\index{$\Sim(A, \sigma)$}
$$\Sim(A, \sigma) = \{a \in A^\times : a\sigma(a) \in Z(A)\}.$$
For a non-degenerate quadratic space $(V, b)$ over $\mathbb{Q}$ with associated involution $\sigma_b$, we know that $\Sim(\End_\mathbb{Q}(V), \sigma_b)$ is the group $\Sim(V, b)$ of similitudes of the quadratic space $(V, b)$, i.e. the set of maps $f \in \GL(V)$ for which there is a scalar $\lambda \in \mathbb{Q}^*$ such that
$$b(v, w) = \lambda b(f(v), f(w))$$
for all $v, w \in V$, see \cite[Section 12.B]{BookOfInvolutions}. Now let $A$ be an algebra with involution and assume that $(V, b)$ is an $A$-equivariant quadratic space. Then $Z(A)$ can be larger than $\mathbb{Q}$, so if we write $\Sim_A(V, b) = \End_A(V) \cap \Sim(V, b)$, we only have an inclusion
$$\Sim_A(V, b) \subseteq \Sim(\End_A(V), \sigma_b),$$
and this is not always an equality.

\begin{lemma}
    \label{lem: abelian factor similitude psi}
    Assume that we are in the setting of Lemma \ref{lem: f : widetilde V to widetilde W respects HS}. Furthermore, assume that $A$ is a simple factor of the group algebra $\mathbb{Q}[G]$ of some finite group $G$. Then there exists an equivariant isometry $\widetilde V \to \widetilde W$ that preserves the Hodge bigrading. In particular, there exists an element $\psi \in \Sim(\widetilde B, \widetilde \tau)$ so that $\psi \circ f$ preserves the Hodge bigrading.
\end{lemma}
\begin{proof}
    Since $A$ is a simple factor of $\mathbb{Q}[G]$, there is a multiplicative map $G \to A$. By composing with the isomorphism $\beta \colon A \to B$, we also obtain a multiplicative map $G \to B$, so $W$ obtains the structure of a $G$-representation. Hence, $V \oplus V^*(-1)$ and $W \oplus W^*(-1)$ are both objects of the category of polarizable $G$-equivariant weight $1$ Hodge structures. By Proposition \ref{prop: category of polarizable G-equiv HS is semisimple}, this category is semisimple, and we have isomorphisms $W \cong W^*(-1)$ and $V \cong V^*(-1)$.

    Since the isomorphism $f \colon \widetilde V \to \widetilde W$ from Lemma \ref{lem: f : widetilde V to widetilde W respects HS} is equivariant and respects the Hodge grading, it gives an isomorphism of equivariant weight $1$ Hodge structures $V \oplus V^*(-1) \to W \oplus W^*(-1)$. Therefore, as $G$-equivariant weight $1$ Hodge structures, $V$ and $W$ must have the same simple factors with the same multiplicities, so there exists a $G$-equivariant isomorphism of Hodge structures $g \colon V \to W$. Let $g^* \colon W^* \to V^*$ be the dual morphism of $g$, and let $\widetilde g = g \oplus (g^*)^{-1}$. Then $\widetilde g \colon \widetilde V \to \widetilde W$ is the desired morphism, and we can take $\psi = \widetilde g \circ f^{-1}$. Then $\psi \circ f$ preserves the Hodge bigrading.
    
    It remains to show that $\psi \in \Sim(\widetilde B, \widetilde \tau)$. Since $\widetilde g$ is an isometry, we have $\Ad(\widetilde g) \circ \widetilde \sigma = \widetilde \tau \circ \Ad(\widetilde g)$. Therefore also $\Ad(\psi) \circ \widetilde \tau = \widetilde \tau \circ \Ad(\psi)$, so for all $b \in \widetilde B$ we have 
    $$\psi \circ \widetilde \tau(b) \circ \psi^{-1} = \widetilde \tau(\psi^{-1}) \circ \widetilde \tau(b) \circ \widetilde \tau(\psi).$$
    Since $\widetilde \tau \colon \widetilde B \to \widetilde B$ is surjective, this implies that $\widetilde \tau(\psi) \psi$ commutes with all $b \in \widetilde B$, so $\widetilde \tau(\psi)\psi \in Z(\widetilde B)$, and therefore $\psi \in \Sim(\widetilde B, \widetilde \tau)$.
\end{proof}

Recall that the group of isometries of an algebra with involution $(A, \sigma)$ is defined by
\[ \Iso(A, \sigma) := \{a \in A^\times : a\sigma(a) = 1\}. \]
Equivalently, one could define this as the set of similitudes $a \in \Sim(A, \sigma)$ with $a\sigma(a) = 1$. In particular, there is an inclusion $\Iso(A, \sigma) \subseteq \Sim(A, \sigma)$.

\begin{lemma}
    \label{lem: abelian factor isometry psi}
    Assume that we are in the same setting as in Lemma \ref{lem: abelian factor similitude psi}. Then there also exists a map $\psi' \in \Iso(\widetilde B, \widetilde \tau)$ with the property that $\psi' \circ f$ preserves the Hodge bigrading. 
\end{lemma}
\begin{proof}
    For every $\lambda \in Z(B)^\times$ with $\tau(\lambda) = \lambda$, let $s_\lambda \colon \widetilde W \to \widetilde W$ be the map defined by
    $$(w, \eta) \mapsto (\lambda w, \eta)$$
    for $w \in W$ and $\eta \in W^*$. We claim that $s_\lambda$ is $B$-equivariant, preserves the Hodge bigrading and satisfies $\widetilde \tau(s_\lambda) s_\lambda = \lambda$ (in particular, we have $s_\lambda \in \Sim(\widetilde B, \widetilde \tau)$). Lemma \ref{lem: abelian factor similitude psi} gives us a $\psi \in \Sim(\widetilde B, \widetilde \tau)$ with the property that $\psi \circ f$ preserves the Hodge bigrading. If we then take $\lambda = (\widetilde \tau(\psi)\psi)^{-1} \in Z(\widetilde B)^\times = Z(B)^\times$ (note that indeed $\tau(\lambda) = \widetilde \tau(\lambda) = \lambda$), the map $\psi' = s_{\lambda}\circ \psi$ is an isometry with the property that $\psi' \circ f$ preserves the Hodge bigrading.

    It is clear from the definition that $s_\lambda$ preserves the grading on $\widetilde W$, and it is $B$-equivariant since $\lambda \in Z(B)$. It preserves the Hodge grading because we assume that the action of $B$ on $W$ preserves the Hodge grading, and $\lambda \in Z(B) \subseteq B$.

    To see that $\widetilde \tau(s_\lambda) s_\lambda = \lambda$, we first compute $\widetilde \tau(s_\lambda)$. Let $s'_\lambda \colon \widetilde W \to \widetilde W$ be the map defined by
    $$(w, \eta) \mapsto (w, \lambda \eta),$$
    for $w \in W$ and $\eta \in W^*$. We claim that $\widetilde \tau(s_\lambda) = s'_\lambda$. The map $\widetilde \tau(s_\lambda)$ is uniquely determined by requiring that
    $$b(s_\lambda(w, \eta), (w', \eta')) = b((w, \eta), \widetilde \tau(s_\lambda)(w', \eta'))$$
    for all $(w, \eta), (w', \eta') \in \widetilde W$. A quick computation shows that 
    $$b(s_\lambda(w, \eta), (w', \eta')) = b((w, \eta), s'_{\tau(\lambda)}(w', \eta'))$$
    for all $(w, \eta), (w', \eta') \in \widetilde W$, so $\widetilde \tau(s_\lambda) = s'_{\tau(\lambda)} = s'_\lambda$. Hence we have $\widetilde \tau(s_\lambda)s_\lambda = s'_{\lambda}s_\lambda$, and this is equal to multiplication by $\lambda$.
\end{proof}

If we summarize the results so far for the case of an isomorphism between factors of abelian origin induced by a derived equivalence, we get the following:
\begin{corollary}
    \label{cor: psi_i for ab var}
    Let $X = (\prod_{i=0}^k X_i)/G$ and $Y = (\prod_{j=0}^l Y_j)/H$ be smooth projective varieties over $\mathbb{C}$ as in Setup \ref{stp: symplectic var}, admitting a holomorphic symplectic form and satisfying condition \ref{condition *}. Let $\Phi \colon D^b(X) \to D^b(Y)$ be a triangulated equivalence. For every $\alpha$, let $\Phi_\alpha^{\llv} \colon \mathfrak{g}_\alpha \to \mathfrak{h}_\alpha$ be the induced isomorphism between the $\alpha$-th simple factors of the LLV algebras. Take an $\alpha$ such such that $\mathfrak{g}_\alpha \cong \mathfrak{s}(\End_A(\widetilde V), \widetilde \sigma)$ and $\mathfrak{h}_\alpha \cong \mathfrak{s}(\End_B(\widetilde W), \widetilde \tau)$ for certain isotypic components $V \subseteq \sH^1(X_0; \mathbb{Q})$ and $W \subseteq \sH^1(Y_0; \mathbb{Q})$ of the representations of $G$ and $H$, and simple algebras with involution $A$ and $B$ (where $A$ and $B$ are simple factors of $\mathbb{Q}[G]$ and $\mathbb{Q}[H]$ respectively), and assume that $\mathfrak{g}_\alpha$ is not of type $A_1$. Then there exists an isomorphism $\beta \colon A \to B$ respecting the involution and an isomorphism of vector spaces $f \colon \widetilde V \to \widetilde W$, equivariant with respect to $\beta$, such that $f$ preserves the Hodge structure and satisfies $\Phi_\alpha^{\llv} = \Ad(f)$. Furthermore, there exists an element $\psi_\alpha \in \Iso(\widetilde B, \widetilde \tau)$ (where $\widetilde B = \End_B(\widetilde W)$) such that $\psi_\alpha \circ f$ preserves the Hodge bigrading.
\end{corollary}
\begin{proof}
    It suffices to show that the conditions of Theorem \ref{thm: involution morphism from LLV isom}, Lemma \ref{lem: f : widetilde V to widetilde W respects HS} and Lemma \ref{lem: abelian factor isometry psi} are satisfied. The Lie algebra $\mathfrak{s}(\widetilde A, \widetilde \sigma)$ is of the right type because $X$ satisfies condition \ref{condition *}. To see that $\dim_\mathbb{Q} V = \dim_\mathbb{Q} W$, we will look at the Verbitsky component.

    The Verbitsky component $\SH_\alpha(X)$ is isomorphic to a subrepresentation of $(\bigwedge^\bullet V)^G$, and is supported in degrees $0$ up to $\dim_\mathbb{Q}(V)$. Hence its Hodge structure also has width $\dim_\mathbb{Q} V$ by Proposition \ref{prop: HS of SH_i is wide}. Similarly, the Hodge structure on $\SH_\alpha(Y)$ has width $\dim_\mathbb{Q} W$. Since the isomorphism $\Phi_\alpha^{\SH} \colon \SH_\alpha(X) \to \SH_\alpha(Y)$ respects the columns of the Hodge diamond, we conclude that $V$ and $W$ have the same dimension.

    To see that $\Phi_\alpha^{\llv}$ sends $h'_{\widetilde V}$ to $h'_{\widetilde W}$, observe that $h'_{\widetilde V} = h'_{X, \alpha}$ and $h'_{\widetilde W} = h'_{Y, \alpha}$ by Lemma \ref{lem: equality of Hodge operators}. The fact that $\Phi^{\sH}$ respects the Hodge structure implies that $\Phi_\alpha^{\llv}(h'_{X,\alpha})= h'_{Y,\alpha}$, so we indeed have $\Phi_\alpha^{\llv}(h'_{\widetilde V}) = h'_{\widetilde W}$. Therefore we can indeed apply Lemma~\ref{lem: abelian factor isometry psi} to obtain the desired map $\psi_\alpha$.
\end{proof}

\subsection{Factors of type $A_1$}
\label{sec: factorwise psi A1}
In this section, we prove Theorem \ref{thm: construction of psi} for simple factors of $\llv(X; \mathbb{Q})$ of type $A_1$. We will first see that the component of the Hodge operator in such a simple factor vanishes. We then only have to deal with the grading, and this can be done with some Lie theory.

\begin{lemma}
    \label{lem: simple factor A1 h' = 0}
    Let $X$ be a smooth projective variety as in Setup \ref{stp: symplectic var}, admitting a holomorphic symplectic form, and let $\mathfrak{g}_\alpha$ be a simple factor of $\llv(X; \mathbb{Q})$ of type $A_1$. Let $h' \in \mathfrak{g}_\alpha$ be the component in $\mathfrak{g}_\alpha$ of the Hodge operator of $X$. Then $h' = 0$.
\end{lemma}
\begin{proof}
    Write $X = (\prod_{i=0}^k X_i)/G$, and let $F$ be the centroid of $\mathfrak{g}_\alpha$. Then $\mathfrak{g}_\alpha \cong \mathfrak{sl}_2(F)$ by Lemma~\ref{lem: simple factor A1 of llv is sl2(F)}. Any simple factor of $\llv(X; \mathbb{Q})$ coming from a hyperkähler factor $X_i$ is of the form $\mathfrak{so}(\widetilde \sH(X_i; \mathbb{Q})^{G_i})$, and the dimension of $\widetilde \sH(X_i; \mathbb{Q})^{G_i}$ is at least $5$ by Lemma~\ref{lem: HK H^2 invariants dimension}. This implies that our simple factor $\mathfrak{g}_\alpha$ must be of abelian origin.

    Theorem \ref{thm: summary LLV computation} implies that $\mathfrak{g}_\alpha$ is of the form $\llv_G(V; \mathbb{Q}) := \llv((\extp^\bullet V)^G; \mathbb{Q})$ for some isotypic component $V \subseteq \sH^1(X_0; \mathbb{Q})$. Since $X$ is projective, $X_0$ is also projective, and therefore $\sH^1(X_0; \mathbb{Q})$ is a polarizable Hodge structure of weight $1$. Since $V \subseteq \sH^1(X_0; \mathbb{Q})$ is a sub Hodge structure, it is polarizable as well. In Proposition \ref{prop: category of polarizable G-equiv HS is semisimple}, we saw that $V$ therefore has a $G$-equivariant polarization, and this gives an algebraic element of $(\extp^2 V)^G$. In particular, we see that $(\extp^2 V)^G \cap (\extp^2 V)^{1,1} \neq 0$. Let $y \in (\extp^2 V)^G \cap (\extp^2 V)^{1,1}$ be such a nonzero element.

    Let $\mathfrak{g}_{\alpha, 2}$ be the degree $2$ part of $\mathfrak{g}_\alpha$, then there is an isomorphism $\mathfrak{g}_{\alpha, 2} \cong (\extp^2 V)^G$ which preserves the Hodge structure. Moreover, since $\mathfrak{g}_\alpha$ is of type $A_1$ with centroid $F$, every element $x \in \mathfrak{g}_{\alpha, 2}$ is of the form $x = \lambda y$ for some $\lambda \in F$.

    We now claim that for every $\lambda \in F$, multiplication by $\lambda$ as a map $\mathfrak{g}_\alpha \to \mathfrak{g}_\alpha$ preserves the Hodge structure. To see why, recall that $\mathfrak{g}_\alpha$ is isomorphic to a Lie subalgebra of $\mathfrak{so}(\widetilde V)$, see diagram \eqref{eq: j pi diagram} and Proposition \ref{prop: inclusion llv_Gpre in s(widetilde A)}. We know that $V \subseteq \sH^1(X_0; \mathbb{Q})$ is an isotypic component. Let $A$ be the simple factor of $\mathbb{Q}[G]$ associated to this isotypic component. The action of $G$ on $V$ gives $\widetilde V$ the structure of an $A$-module, and through this $\End_\mathbb{Q}(\widetilde V)$ also becomes an $A$-module (defined for $f \in \End_\mathbb{Q}(\widetilde V)$ and $a \in A$ by $(a \cdot f)(v) = af(\sigma(a)v)$ for all $v \in \widetilde V$, where $\sigma$ is the involution on $A$). Since the action of $G$ on $V$ preserves the Hodge structure, the action of every $a \in A$ on $V$ preserves the Hodge structure. In particular, the action of $F = Z(A)^\sigma$ on $\End_\mathbb{Q}(\widetilde V)$ preserves the Hodge structure, and this implies that for every $\lambda \in F$, multiplication by $\lambda$ as a map $\mathfrak{g}_\alpha \to \mathfrak{g}_\alpha$ preserves the Hodge structure, which proves the claim.
    
    It now follows that all elements of $\mathfrak{g}_{\alpha, 2}$ are algebraic, and therefore all elements of $\mathfrak{g}_\alpha$ are algebraic, which implies that $h' = 0$.
\end{proof}

\begin{lemma}
    \label{lem: automorphism of sl2 triples}
    Let $F$ be a field of characteristic $0$, let $\mathfrak{g}_\alpha \cong \mathfrak{sl}_2(F)$ and let $\langle e_1, h_1, f_1\rangle$ and $\langle e_2, h_2, f_2\rangle$ be two $\mathfrak{sl}_2$-triples in $\mathfrak{g}_\alpha$. Then the $F$-linear map $\mathfrak{g}_\alpha \to \mathfrak{g}_\alpha$ which sends $e_1$ to $e_2$, sends $h_1$ to $h_2$ and sends $f_1$ to $f_2$ is a Lie algebra automorphism of $\mathfrak{g}_\alpha$.
\end{lemma}
\begin{proof}
    This follows immediately from the definition of an $\mathfrak{sl}_2$-triple.
\end{proof}

By using these two results, we can prove Theorem \ref{thm: construction of psi} for simple factors of type $A_1$.
\begin{corollary}
    \label{cor: psi_i for type A1}
    Let $X$ and $Y$ be two smooth projective varieties over $\mathbb{C}$, satisfying condition \ref{condition *}. Let $\Phi \colon D^b(X) \to D^b(Y)$ be a triangulated equivalence, let $\mathfrak{g}_\alpha$ be a simple factor of $\llv(X; \mathbb{Q})$ of type $A_1$ and let $\mathfrak{h}_\alpha$ be the simple factor of $\llv(Y; \mathbb{Q})$ such that $\Phi^{\llv}$ induces an isomorphism $\Phi^{\llv}_\alpha \colon \mathfrak{g}_\alpha \to \mathfrak{h}_\alpha$. Let $F$ be the centroid of $\mathfrak{g}_\alpha$. Then there is a $\psi_\alpha \in \GL_2(F)$ such that $\Ad(\psi_\alpha) \circ \Phi^{\llv}_\alpha \colon \mathfrak{g}_\alpha \to \mathfrak{h}_\alpha$ preserves the Hodge bigrading.
\end{corollary}
\begin{proof}
    First observe that $\mathfrak{h}_\alpha \cong \mathfrak{sl}_2(F)$ by Lemma \ref{lem: simple factor A1 of llv is sl2(F)}. Moreover, the Hodge operators in $\mathfrak{g}_\alpha$ and $\mathfrak{h}_\alpha$ are $0$ by Lemma \ref{lem: simple factor A1 h' = 0}.

    Let $e_X \in \mathfrak{g}_{\alpha, 2}$ be an element with the hard Lefschetz property, which gives an $\mathfrak{sl}_2$-triple $\langle e_X, h_X, f_X \rangle$ over $F$, where $h_X$ is the grading operator on $\mathfrak{g}_\alpha$. Similarly, an element $e_Y \in \mathfrak{h}_{\alpha, 2}$ with the hard Lefschetz property gives an $\mathfrak{sl}_2$-triple $\langle e_Y, h_Y, f_Y\rangle$ which is an $F$-basis of $\mathfrak{h}_\alpha$. By Lemma \ref{lem: automorphism of sl2 triples}, there exists an automorphism $\Psi$ of $\mathfrak{h}_\alpha \cong \mathfrak{sl}_2(F)$ which sends the $\mathfrak{sl}_2$-triple $\langle \Phi^{\llv}_\alpha(e_X), \Phi^{\llv}_\alpha(h_X), \Phi^{\llv}_\alpha(f_X)\rangle$ to $\langle e_Y, h_Y, f_Y\rangle$. It follows that $\Psi \circ \Phi^{\llv}_\alpha \colon \mathfrak{g}_\alpha \to \mathfrak{h}_\alpha$ respects the Hodge bigrading.

    Since the Dynkin diagram of $\mathfrak{sl}_2(F)$ is trivial, its group of automorphisms is isomorphic to $\PGL_2(F)$. There is a surjection $\GL_2(F) \to \PGL_2(F)$, so there is a $\psi_\alpha \in \GL_2(F)$ such that $\Psi = \Ad(\psi_\alpha)$. 
\end{proof}

\subsection{Factors of hyperkähler origin}
\label{sec: factorwise psi HK}
Now suppose $\Phi_\alpha^{\llv} \colon \mathfrak{g}_\alpha \to \mathfrak{h}_\alpha$ is an isomorphism of simple factors of hyperkähler origin. Recall that $X = (\prod_{i=0}^l X_i)/G$ and $Y = (\prod_{j=0}^s Y_j)/H$. There are hyperkähler varieties $X_j$ and $Y_k$ such that $\mathfrak{g}_\alpha \cong \mathfrak{so}(\widetilde \sH(X_j; \mathbb{Q})^{G_j})$ and $\mathfrak{h}_\alpha \cong \mathfrak{so}(\widetilde \sH(Y_k; \mathbb{Q})^{H_k})$. Let $[j] \subseteq \{1, \dots, l\}$ be the $G$-orbit of $j$ and let $[k] \subseteq \{1, \dots, s\}$ be the $H$-orbit of $k$. Recall that the decomposition $\llv(X; \mathbb{Q}) = \prod_\beta \mathfrak{g}_\beta$ induces a factorization $\SH(X; \mathbb{Q}) = \bigotimes_\beta \SH_\beta(X)$, where each $\SH_\beta(X)$ is a representation of $\mathfrak{g}_\beta$.

The statement and proof of the following theorem are analogous to \cite[Proposition~4.9]{TaelmanDerivedEquivalences}.

\begin{theorem}
    \label{thm: hyperkahler factor similitude}
    Suppose we are in the setting as described above. Let $d_j, d_k \in \mathbb{N}$ be defined by $\dim X_j = 2d_j$ and $\dim Y_k = 2d_k$. We then have an equality $\#[j]d_j = \#[k]d_k$ and there exist an element $\lambda \in \mathbb{Q}^*$ and a similitude $\varphi \colon \widetilde \sH(X_j;\mathbb{Q})^{G_j} \to \widetilde \sH(Y_k;\mathbb{Q})^{H_k}$ preserving the Hodge grading such that the diagram
    \begin{equation}
    \label{eq: SH alpha diagram}
        \begin{tikzcd}
            \SH_\alpha(X) \arrow[rr, "\Phi_\alpha^{\SH}"] \arrow[d] & & \SH_\alpha(Y) \arrow[d]\\
            \Sym^{\#[j]\cdot d_j} \widetilde \sH(X_j; \mathbb{Q})^{G_j} \arrow[rr, "\lambda \Sym^{\#[j]d_j}\varphi"] & & \Sym^{\#[k]\cdot d_k} \widetilde \sH(Y_k; \mathbb{Q})^{H_k}
        \end{tikzcd}
    \end{equation}
    commutes, where the vertical arrows are the inclusions from Theorem \ref{thm: computation of SH of HK factor}.
\end{theorem}
Before we prove this, we recall a bit of notation from \cite[Section 4]{TaelmanDerivedEquivalences}. Given a quadratic space $V$ and an integer $d$, there is a contraction operator $\Delta \colon \Sym^d V \to \Sym^{d-2}V$, and this is a morphism of $\mathfrak{so}(V)$-representations with irreducible kernel (see Lemma \ref{lem: SH is S_{[d]}}). We write $S_d V$ for the kernel of $\Delta \colon \Sym^d V \to \Sym^{d-2}V$.

\begin{proof}[{Proof of Theorem \ref{thm: hyperkahler factor similitude}}]
    With the notation just introduced, Lemma \ref{lem: SH is S_{[d]}} says that $\SH_\alpha(X) \cong S_{\#[j]d_j} \widetilde \sH(X_j;\mathbb{Q})^{G_j}$ and $\SH_\alpha(Y) \cong S_{\#[k]d_k} \widetilde \sH(Y_k; \mathbb{Q})^{H_k}$. The space $S_{\#[k]d_k} \widetilde \sH(Y_k;\mathbb{Q})^{H_k}$ becomes a representation of $\mathfrak{g}_\alpha$ via $\Phi^{\llv}_\alpha$, and then $\Phi^{\SH}_\alpha$ gives an isomorphism of $\mathfrak{g}_\alpha$-representations 
    \[S_{\#[j]d_j}\widetilde \sH(X_j;\mathbb{Q})^{G_j} \isomto S_{\#[k]d_k} \widetilde \sH(Y_k;\mathbb{Q})^{H_k}.\]
    The spaces $\widetilde \sH(X_j; \mathbb{Q})^{G_j}$ and $\widetilde \sH(Y_k; \mathbb{Q})^{H_k}$ have the same dimension, because $\Phi^{\llv}_\alpha$ is an isomorphism $\mathfrak{so}( \widetilde \sH(X_j; \mathbb{Q})^{G_j}) \to \mathfrak{so}(\widetilde \sH(Y_k; \mathbb{Q})^{H_k})$. Therefore we must have $\#[j]d_j = \#[k]d_k$. 
    
    Since $\Phi^{\SH}_\alpha$ is equivariant with respect to $\Phi^{\llv}_\alpha$, we have 
    $$\Phi^{\SH}_\alpha \mathfrak{so}(\widetilde \sH(X_j;\mathbb{Q})^{G_j}) (\Phi^{\SH}_\alpha)^{-1} = \mathfrak{so}(\widetilde \sH(Y_k;\mathbb{Q})^{H_k})$$
    as subspaces of $\End(\SH_\alpha(Y))$, and then \cite[Proposition 4.3]{TaelmanDerivedEquivalences} gives the desired similitude $\varphi$ and scalar $\lambda$. 

    It remains to show that $\varphi$ respects the Hodge grading. For this, we can use the same strategy as in \cite[Proposition 4.9]{TaelmanDerivedEquivalences}. By using the Deligne torus, the Hodge structure on $\widetilde\sH(X_j; \mathbb{Q})^{G_j}$ is determined by a homomorphism $\mathbb{C}^\times \to \sO(\widetilde\sH(X_j; \mathbb{R})^{G_j})$, while the Hodge structure on $\SH_\alpha(X)$ is determined by a homomorphism $\mathbb{C}^\times \to \GL(\SH_\alpha(X) \otimes_\mathbb{Q} \mathbb{R})$, and similarly for $Y$. The Hodge gradings are then determined by the restrictions of these morphisms to the unit circle $S^1 \subseteq \mathbb{C}^\times$. Using Lemma \ref{lem: SH is S_{[d]}}, every isometry of $\widetilde \sH(X_j; \mathbb{Q})^{G_j}$ gives an automorphism of $\SH_\alpha(X)$, so we have a group homomorphism $\sO(\widetilde\sH(X_j; \mathbb{R})^{G_j}) \to \GL(\SH_\alpha(X) \otimes_\mathbb{Q} \mathbb{R})$, and this commutes with the maps from $S^1$ by a proof analogous to \cite[Lemma 4.8]{TaelmanDerivedEquivalences}. The kernel of $\sO(\widetilde \sH(X_j; \mathbb{R})^{G_j}) \to \GL(\SH_\alpha(X) \otimes_\mathbb{Q} \mathbb{R})$ is either trivial or $\{\pm 1\}$, depending on whether $\#[j]d_j$ is odd or even. Putting everything together, we have a diagram
    \begin{equation*}
        \begin{tikzcd}
            & S^1 \arrow[ld] \arrow[rd] & \\
            \sO(\widetilde \sH(X_j; \mathbb{R})^{G_j}) \arrow[rr, "\Ad(\varphi)"] \arrow{d} && \sO(\widetilde \sH(Y_k; \mathbb{R})^{H_k}) \arrow{d}\\
            \GL(\SH_\alpha(X) \otimes_\mathbb{Q} \mathbb{R}) \arrow[rr, "\Ad(\Phi_\alpha^{\SH})"] && \GL(\SH_\alpha(Y) \otimes_\mathbb{Q} \mathbb{R}).
        \end{tikzcd}
    \end{equation*}
    The outer diagram commutes, since $\Phi_\alpha^{\SH}$ respects the Hodge grading. Using that~\eqref{eq: SH alpha diagram} commutes, we see that the bottom square also commutes. Therefore the top triangle commutes up to a map $S^1 \to \ker(\sO(\widetilde \sH(Y_k;\mathbb{R})) \to \GL(\SH_\alpha(Y) \otimes_\mathbb{Q} \mathbb{R}))$. Since this kernel is discrete, the top triangle must commute, which means that $\varphi$ respects the Hodge grading.
\end{proof}

\begin{corollary}
    \label{cor: psi_i for HK}
    Under the same assumptions as above, there exists an isometry $\psi_\alpha \in \SO(\widetilde \sH(Y_k; \mathbb{Q})^{H_k})$ preserving the Hodge grading such that $\psi_\alpha\varphi \colon \widetilde \sH(X_i; \mathbb{Q})^{G_i} \to \widetilde \sH(Y_k; \mathbb{Q})^{H_k}$ is a similitude which respects the Hodge bigrading.
\end{corollary}
\begin{proof}
    This is exactly the same as the first part of the proof of \cite[Theorem~5.3]{TaelmanDerivedEquivalences}. 
\end{proof}

\subsection{Summary}
\label{sec: factorwise psi summary}
In this section, we give the proof of Theorem \ref{thm: construction of psi}, summarizing the results obtained in this chapter. Recall from Theorem \ref{thm: summary LLV computation} that we have a decomposition $\llv(Y; \mathbb{Q}) \cong \prod_\alpha \mathfrak{h}_\alpha$ into simple factors. Every simple factor is of the form $\mathfrak{so}(\widetilde \sH(Y_j; \mathbb{Q})^{H_j})$ for some hyperkähler factor $Y_j$, of the form $\mathfrak{s}(\widetilde B, \widetilde \tau)$ for an algebra with involution $(\widetilde B, \widetilde \tau)$, or of type $A_1$.

\begin{proof}[{Proof of Theorem \ref{thm: construction of psi}}]
    The equivalence $\Phi \colon D^b(X) \to D^b(Y)$ induces an isomorphism $\Phi^{\llv} \colon \llv(X; \mathbb{Q}) \to \llv(Y; \mathbb{Q})$ by \cite[Theorem A]{TaelmanDerivedEquivalences}. Theorem~\ref{thm: split Lie alg reps} implies that, after reordering the simple factors of $\llv(Y; \mathbb{Q}) \cong \prod_\beta \mathfrak{h}_\beta$, the map $\Phi^{\llv}$ is a direct product of isomorphisms $\Phi^{\llv}_\alpha \colon \mathfrak{g}_\alpha \to \mathfrak{h}_\alpha$. By Proposition~\ref{prop: matching of LLV factors}, for every $\alpha$ the factors $\mathfrak{g}_\alpha$ and $\mathfrak{h}_\alpha$ are both of hyperkähler or both of abelian origin. If they are of abelian origin they are either of type $A_1$ or both come from an algebra with involution, and Proposition~\ref{prop: matching of LLV factors} also says that in the latter case those involutions have the same type.
    
    For every $\alpha$, define the group $\mathcal{H}_\alpha$ by:
    \begin{equation*}
        \mathcal{H}_\alpha = \begin{cases}
            \uSO(\widetilde \sH(Y_j; \mathbb{Q})^{H_j}) & \text{if } \mathfrak{h}_\alpha \cong \mathfrak{so}(\widetilde \sH(Y_j; \mathbb{Q})^{H_j})\\
            \sR_{F/\mathbb{Q}}\uGL_2(F)  & \text{if }\mathfrak{h}_\alpha \cong \mathfrak{sl}_2(F)\\
            \sR_{F/\mathbb{Q}}\uIso(\widetilde B, \widetilde \tau) & \text{if } \mathfrak{h}_\alpha \cong \mathfrak{s}(\widetilde B, \widetilde \tau),
        \end{cases}
    \end{equation*}
    where $F = Z(\widetilde B)^{\widetilde \tau}$ in the third case. Then Corollaries~\ref{cor: psi_i for ab var}, \ref{cor: psi_i for type A1} and~\ref{cor: psi_i for HK} give the desired elements $\psi_\alpha \in \mathcal{H}_\alpha(\mathbb{Q})$.
\end{proof}

%% file: chapters/chapter7.tex
As in Chapter \ref{chap: fatcorwise psi}, we continue with a derived equivalence $\Phi \colon D^b(X) \to D^b(Y)$, where $X$ and $Y$ are smooth projective varieties over $\mathbb{C}$ admitting a holomorphic symplectic form and satisfying condition \ref{condition *} from the introduction. Recall from Theorem \ref{thm: LLV alg derived invariant} that the equivalence $\Phi$ induces an isomorphism of LLV algebras $\Phi^{\llv} \colon \llv(X; \mathbb{Q}) \to \llv(Y; \mathbb{Q})$. In Theorem \ref{thm: construction of psi}, we constructed an algebraic group $\mathcal{H}$ over $\mathbb{Q}$ and an element $\psi \in \mathcal{H}(\mathbb{Q})$ such that $\Ad(\psi)$ induces an automorphism of $\llv(Y; \mathbb{Q})$ with the property that $\Ad(\psi) \circ \Phi^{\llv} \colon \llv(X; \mathbb{Q}) \to \llv(Y; \mathbb{Q})$ preserves the Hodge bigrading.

In this chapter, we will lift the action of $\psi$ on $\llv(Y; \mathbb{Q})$ to an automorphism $\Psi^{\sH}$ of $\sH^\bullet(Y; \mathbb{Q})$, such that the composition $\Psi^{\sH} \circ \Phi^{\sH} \colon \sH^\bullet(X; \mathbb{Q}) \to \sH^\bullet(Y; \mathbb{Q})$ preserves the Hodge bigrading, which proves Theorem \ref{thm: main theorem}. That is, we now have all solid arrows in the diagram below, and will construct the dashed arrow:
\begin{equation*}
    \begin{tikzcd}
        \sH^\bullet(X; \mathbb{Q}) \arrow[r, "\Phi^{\sH}"] \arrow[loop below] & \sH^\bullet(Y; \mathbb{Q})  \arrow[loop below] \arrow[r, dashed, "\Psi^{\sH}"] & \sH^\bullet(Y; \mathbb{Q}) \arrow[loop below]\\
        \llv(X; \mathbb{Q}) \arrow[r, "\Phi^{\llv}"] & \llv(Y; \mathbb{Q}) \arrow[r, "\Ad(\psi)"] & \llv(Y; \mathbb{Q}). 
    \end{tikzcd}
\end{equation*}

If $\mathcal{H}$ would be simply connected with $\Lie(\mathcal{H}) = \llv(Y; \mathbb{Q})$, then we could use Theorem \ref{thm: integrate Lie algebra action} to integrate the representation of $\llv(Y; \mathbb{Q})$ on $\sH^\bullet(Y; \mathbb{Q})$ to a representation of $\mathcal{H}$ on $\sH^\bullet(Y; \mathbb{Q})$. The induced action of $\psi \in \mathcal{H}(\mathbb{Q})$ on $\sH^\bullet(Y; \mathbb{Q})$ would then give the desired automorphism $\Psi^{\sH}$ of $\sH^\bullet(Y; \mathbb{Q})$. Unfortunately, the above idea will not work for two reasons. On the one hand, $\Lie(\mathcal{H})$ can be bigger than $\llv(Y; \mathbb{Q})$. The other problem is that the semisimple part of $\mathcal{H}$ is not simply connected, so we have to pass to a finite covering of the semisimple part to make this idea work.

To fix this, we will define an auxiliary group $\Gamma\mathcal{H}$ acting on $\sH^\bullet(Y; \mathbb{Q})$ and an element $\psi' \in \Gamma\mathcal{H}(\mathbb{Q})$ constructed from $\psi$. We will then define $\Psi^{\sH} \colon \sH^\bullet(Y; \mathbb{Q}) \to \sH^\bullet(Y; \mathbb{Q})$ to be the action of $\psi'$.

\begin{example}
    Let $X$ be a hyperkähler variety. Then $\llv(X; \mathbb{Q}) \cong \mathfrak{so}(\widetilde \sH(X; \mathbb{Q}))$. The associated group $\mathcal{H}$ is $\uSO(\widetilde \sH(X; \mathbb{Q}))$, which is not simply connected. Let $\Gamma\mathcal{H} = \uGSpin(\widetilde \sH(X; \mathbb{Q}))$. Then $\Gamma\mathcal{H}$ is a reductive algebraic group over $\mathbb{Q}$, and its semisimple part is a simply connected covering of $\mathcal{H}$. These objects are already used in \cite[Chapter~5]{TaelmanDerivedEquivalences}. See Section~\ref{subsec: HK lifting} for more details, where we extend this to the simple factors of hyperkähler origin of the LLV algebra of a variety admitting a holomorphic symplectic form.
\end{example}

Recall the decomposition $\llv(Y; \mathbb{Q}) = \prod_\alpha \mathfrak{h}_\alpha$ as a product of simple factors from Theorem~\ref{thm: summary LLV computation}. In Theorem~\ref{thm: construction of psi}, we constructed $\mathcal{H}$ and $\psi$ by constructing a group $\mathcal{H}_\alpha$ and element $\psi_\alpha \in \mathcal{H}_\alpha(\mathbb{Q})$ for every $\alpha$, and then defined $\mathcal{H} = \prod_\alpha \mathcal{H}_\alpha$ and $\psi = (\psi_\alpha)_\alpha$. We will also construct $\Gamma\mathcal{H}$ and $\psi'$ factorwise, and then let $\Gamma\mathcal{H} = \prod_\alpha \Gamma\mathcal{H}_\alpha$ and $\psi' = (\psi'_\alpha)_\alpha$. To be precise, in this chapter we will prove the following theorem:

\begin{theorem}
    \label{thm: properties of GammaH_i}
    Let $X$ and $Y$ be smooth connected projective varieties over $\mathbb{C}$, admitting a holomorphic symplectic form and satisfying condition \ref{condition *}, and let $\Phi \colon D^b(X) \to D^b(Y)$ be an equivalence of triangulated categories. Let $\psi = (\psi_\alpha)_\alpha \in \mathcal{H}(\mathbb{Q}) = \prod_\alpha \mathcal{H}_\alpha(\mathbb{Q})$ be the element constructed in Theorem \ref{thm: construction of psi}. Then for every~$\alpha$, there exist:
    \begin{enumerate}
        \item A reductive group $\Gamma\mathcal{H}_\alpha$ over $\mathbb{Q}$.
        \item Inclusions $\mathfrak{h}_\alpha \subseteq \Lie(\mathcal{H}_\alpha) \subseteq \Lie(\Gamma\mathcal{H}_\alpha)$.
        \item A representation of $\Gamma\mathcal{H}_\alpha$ on $\sH^\bullet(Y; \mathbb{Q})$.
        \item An element $\psi'_\alpha \in \Gamma\mathcal{H}_\alpha(\mathbb{Q})$.
    \end{enumerate}
    These satisfy:
    \begin{enumerate}
        \item For every $\alpha$, the induced representation of $\Lie(\Gamma\mathcal{H}_\alpha)$ on $\sH^\bullet(Y; \mathbb{Q})$ extends the representation of $\mathfrak{h}_\alpha$ on $\sH^\bullet(Y; \mathbb{Q})$.
        \item The representations of the groups $\Gamma\mathcal{H}_\alpha$ on $\sH^\bullet(Y; \mathbb{Q})$ mutually commute.
        \item For every $\alpha$, the automorphism $\Ad(\psi'_\alpha)$ of $\Lie(\Gamma\mathcal{H}_\alpha)$ restricts to the automorphism $\Ad(\psi_\alpha)$ on $\mathfrak{h}_\alpha$.
    \end{enumerate}
\end{theorem}

In particular, if we denote the inclusion $\llv(Y; \mathbb{Q}) \to \mathfrak{gl}(\sH^\bullet(Y; \mathbb{Q}))$ by $\iota$ and write $j_\alpha \colon \mathfrak{h}_\alpha \to \Lie(\Gamma\mathcal{H}_\alpha)$ for the inclusion from the second part of the theorem, we obtain a commutative diagram
\begin{equation*}
    \begin{tikzcd}
        \llv(Y; \mathbb{Q}) = \prod_\alpha \mathfrak{h}_\alpha \arrow[rr, "\iota"] \arrow[dr, "\prod_\alpha j_\alpha"'] &&\mathfrak{gl}(\sH^\bullet(Y; \mathbb{Q})).\\
        &\prod_{\alpha} \Lie(\Gamma\mathcal{H}_\alpha) \arrow[ur]
    \end{tikzcd}
\end{equation*}
For every $\alpha$ the Lie algebra $\Lie(\Gamma\mathcal{H}_\alpha)$ will be the direct product of $\mathfrak{h}_\alpha$ and an abelian Lie algebra. In particular, $\mathfrak{h}_\alpha$ is the semisimple part of $\Lie(\Gamma\mathcal{H}_\alpha)$. We will see that for all $\alpha$, the Lie algebra $\Lie(\mathcal{H}_\alpha)$ will either be equal to $\mathfrak{h}_\alpha$ or $\Lie(\Gamma\mathcal{H}_\alpha)$.

Before we prove the above theorem, let us see why it implies Theorem \ref{thm: main theorem}:

\begin{proof}[{Proof of Theorem \ref{thm: main theorem}}]
    \label{pf: main theorem}
    By \cite[Theorem A]{TaelmanDerivedEquivalences}, the equivalence $\Phi \colon D^b(X) \to D^b(Y)$ induces an isomorphism $\Phi^{\llv} \colon \llv(X; \mathbb{Q}) \to \llv(Y; \mathbb{Q})$. Theorem~\ref{thm: split Lie alg reps} implies that $\Phi^{\llv}$ splits as a direct product of isomorphisms $\Phi^{\llv}_\alpha \colon \mathfrak{g}_\alpha \to \mathfrak{h}_\alpha$ between the simple factors of $\llv(X; \mathbb{Q})$ and $\llv(Y; \mathbb{Q})$.

    Let $\Gamma\mathcal{H}_\alpha$ and $\psi_\alpha' \in \Gamma\mathcal{H}_\alpha(\mathbb{Q})$ be the groups and elements constructed in Theorem~\ref{thm: properties of GammaH_i}, and define $\Gamma\mathcal{H} = \prod_\alpha \Gamma\mathcal{H}_\alpha$ and $\psi' = (\psi'_\alpha)_\alpha \in \Gamma\mathcal{H}(\mathbb{Q})$. The representations of the groups $\Gamma\mathcal{H}_\alpha$ on $\sH^\bullet(Y; \mathbb{Q})$ mutually commute, so they give a representation of $\Gamma\mathcal{H}$ on $\sH^\bullet(Y; \mathbb{Q})$. Write $\Psi^{\sH} \colon \sH^\bullet(Y; \mathbb{Q}) \to \sH^\bullet(Y; \mathbb{Q})$ for the automorphism of $\sH^\bullet(Y; \mathbb{Q})$ induced by the action of the element $\psi' \in \Gamma\mathcal{H}(\mathbb{Q})$ under the representation of $\Gamma\mathcal{H}$ on $\sH^\bullet(Y;\mathbb{Q})$. We are done once we show that $\Psi^{\sH} \circ \Phi^{\sH}$ preserves the Hodge bigrading. To do this, we will show that the map 
    \[ \Ad(\Psi^{\sH} \circ \Phi^{\sH}) \colon \End(\sH^\bullet(X; \mathbb{Q})) \to \End(\sH^\bullet(Y; \mathbb{Q})) \]
    sends $h_X$ to $h_Y$ and $h'_X$ to $h'_Y$.

    We first claim that $\Ad(\Psi^{\sH} \circ \Phi^{\sH})$ restricts to the isomorphism 
    \[\Ad(\psi) \circ \Phi^{\llv} \colon \llv(X; \mathbb{Q}) \to \llv(Y; \mathbb{Q}),\]
    where $\psi = (\psi_\alpha)_\alpha \in \prod \mathcal{H}_\alpha(\mathbb{Q})$ is the element constructed in Theorem \ref{thm: construction of psi}. By \cite[Theorem~A]{TaelmanDerivedEquivalences}, the map $\Ad(\Phi^{\sH})$ restricts to the isomorphism of LLV algebras $\Phi^{\llv} \colon \llv(X; \mathbb{Q}) \to \llv(Y; \mathbb{Q})$, so it suffices to show that $\Ad(\Psi^{\sH})$ restricts to the automorphism $\Ad(\psi)$ of $\llv(Y; \mathbb{Q})$. Theorem~\ref{thm: properties of GammaH_i} gives the following commutative diagram:
    \begin{equation*}
        \begin{tikzcd}
            \llv(Y; \mathbb{Q}) = \prod_\alpha \mathfrak{h}_\alpha \arrow[d, "\Ad((\psi_\alpha)_\alpha)"] \arrow[r] & \Lie(\prod_\alpha \Gamma\mathcal{H}_\alpha) \arrow[r] \arrow[d, "\Ad((\psi_\alpha')_\alpha)"] & \mathfrak{gl}(\sH^\bullet(Y; \mathbb{Q})) \arrow[d, "\Ad(\Psi^{\sH})"]\\
            \llv(Y; \mathbb{Q}) = \prod_i \mathfrak{h}_\alpha \arrow[r] & \Lie(\prod_\alpha \Gamma\mathcal{H}_\alpha) \arrow[r] & \mathfrak{gl}(\sH^\bullet(Y; \mathbb{Q})),
        \end{tikzcd}
    \end{equation*}
    where the horizontal compositions are the inclusion $\iota \colon \llv(Y; \mathbb{Q}) \to \mathfrak{gl}(\sH^\bullet(Y; \mathbb{Q}))$. Therefore $\Ad(\Psi^{\sH})$ indeed restricts to $\Ad(\psi) = (\Ad(\psi_\alpha))_\alpha$.
    
    Theorem \ref{thm: construction of psi} then gives that $\Ad(\psi) \circ \Phi^{\llv}$ sends $h_X$ to $h_Y$ and $h'_X$ to $h'_Y$, which finishes the proof.
\end{proof}

The remainder of this chapter is dedicated to proving Theorem \ref{thm: properties of GammaH_i}. In sections~\ref{subsec: HK lifting} and~\ref{sec: abelian type groups}, we will define and study the group $\Gamma\mathcal{H}_\alpha$ for the different types of simple factors $\mathfrak{h}_\alpha$ that occur in $\llv(Y; \mathbb{Q})$. Throughout these sections, we will also provide the elements $\psi_\alpha' \in \Gamma\mathcal{H}_\alpha(\mathbb{Q})$, and show that the automorphism $\Ad(\psi'_\alpha)$ of $\Lie(\Gamma\mathcal{H}_\alpha)$ restricts to the automorphism $\Ad(\psi_\alpha)$ on $\mathfrak{h}_\alpha$. In section~\ref{subsec: proof of thm properties of GammaH_i}, we will then construct the representations of the $\Gamma\mathcal{H}_\alpha$ on $\sH^\bullet(Y; \mathbb{Q})$ and show that they mutually commute.

\section{Lifting for factors of hyperkähler origin}
\label{subsec: HK lifting}
Throughout this section, we assume that $\mathfrak{h}_\alpha \cong \mathfrak{so}(\widetilde \sH(Y_k; \mathbb{Q})^{H_k})$ for some hyperkähler factor $Y_k$ of $Y = (\prod_j Y_j)/H$. In Corollary \ref{cor: psi_i for HK}, we constructed a Hodge isometry $\psi_\alpha \in \SO(\widetilde \sH(Y_k; \mathbb{Q})^{H_k})$ such that $\Ad(\psi_\alpha) \circ \Phi^{\llv}_\alpha \colon \mathfrak{g}_\alpha \to \mathfrak{h}_\alpha$ respects the Hodge bigrading. We have $\mathcal{H}_\alpha = \uSO(\widetilde \sH(Y_k; \mathbb{Q})^{H_k})$, and $\Lie(\mathcal{H}_\alpha) = \mathfrak{h}_\alpha$.
\begin{definition}
    Assume that $\mathfrak{h}_\alpha \cong \mathfrak{so}(\widetilde \sH(Y_k; \mathbb{Q})^{H_k})$. We then define $\Gamma\mathcal{H}_\alpha = \uGSpin(\widetilde \sH(Y_k;\mathbb{Q})^{H_k})$. 
\end{definition}

\begin{lemma}
    \label{lem: HK Spin short exact sequences}
     There exist short exact sequences of algebraic groups
     \begin{equation*}
        1 \to \mu_{2, \mathbb{Q}} \to \mathbb{G}_{m, \mathbb{Q}} \times \uSpin(\widetilde \sH(Y_k;\mathbb{Q})^{H_k}) \to \uGSpin(\widetilde \sH(Y_k;\mathbb{Q})^{H_k}) \to 1
    \end{equation*}
    and
    \begin{equation*}
        1 \to \mathbb{G}_{m, \mathbb{Q}} \to \uGSpin(\widetilde \sH(Y_k;\mathbb{Q})^{H_k}) \to \uSO(\widetilde \sH(Y_k;\mathbb{Q})^{H_k}) \to 1.
    \end{equation*}
\end{lemma}
\begin{proof}
    By Lemma \ref{lem: exactness algebraic closure} it suffices to check the exactness for the $\overline{\mathbb{Q}}$-points, where this is well-known.
\end{proof}

Using this lemma, it is not too hard to construct for the simple factors of hyperkähler origin an element $\psi'_\alpha \in \Gamma\mathcal{H}_\alpha(\mathbb{Q})$ which satisfies the third condition of Theorem \ref{thm: properties of GammaH_i}:
\begin{proposition}
    \label{prop: GammaH_i for HK}
    Assume that we are in the setting of Theorem \ref{thm: properties of GammaH_i}, and suppose that $\mathfrak{h}_\alpha \cong \mathfrak{so}(\widetilde \sH(Y_k; \mathbb{Q})^{H_k})$ is a simple factor of hyperkähler origin of $\llv(Y; \mathbb{Q})$. Then:
    \begin{enumerate}
        \item There is a natural inclusion $\mathfrak{h}_\alpha \subseteq \Lie(\Gamma\mathcal{H}_\alpha)$.
        \item There is an element $\psi'_\alpha \in \Gamma\mathcal{H}_\alpha(\mathbb{Q})$ with the property that $\Ad(\psi'_\alpha)$ restricts to the automorphism $\Ad(\psi_\alpha)$ of $\mathfrak{h}_\alpha$.
    \end{enumerate}
\end{proposition}

\begin{proof}
    By taking Lie algebras in the first sequence of Lemma \ref{lem: HK Spin short exact sequences}, we see that 
    \[\Lie(\uGSpin(\widetilde \sH(Y_k;\mathbb{Q})^{H_k})) \cong \mathbb{Q} \times \mathfrak{so}(\widetilde \sH(Y_k;\mathbb{Q})^{H_k}).\]
    This gives a natural inclusion $\mathfrak{h}_\alpha \to \Lie(\Gamma\mathcal{H}_\alpha)$. Moreover, Lemma~\ref{lem: HK H^2 invariants dimension} gives $\dim \widetilde \sH(Y_k; \mathbb{Q})^{H_k} \geq 5$, so the Lie algebra $\mathfrak{so}(\widetilde \sH(Y_k;\mathbb{Q})^{H_k})$ is simple. Therefore, the derived Lie algebra of $\Lie(\uGSpin(\widetilde \sH(Y_k;\mathbb{Q})^{H_k}))$ is isomorphic to $\mathfrak{so}(\sH(Y_k;\mathbb{Q})^{H_k})$. It follows that for every element $\varphi \in \Gamma\mathcal{H}_\alpha(\mathbb{Q})$, the automorphism $\Ad(\varphi)$ of $\Lie(\Gamma\mathcal{H}_\alpha)$ restricts to an automorphism of $\mathfrak{h}_\alpha$.
    
    Hilbert 90 implies that $\sH^1(\mathbb{Q}, \mathbb{G}_{m, \mathbb{Q}}) = 0$, so the second sequence in Lemma \ref{lem: HK Spin short exact sequences} gives a surjective group homomorphism $\GSpin(\widetilde \sH(Y_k; \mathbb{Q})^{H_k}) \to \SO(\widetilde \sH(Y_k; \mathbb{Q})^{H_k})$. Let $\psi'_\alpha$ be a preimage of $\psi_\alpha$ under this surjection, then this map $\psi_\alpha'$ has the desired properties.
\end{proof}

\section{Lifting for factors of abelian origin}
\label{sec: abelian type groups}
We now assume that $\mathfrak{h}_\alpha$ is a simple factor of abelian origin. Then $\mathfrak{h}_\alpha$ is of type $A_1$ or $\mathfrak{h}_\alpha \cong \mathfrak{s}(\widetilde B, \widetilde \tau)$, where $W \subseteq \sH^1(Y_0; \mathbb{Q})$ is an isotypic component for the representation of $H$ on $\sH^1(Y_0; \mathbb{Q})$ and $\widetilde B = \End_B(\widetilde W)$, equipped with the involution $\widetilde \tau$ coming from the bilinear form from Lemma~\ref{lem: bilin form on V tilde is equivariant}.

In Corollary \ref{cor: psi_i for ab var}, we constructed for $\mathfrak{h}_\alpha \cong \mathfrak{s}(\widetilde B, \widetilde \tau)$ a $\psi_\alpha \in \Iso(\widetilde B, \widetilde \tau)$ such that the map 
\[\Ad(\psi_\alpha) \circ \Phi^{\llv}_\alpha \colon \mathfrak{g}_\alpha \to \mathfrak{h}_\alpha\]
respects the Hodge bigrading. For $\mathfrak{h}_\alpha$ of type $A_1$ with centroid $F$, we constructed a $\psi_\alpha \in \GL_2(F)$ such that $\Ad(\psi_\alpha) \circ \Phi^{\llv}_\alpha$ respects the Hodge bigrading in Corollary~\ref{cor: psi_i for type A1}. We will first discuss some properties of groups of isometries of algebras with involution. After that, we will partially prove Theorem~\ref{thm: properties of GammaH_i} for each type of involution by constructing $\psi'_\alpha$ and proving the third property. In Section~\ref{sec: A1 case} we treat the case where $\mathfrak{h}_\alpha$ is of type $A_1$.

\subsection{Algebraic groups associated to algebras with involutions}
\label{subsubsec: alg groups of involutions}
Let $(A, \sigma)$ be a simple algebra with involution over a field $k$ and let $\mathfrak{g} = \mathfrak{s}(A, \sigma)$. Except for some low-degree exceptions, $\mathfrak{g}$ is a semisimple Lie algebra over $k^\sigma$, see Theorem \ref{Thm: computation of s(A, sigma)}. Hence there exists a unique simply connected linear algebraic group $\Gamma$ over $k^\sigma$ with $\Lie(\Gamma) \cong \mathfrak{g}$. In this section, we will explicitly construct $\Gamma$.

We will start with some generalities, following \cite[Section 23]{BookOfInvolutions}. Let $(A, \sigma)$ be a simple algebra with involution over a field $F$, and assume that the center of $A$ is either $F$ (if $(A, \sigma)$ is of the first kind) or an étale quadratic extension $Z(A)/F$ with $F = Z(A)^\sigma$ (if $(A, \sigma)$ is of the second kind). 

Associated to such a pair $(A, \sigma)$, one can associate the group of isometries:
\begin{definition}
    The algebraic group \myindex{$\uIso(A, \sigma)$} of \textbf{isometries} of $(A, \sigma)$ is the algebraic group over $F = Z(A)^\sigma$, whose $R$-points for an $F$-algebra $R$ are given by:
     $$\uIso(A, \sigma)(R) = \{a \in A_R^\times : a\sigma(a) = 1\},$$
     where $A_R = A \otimes_{F} R$.
\end{definition}

Recall from Lemma \ref{lem: isometries algebraically closed field} that the group of isometries over an algebraically closed field is either the orthogonal group, the symplectic group or the general linear group, depending on whether $\sigma$ is orthogonal, symplectic or unitary.

The Lie algebra of the group of isometries is computed in \cite[Section~23.A]{BookOfInvolutions}:
\begin{lemma}
    \label{lem: Lie algebra of iso is skew}
    Let $(A, \sigma)$ be a simple algebra with involution, and let $F = Z(A)^\sigma$. Then there is an isomorphism $\Lie(\uIso(A, \sigma)) \cong \Skew(A,\sigma)$ of Lie algebras over~$F$. \qed
\end{lemma}
This implies the following (recall that $\mathfrak{s}(A, \sigma)$ is the derived Lie algebra of $\Skew(A, \sigma)$):
\begin{corollary}
    \label{cor: Lie iso is s(A, sigma)}
    Let $(A, \sigma)$ be a finite-dimensional simple algebra with involution over a field $k$ of characteristic zero, and assume that $k = k^\sigma$. Let $F = Z(A)^\sigma$, let the algebraic $F$-group $\uIso(A, \sigma)'$ be the commutator subgroup of $\uIso(A, \sigma)$, and let $\Gamma = \sR_{F/k} \uIso(A, \sigma)'$. Then there is an isomorphism 
    $$\Lie(\Gamma) \cong \mathfrak{s}(A, \sigma)$$
    of Lie algebras over $k$, compatible with the isomorphism from Lemma \ref{lem: Lie algebra of iso is skew}.
\end{corollary}
\begin{proof}
    By Lemma \ref{lem: Lie algebra of iso is skew}, there is an isomorphism $\Lie(\uIso(A, \sigma)) \cong \Skew(A, \sigma)$ of Lie algebras over $F$. Taking the commutator subgroup corresponds to taking the derived Lie algebra, so this induces an isomorphism $\Lie(\uIso(A, \sigma)') \cong \mathfrak{s}(A, \sigma)$ of Lie algebras over $F$. The field extension $F/k$ is finite because $A$ is finite-dimensional over $k$, and it is separable because $k$ has characteristic $0$, so the result now follows from Lemma \ref{lem: Lie algebra of Weil restriction}.
\end{proof}

\subsection{Orthogonal involutions}
\label{sec: orthogonal case}
In this section we partially prove Theorem \ref{thm: properties of GammaH_i} for $\mathfrak{h}_\alpha \cong \mathfrak{s}(\widetilde B, \widetilde \tau)$ with $\widetilde \tau$ an orthogonal involution. Before we do this, we summarize some properties of spin groups associated to algebras with orthogonal involutions.

\subsubsection{Spin groups}
Let $A$ be a finite-dimensional central simple algebra over a field $k$ of characteristic~$0$, equipped with an orthogonal involution $\sigma$. We write \myindex{$\uO(A, \sigma)$} for $\uIso(A, \sigma)$, this is an algebraic group over $k$. This algebraic group is disconnected, and we denote the connected component of the identity by \myindex{$\uSO(A, \sigma)$}. Since every connected variety with a rational point is geometrically connected \cite[\href{https://stacks.math.columbia.edu/tag/04KV}{Lemma~04KV}]{stacks}, we see that $\uSO(A, \sigma)(\overline{k}) \cong \uSO_{n}(\overline{k})$ for some integer $n$.

The group $\uSO(A, \sigma)$ is not simply connected. However, $\uSO(A, \sigma)$ admits a double covering by a simply connected algebraic group over $k$, denoted by $\uSpin(A, \sigma)$. Similar to the case of the classical spin group, there is also a version of the general spin group, which we denote by $\uGSpin(A, \sigma)$. We now summarize the most important properties of $\uSpin(A, \sigma)$ and $\uGSpin(A, \sigma)$.

\begin{definition}
    Let $A$ be a central simple algebra over a field $k$ of characteristic~$0$, equipped with an orthogonal involution $\sigma$. We define \myindex{$\uSpin(A, \sigma)$} to be the simply connected cover of $\uSO(A, \sigma)$, and call $\uSpin(A, \sigma)$ the \textbf{spin group} of $(A, \sigma)$. It is an algebraic group over $k = Z(A)$.
\end{definition}
\begin{remark}
    For an explicit construction of $\uSpin(A, \sigma)$ in terms of Clifford algebras, see \cite[Section~23.B]{BookOfInvolutions}.
\end{remark}

The essential properties of $\uSpin(A, \sigma)$ that we will need are summarized in the following theorem. The following properties (and a lot more about spin groups) are discussed in the Book of Involutions \cite[Section~23.B]{BookOfInvolutions}.
\begin{theorem}
    \label{thm: Spin(A,sigma) properties}
    Let $A$ be a central simple algebra over a field $k$ of characteristic~$0$, equipped with an orthogonal involution $\sigma$. Then there exists a short exact sequence of algebraic groups over $k$
    \begin{equation*}
        1 \to \mu_{2, k} \to \uSpin(A, \sigma) \to \uSO(A, \sigma) \to 1.
    \end{equation*}
    Moreover, this short exact sequence has the following properties:
    \begin{enumerate}
        \item The formation of this short exact sequence commutes with base change.
        \item If $A \cong \End_k(V)$ for a finite-dimensional $k$-vector space $V$, and $V$ is equipped with a nondegenerate symmetric bilinear form $b$ such that $\sigma$ is the adjoint involution of $b$, then the short exact sequence is isomorphic to the short exact sequence
        \begin{equation*}
            1 \to \mu_{2, k} \to \uSpin(V, b) \to \uSO(V, b) \to 1.
        \end{equation*}
    \end{enumerate}\qed
\end{theorem}

Now that we have the group $\uSpin(A, \sigma)$, we can define $\uGSpin(A, \sigma)$ as follows:
\begin{definition}
    Let $A$ be a central simple algebra over a field $k$ of characteristic~$0$, equipped with an orthogonal involution $\sigma$. Let \myindex{$\uGSpin(A, \sigma)$} be the algebraic group over $k$ defined by 
    \[\uGSpin(A, \sigma) := (\mathbb{G}_{m, k} \times \uSpin(A, \sigma))/\mu_{2, k},\]
    where the inclusion $\mu_2 \to \mathbb{G}_m \times \uSpin(A, \sigma)$ is the product of the standard inclusion $\mu_2 \to \mathbb{G}_m$ and the inclusion $\mu_2 \to \uSpin(A, \sigma)$ from the above theorem. We call $\uGSpin(A, \sigma)$ the \textbf{general spin group} of $A$.
\end{definition}

\begin{remark}
    The groups $\uSpin(A, \sigma)$ and $\uGSpin(A, \sigma)$ are algebraic groups over the center of $A$. Often, we will consider $A$ as an algebra over a field $k$, with possibly a bigger center $F = Z(A)$. In those cases, the group $\uSpin(A, \sigma)$ is defined over $F$, and if we want to obtain a group over $k$, we take the Weil restriction $\sR_{F/k}\uSpin(A, \sigma)$.
\end{remark}

The group $\uGSpin(A, \sigma)$ comes with a map to $\uSO(A, \sigma)$. Indeed, consider the trivial map $\mathbb{G}_m \to \uSO(A, \sigma)$ which sends everything to the identity, and the map $\uSpin(A, \sigma) \to \uSO(A, \sigma)$. The product of these two maps contains $\mu_2$ in its kernel, and hence induces a map $\uGSpin(A, \sigma) \to \uSO(A, \sigma)$. The group $\uGSpin(A, \sigma)$ interacts just as expected with the groups $\uSO(A, \sigma)$ and $\uSpin(A, \sigma)$:

\begin{lemma}
    \label{lem: GSpin SO short exact sequences}
    Let $(A, \sigma)$ be a finite-dimensional simple algebra with orthogonal involution over a field~$k$, and let $F = Z(A)$. Then the group $\sR_{F/k}\uGSpin(A, \sigma)$ fits in the short exact sequences:
    \[ 1 \to \sR_{F/k}\mu_{2, F} \to \sR_{F/k}\mathbb{G}_{m, F} \times \sR_{F/k}\uSpin(A, \sigma) \to \sR_{F/k}\uGSpin(A, \sigma) \to 1 \]
    and
    \[1 \to \sR_{F/k}\mathbb{G}_{m, F} \to \sR_{F/k}\uGSpin(A, \sigma) \to \sR_{F/k}\uSO(A, \sigma) \to 1.\]
\end{lemma}
\begin{proof}
    Since Weil restriction commutes with products and preserves exactness (Lemmas \ref{lem: weil restriction product} and \ref{lem: weil restriction exact}), it suffices to show that the sequences without Weil restriction are exact. For the first one, this is just a rewriting of the definition of $\uGSpin(A, \sigma)$. 

    The exactness of the second sequence then follows from the definition of the map $\uGSpin(A, \sigma) \to \uSO(A, \sigma)$.
\end{proof}

\subsubsection{Construction of $\Gamma\mathcal{H}_\alpha$ and $\psi_\alpha'$ for orthogonal involutions}
Let $(B, \tau)$ be a simple finite-dimensional $\mathbb{Q}$-algebra with involution, with center $F = Z(B)$. Let $W$ be a $B$-module and let $\widetilde B = \End_B(\widetilde W)$. Assume that the induced involution $\widetilde \tau$ on $\widetilde B$ is orthogonal. Then the $\mathbb{Q}$-points of $\sR_{F/\mathbb{Q}}\uO(\widetilde B, \widetilde \tau)$ are equal to $\sO_B(\widetilde W)$ by Lemma~\ref{lem: SO_A is Iso(A, sigma)}. We will denote the $\mathbb{Q}$-points of $\sR_{F/\mathbb{Q}}\uSO(\widetilde B, \widetilde \tau)$ by \myindex{$\SO_B(\widetilde W)$}.

Assume that we have an isomorphism $\Phi^{\llv}_\alpha \colon \mathfrak{g}_\alpha \to \mathfrak{h}_\alpha$ between simple factors of the LLV algebras of two varieties $X$ and $Y$ admitting a holomorphic symplectic form and satisfying condition \ref{condition *}. Assume that $\mathfrak{g}_\alpha \cong \mathfrak{s}(\widetilde A, \widetilde \sigma)$ and $\mathfrak{h}_\alpha \cong \mathfrak{s}(\widetilde B, \widetilde \tau)$, with $\widetilde \sigma$ and $\widetilde \tau$ orthogonal involutions. By Corollary \ref{cor: psi_i for ab var}, we obtain a $\psi_\alpha \in \Iso(\widetilde B, \widetilde \tau)$ with the property that $\psi_\alpha \circ \Phi_\alpha^{\llv} \colon \mathfrak{g}_\alpha \to \mathfrak{h}_\alpha$ preserves the Hodge bigrading.

Analogous to the proof of Proposition~\ref{prop: GammaH_i for HK}, we will want to lift this map $\psi_\alpha$ to an element of the group $\GSpin(\widetilde B, \widetilde \tau)$. However, the element $\psi_\alpha$ from Corollary~\ref{cor: psi_i for ab var} is only in $\Iso(\widetilde B, \widetilde \tau) = \sO_B(\widetilde W)$, and to lift it, it needs to be in $\SO_B(\widetilde W)$. Before we prove that $\psi_\alpha$ is actually in $\SO_B(\widetilde W)$, we need the following technical result.

Let $(A, \sigma)$ be an algebra with orthogonal involution over $\mathbb{Q}$, let $F = Z(A)$ and let $k = [F : \mathbb{Q}]$. Assume that the degree of $A$ is even, so we can write $\deg(A) = 2n$, and $\mathfrak{s}(A, \sigma)$ is of type $D_n$. Then there is an isomorphism
\[ \mathfrak{s}(A, \sigma) \otimes_\mathbb{Q} \mathbb{C} \cong \prod_{i=1}^k \mathfrak{so}_{2n}(\mathbb{C}),\]
so the Dynkin diagram of $\mathfrak{s}(A, \sigma)$ as a Lie algebra over $\mathbb{Q}$ is isomorphic to the disjoint union $\coprod_{i=1}^k D_n$. The group $\Iso(A, \sigma)$ acts on $\mathfrak{s}(A, \sigma)$ by conjugation, so there is an induced action of $\Iso(A, \sigma)$ on $\coprod_{i=1}^k D_n$.

\begin{lemma}
    \label{lem: kernel Iso to Aut(Dn)}
    Let $(A, \sigma)$ be an algebra with orthogonal involution over $\mathbb{Q}$. Let $F = Z(A)$, and let $k = [F : \mathbb{Q}]$. Assume that there is an $n$ such that $\mathfrak{s}(A, \sigma)$ is of type $D_n$. Then the kernel of the map
    \[(\sR_{F/ \mathbb{Q}}\uIso(A, \sigma))(\mathbb{C}) \to \Aut\left(\coprod_{i=1}^k D_n\right)\]
    is equal to the connected component of the identity of $(\sR_{F/ \mathbb{Q}}\uIso(A, \sigma))(\mathbb{C})$.
\end{lemma}
\begin{proof}
    Since $\sigma$ is orthogonal, Lemmas~\ref{lem: Weil restriction F-bar} and~\ref{lem: isometries algebraically closed field} give an isomorphism
    \[ (\sR_{F/\mathbb{Q}} \uIso(A,\sigma))(\mathbb{C}) \cong \prod_{i=1}^k \sO_{2n}(\mathbb{C}).\]
    The conjugation action of $\prod_i \sO_{2n}(\mathbb{C})$ on $\mathfrak{s}(A, \sigma) \otimes_\mathbb{Q} \mathbb{C} \cong \prod_i \mathfrak{so}_{2n}(\mathbb{C})$ does not permute the different factors $\mathfrak{so}_{2n}(\mathbb{C})$, so the action of $\prod_i \sO_{2n}(\mathbb{C})$ on $\coprod_i D_n$ does not permute the components of the Dynkin diagram. The result now follows from the fact that the kernel of the map $\sO_{2n}(\mathbb{C}) \to \Aut(D_n)$ is equal to $\SO_{2n}(\mathbb{C})$.
\end{proof}

\begin{lemma}
    \label{lem: psi_i orthogonal in SO_B}
    Let $(\widetilde B, \widetilde \tau)$ be as in Corollary~\ref{cor: psi_i for ab var}, and assume that $\widetilde \tau$ is an orthogonal involution. Let $\psi_\alpha \in \Iso(\widetilde B, \widetilde \tau)$ be the isometry constructed in Corollary~\ref{cor: psi_i for ab var}. Then $\psi_\alpha \in \SO_B(\widetilde W) \subseteq \Iso(\widetilde B, \widetilde \tau)$.
\end{lemma}
\begin{proof}
    Let $F = Z(\widetilde B)$, let $k = [F : \mathbb{Q}]$ be the degree of $F$ over $\mathbb{Q}$, and let $n$ be the integer such that $\mathfrak{s}(\widetilde B, \widetilde \tau)$ is of type $D_n$. We saw above that the Dynkin diagram of $\mathfrak{s}(\widetilde B, \widetilde \tau)$ as a Lie algebra over $\mathbb{Q}$ is equal to the disjoint union $\coprod_{i=1}^k D_n$. Lemma~\ref{lem: kernel Iso to Aut(Dn)} implies that the kernel of the composition $\sO_B(\widetilde W) \to \Aut(\mathfrak{s}(\widetilde B, \widetilde \tau)) \to \Aut(\coprod_{i=1}^k D_n)$ is equal to $\SO_B(\widetilde W)$. Therefore, we are done once we show that $\psi_\alpha$ acts trivially on the Dynkin diagram $\coprod_{i=1}^k D_n$. 

    We saw in the proof of Lemma~\ref{lem: kernel Iso to Aut(Dn)} that the action of $\Iso(\widetilde B, \widetilde \tau)$ on $\coprod_i D_n$ preserves the connected components. Since $\psi_\alpha \in \Iso(\widetilde B, \widetilde \tau)$, it remains to show that $\psi_\alpha$ acts trivially on each component $D_n$ of $\coprod_i D_n$.

    Recall that $\psi_\alpha = s_\lambda \circ \widetilde g \circ f^{-1}$ as constructed in Lemmas \ref{lem: abelian factor similitude psi} and \ref{lem: abelian factor isometry psi}. We know that $\Ad(f) = \Phi^{\llv}_\alpha \colon \mathfrak{g}_\alpha \to \mathfrak{h}_\alpha$. There is a decomposition
    \[\SH_\alpha(Y) \otimes_{\mathbb{Q}} \mathbb{C} \cong \bigotimes_{i=1}^k \SH_{\alpha, i}(Y),\] 
    where each $\SH_{\alpha, i}(Y)$ is an irreducible representation of the \mbox{$i$-th} factor $\mathfrak{so}_{2n}(\mathbb{C})$ in the decomposition $\mathfrak{s}(\widetilde B, \widetilde \tau) \otimes_\mathbb{Q} \mathbb{C} \cong \prod_{i=1}^k \mathfrak{so}_{2n}(\mathbb{C})$. Similarly, we have a decomposition $\SH_\alpha(X) \otimes_\mathbb{Q} \mathbb{C} \cong \bigotimes_{i=1}^k \SH_{\alpha, i}(X)$. By Theorem \ref{thm: split Lie alg reps}, the isomorphism $\Phi^{\SH}_\alpha \colon \SH_\alpha(X) \to \SH_\alpha(Y)$ is equivariant with respect to $\Phi^{\llv}_\alpha$. From Proposition~\ref{prop: weight of Dn SH} we know that the $\SH_{\alpha, i}(X)$ and $\SH_{\alpha, i}(Y)$ are representations with highest weight $k\omega'$ for some integer $k > 0$, where $\omega'$ corresponds to the even spinor representation.

    Next, we claim that $\Ad(\widetilde g) \colon \mathfrak{g}_\alpha \to \mathfrak{h}_\alpha$ also preserves the highest weight of the even spinor representation in each component $D_n$ of $\coprod_i D_n$. To see this, recall from Lemma \ref{lem: abelian factor similitude psi} that $\widetilde g = g \oplus (g^*)^{-1} \colon \widetilde V \to \widetilde W$ for some isomorphism $g \colon V \to W$, equivariant with respect to the isomorphism $A \to B$. Then the degree-preserving linear isomorphism $\wedge g \colon \extp^\bullet V \to \extp^\bullet W$ is equivariant with respect to $\Ad(g)$. Let $\SH(\extp^\bullet V, \mathfrak{g}_\alpha) \subseteq \extp^\bullet V$ be the smallest $\mathfrak{g}_\alpha$-subrepresentation containing $\extp^0 V$, and define $\SH(W, \mathfrak{h}_\alpha) \subseteq \extp^\bullet W$ similarly. Then $\wedge g$ restricts to an isomorphism $\SH(V, \mathfrak{g}_\alpha) \to \SH(W, \mathfrak{h}_\alpha)$.

    There is a decomposition $B \otimes_\mathbb{Q} \mathbb{C} \cong \prod_{i=1}^k B_i$, where each $B_i$ is a central simple $\mathbb{C}$-algebra. Associated to this, there is a decomposition $W \otimes_\mathbb{Q} \mathbb{C} \cong \bigoplus_i W_i$, where each $W_i$ is a $B_i$-module. Similarly, we have a decomposition $A \otimes_\mathbb{Q} \mathbb{C} \cong \prod_{i=1}^k A_i$, where the $A_i$ are central simple $\mathbb{C}$-algebras. We also obtain a decomposition $V \otimes_\mathbb{Q} \mathbb{C} \cong \bigoplus_{i=1}^k V_i$, where each $V_i$ is an $A_i$-module. This induces decompositions $\mathfrak{g}_\alpha \otimes_\mathbb{Q} \mathbb{C} \cong \prod_{i=1}^k \mathfrak{g}_{\alpha, i}$ and $\mathfrak{h}_\alpha \otimes_\mathbb{Q} \mathbb{C} \cong \prod_{i=1}^k \mathfrak{h}_{\alpha, i}$ with the $\mathfrak{g}_{\alpha, i}$ and $\mathfrak{h}_{\alpha, i}$ simple Lie algebras over $\mathbb{C}$ (they are isomorphic to $\mathfrak{so}_{2n}(\mathbb{C})$). From this, we obtain factorizations $\SH(V, \mathfrak{g}_\alpha) \cong \bigotimes_i \SH(V_i, \mathfrak{g}_{\alpha, i})$ and $\SH(W, \mathfrak{h}_\alpha) \cong \bigotimes_i \SH(W_i, \mathfrak{h}_{\alpha, i})$. Proposition \ref{prop: weight of Dn SH} implies that the $\SH(V_i, \mathfrak{g}_{\alpha, i})$ and $\SH(W_i, \mathfrak{h}_{\alpha, i})$ correspond to the representations of highest weight $k\omega'$ for some integer $k' > 0$, where the weight $\omega'$ is the highest weight associated to the even spinor representation.

    Therefore, it follows that the composition $\Ad(\widetilde g) \circ \Ad(f^{-1}) \in \Aut(\mathfrak{h}_\alpha)$ fixes the weight~$\omega'$ in each component $D_n$ of $\coprod_i D_n$, so $\Ad(\widetilde g) \circ \Ad(f^{-1})$ induces the identity on the Dynkin diagram $\coprod_i D_n$.
    
    To see that $s_\lambda$ acts trivially on $\coprod_i D_n$, we again pass to $\mathbb{C}$. The adjoint map $\Ad \colon \Iso(\widetilde B, \widetilde \tau) \to \Aut(\mathfrak{s}(\widetilde B, \widetilde \tau))$ extends to an adjoint map 
    \[\Ad \colon \Sim(\widetilde B, \widetilde \tau) \to \Aut(\mathfrak{s}(\widetilde B, \widetilde \tau)),\]
    also given by conjugation. By choosing a path from $\lambda \otimes 1$ to $1$ in $Z(B) \otimes_\mathbb{Q} \mathbb{C}^\times$, we see that there is a continuous path from $s_\lambda$ to the identity in $\Sim(\widetilde B \otimes_\mathbb{Q} \mathbb{C}, \widetilde \tau)$, so $s_\lambda$ is in the identity component of $\Sim(\widetilde B, \widetilde \tau)$. Hence $s_\lambda$ must act trivially on $\coprod_i D_n$, and therefore $\psi_\alpha = s_\lambda \circ \widetilde g \circ f^{-1}$ acts trivially on $\coprod_i D_n$.
\end{proof}

For $\mathfrak{h}_\alpha$ coming from an algebra with an orthogonal involution, we can now construct $\Gamma\mathcal{H}_\alpha$ and $\psi_\alpha'$ in Theorem \ref{thm: properties of GammaH_i} and show that $\Ad(\psi_\alpha')$ restricts to the automorphism $\Ad(\psi_\alpha)$ of $\mathfrak{h}_\alpha$.
\begin{proposition}
    \label{prop: GammaHi orthogonal}
    Assume that we are in the setting of Theorem \ref{thm: properties of GammaH_i}, and let $\mathfrak{h}_\alpha \cong \mathfrak{s}(\widetilde B, \widetilde \tau)$ be a simple factor of abelian origin of $\llv(Y; \mathbb{Q})$ with $\widetilde \tau$ an orthogonal involution. Assume that $\mathfrak{h}_\alpha$ is of type $D_n$ with $n = 3$ or $n \geq 5$. Let $F = Z(B)$ and let $\Gamma\mathcal{H}_\alpha = \sR_{F/\mathbb{Q}}\uGSpin(\widetilde B, \widetilde \tau)$. Then there is a natural inclusion $\mathfrak{h}_\alpha \subseteq \Lie(\Gamma\mathcal{H}_\alpha)$ and there exists a $\psi'_\alpha \in \Gamma\mathcal{H}_\alpha(\mathbb{Q}) = \GSpin(\widetilde B, \widetilde \tau)$ with the property that $\Ad(\psi'_\alpha) = \Ad(\psi_\alpha)$ as automorphisms of $\mathfrak{h}_\alpha$.
\end{proposition}
\begin{proof}
    The proof is similar to the proof of Proposition \ref{prop: GammaH_i for HK}. Part (2) of Theorem~\ref{thm: Spin(A,sigma) properties} implies that $\Lie(\uSpin(\widetilde B, \widetilde \tau)) \cong \mathfrak{s}(\widetilde B, \widetilde \tau)$. The first short exact sequence from Lemma \ref{lem: GSpin SO short exact sequences} then implies that 
    \[\Lie(\sR_{F/\mathbb{Q}}\uGSpin(\widetilde B, \widetilde \tau)) \cong Z(B) \times \mathfrak{s}(\widetilde B, \widetilde \tau),\]
    with $Z(B)$ an abelian Lie algebra. This gives the inclusion $\mathfrak{h}_\alpha \subseteq \Lie(\Gamma\mathcal{H}_\alpha)$. Since the variety $Y$ satisfies condition \ref{condition *}, Theorem \ref{thm: summary LLV computation} implies that $\mathfrak{s}(\widetilde B, \widetilde \tau)$ is simple, so it follows that $\mathfrak{s}(\widetilde B, \widetilde \tau)$ is the derived Lie algebra of $\Lie(\sR_{F/\mathbb{Q}}\uGSpin(\widetilde B, \widetilde \tau))$. This implies that for every $\varphi \in \Gamma\mathcal{H}_\alpha(\mathbb{Q})$ the automorphism $\Ad(\varphi)$ of $\Lie(\Gamma\mathcal{H}_\alpha)$ restricts to an automorphism of $\mathfrak{h}_\alpha$.

    For the construction of $\psi'_\alpha$, note that $\psi_\alpha \in \SO_B(\widetilde W)$ by Lemma \ref{lem: psi_i orthogonal in SO_B}. By taking $\mathbb{Q}$-points in the second short exact sequence from Lemma~\ref{lem: GSpin SO short exact sequences} and applying Lemma~\ref{lem: Hilbert 90}, we see that there is a surjection 
    \[\GSpin(\widetilde B, \widetilde \tau) \to \SO_B(\widetilde W).\]
    Let $\psi'_\alpha$ be a preimage of $\psi_\alpha$ under this surjection, then $\psi'_\alpha$ satisfies $\Ad(\psi'_\alpha) = \Ad(\psi_\alpha)$ as automorphisms of $\mathfrak{s}(\widetilde B, \widetilde \tau)$.
\end{proof}

\subsection{Symplectic involutions}
Now assume that $\mathfrak{h}_\alpha \cong \mathfrak{s}(\widetilde B, \widetilde \tau)$ with $\widetilde \tau$ a symplectic involution. If we take the $\overline{\mathbb{Q}}$-points of $\uIso(\widetilde B, \widetilde \tau)$, we get the symplectic group $\Sp(n, \overline{\mathbb{Q}})$ for some integer $n$. Lemma \ref{lem: simply connected iff} then implies that $\uIso(\widetilde B, \widetilde \tau)$ is simply connected, and it has Lie algebra equal to $\mathfrak{g} = \mathfrak{s}(\widetilde B, \widetilde \tau)$. Then Lemma \ref{lem: Lie algebra of Weil restriction} gives that $\sR_{F/\mathbb{Q}} \uIso(\widetilde B, \widetilde \tau)$ is the simply connected algebraic group over $\mathbb{Q}$ with Lie algebra $\mathfrak{h}_\alpha$.

The proof of Theorem \ref{thm: properties of GammaH_i} is straightforward in this case, since we take $\Gamma\mathcal{H}_\alpha = \mathcal{H}_\alpha$ and $\psi_\alpha' = \psi_\alpha$.
\begin{proposition}
    \label{prop: GammaHi symplectic}
    Assume that we are in the setting of Theorem \ref{thm: properties of GammaH_i}, and let $\mathfrak{h}_\alpha \cong \mathfrak{s}(\widetilde B, \widetilde \tau)$ be a simple factor of abelian origin of $\llv(Y; \mathbb{Q})$ with $\widetilde \tau$ a symplectic involution, and assume that $\mathfrak{h}_\alpha$ is of type $C_l$ with $l \geq 2$. Let $\Gamma\mathcal{H}_\alpha = \sR_{K/\mathbb{Q}}\uIso(\widetilde B, \widetilde \tau)$, where $K = Z(B)$, and let $\psi'_\alpha = \psi_\alpha$. Then the third property of Theorem \ref{thm: properties of GammaH_i} holds.
\end{proposition}
\begin{proof}
    We have $\Lie(\Gamma\mathcal{H}_\alpha) \cong \mathfrak{s}(\widetilde B, \widetilde \tau)$ by Corollary \ref{cor: Lie iso is s(A, sigma)}, so indeed $\mathfrak{h}_\alpha \cong \Lie(\mathcal{H}_\alpha) \cong \Lie(\Gamma\mathcal{H}_\alpha)$. Since $\psi_\alpha = \psi'_\alpha$, we indeed have that $\Ad(\psi_\alpha') = \Ad(\psi_\alpha)$ as automorphisms of $\mathfrak{h}_\alpha$.
\end{proof}

\subsection{Unitary involutions}
We will now treat the unitary case, so assume that $\mathfrak{h}_\alpha \cong \mathfrak{s}(\widetilde B, \widetilde \tau)$ with $\widetilde \tau$ a unitary involution, and assume that $\mathfrak{h}_\alpha$ is not of type $A_1$. Let $K = Z(\widetilde B)$, so $\tau$ acts nontrivially on $K$, and we have $[K : K^\tau] = 2$. Define $\Gamma\mathcal{H}_\alpha = \sR_{K^{\widetilde \tau}/\mathbb{Q}}\uIso(\widetilde B, \widetilde \tau)$ as an algebraic group over $\mathbb{Q}$. We then have $\psi_\alpha \in \Gamma\mathcal{H}_\alpha(\mathbb{Q})$ for the element $\psi_\alpha$ constructed in Corollary \ref{cor: psi_i for ab var}. The difficulty in proving Theorem \ref{thm: properties of GammaH_i} in this case is the construction of a representation of $\Gamma\mathcal{H}_\alpha$ on $\sH^\bullet(Y; \mathbb{Q})$.

If $\widetilde \tau$ is unitary, recall from Lemma~\ref{lem: isometries algebraically closed field} that the algebraic group $\uIso(\widetilde B, \widetilde \tau)$ is a form of the general linear group $\uGL_n$ for some $n \geq 1$. In particular, $\uIso(\widetilde B, \widetilde \tau)$ is reductive but not semisimple. Its Lie algebra is $\Skew(\widetilde B, \widetilde \tau)$, and we have a strict inclusion $\mathfrak{s}(\widetilde B, \widetilde \tau) \subset \Skew(\widetilde B, \widetilde \tau)$. Only the Lie subalgebra $\mathfrak{s}(\widetilde B, \widetilde \tau)$ occurs in the LLV algebra, so by integrating the action of the LLV algebra on the cohomology we only get a representation of $\sR_{K^{\widetilde \tau}/\mathbb{Q}}\uIso(\widetilde B, \widetilde \tau)'$ on $\sH^\bullet(Y; \mathbb{Q})$ (note that $\Lie(\sR_{K^{\widetilde \tau}/\mathbb{Q}}\uIso(\widetilde B, \widetilde \tau)') \cong \mathfrak{s}(\widetilde B, \widetilde \tau)$ by Corollary~\ref{cor: Lie iso is s(A, sigma)}). Since $\psi_\alpha$ need not be in $\Iso(\widetilde B, \widetilde \tau)'$, we wish to extend the representation of $\sR_{K^{\widetilde \tau}/\mathbb{Q}}\uIso(\widetilde B, \widetilde \tau)'$ on $\sH^\bullet(Y; \mathbb{Q})$ to a representation of $\sR_{K^{\widetilde \tau}/\mathbb{Q}}\uIso(\widetilde B, \widetilde \tau)$. 

First, let us set the notation straight for this section:
\begin{setup}
    Suppose we have the following data:
    \begin{itemize}
        \item A finite group $H$.
        \item A $\rho \in \Irr_\mathbb{Q}(H)$ which splits into non-self-dual representations over $\overline{\mathbb{Q}}$.
        \item An $H$-representation $W$ over $\mathbb{Q}$, and assume that $W$ is $\rho$-isotypic, i.e. there is an $m \in \mathbb{N}$ with $W \cong \rho^m$ as $H$-representations.
    \end{itemize}
    We obtain the following objects:
    \begin{itemize}
        \item The simple factor $B$ of $\mathbb{Q}[H]$ associated to $\rho$.
        \item $\widetilde W = W \oplus W^*$ as an $H$-representation. This carries the nondegenerate bilinear form $b$ described above Lemma \ref{lem: bilin form on V tilde is equivariant}.
        \item $\widetilde B = \End_H(\widetilde W)$, equipped with the adjoint involution $\widetilde \tau$ corresponding to the bilinear form $b$.
        \item $K = Z(\widetilde B)$, and the restriction $\tau = \widetilde \tau|_K$.
        \item There is an integer $n$ such that $\widetilde B \otimes_{K^\tau} \overline{K} \cong M_n(\overline{K}) \times M_n(\overline{K})$, and we fix this integer $n$ throughout this section.
    \end{itemize}
    Note that $\widetilde \tau$ is an involution of the second kind. Its restriction $\tau$ to $K$ is nontrivial, so $K/K^\tau$ is a degree $2$ field extension. 
\end{setup}

The inclusion $\widetilde B \subseteq \End_\mathbb{Q}(\widetilde W)$ induces and inclusion $\Skew(\widetilde B, \widetilde \tau) \subseteq \mathfrak{so}(\widetilde W)$, and $\mathfrak{so}(\widetilde W)$ acts on $\bigwedge^\bullet W$ via the spinor representation as described in \cite[Chapter 3]{LooijengaLunts}. From this we obtain a map $\mathfrak{s}(\widetilde B, \widetilde \tau) \to \mathfrak{gl}(\extp^\bullet W)$. Since $\uIso(\widetilde B, \widetilde \tau)'$ is simply connected, the map $\mathfrak{s}(\widetilde B, \widetilde \tau) \to \mathfrak{gl}(\bigwedge^\bullet W)$ integrates to a representation of $\sR_{K^\tau/\mathbb{Q}}\uIso(\widetilde B, \widetilde \tau)'$ on $\bigwedge^\bullet W$. We wish to extend this to a representation of $\sR_{K^\tau/\mathbb{Q}}\uIso(\widetilde B, \widetilde \tau)$. In the proof of Theorem~\ref{thm: properties of GammaH_i} in Section~\ref{subsec: proof of thm properties of GammaH_i} we will see why this suffices for extending the representation of $\sR_{K^\tau/\mathbb{Q}}\uIso(\widetilde B, \widetilde \tau)'$ on $\sH^\bullet(Y; \mathbb{Q})$ to a representation of $\sR_{K^\tau/\mathbb{Q}}\uIso(\widetilde B, \widetilde \tau)$.

\subsubsection{The norm homomorphism}
We will now discuss a generalization of the norm homomorphism to algebraic groups, and study its kernel. For a detailed reference, see \cite{NormTorus}.

For a finite field extension $K/F$, the norm homomorphism $N_{K/F} \colon K^\times \to F^\times$ induces a homomorphism
\[N_{K/F} \colon \sR_{K/F}\mathbb{G}_{m, K} \to \mathbb{G}_{m,F} \]
of algebraic groups over $F$; the induced map on $F$-points is precisely the norm homomorphism $N_{K/F} \colon K^\times \to F^\times$.
\begin{definition}
    Let $K/F$ be a finite separable field extension. We write $\mathbb{S}_{K/F}$ for the algebraic group over $F$ defined by
    \[ \mathbb{S}_{K/F} := \ker(N_{K/F} \colon \sR_{K/F} \mathbb{G}_{m, K} \to \mathbb{G}_{m, F}). \]
\end{definition}
For $t \in \mathbb{N}$, we will also look at the subgroup $\mathbb{S}_{K/F}[t] \subseteq \mathbb{S}_{K/F}$ of $t$-torsion elements. The norm homomorphism restricts to a homomorphism $\sR_{K/F} \mu_{t, K} \to \mu_{t, F}$, and this gives an isomorphism 
$$\mathbb{S}_{K/F}[t] \cong \ker(N_{K/F} \colon \sR_{K/F} \mu_{t, K} \to \mu_{t, F}).$$

During the rest of this section, the field $K$ will be $Z(\widetilde B)$, and we will always take $F = K^\tau$. For this setting, we will just write \myindex{$\mathbb{S}$} instead of $\mathbb{S}_{K/K^\tau}$, this is an algebraic group over $K^\tau$. We will denote the Weil restriction $\sR_{K^{\tau}/\mathbb{Q}} \mathbb{S}$ by $\sRS$, this is an algebraic group over $\mathbb{Q}$. For the extension $K/K^\tau$, the norm homomorphism is given by sending $x \in K^\times$ to $x\tau(x)$, so 
\[\mathbb{S}(K^{\tau}) = \{x \in K^\times : x\tau(x)=1\}.\]

This description of the norm homomorphism implies that we have a homomorphism $\mathbb{S} \to \uIso(\widetilde B, \widetilde \tau)$. By computing the base change to $\overline{K}$, we see that this identifies $\mathbb{S}$ with the center of $\uIso(\widetilde B, \widetilde \tau)$.

Recall that $n$ is the integer such that $\widetilde B \otimes_{K^\tau} \overline{K} \cong M_n(\overline{K}) \times M_n(\overline{K})^{\op}$, and that $\uIso(\widetilde B, \widetilde \tau)'$ denotes the commutator subgroup of $\uIso(\widetilde B, \widetilde \tau)$.

\begin{lemma}
    \label{lem: S[n] S Iso short exact sequence}
    There is a short exact sequence of algebraic groups over $K^\tau$
    \begin{equation}
        \label{eq: S[n] to S times Iso' ses}
        1 \to \mathbb{S}[n] \to \mathbb{S} \times \uIso(\widetilde B, \widetilde \tau)' \to \uIso(\widetilde B, \widetilde \tau) \to 1,
    \end{equation}
    where the first map is induced by sending $x \in \mathbb{S}[n]$ to $(x, x^{-1})$.
\end{lemma}

\begin{proof}
    We saw above that there is an inclusion $\mathbb{S} \to \uIso(\widetilde  B, \widetilde \tau)$, and there is also an inclusion $\uIso(\widetilde B, \widetilde \tau)' \subseteq \uIso(\widetilde B, \widetilde \tau)$. We then obtain a map $\mathbb{S} \times \uIso(\widetilde B, \widetilde \tau)' \to \uIso(\widetilde B, \widetilde \tau)$ by sending $(\lambda, x)$ to $\lambda x$. Note that the images of $\mathbb{S}$ and $\uIso(\widetilde B, \widetilde \tau)'$ commute, since $\mathbb{S}$ maps onto the center of $\uIso(\widetilde B, \widetilde \tau)$. Under the isomorphism $\widetilde B \otimes_{K^\tau} \overline{K} \cong M_n(\overline{K}) \times M_n(\overline{K})^{\op}$, the involution $\widetilde \tau$ corresponds with the swapping involution. From this, we see that $\Iso(\widetilde B, \widetilde \tau)_{\overline{K}} \cong \GL_n(\overline{K})$, and $\Iso(\widetilde B, \widetilde \tau)'_{\overline{K}} \cong \SL_n(\overline{K})$. This implies that $\mathbb{S}[n]$ maps to $\uIso(\widetilde B, \widetilde \tau)'$ under the inclusion $\mathbb{S} \to \uIso(\widetilde B, \widetilde \tau)$. Therefore, the maps in sequence \eqref{eq: S[n] to S times Iso' ses} are indeed defined.
    
     By Lemma \ref{lem: exactness algebraic closure}, it suffices to check the exactness of sequence \eqref{eq: S[n] to S times Iso' ses} on the $\overline{K^\tau}$-points, where it becomes the short exact sequence
     \[1 \to \mu_{n}(\overline{K^\tau}) \to \overline{K^\tau}^\times \times \SL_n(\overline{K^\tau}) \to \GL_n(\overline{K^\tau}) \to 1.\]
\end{proof}

Using the short exact sequence of the above lemma, we can construct $\Gamma\mathcal{H}_\alpha$ and $\psi'_\alpha$ in Theorem \ref{thm: properties of GammaH_i} in the unitary case:
\begin{proposition}
    \label{prop: GammaHi unitary}
    Assume that we are in the setting of Theorem \ref{thm: properties of GammaH_i}, and let $\mathfrak{h}_\alpha \cong \mathfrak{s}(\widetilde B, \widetilde \tau)$ be a simple factor of abelian origin of $\llv(Y; \mathbb{Q})$ with $\widetilde \tau$ a unitary involution, and assume that $\mathfrak{h}_\alpha$ is not of type $A_1$. Let $\Gamma\mathcal{H}_\alpha = \sR_{K^\tau/\mathbb{Q}}\uIso(\widetilde B, \widetilde \tau)$, and let $\psi'_\alpha = \psi_\alpha$. Then there is a natural inclusion $\mathfrak{h}_\alpha \subseteq \Lie(\Gamma\mathcal{H}_\alpha)$, and the restrictions of $\Ad(\psi'_\alpha)$ and $\Ad(\psi_\alpha)$ to $\mathfrak{h}_\alpha$ coincide.
\end{proposition}
\begin{proof}
    Recall that we defined $\mathcal{H}_\alpha = \sR_{K^\tau/\mathbb{Q}}\uIso(\widetilde B, \widetilde \tau)$ in this case, so $\Gamma\mathcal{H}_\alpha = \mathcal{H}_\alpha$. Similar to the computation in the proof of Proposition \ref{prop: GammaH_i for HK}, Lemma \ref{lem: S[n] S Iso short exact sequence} implies that the derived Lie algebra of $\Lie(\Gamma\mathcal{H}_\alpha)$ is $\mathfrak{h}_\alpha$. Therefore, we have a natural inclusion $\mathfrak{h}_\alpha \subseteq \Lie(\mathcal{H}_\alpha) = \Lie(\Gamma\mathcal{H}_\alpha)$. Since $\psi_\alpha \in \Gamma\mathcal{H}_\alpha(\mathbb{Q})$, we can take $\psi'_\alpha = \psi_\alpha$, and the equality $\Ad(\psi'_\alpha) = \Ad(\psi_\alpha)$ then automatically holds.
\end{proof}

\subsubsection{Computing the action of $\mathbb{S}[n]$}
For ease of notation, we drop the subscript $K^\tau/\mathbb{Q}$ in the Weil restrictions; the algebraic groups obtained by this are algebraic groups over $\mathbb{Q}$. We will use the rest of this section to prove Proposition \ref{prop: extend iso' to iso}, which we will use to extend the representation of $\sR\uIso(\widetilde B, \widetilde \tau)'$ on $\sH^\bullet(Y; \mathbb{Q})$ to a representation of $\sR\uIso(\widetilde B, \widetilde \tau)$. To do this, we will use the short exact sequence from Lemma~\ref{lem: S[n] S Iso short exact sequence}. We will first compute the action of the algebraic subgroup $\sRS[n] \subseteq \sR\uIso(\widetilde B, \widetilde \tau)'$ on $\extp^\bullet W$, and then extend this action to an action of $\sRS$ that commutes with the action of $\sR \uIso(\widetilde B, \widetilde \tau)'$.

First observe that there is an injective homomorphism of Lie algebras $\mathfrak{gl}(W) \to \mathfrak{so}(\widetilde W)$, sending a linear operator $A \colon W \to W$ to the map $\widetilde A \colon \widetilde W \to \widetilde W$ given by $(w, \eta) \mapsto (Aw, -\eta \circ A)$. By composing with the spinor representation $\mathfrak{so}(\widetilde W) \to \mathfrak{gl}(\bigwedge^\bullet W)$, we obtain a map $\mathfrak{gl}(W) \to \mathfrak{gl}(\bigwedge^\bullet W)$.
\begin{lemma}
    \label{lem: degree 0 of spin representation}
    The composition $\mathfrak{gl}(W) \to \mathfrak{so}(\widetilde W) \to \mathfrak{gl}(\bigwedge^\bullet W)$ sends an endomorphism $A$ of $W$ to the endomorphism of $\bigwedge^\bullet W$ which, for every $k$, sends $w_1 \wedge \dots \wedge w_k \in \bigwedge^k W$ to
    \[-\frac{1}{2}\Tr(A) w_1 \wedge \dots \wedge w_k + \sum_{i=1}^k w_1 \wedge \dots \wedge Aw_i \wedge \dots \wedge w_k.\]
\end{lemma}
\begin{proof}
    This is \cite[Proposition 3.2]{LooijengaLunts}. 
\end{proof}

The Lie algebra $\mathfrak{sl}(W) \subseteq \mathfrak{gl}(W)$ is semisimple, and it has the associated simply connected algebraic group $\uSL(W)$. Therefore the inclusion $\mathfrak{sl}(W) \to \mathfrak{so}(\widetilde W)$ integrates to a morphism of algebraic groups $\uSL(W) \to \uSpin(\widetilde W)$, and by composing with the spinor representation we obtain a representation $\uSL(W) \to \uGL(\bigwedge^\bullet W)$.

\begin{corollary}
    \label{cor: SL(W) spinor rep computation}
    For every commutative $\mathbb{Q}$-algebra $R$, the morphism $\uSL(W) \to \uGL(\bigwedge^\bullet W)$ is given on $R$-points by the map $\SL(W \otimes R) \to \GL(\extp^\bullet_R(W \otimes R))$ which sends the element $A \in \SL(W \otimes R)$ to the automorphism of $\bigwedge_R^\bullet (W \otimes R)$ which, for every $k$, sends $w_1 \wedge \dots \wedge w_k \in \bigwedge_R^k (W \otimes R)$ to
    \begin{equation}
        \label{eq: action of SL on spinor rep}
        Aw_1 \wedge \dots \wedge Aw_k.
    \end{equation}
\end{corollary}
\begin{proof}
    This follows directly by taking derivatives and seeing that we get the action described in Lemma \ref{lem: degree 0 of spin representation}.
\end{proof}

Recall that $\widetilde B = \End_B(\widetilde W)$, and that $\widetilde W$ carries a grading where $W$ sits in degree $1$ and $W^*$ in degree $-1$. The group $\uIso(\widetilde B, \widetilde \tau)'$ acts on $\widetilde W$, so we can look at the subgroup of degree-preserving elements:
\begin{definition}
    Let \myindex{$\uIso_0$} be the algebraic subgroup over $K^\tau$ of $\uIso(\widetilde B, \widetilde \tau)'$ consisting of the degree-preserving elements for the action of $\uIso(\widetilde B, \widetilde \tau)'$ on $\widetilde W$.
\end{definition}

The grading operator $h \colon \widetilde W \to \widetilde W$, given by multiplication by $1$ on $W$ and multiplication by $-1$ on $W^*$, is an element of $\Iso(\widetilde B, \widetilde \tau)'$. The algebraic group $\uIso_0$ is equal to the centralizer subgroup of $h$ in $\uIso(\widetilde B, \widetilde \tau)'$.

\begin{lemma}
    \label{lem: inclusion S[n] into Iso0}
    The algebraic subgroup $\mathbb{S}[n] \subseteq \uIso(\widetilde B, \widetilde \tau)'$ is contained in $\uIso_0$.
\end{lemma}
\begin{proof}
    Since $\mathbb{S}[n]$ is contained in the center of $\uIso(\widetilde B, \widetilde \tau)'$, it automatically commutes with the grading operator $h$, so $\mathbb{S}[n] \subseteq \uIso_0$.
\end{proof}

\begin{lemma}
    \label{lem: Iso0' connected}
    The algebraic $K^\tau$-group $\uIso_0$ is connected.
\end{lemma}
\begin{proof}
    Choose an embedding $K^\tau \injto \mathbb{C}$. It suffices to show that $\uIso_0(\mathbb{C})$ is connected in the analytic topology. Since~$\tau$ is an involution of the second kind, we know that $K/K^\tau$ is a field extension of degree $2$, so if we choose an embedding $K^\tau \to \mathbb{C}$, we have an isomorphism $B \otimes_{K^\tau} \mathbb{C} \cong M_l(\mathbb{C}) \times M_l(\mathbb{C})$ for some integer $l \geq 1$. Let $W_0$ and $W_1$ be the two (up to isomorphism) unique simple $B \otimes_{K^\tau} \mathbb{C}$-modules. We have isomorphisms $W_0^* \cong W_1$ and $W_1^* \cong W_0$, where the duals are computed via the involution (see Definition \ref{def: dual of module over algebra with involution}). Then we have $W \otimes_{K^\tau} \mathbb{C} \cong W_0^m \oplus W_1^m$ for some integer $m$. 

    Since there are no $B$-equivariant morphisms between $W_0$ and $W_1$, we have an isomorphism
    \[ \widetilde B \otimes_{K^\tau} \mathbb{C} \cong \End_B(W_0^m \oplus (W_1^*)^m) \times \End_B(W_1^m \oplus (W_0^*)^m) \cong M_{2m}(\mathbb{C}) \times M_{2m}(\mathbb{C})^{\op}.\]
    Note that $2m = n$. Under this identification, the involution $\tau$ corresponds to the swapping involution. We then see that
    \[\Iso(\widetilde B \otimes_{K^{\tau}} \mathbb{C}, \widetilde \tau) \cong \GL_B(W^m_0 \oplus (W_1^*)^m),\]
    and therefore $\uIso(\widetilde B \otimes_{K^{\tau}} \mathbb{C}, \widetilde \tau)' \cong \SL_B(W_0^m \oplus (W_1^*)^m) \cong \SL_{2m}(\mathbb{C})$.

    The subgroup $\uIso_0(\mathbb{C})$ is then given by those elements which preserve the subspaces $W_0^m$ and $(W_1^*)^m$ of $W^m_0 \oplus (W_1^*)^m$. So, under the identification with $\SL_{2m}(\mathbb{C})$, we see that $\uIso_0(\mathbb{C})$ is given by:
    \begin{equation*}
        \uIso_0(\mathbb{C}) = \left\{
        \begin{pmatrix}
        A & 0\\
        0 & B
        \end{pmatrix}: A, B \in \GL_m(\mathbb{C}) \text{ and } \det(A)\det(B) = 1 \right\}\lowdot
    \end{equation*}
    This description implies that there is a surjective homomorphism of groups $\mathbb{C}^\times \times \SL_m(\mathbb{C}) \times \SL_m(\mathbb{C}) \to \uIso_0(\mathbb{C})$, given by
    \[ (\lambda, A, B) \mapsto \begin{pmatrix} \lambda A & 0\\ 0 & \lambda^{-1}B\end{pmatrix},\]
    and this implies that $\uIso_0(\mathbb{C})$ is connected.
\end{proof}

Since $W$ is an $H$-representation, it is a $B$-module, and hence a $K$-vector space (recall that $K = Z(B)$). In particular, we can write down an action of $\sRS$ on $W$ by scalar multiplication via the inclusion $\sRS \subseteq \sR_{K/\mathbb{Q}}\mathbb{G}_{m, K}$. This gives a map $\sRS \to \uGL(W)$ of algebraic groups over $\mathbb{Q}$.

\begin{lemma}
    \label{lem: sRS to SL(W)}
    The image of the map $\sRS \to \uGL(W)$ is contained in $\uSL(W)$.
\end{lemma}
\begin{proof}
    It suffices to show that for $\lambda \in \sRS(\mathbb{Q})$ the determinant of $\lambda \colon W \to W$ as a $\mathbb{Q}$-linear map is equal to $1$. Since $W$ has the structure of a $K$-vector space, we can write $W \cong K^m$ for an integer $m$, so we have
    \[\det(\lambda \colon  W \to W) = \det(\lambda \colon K \to K)^m.\]
    By \cite[Corollary I.2.7]{Neukirch}, we have
    \[\det(\lambda \colon K \to K) = N_{K/\mathbb{Q}}(\lambda) = N_{K^{\tau}/\mathbb{Q}}(N_{K/K^{\tau}}(\lambda)),\]
    and this is equal to $1$ because $N_{K/K^{\tau}}(\lambda) = 1$, since $\lambda \in \sRS(\mathbb{Q})$.
\end{proof}

By restricting the $\mathbb{Q}$-linear action of $\sRS$ on $W$ to $\sRS[n]$, there is an induced action of $\sRS[n]$ on $W$ by $K$-scalar multiplication, since $(\sRS[n])(\mathbb{Q}) \subseteq K^*$. By Lemma~\ref{lem: inclusion S[n] into Iso0}, we have an inclusion $\mathbb{S}[n] \subseteq \uIso_0$. Since the group $\uIso_0$ consists of the degree-preserving elements, there is a morphism $\sR\uIso_0 \to \uGL(W)$ of algebraic groups over $\mathbb{Q}$. Under the morphism $\sR \uIso_0 \to \uGL(W)$, the group $\sRS[n] \subseteq \sR \uIso_0$ acts on $W$ by scalar multiplication.

By composing the inclusion $\sR\uIso_0 \to \uGL(W)$ with the determinant map $\det \colon \uGL(W) \to \mathbb{G}_{m, \mathbb{Q}}$, we get a determinant map $\det \colon \sR\uIso_0 \to \mathbb{G}_{m, \mathbb{Q}}$. 
\begin{definition}
    Let \myindex{$\uSIso_0$} be the algebgraic group over $\mathbb{Q}$ defined by 
    \[\uSIso_0 = \ker(\det \colon \sR\uIso_0 \to \mathbb{G}_m).\]
\end{definition}
Lemma \ref{lem: sRS to SL(W)} implies that $\sRS[n] \subseteq \uSIso_0$.

The group $\sR \uIso(\widetilde B, \widetilde \tau)'$ is a semisimple simply connected algebraic group, because $\uIso(\widetilde B, \widetilde \tau)'(\mathbb{C}) \cong \SL_{n}(\mathbb{C})$. The inclusion $\mathfrak{s}(\widetilde B, \widetilde \tau) \to \mathfrak{so}(\widetilde W)$ integrates to a morphism 
\begin{equation}
    \label{eq: Iso' to Spin}
    \sR\uIso(\widetilde B, \widetilde \tau)' \to \uSpin(\widetilde W)
\end{equation}
of algebraic groups over $\mathbb{Q}$. Moreover, since $\uSIso_0 \subseteq \sR\uIso_0$, and there is a map $\sR\uIso_0 \to \uGL(W)$, there is a map $\uSIso_0 \to \uSL(W)$ of algebraic groups over $\mathbb{Q}$. There is an inclusion $\uSIso_0 \to \sR\uIso'$ obtained by composing the two inclusions $\uSIso_0 \to \sR\uIso_0$ and $\sR\uIso_0 \to \sR\uIso'$.

These maps give a diagram of algebraic groups over $\mathbb{Q}$:
\begin{equation}
    \label{eq: SIso_0 to GL(wedge W) diagram}
    \begin{tikzcd}
        \uSIso_0 \arrow[d] \arrow[r, hookrightarrow] & \sR\uIso' \arrow[r, "\eqref{eq: Iso' to Spin}"] & \uSpin(\widetilde W) \arrow[d, "\substack{\text{Spinor}\\ \text{representation}}"]\\
        \uSL(W) \arrow[rr, "\eqref{eq: action of SL on spinor rep}"] & & \uGL(\extp^\bullet W),
    \end{tikzcd}
\end{equation}

\begin{remark}
    One could wonder why we do not put $\sR\uIso_0$ in the top left spot in the above diagram. The reason for this is that this group only maps to $\uGL(W)$, and not to $\uSL(W)$. So then we would have to understand the map $\uGL(W) \to \uGL(\extp^\bullet W)$ instead of just the map $\uSL(W) \to \uGL(\extp^\bullet W)$, and this is more difficult due to the extra factor with the trace in Lemma~\ref{lem: degree 0 of spin representation}.
\end{remark}

\begin{lemma}
    \label{lem: SIso action on extp W}
    Diagram \eqref{eq: SIso_0 to GL(wedge W) diagram} commutes.
\end{lemma}

\begin{proof}
    Recall that we have a double cover $\pi \colon \uSpin(\widetilde W) \to \uSO(\widetilde W)$ of algebraic groups over $\mathbb{Q}$. We have the subgroup $\uGL(W) \subseteq \uSO(\widetilde W)$, and can consider its preimage $\pi^{-1}(\uGL(W)) \subseteq \uSpin(\widetilde W)$. The image of the composition $\sR\uIso_0 \to \sR\uIso \to \uSpin(\widetilde W) \to \uSO(\widetilde W)$ is contained in $\uGL(W)$, because $\uGL(W) \subseteq \uSO(\widetilde W)$ is the centralizer of the grading operator. However, we also have the map $\sR\uIso_0 \to \uGL(W)$ coming from the fact that $\uIso_0$ consists of the degree-preserving elements. Hence we obtain a diagram of algebraic groups
    \begin{equation}
        \label{eq: Iso0 GL and Spin}
        \begin{tikzcd}
            \sR\uIso_0 \arrow[r] \arrow[d] & \uSpin(\widetilde W) \arrow[d, "\pi"]\\
            \uGL(W) \arrow[r] & \uSO(\widetilde W).
        \end{tikzcd}
    \end{equation}
    The induced diagram of Lie algebras commutes, and since  $\uIso_0$ is connected (see Lemma~\ref{lem: Iso0' connected}), Lemma \ref{lem: morphism on algebraic groups determined by Lie alg} implies that diagram \eqref{eq: Iso0 GL and Spin} itself commutes.

    It follows that the image of $\sR\uIso_0 \to \uSpin(\widetilde W)$ is contained in $\pi^{-1}(\uGL(W))$. Since $\uSL(W)$ is simply connected and the map $\pi \colon \pi^{-1}(\uSL(W)) \to \uSL(W)$ is a double covering, the group $\pi^{-1}(\uSL(W))$ is disconnected, and its identity component $\pi^{-1}(\uSL(W))^o$ is isomorphic to $\uSL(W)$. The algebraic subgroup $\uSIso_0 \subseteq \sR\uIso_0$ maps to $\uSL(W)$, so from~\eqref{eq: Iso0 GL and Spin} we obtain a commutative diagram of algebraic groups over $\mathbb{Q}$:

    \begin{equation*}
        \begin{tikzcd}
            \uSIso_0 \arrow[r] \arrow[rd] & \pi^{-1}(\uSL(W))^o \arrow[d, "\pi"] \arrow[r] & \uSpin(\widetilde W) \arrow[d, "\pi"]\\
            & \uSL(W) \arrow[r] & \uSO(\widetilde W),
        \end{tikzcd}
    \end{equation*}
    where the vertical map $\pi \colon \pi^{-1}(\uSL(W))^o \to \uSL(W)$ is an isomorphism. It follows that in diagram \eqref{eq: SIso_0 to GL(wedge W) diagram}, by using the inverse of $\pi$, we obtain a morphism $\uSL(W) \to \uSpin(\widetilde W)$:
    \begin{equation*}
        \begin{tikzcd}
            \uSIso_0 \arrow[d] \arrow[r] & \sR\uIso' \arrow[r] & \uSpin(\widetilde W) \arrow[d]\\
            \uSL(W) \arrow[rr] \arrow[rru] & & \uGL(\extp^\bullet W).
        \end{tikzcd}
    \end{equation*}
    We have just shown that the upper triangle commutes, and Corollary \ref{cor: SL(W) spinor rep computation} implies that the bottom triangle commutes.
\end{proof}

We now let the algebraic $\mathbb{Q}$-group $\sRS$ act $\mathbb{Q}$-linearly on $\bigwedge^\bullet W$ by sending, for any commutative $\mathbb{Q}$-algebra $R$, the element $\lambda \in \mathbb{S}(K^\tau \otimes_\mathbb{Q} R) \subseteq (K^\tau \otimes_\mathbb{Q} R)^\times$ to the automorphism of $\bigwedge_R^\bullet W_R$ defined by
\begin{equation}
    \label{eq: action of RS on wedge W}
    w_1 \wedge \dots \wedge w_k \mapsto \lambda w_1 \wedge \dots \wedge \lambda w_k.
\end{equation}
This defines a morphism $\sRS \to \uGL(\extp^\bullet W)$ of algebraic groups over $\mathbb{Q}$.

\begin{lemma}
    \label{lem: Computation of RS[n]}
    The above representation of $\sR_{K^\tau/\mathbb{Q}} \mathbb{S}$ on $\bigwedge^\bullet W$ extends the representation of $\sR_{K^\tau/\mathbb{Q}} \mathbb{S}[n]$ on $\bigwedge^\bullet W$ coming from the representation of $\sR\uIso(\widetilde B, \widetilde \tau)'$, i.e. the following diagram of algebraic groups over $\mathbb{Q}$ commutes:
    \begin{equation*}
        \begin{tikzcd}
            \sRS[n] \arrow[r] \arrow[d] & \sRS \arrow[d, "\eqref{eq: action of RS on wedge W}"]\\
            \sR\uIso(\widetilde B, \widetilde \tau)'_0 \arrow[r] & \uGL(\extp^\bullet W).
        \end{tikzcd}
    \end{equation*}
\end{lemma}

\begin{proof}
    We need to show that for every commutative $\mathbb{Q}$-algebra $R$, every element $\lambda \in \sR\mathbb{S}[n](R)$ acts on $\bigwedge_R^k W_R$ by
    \begin{equation}
        \label{eq: action of S[n] on spinor rep}
        w_1 \wedge \dots \wedge w_k \mapsto \lambda w_1 \wedge \dots \wedge \lambda w_k.
    \end{equation}
    By Lemma \ref{lem: inclusion S[n] into Iso0}, we have an inclusion $\sRS[n] \subseteq \sR\uIso_0$, and Lemma \ref{lem: sRS to SL(W)} implies that $\sRS[n]$ is in fact contained in the subgroup $\uSIso_0$. By Lemma \ref{lem: SIso action on extp W}, we may compute the action of the subgroup $\uSIso_0 \subseteq \sR\uIso(\widetilde B, \widetilde \tau)$ on $\bigwedge^\bullet W$ via the action of $\uSL(W)$ on $\bigwedge^\bullet W$ from Corollary \ref{cor: SL(W) spinor rep computation}, and therefore $\sRS[n]$ indeed acts on $\extp^\bullet W$ as described in equation~\eqref{eq: action of S[n] on spinor rep}.
\end{proof}

\subsubsection{Extending the representation of $\sR_{K^\tau/\mathbb{Q}}\uIso(\widetilde B, \widetilde \tau)'$}
Lastly, we will need to check that the representation of $\sR_{K^\tau/\mathbb{Q}} \mathbb{S}$ on $\extp^\bullet W$ described in equation \eqref{eq: action of RS on wedge W} commutes with the representation of $\sR_{K^\tau/\mathbb{Q}}\uIso(\widetilde B, \widetilde \tau)'$ on $\bigwedge^\bullet W$, so that we actually get a representation of $\sRS \times \sR \uIso(\widetilde B, \widetilde \tau)'$. First, we will show that the maps from $\sR \mathbb{S}$ and $\sR\uIso(\widetilde B, \widetilde \tau)'$ to $\uGL(\extp^\bullet W)$ both factor via the spinor representation $\uSpin(\widetilde W) \to \uGL(\extp^\bullet W)$. Once we have this, it suffices to show that the images of $\sR \mathbb{S}$ and $\sR\uIso(\widetilde B, \widetilde \tau)'$ in $\uSpin(\widetilde W)$ commute, which we check by computing the images of their Lie algebras in $\mathfrak{so}(\widetilde W)$.

In Corollary \ref{cor: SL(W) spinor rep computation}, we described the morphism $\uSL(W) \to \uGL(\extp^\bullet W)$ integrating the composition $\mathfrak{sl}(W) \to \mathfrak{so}(\widetilde W) \to \mathfrak{gl}(\extp^\bullet W)$, where the second map is the spinor representation. 
\begin{lemma}
    \label{lem: image of sRS in SL(W)}
    The image of the map $\sR_{K^\tau/\mathbb{Q}}\mathbb{S} \to \uGL(\extp^\bullet W)$ is contained in the image of the map $\uSL(W) \to \uGL(\extp^\bullet W)$ from Corollary \ref{cor: SL(W) spinor rep computation}.
\end{lemma}
\begin{proof}
    Since $\sRS$ acts on $W$ by scalar multiplication, we have a map $\sRS \to \uGL(W)$. We know from Lemma \ref{lem: sRS to SL(W)} that the image of this map is actually contained in $\uSL(W)$. It follows from Corollary \ref{cor: SL(W) spinor rep computation} that the composition $\sRS \to \uSL(W) \to \uGL(\extp^\bullet W)$ is equal to the map $\sRS \to \uGL(\extp^\bullet W)$ from Equation~\eqref{eq: action of RS on wedge W}.
\end{proof}

\begin{lemma}
    The Lie algebra of $\mathbb{S}$ is isomorphic to $\Skew(K, \tau)$.
\end{lemma}
\begin{proof}
    Note that $\mathbb{S}$ is the algebraic subgroup of $\sR_{K/K^\tau}\mathbb{G}_{m, K}$ given by the equation $\lambda \tau(\lambda) = 1$. Since $\Lie(\sR_{K/K^\tau}\mathbb{G}_{m, K}) = K$, we see that $\Lie(\mathbb{S})$ is the Lie subalgebra of $K$ defined by the equation $x + \tau(x) = 0$, which is exactly $\Skew(K, \tau)$.
\end{proof}

\begin{lemma}
    \label{lem: Skew of unitary involution}
    Let $(\widetilde B, \widetilde \tau)$ be an algebra with involution of the second kind, and let $K = Z(\widetilde B)$. Then we have a direct sum decomposition of Lie algebras
    \[ \Skew(\widetilde B, \widetilde \tau) \cong \mathfrak{s}(\widetilde B, \widetilde \tau) \oplus \Skew(K, \tau). \]
\end{lemma}

\begin{proof}
    By definition we have an inclusion $\mathfrak{s}(\widetilde B, \widetilde \tau) \subseteq \Skew(\widetilde B, \widetilde \tau)$, and the inclusion $K \subseteq B$ induces an inclusion $\Skew(K, \tau) \subseteq \Skew(\widetilde B, \widetilde \tau)$. Hence we obtain a map
    \begin{equation}
        \label{eq: s oplus Skew to Skew}
        \mathfrak{s}(\widetilde B, \widetilde \tau) \oplus \Skew(K, \tau) \to \Skew(\widetilde B, \widetilde \tau).
    \end{equation}
    To check that this is an isomorphism, it suffices to do this after taking the tensor product $- \otimes_{K^\tau} \overline{K}$. Let $n$ be the integer such that $\widetilde B \otimes_{K^\tau} \overline K \cong M_n(\overline{K}) \times M_n(\overline K)$, then the map in Equation \eqref{eq: s oplus Skew to Skew} becomes the isomorphism
    \[\mathfrak{sl}_n(\overline{K}) \oplus \overline K \cong \mathfrak{gl}_n(\overline K).\]
\end{proof}

\begin{lemma}
    \label{lem: actions of S and Iso' commute}
    The actions of $\sR_{K^\tau/\mathbb{Q}} \mathbb{S}$ and $\sR_{K^\tau/\mathbb{Q}} \uIso(\widetilde B, \widetilde \tau)'$ on $\bigwedge^\bullet W$ commute.
\end{lemma}
\begin{proof}
    Lemma \ref{lem: image of sRS in SL(W)} implies that the map $\sRS \to \uGL(\extp^\bullet W)$ factorizes via the group $\uSpin(\widetilde W)$. On the level of Lie algebras this map $\sRS \to \uSpin(\widetilde W)$ induces the map $\Skew(K, \tau) \to \mathfrak{so}(\widetilde W)$, which is the restriction of the map $K \to \End_\mathbb{Q}(\widetilde W)$ (with $W^*$ a $K$-module via $\tau$).

    The map $\sR \uIso(\widetilde B, \widetilde \tau)' \to \uGL(\extp^\bullet W)$ was defined as the map integrating the composition of Lie algebra homomorphisms $\mathfrak{s}(\widetilde B, \widetilde \tau) \to \mathfrak{so}(\widetilde W) \to \mathfrak{gl}(\extp^\bullet W)$, so it is clear that the map $\sR \uIso(\widetilde B, \widetilde \tau)' \to \uGL(\extp^\bullet W)$ also factorizes via $\uSpin(\widetilde W)$, and that the map on Lie algebras induced by $\sR \uIso(\widetilde B, \widetilde \tau)' \to \uSpin(\widetilde W)$ is just the inclusion $\mathfrak{s}(\widetilde B, \widetilde \tau) \to \mathfrak{so}(\widetilde W)$.

    The inclusion $\widetilde B := \End_B(\widetilde W) \subseteq \mathfrak{gl}(\widetilde W)$ gives an inclusion $\Skew(\widetilde B, \widetilde \tau) \subseteq \mathfrak{so}(\widetilde W)$, and the inclusions $\Skew(K, \tau) \to \mathfrak{so}(\widetilde W)$ and $\mathfrak{s}(\widetilde B, \widetilde \tau) \to \mathfrak{so}(\widetilde W)$ from above induce the isomorphism $\mathfrak{s}(\widetilde B, \widetilde \tau) \oplus \Skew(K ,\tau) \cong \Skew(\widetilde B, \widetilde \tau)$ from Lemma~\ref{lem: Skew of unitary involution}. The Lie algebas $\Skew(K, \tau)$ and $\mathfrak{s}(\widetilde B, \widetilde \tau)$ commute as Lie subalgebras of $\Skew(\widetilde B, \widetilde \tau)$, and hence the actions of $\sRS$ and $\sR \uIso(\widetilde B, \widetilde \tau)$ on $\extp^\bullet W$ also commute.
\end{proof}

We can now prove the extension result that we are really interested in:
\begin{proposition}
    \label{prop: extend iso' to iso}
    The representation of $\sR_{K^\tau/\mathbb{Q}} \uIso(\widetilde B, \widetilde \tau)'$ on $\extp^\bullet W$ obtained by integrating the representation of $\mathfrak{s}(\widetilde B, \widetilde \tau)$ extends to a representation of the group $\sR_{K^\tau/\mathbb{Q}} \uIso(\widetilde B, \widetilde \tau)$.
\end{proposition}
\begin{proof}
    By Lemma \ref{lem: S[n] S Iso short exact sequence}, it suffices to write down an action of $\sRS$ on $\extp^\bullet W$ which commutes with the action of $\sR \uIso(\widetilde B, \widetilde \tau)'$ and extends the action of $\sRS[n]$. For this we use the action given in Equation \eqref{eq: action of RS on wedge W}, which indeed commutes with the action of $\sR\uIso(\widetilde B, \widetilde \tau)'$ by Lemma \ref{lem: actions of S and Iso' commute}, and it extends the action of $\sRS[n]$ by Lemma~\ref{lem: Computation of RS[n]}.
\end{proof}

Recall that we defined $\Gamma\mathcal{H}_\alpha = \sR_{K^\tau/\mathbb{Q}} \uIso(\widetilde B, \widetilde \tau)$. The induced representation of $\Lie(\Gamma\mathcal{H}_\alpha)$ on $\extp^\bullet W$ is not too hard to compute now:
\begin{lemma}
    \label{lem: representation of Lie(Gamma H) unitary case}
    There is an isomorphism $\Lie(\Gamma\mathcal{H}_\alpha) \cong \Skew(\widetilde B, \widetilde \tau)$, and under this isomorphism the induced representation of $\Lie(\Gamma\mathcal{H}_\alpha)$ on $\extp^\bullet W$ factors as
    \[\Skew(\widetilde B, \widetilde \tau) \to \mathfrak{so}(\widetilde W) \to \mathfrak{gl}\left(\extp^\bullet W\right),\]
    where the first map comes from the inclusion $\widetilde B \subseteq \End_\mathbb{Q}(\widetilde W)$, and the second map is the spinor representation.
\end{lemma}
\begin{proof}
    The isomorphism $\Lie(\Gamma\mathcal{H}_\alpha) \cong \Skew(\widetilde B, \widetilde \tau)$ follows from Lemma~\ref{lem: Lie algebra of Weil restriction} and Lemma~\ref{lem: Lie algebra of iso is skew}. To prove that $\Skew(\widetilde B, \widetilde \tau) \to \mathfrak{gl}(\extp^\bullet W)$ factors over the spinor representation, it suffices to show this for the factors $\Skew(K, \tau)$ and $\mathfrak{s}(\widetilde B, \widetilde \tau)$ in the direct sum decomposition from Lemma~\ref{lem: Skew of unitary involution}. 

    For the factor $\mathfrak{s}(\widetilde B, \widetilde \tau)$, this follows from the construction of the representation of $\uIso(\widetilde B, \widetilde \tau)'$ on $\extp^\bullet W$ by integrating the spinor representation. For the factor $\Skew(K, \tau)$, this follows from Lemma~\ref{lem: image of sRS in SL(W)}.
\end{proof}

\begin{remark}
    One could wonder why we did not use a strategy similar to what was used in Corollary \ref{cor: extending mu_n-rep to G_m-rep} to extend the representation from $\sRS[n]$ to $\sR \mathbb{S}$. The reason for this is that, in the proof of Proposition \ref{prop: Rmu to RG extension}, we used in an essential way that for every Galois orbit $T \subseteq X(\sR \mu_n)$ there is a Galois orbit $T' \subseteq X(\sR\mathbb{G}_m)$ mapping bijectively onto $T$. Unfortunately, the inclusion $\sRS[n] \subseteq \sR\mathbb{S}$ does not satisfy the analogous property. To see why, assume for simplicity that we are in the situation where $K^\tau = \mathbb{Q}$ (and therefore $K/\mathbb{Q}$ is a quadratic extension). Then $X(\mathbb{S}) \cong \mathbb{Z}$, and the action of $\Gal(\overline{\mathbb{Q}} / \mathbb{Q})$ on $X(\mathbb{S})$ comes from the action of $\Gal(K/\mathbb{Q}) \cong \mathbb{Z}/2\mathbb{Z}$ on $\mathbb{Z}$, where the nontrivial element acts as multiplication by $-1$. Similarly, we have $X(\mathbb{S}[n]) \cong \mathbb{Z}/n\mathbb{Z}$, and the nontrivial element of $\Gal(K/\mathbb{Q})$ acts on this as multiplication by $-1$. Since $n$ is even (see the proof of Lemma~\ref{lem: Iso0' connected}), we see that $\{n/2\} \subseteq \mathbb{Z}/n\mathbb{Z}$ is a Galois orbit, but there is no Galois orbit in $\mathbb{Z}$ which maps bijectively onto $\{n/2\}$.
\end{remark}

\subsection{Simple factors of type $A_1$}
\label{sec: A1 case}
Lastly, suppose that $\mathfrak{h}_\alpha$ is a simple factor of $\llv(Y; \mathbb{Q})$ of type $A_1$ and let $F$ be its centroid. Recall that Lemma~\ref{lem: simple factor A1 of llv is sl2(F)} then implies that $\mathfrak{h}_\alpha \cong \mathfrak{sl}_2(F)$. We define $\Gamma\mathcal{H}_\alpha = \sR_{F/\mathbb{Q}} \uGL_2(F)$, recall that this is equal to $\mathcal{H}_\alpha$. 

\begin{lemma}
    \label{lem: GL_n ses}
    There is a short exact sequence of algebraic groups over $\mathbb{Q}$
    \[ 1 \to \sR_{F/\mathbb{Q}} \mu_{2, F} \to \sR_{F/\mathbb{Q}} \mathbb{G}_{m, F} \times \sR_{F/\mathbb{Q}} \uSL_2(F) \to \sR_{F/\mathbb{Q}} \uGL_2(F) \to 1.\]
\end{lemma}
\begin{proof}
    This follows from Lemmas \ref{lem: weil restriction exact} and \ref{lem: weil restriction product}.
\end{proof}

Just as in the symplectic case, the proof of the third property of Theorem \ref{thm: properties of GammaH_i} is straightforward in the $A_1$ case:
\begin{proposition}
    \label{prop: GammaH_i type A1}
    Assume that we are in the setting of Theorem \ref{thm: properties of GammaH_i}, and let $\mathfrak{h}_\alpha$ be a simple factor of $\llv(Y; \mathbb{Q})$ of type $A_1$. Then:
    \begin{enumerate}
        \item There is a natural inclusion $\mathfrak{h}_\alpha \subseteq \Lie(\Gamma\mathcal{H}_\alpha)$.
        \item There is a $\psi_\alpha' \in \Gamma\mathcal{H}_\alpha(\mathbb{Q})$ such that $\Ad(\psi_\alpha) = \Ad(\psi'_\alpha)$ as automorphisms of $\mathfrak{h}_\alpha$.
    \end{enumerate}
\end{proposition}
\begin{proof}
    Lemma \ref{lem: simple factor A1 of llv is sl2(F)} implies that $\mathfrak{h}_\alpha \cong \mathfrak{sl}_2(F)$. Since $\Lie(\Gamma\mathcal{H}_\alpha) \cong F \times \mathfrak{sl}_2(F)$, we obtain an inclusion $\mathfrak{h}_\alpha \subseteq \Lie(\Gamma\mathcal{H}_\alpha)$.

    For the second statement, observe that $\Gamma\mathcal{H}_\alpha = \mathcal{H}_\alpha$, so we can just take $\psi_\alpha' = \psi_\alpha$.
\end{proof}

\section{Finishing the proof of Theorem \ref{thm: properties of GammaH_i}}
\label{subsec: proof of thm properties of GammaH_i}
In this section we will combine the knowledge gained about the groups $\Gamma\mathcal{H}_\alpha$ for the different types of $\mathfrak{h}_\alpha$ to finish the proof of Theorem \ref{thm: properties of GammaH_i}. As seen in the beginning of this chapter, this finishes the proof of Theorem \ref{thm: main theorem}.

In the previous sections, we defined for every simple factor $\mathfrak{h}_\alpha$ of $\llv(Y; \mathbb{Q})$ an algebraic group $\Gamma\mathcal{H}_\alpha$ over $\mathbb{Q}$ by:
\begin{equation*}
    \Gamma\mathcal{H}_\alpha = \begin{cases}
         \sR_{F/\mathbb{Q}} \uIso(\widetilde B, \widetilde \tau) &\text{ if } \mathfrak{h}_\alpha \cong \mathfrak{s}(\widetilde B, \widetilde \tau) \text{ and } \widetilde \tau \text{ is symplectic}\\
            \sR_{F/\mathbb{Q}} \uIso(\widetilde B, \widetilde \tau) &\text{ if } \mathfrak{h}_\alpha \cong \mathfrak{s}(\widetilde B, \widetilde \tau) \text{ and } \widetilde \tau \text{ is unitary}\\ 
            \sR_{F/\mathbb{Q}} \uGSpin(\widetilde B, \widetilde \tau) &\text{ if } \mathfrak{h}_\alpha \cong \mathfrak{s}(\widetilde B, \widetilde \tau) \text{ and } \widetilde \tau \text{ is orthogonal}\\
            \sR_{F/\mathbb{Q}}\uGL_2(F) & \text{ if  } \mathfrak{h}_\alpha \cong \mathfrak{sl}_2(F)\\
            \uGSpin(\widetilde \sH(Y_j; \mathbb{Q})^{H_j}) &\text{ if }\mathfrak{h}_\alpha \cong \mathfrak{so}(\widetilde \sH(Y_j; \mathbb{Q})^{H_j}),
    \end{cases}
\end{equation*}
where in the first three cases we let $F = Z(\widetilde B)^{\widetilde \tau}$ and assume that $\mathfrak{h}_\alpha$ is not of type $A_1$.
In order to prove Theorem~\ref{thm: properties of GammaH_i} it remains to construct representations of the groups $\Gamma\mathcal{H}_\alpha$ on $\sH^\bullet(Y; \mathbb{Q})$ satisfying properties (1) and (2) of the theorem. To do this, we will first need the simply connected algebraic groups with Lie algebra $\mathfrak{h}_\alpha$.

\begin{definition}
    \label{def: widetilde Hi}
    For every $\alpha$, let $\widetilde{\mathcal{H}}_\alpha$ be the algebraic group over $\mathbb{Q}$ defined by
    \begin{equation*}
        \widetilde{\mathcal{H}}_\alpha = \begin{cases}
            \sR_{F/\mathbb{Q}} \uIso(\widetilde B, \widetilde \tau) &\text{ if } \mathfrak{h}_\alpha \cong \mathfrak{s}(\widetilde B, \widetilde \tau) \text{ and } \widetilde \tau \text{ is symplectic}\\
            \sR_{F/\mathbb{Q}} \uIso(\widetilde B, \widetilde \tau)' &\text{ if } \mathfrak{h}_\alpha \cong \mathfrak{s}(\widetilde B, \widetilde \tau) \text{ and } \widetilde \tau \text{ is unitary}\\ 
            \sR_{F/\mathbb{Q}} \uSpin(\widetilde B, \widetilde \tau) &\text{ if } \mathfrak{h}_\alpha \cong \mathfrak{s}(\widetilde B, \widetilde \tau) \text{ and } \widetilde \tau \text{ is orthogonal}\\
            \sR_{F/\mathbb{Q}} \uSL_2(F) & \text{ if } \mathfrak{h}_\alpha \cong \mathfrak{sl}_2(F)\\
            \uSpin(\widetilde \sH(Y_j; \mathbb{Q})^{H_j}) &\text{ if }\mathfrak{h}_\alpha \cong \mathfrak{so}(\widetilde \sH(Y_j; \mathbb{Q})^{H_j}),
        \end{cases}
    \end{equation*}
    where in the first three cases we let $F = Z(\widetilde B)^{\widetilde \tau}$ and assume that $\mathfrak{h}_\alpha$ is not of type $A_1$.
\end{definition}
By passing to $\mathbb{C}$ and using Lemma \ref{lem: simply connected iff}, we indeed see for every $\alpha$ that $\widetilde{\mathcal{H}}_\alpha$ is the simply connected algebraic group over $\mathbb{Q}$ with Lie algebra $\mathfrak{h}_\alpha$.

In the proof, we will also need the following notions:
\begin{definition}
    Let $k$ be a field, let $\mathfrak{g}$ be a Lie algebra over $k$, let $H$ be a finite group and let $V$ be an $H$-representation over $k$, and assume that $V$ is also a $\mathfrak{g}$-representation. We say that $\mathfrak{g}$ acts $H$\textbf{-equivariantly} on $V$ if for all $h \in H$ and $x \in \mathfrak{g}$ the actions of $h$ and $x$ on $V$ commute. 
\end{definition}
For a field $k$ and a finite group $H$, if $V$ is an $H$-representation over $k$, the map $\rho \colon H \to \GL(V)$ induces a map $\underline{\rho} \colon H \to \uGL(V)$ of algebraic groups over $k$, where $H$ is seen as a constant algebraic group over $k$. 

\begin{definition}
    Let $\mathcal{G}$ be an algebraic group over a field $k$, let $H$ be a finite group, and let $\rho \colon H \to \GL(V)$ be an $H$-representation over $k$. Assume that $V$ is also a $\mathcal{G}$-representation. We say that $\mathcal{G}$ acts $H$\textbf{-equivariantly} on $V$ if the image of $\mathcal{G} \to \uGL(V)$ is contained in the centralizer of $\underline{\rho}(H) \subseteq \uGL(V)$.
\end{definition}
If this is the case, then \cite[Proposition 1.92]{milneiAG} implies that for every $k$-algebra $R$, the actions of $\mathcal{G}(R)$ and $H$ on $V \otimes_k R$ commute.
\begin{lemma}
    \label{lem: equivariant group action checked on Lie algebra}
    Let $k$ be a field of characteristic $0$, let $\mathcal{G}$ be a connected algebraic group over $k$, let $H$ be a finite group, and let $V$ be an $H$-representation over $k$. Assume that $V$ is also a $\mathcal{G}$-representation. If the induced action of $\Lie(\mathcal{G})$ on $V$ is $H$-equivariant, then the action of $\mathcal{G}$ on $V$ is also $H$-equivariant.
\end{lemma}
\begin{proof}
    Take $h \in H$. We have to show that the two maps $\rho, h\rho h^{-1} \colon \mathcal{G} \to \uGL(V)$ coincide. Since $\Lie(\mathcal{G})$ acts $H$-equivariantly on $V$, we know that the induced maps on the level of Lie algebras $\Lie(\rho), \Lie(h\rho h^{-1}) \colon \Lie(\mathcal{G}) \to \mathfrak{gl}(V)$ coincide. Since $\mathcal{G}$ is connected, the result follows from Lemma \ref{lem: morphism on algebraic groups determined by Lie alg}.
\end{proof}

\begin{proof}[{Proof of Theorem \ref{thm: properties of GammaH_i}}]
    For every $\alpha$ we have constructed groups $\Gamma\mathcal{H}_\alpha$ and elements $\psi_\alpha' \in \Gamma\mathcal{H}_\alpha(\mathbb{Q})$ satisfying condition (3) in Propositions~\ref{prop: GammaH_i for HK}, \ref{prop: GammaHi orthogonal}, \ref{prop: GammaHi symplectic}, \ref{prop: GammaHi unitary} and~\ref{prop: GammaH_i type A1}. So to finish the proof of the theorem, we need to construct mutually commuting representations of the groups $\Gamma\mathcal{H}_\alpha$ on $\sH^\bullet(Y; \mathbb{Q})$, satisfying condition~(1).

    Write $Y = (\prod_{j=0}^l Y_j)/H$ as in Setup \ref{stp: symplectic var}, then
    \[\sH^\bullet(Y; \mathbb{Q}) \cong \sH^\bullet\Bigl (\prod_j Y_j; \mathbb{Q} \Bigl)^H.\]
    Every simple factor $\mathfrak{h}_\alpha$ of $\llv(Y; \mathbb{Q})$ is either isomorphic to $\mathfrak{so}(\widetilde \sH(Y_i; \mathbb{Q})^{H_i})$ for some hyperkähler factor $Y_i$, or $\mathfrak{h}_\alpha$ is isomorphic to (a Lie subalgebra of) $\mathfrak{s}(\End_B(\widetilde W_\rho), \widetilde \tau)$ for some $\rho \in \Irr_\mathbb{Q}(H)$, with $W_\rho \subseteq \sH^1(Y_0; \mathbb{Q})$ the $\rho$-isotypic component. We then define for every $\alpha$ the subring $\sH^\bullet_\alpha(Y; \mathbb{Q})$ of the cohomology $\sH^\bullet(\prod_j Y_j; \mathbb{Q})$ as follows:
    \begin{equation*}
        \sH^\bullet_\alpha(Y; \mathbb{Q}) = \begin{cases}
            \sH^\bullet(\prod_{j \in [i]} Y_j; \mathbb{Q}) & \text{if } \mathfrak{h}_\alpha \cong \mathfrak{so}(\widetilde \sH(Y_i; \mathbb{Q})^{H_i})\\
            \extp^\bullet W_\rho & \text{if } \mathfrak{h}_\alpha \injto \mathfrak{s}(\End_B(\widetilde W_\rho), \widetilde \tau).
        \end{cases}
    \end{equation*}
    For every $\alpha$, the group $H$ acts on $\sH^\bullet_\alpha(Y; \mathbb{Q})$, and we have an isomorphism 
    \[\sH^\bullet(Y; \mathbb{Q}) \cong \left(\bigotimes_\alpha \sH_\alpha^\bullet(Y; \mathbb{Q})\right)^H.\]
    We will first construct for every simple factor $\mathfrak{h}_\alpha$ an $H$-equivariant representation of $\mathfrak{h}_\alpha$ on $\sH^\bullet_\alpha(Y; \mathbb{Q})$.
    
    For the factors of hyperkähler origin, we can use the inclusion 
    \[ \mathfrak{h}_\alpha \cong \mathfrak{so}(\widetilde \sH(Y_i; \mathbb{Q})^{H_i}) \subseteq \prod_{j \in [i]} \mathfrak{so}(\widetilde \sH(Y_j; \mathbb{Q}))\]
    from the proof of Proposition \ref{prop: LLV of HK}, and let the product $\prod_{j \in [i]} \mathfrak{so}(\widetilde \sH(Y_j; \mathbb{Q}))$ act on $\sH^\bullet_\alpha(Y; \mathbb{Q}) \cong \bigotimes_{j \in [i]} \sH^\bullet(Y_j; \mathbb{Q})$ by taking the tensor product of the representations of the Lie algebras $\mathfrak{so}(\widetilde \sH(Y_j; \mathbb{Q}))$ on $\sH^\bullet(Y_j; \mathbb{Q})$. Since all elements of $\mathfrak{h}_\alpha$ are $H$-invariant, the Lie algebra $\mathfrak{h}_\alpha$ acts $H$-equivariantly on $\sH^\bullet_\alpha(Y; \mathbb{Q})$.
    
    If $\mathfrak{h}_\alpha \injto \mathfrak{s}(\widetilde B, \widetilde \tau)$ is of abelian origin, then we have an injective map $\mathfrak{h}_\alpha \injto \mathfrak{so}(\widetilde{W}_\rho)$ which is the composition of the map $\mathfrak{h}_\alpha \injto \mathfrak{s}(\widetilde B, \widetilde \tau)$ (which is an isomorphism if $\mathfrak{h}_\alpha$ is not of type $A_1$) and the inclusion $\mathfrak{s}(\widetilde B, \widetilde \tau) \subseteq \mathfrak{so}(\widetilde W_{\rho})$. Since $\mathfrak{so}(\widetilde W_\rho)$ acts on $\extp^\bullet W_\rho$ via the spinor representation, we obtain an action of $\mathfrak{h}_\alpha$ on $\extp^\bullet W_\rho$. Since we even have $\mathfrak{h}_\alpha \subseteq \mathfrak{so}(\widetilde{W}_\rho)^H$, we see that $\mathfrak{h}_\alpha$ acts $H$-equivariantly on $\sH^\bullet_\alpha(Y; \mathbb{Q})$.
    
    By taking the tensor product of these representations of the $\mathfrak{h}_\alpha$, we obtain an $H$-equivariant action of $\llv(Y; \mathbb{Q})$ on $\bigotimes_\alpha \sH^\bullet_\alpha(Y; \mathbb{Q})$. Since this action is $H$-equivariant, it restricts to a representation of $\llv(Y; \mathbb{Q})$ on $(\bigotimes_\alpha \sH^\bullet_\alpha(Y; \mathbb{Q}))^H = \sH^\bullet(Y; \mathbb{Q})$, and this recovers the action of $\llv(Y; \mathbb{Q})$ on $\sH^\bullet(Y; \mathbb{Q})$.

    For every $\alpha$, the representation of $\mathfrak{h}_\alpha$ on $\sH^\bullet_\alpha(Y; \mathbb{Q})$ integrates to a representation of $\widetilde{\mathcal{H}}_\alpha$ by Theorem~\ref{thm: integrate Lie algebra action}. Since the elements of $\mathfrak{h}_\alpha$ act $H$-equivariantly on $\sH^\bullet_\alpha(Y; \mathbb{Q})$, the action of $\widetilde{\mathcal{H}}_\alpha$ on $\sH^\bullet_\alpha(Y; \mathbb{Q})$ is also $H$-equivariant by Lemma~\ref{lem: equivariant group action checked on Lie algebra}. We now claim that these representations of the $\widetilde{\mathcal{H}}_\alpha$ extend to representations of the groups $\Gamma\mathcal{H}_\alpha$, in such a way that all elements of $\Gamma\mathcal{H}_\alpha(\mathbb{Q})$ act $H$-equivariantly on $\sH^\bullet_\alpha(Y; \mathbb{Q})$. By then taking the tensor product of the representations of the $\Gamma\mathcal{H}_\alpha$ on the $\sH^\bullet_\alpha(Y; \mathbb{Q})$, we get an $H$-equivariant representation of $\prod_\alpha \Gamma\mathcal{H}_\alpha$ on $\bigotimes_\alpha \sH^\bullet_\alpha(Y; \mathbb{Q})$. By taking the $H$-invariants, this gives a representation of $\prod_\alpha \Gamma\mathcal{H}_\alpha$ on $\sH^\bullet(Y; \mathbb{Q})$, which shows that the representations of the $\Gamma\mathcal{H}_\alpha$ indeed commute. Hence, we are done once we have extended the representations from the groups $\widetilde{\mathcal{H}}_\alpha$ to $\Gamma\mathcal{H}_\alpha$. There are five different cases to cover.

    First assume that $\mathfrak{h}_\alpha \cong \mathfrak{so}(\widetilde \sH(Y_i; \mathbb{Q})^{H_i})$ is of hyperkähler origin. Lemma~\ref{lem: HK Spin short exact sequences} gives a short exact sequence
    \[ 1 \to \mu_{2, \mathbb{Q}} \to \mathbb{G}_{m, \mathbb{Q}} \times \widetilde{\mathcal{H}}_\alpha \to \Gamma\mathcal{H}_\alpha \to 1.\]
    By Corollary \ref{cor: extending mu_n-rep to G_m-rep}, the representation of $\widetilde{\mathcal{H}}_\alpha$ on $\sH^\bullet_\alpha(Y; \mathbb{Q})$ extends to a representation of $\Gamma\mathcal{H}_\alpha$, and the centralizers of $\widetilde{\mathcal{H}}_\alpha$ and $\Gamma\mathcal{H}_\alpha$ in $\uGL(\sH^\bullet(Y; {\mathbb{Q}}))$ coincide. Since the action of $\widetilde{\mathcal{H}}_\alpha$ on $\sH^\bullet(Y; \mathbb{Q})$ is $H$-equivariant, we have an inclusion 
    \[ H \subseteq C_{\uGL(\sH^\bullet_\alpha(Y; {\mathbb{Q}}))}(\widetilde{\mathcal{H}}_\alpha({\mathbb{Q}})) = C_{\uGL(\sH^\bullet_\alpha(Y; {\mathbb{Q}}))}(\Gamma\mathcal{H}_\alpha({\mathbb{Q}})),\]
    which implies that $\Gamma\mathcal{H}_\alpha$ acts $H$-equivariantly on $\sH^\bullet_\alpha(Y; \mathbb{Q})$.
    
    If $\mathfrak{h}_\alpha \cong \mathfrak{s}(\widetilde B, \widetilde \tau)$ with $\widetilde \tau$ an orthogonal involution and $\mathfrak{h}_\alpha$ is not of type $A_1$, then we can argue similarly as above. Lemma \ref{lem: GSpin SO short exact sequences} provides a short exact sequence
    \[ 1 \to \sR_{F/k}\mu_{2, F} \to \sR_{F/k}\mathbb{G}_{m, F} \times \widetilde{\mathcal{H}}_\alpha \to \Gamma\mathcal{H}_\alpha \to 1,\]
    where $F = Z(\widetilde B)$, and we can again apply Corollary \ref{cor: extending mu_n-rep to G_m-rep} to extend the representation of $\widetilde{\mathcal{H}}_\alpha$ to an $H$-equivariant representation of $\Gamma\mathcal{H}_\alpha$.

    For the third case, if $\mathfrak{h}_\alpha \cong \mathfrak{s}(\widetilde B, \widetilde \tau)$ with $\widetilde \tau$ a symplectic involution and $\mathfrak{h}_\alpha$ is not of type $A_1$, then we are directly done, since $\widetilde{\mathcal{H}}_\alpha = \Gamma\mathcal{H}_\alpha$.

    For the fourth case, we have $\mathfrak{h}_\alpha \cong \mathfrak{s}(\widetilde B, \widetilde \tau)$ with $\widetilde \tau$ a unitary involution and $\mathfrak{h}_\alpha$ not of type $A_1$. Proposition~\ref{prop: extend iso' to iso} implies that the representation of $\widetilde{\mathcal{H}}_\alpha$ on $\extp^\bullet W_\rho$ extends to a representation of $\Gamma\mathcal{H}_\alpha$. It remains to see that $\Gamma\mathcal{H}_\alpha$ acts $H$-equivariantly on $\extp^\bullet W_\rho$. We have $\Lie(\Gamma\mathcal{H}_\alpha) = \Skew(\widetilde B, \widetilde \tau)$, and this is contained in $\mathfrak{so}(\widetilde{W}_\rho)^H$. Lemma~\ref{lem: representation of Lie(Gamma H) unitary case} then implies that the map on Lie algebras induced by the morphism $\Gamma\mathcal{H}_\alpha \to \uGL(\extp^\bullet W_\rho)$ of algebraic groups over $\mathbb{Q}$ factorizes as
    \[ \Skew(\widetilde B, \widetilde \tau) \to \mathfrak{so}(\widetilde W_{\rho})^H \to \mathfrak{gl}\Bigl(\extp^\bullet W_\rho\Bigl),\]
    where the second map is the restriction of the spinor representation. It then follows from Lemma \ref{lem: equivariant group action checked on Lie algebra} that $\Gamma\mathcal{H}_\alpha$ acts $H$-equivariantly on $\extp^\bullet W_\rho$.

    Finally, suppose that $\mathfrak{h}_\alpha \cong \mathfrak{sl}_2(F)$. Lemma \ref{lem: GL_n ses} gives a short exact sequence
    \[1 \to \sR_{F/\mathbb{Q}} \mu_{2, F} \to \sR_{F/\mathbb{Q}} \mathbb{G}_{m, \mathbb{Q}} \times \widetilde{\mathcal{H}}_\alpha \to \Gamma\mathcal{H}_\alpha \to 1.\]
    As above, Corollary \ref{cor: extending mu_n-rep to G_m-rep} implies that the $H$-equivariant representation of $\widetilde{\mathcal{H}}_\alpha$ on $\sH^\bullet_\alpha(Y; \mathbb{Q})$ extends to an $H$-equivariant representation of $\Gamma\mathcal{H}_\alpha$, and this finishes the proof.

\end{proof}

%% file: chapters/chapter8.tex
During the proof of Theorem \ref{thm: main theorem}, we had to assume that the varieties $X$ and $Y$ satisfy condition \ref{condition *}. In this chapter, we will discuss some ways in which condition~\ref{condition *} can be relaxed in Theorem \ref{thm: main theorem}.

\section{Restrictions of simple factors}
Let $X$ be a smooth projective variety over $\mathbb{C}$ admitting a holomorphic symplectic form. By Theorem \ref{thm: symplectic variety is quotient}, we can write $X$ as a quotient $X \cong (\prod_{i=0}^k X_i)/G$, where $X_0$ is an abelian variety, and the $X_i$ for $i = 1, \dots, k$ are hyperkähler varieties. Using this quotient presentation, we computed $\llv(X; \mathbb{Q})$ in Theorem \ref{thm: summary LLV computation}, and we saw that each simple factor of $\llv(X; \mathbb{Q})$ either comes from the abelian part or the hyperkähler part. Moreover, according to Lemma \ref{lem: llvAb and llvHK independent}, this distinction is independent of the chosen quotient presentation.

Each simple factor of hyperkähler origin is of type $B_n$ or $D_n$ for some $n \geq 2$, while each simple factor of abelian origin is of type $A_n, C_n$ or $D_n$ for some $n \geq 1$ (see Theorem \ref{thm: summary LLV computation}). In Proposition \ref{prop: matching of LLV factors} we saw that, if $X$ satisfies condition~\ref{condition *}, the isomorphism class of a simple factor of $\llv(X; \mathbb{Q})$ determines whether it is of hyperkähler or of abelian origin. However, this is not true if $X$ does not satisfy condition~\ref{condition *} anymore. This is caused by the existence of some exceptional isomorphisms of Lie algebras (notably, $B_2 \cong C_2$ and $A_3 \cong D_3$). For example, a simple factor of type $C_2$ does not necessarily have to be of abelian origin, it can also be of hyperkähler origin.

This means that the proof of Proposition \ref{prop: matching of LLV factors} breaks down without condition~\ref{condition *}. However, we can slightly weaken this condition by requiring that (some of) the simple factors from condition \ref{condition *} do not occur both of hyperkähler and abelian origin at the same time.

\subsection{No simple factors of abelian origin of type $C_2, A_3$ or $D_4$}
The proof of Corollary \ref{cor: psi_i for HK} does not use condition \ref{condition *} in any way. This condition is only used in the proof of Theorem \ref{thm: main theorem} in the following two points:
\begin{enumerate}
    \item To show that if two simple factors $\mathfrak{g}_\alpha$ and $\mathfrak{h}_\alpha$ are isomorphic, then either they are both of abelian origin, or both of hyperkähler origin.
    \item In the proof of Corollary \ref{cor: psi_i for ab var}.
\end{enumerate}
Therefore, we can refine our main theorem to the following:
\begin{theorem}
    \label{thm: * restricted to HK}
    Let $X$ and $Y$ be smooth projective varieties over $\mathbb{C}$, both admitting a holomorphic symplectic form. Moreover, assume that all simple factors $\mathfrak{g}_\alpha$ of $\llv(X; \mathbb{Q})$ and $\llv(Y; \mathbb{Q})$ of type $C_2, A_3$ or $D_4$ are of hyperkähler origin. If $\Phi \colon D^b(X) \simeq D^b(Y)$ is a triangulated equivalence, then there is an isomorphism $\sH^\bullet(X; \mathbb{Q}) \cong \sH^\bullet(Y; \mathbb{Q})$ which preserves the Hodge bigrading.\qed
\end{theorem}

\subsection{No hyperkähler factors of types $C_2$ and $A_3$, and no factors of type $D_4$}
Similar to the previous section, we can wonder what happens if we instead restrict the possible simple factors of $\llv(X; \mathbb{Q})$ and $\llv(Y; \mathbb{Q})$ of hyperkähler origin, and loosen the restrictions on the simple factors of abelian origin. Here the answer is unfortunately not as nice as above, but we can still make some progress. For the three simple types of factors excluded by condition~\ref{condition *}, we now only exclude type $D_4$, and assume that the factors of type $C_2$ and $A_3$ are of abelian origin:

\begin{theorem}
    \label{thm: * restricted to abelian}
    Let $X$ and $Y$ be smooth projective varieties over $\mathbb{C}$, both admitting a holomorphic symplectic form. Moreover, assume that $\llv(X; \mathbb{Q})$ and $\llv(Y; \mathbb{Q})$ do not have any simple factors of type $D_4$, and that all simple factors of types $C_2$ and $A_3$ in $\llv(X; \mathbb{Q})$ and $\llv(Y; \mathbb{Q})$ are of abelian origin. Let $\Phi \colon D^b(X) \simeq D^b(Y)$ be a triangulated equivalence. Then there is an isomorphism $\sH^\bullet(X; \mathbb{Q}) \cong \sH^\bullet(Y; \mathbb{Q})$ which preserves the Hodge bigrading.
\end{theorem}

Before we prove this theorem, we need a small refinement of Proposition~\ref{prop: matching of LLV factors}. If $\mathfrak{g}_\alpha$ is a simple factor of $\llv(X; \mathbb{Q})$ of type $A_3$, there are three possible ways in which $\mathfrak{g}_\alpha$ could have been constructed (recall that $D_3 \cong A_3$):
\begin{enumerate}
    \item It is a simple factor of hyperkähler origin, so we have $\mathfrak{g}_\alpha \cong \mathfrak{so}(\widetilde \sH(X_i; \mathbb{Q})^{G_i})$ for some hyperkähler factor $X_i$ in the Beauville-Bogomolov covering of $X$, and $\widetilde \sH(X_i; \mathbb{Q})^{G_i}$ has dimension $6$.
    \item It is a simple factor of abelian origin, coming from an orthogonal involution, so we have $\mathfrak{g}_\alpha \cong \mathfrak{s}(\widetilde A, \widetilde \sigma)$ for some orthogonal involution $\widetilde \sigma$.
    \item It is a simple factor of abelian origin, coming from a unitary involution, so we have $\mathfrak{g}_\alpha \cong \mathfrak{s}(\widetilde A, \widetilde \sigma)$ for some unitary involution $\widetilde \sigma$.
\end{enumerate}
To prove Theorem \ref{thm: * restricted to abelian}, we need to be able to distinguish between the simple factors coming from the unitary and the orthogonal involutions.

\begin{lemma}
    \label{lem: Phi llv orthogonal vs unitary}
    Let $X$ and $Y$ be smooth projective varieties, admitting a holomorphic symplectic form, and let $\Phi \colon D^b(X) \to D^b(Y)$ be an equivalence of triangulated categories. Let $\mathfrak{g}_\alpha$ be a simple factor of $\llv(X; \mathbb{Q})$, and let $\mathfrak{h}_\alpha$ be the corresponding simple factor of $\llv(Y; \mathbb{Q})$ such that $\Phi^{\llv}$ restricts to an isomorphism $\Phi^{\llv}_\alpha \colon \mathfrak{g}_\alpha \to \mathfrak{h}_\alpha$. Assume that $\mathfrak{g}_\alpha$ and $\mathfrak{h}_\alpha$ both are of abelian origin, and both are of type $A_3$. Then $\mathfrak{g}_\alpha$ comes from an orthogonal involution if and only if $\mathfrak{h}_\alpha$ comes from an orthogonal involution.
\end{lemma}
\begin{proof}
    Suppose otherwise, so assume that $\mathfrak{g}_\alpha$ comes from a unitary involution, while $\mathfrak{h}_\alpha$ comes from an orthogonal involution. Let $k$ be the integer such that
    \begin{equation}
        \label{eq: g_alpha complexification}
        \mathfrak{g}_\alpha \otimes_\mathbb{Q} \mathbb{C} \cong \prod_{i=1}^k \mathfrak{sl}_4(\mathbb{C}).
    \end{equation}
    Then factor $\SH_\alpha(X; \mathbb{Q})$ of $\SH(X; \mathbb{Q})$ associated to $\mathfrak{g}_\alpha$ splits, after complexification, as
    \[ \SH_\alpha(X; \mathbb{Q}) \otimes_\mathbb{Q} \mathbb{C} \cong \bigotimes_{i=1}^k \SH_{\alpha, i}(X),\]
    where each component $\SH_{\alpha, i}(X)$ is an irreducible representation of the $i$-th component in decomposition \eqref{eq: g_alpha complexification}. Similarly, we can decompose 
    \begin{equation}
        \label{eq: h_alpha complexification}
        \mathfrak{h}_\alpha \otimes_\mathbb{Q} \mathbb{C} \cong \prod_{i=1}^k \mathfrak{sl}_4(\mathbb{C}),
    \end{equation}
    and there is an associated decomposition
    \begin{equation*}
        \SH_\alpha(Y; \mathbb{Q}) \otimes_\mathbb{Q} \mathbb{C} \cong \bigotimes_{i=1}^k \SH_{\alpha, i}(Y).
    \end{equation*}
    From the equivalence $\Phi$, we obtain an isomorphism of Verbitsky components $\Phi^{\SH}_\alpha \colon \SH_\alpha(X; \mathbb{Q}) \to \SH_\alpha(Y; \mathbb{Q})$, which is equivariant with respect to $\Phi^{\llv}_\alpha$. After complexification, $\Phi^{\llv}_\alpha$ splits as a direct product of isomorphisms between the simple factors in the decompositions~\eqref{eq: g_alpha complexification} and~\eqref{eq: h_alpha complexification}. After reindexing, we may assume that the $i$-th factors are mapped to each other, and we have equivariant isomorphisms $\Phi^{\SH}_{\alpha, i} \colon \SH_{\alpha, i}(X) \to \SH_{\alpha, i}(Y)$ for every $i$. We know from Proposition~\ref{prop: weight of Dn SH} that $\SH_{\alpha, i}(X)$ is an irreducible representation of highest weight $t\omega_2$ for some integer $t \geq 1$, where $\omega_2$ is the weight in the middle of the Dynkin diagram $A_3$. 

    However, if we now look at the $\SH_{\alpha, i}(Y)$, then Proposition~\ref{prop: weight of Dn SH} gives that for all~$i$ the representation $\SH_{\alpha, i}(Y)$ is irreducible of highest weight $s\omega'$ for some integer $s \geq 1$, where $\omega'$ is the weight corresponding to the even spinor representation. Since the weights $\omega_2$ and $\omega'$ can not be mapped to each other by the exceptional isomorphism $A_3 \cong D_3$ (see the diagrams below), we see that the maps $\Phi_{\alpha, i}^{\SH}$ cannot exists, which proves the lemma.

    \begin{equation*}
        \begin{tikzcd}[arrows=-]
            & & & & \bullet \ar[phantom, "\omega'"above=4pt]\\
            \bullet \arrow[r] & \bullet \ar[phantom, "\omega_2"above=4pt] \arrow[r] & \bullet & \bullet \arrow[ur] \arrow[dr]\\
            & & & & \bullet
        \end{tikzcd}
    \end{equation*}
    
\end{proof}

\begin{proof}[Proof of Theorem \ref{thm: * restricted to abelian}]
    The proof is analogous to the proof of Theorem \ref{thm: main theorem}, so we only mention the changes that need to be made. Observe that the outcome of Proposition~\ref{prop: matching of LLV factors} is still true, because we do not assume $X$ and $Y$ to have any simple factor in the LLV algebra of type $C_2, A_3$ or $D_4$ of hyperkähler origin, and Lemma~\ref{lem: Phi llv orthogonal vs unitary} guarantees that also for simple factors of type $A_3$ the type of the involution is preserved by the maps $\Phi^{\llv}_{\alpha}$. Hence we can still construct the maps $\psi_\alpha$ factorwise, separately for factors of hyperkähler or abelian origin. For the factors of hyperkähler origin, nothing changes from the proof of Theorem~\ref{thm: main theorem}. For the factors of abelian origin, we now also allow simple factors of type $C_2$, or of type $A_3$ coming from a unitary or orthogonal involution, compared to the original proof. Since Theorem~\ref{thm: isom lifts from skew to algebra} is also true for factors of type $C_2$ and $A_3$, the rest of the proof of Theorem~\ref{thm: main theorem} still goes through.
\end{proof}

\section{Geometry of symplectic varieties}
\label{sec: geometry of symplectic varieties}
In this section, we find some restrictions on the possible group actions $G \acts \prod_{i=0}^k X_i$, where $X_0$ is an even-dimensional abelian variety, the $X_i$ for $i \geq 1$ are hyperkähler varieties, and $G$ is a finite group acting freely on $\prod_i X_i$ and preserving a holomorphic symplectic form. We will use the understanding of such actions to prove Theorem~\ref{thm: main theorem} in dimension $4$ without requiring the varieties to satisfy condition~\ref{condition *}, see Section~\ref{sec: orlov dimension 4}.

Recall from Theorem \ref{thm: HKA automorphisms} that for every $g \in G$, the automorphism $g$ of $\prod_i X_i$ has the form
\[g(x_0, x_1, \dots, x_k) = (g_0(x_0), g_1(x_{\sigma^{-1}(1)}), \dots, g_k(x_{\sigma^{-1}(k)}))\]
where $x_i \in X_i$ for $i \geq 0$, where $\sigma \in S_k$ a permutation, and where $g_0$ is an automorphism of $X_0$ and the $g_i$ for $i \geq 1$ are isomorphisms $g_i \colon X_{\sigma^{-1}(i)} \to X_i$.

We start by looking at the hyperkähler part, where it turns out that a symplectic group action always has a fixed point. We first treat the case of a single hyperkähler variety. The following result is well-known, but I could not find a precise reference.
\begin{lemma}
    \label{lem: symplectic HK autom fixed points}
    Let $G$ be a finite group acting on a hyperkähler variety $X$ and preserving the holomorphic symplectic form. Then every $g \in G$ has at least one fixed point.
\end{lemma}
\begin{proof}
    Let $2d = \dim(X)$, and let $\omega \in \sH^0(X, \Omega^2_X)$ be a symplectic form (it is unique up to scaling). Since $G$ preserves the symplectic form, $g^*$ acts as the identity on $\sH^{2, 0}(X) = \sH^0(X, \Omega^2_{X})$. Since $g^*$ acts by ring automorphisms on $\sH^\bullet(X; \mathbb{Q})$, and
    \[ \sH^{p, 0}(X) = \begin{cases}
        \mathbb{C}\omega^{p/2} & \text{ if } 2 | p \text{ and } 0 \leq p \leq d/2\\
        0 & \text{ else}
    \end{cases}
    \]
    (see \cite[Proposition 3]{BeauvilleC1Zero}) it follows that $g^*$ acts as the identity on the ring $\bigoplus_{p \geq 0} \sH^{p,0}(X)$. Then $g \colon X \to X$ has a fixed point by the holomorphic Lefschetz theorem \cite[Page~426]{GriffithsHarris}.
\end{proof}

\begin{proposition}
    \label{prop: symplectic action product HK not free}
    Let $X = \prod_{i=1}^k X_i$ be a product of hyperkähler varieties, and let $G$ be a finite group acting on $X$. Assume that $X$ has a $G$-equivariant holomorphic symplectic form. Then the action of $G$ on $X$ is not free.
\end{proposition}
\begin{proof}
    Take a $g \in G$ not equal to the identity. We will show that $g$ has at least one fixed point, so we may replace $G$ by $\langle g \rangle$. As observed in Section~\ref{sec: reduction to quotients}, we have a homomorphism $G \to S_k$, and the set $I = \{1, \dots, k\}$ is a disjoint union of orbits of the induced $G$-action on $I$. This decomposes $X$ into the product $\prod_{[i] \in I/G} X_{[i]}$, where $X_{[i]} = \prod_{j \in [i]} X_j$, and the $G$-action preserves the factors $X_{[i]}$ of this decomposition. If we show that $G$ has a fixed point on each $X_{[i]}$, then the action of $G$ on $X$ also has a fixed point, so we may assume that $G$ acts transitively on $I$.

    Let $\ord(g)$ be the order of $g$. Since $G = \langle g \rangle$ acts transitively on $I = \{1, \dots, k\}$, we have $k | \ord(g)$. Furthermore, $g$ must act as a $k$-cycle on $\{1, \dots, k\}$. Therefore, if $\ord(g) > k$, then $g^k$ is not equal to the identity, but does act trivially on $I$. Therefore, $g^k$ is a product of symplectic automorphisms of the $X_i$. For every $i$, the induced automorphism of $X_i$ has a fixed point by Lemma \ref{lem: symplectic HK autom fixed points}, so $g^k$ has a fixed point, implying that the action of $G$ on $X$ is not free.

    Otherwise, $g$ has order exactly $k$. Since it acts as a $k$-cycle on $\{1, \dots, k\}$, we may assume, after reindexing, that $g$ induces isomorphisms 
    \[X_1 \isomto X_2 \isomto X_3 \dots \isomto X_k \isomto X_1,\] 
    whose composition $X_1 \isomto X_1$ is the identity. Take an element $x \in X_1$. Then the element $(x, gx, g^2x, \dots, g^{k-1}x) \in \prod_i X_i$ is fixed by $g$.
\end{proof}
Consequently, if $G$ acts freely on $\prod_{i=0}^k X_i$, the action of $G$ on $X_0$ cannot have any fixed points. We will now try to understand such actions a bit better. Note that we do not require the automorphisms of $X_0$ to respect the group structure, something which is often required of automorphisms of abelian varieties. Recall that an automorphism of an abelian variety is a group automorphism if and only if it sends $0$ to $0$ \cite[II.4, Corollary 1]{MumfordAbelian}. 

\begin{lemma}
    \label{lem: tilde g modification}
    Let $G$ be a group of automorphisms of an abelian variety $X_0$ over $\mathbb{C}$. For $g \in G$, define the group automorphism $\widetilde g$ of $X_0$ by
    \[ \widetilde g(x) := g(x) - g(0)\]
    Then the assignment $g \mapsto \widetilde g$ gives a group homomorphism $G \to \Aut(X_0)$. Moreover, for every $g \in G$ the two automorphisms $g^*, \widetilde g^*  \in \Aut(\sH^\bullet(X_0; \mathbb{Q}))$ coincide.
\end{lemma}
\begin{proof}
    Take $g, h \in G$, we need to show that $\widetilde g \widetilde h = \widetilde{gh}$ as automorphisms of $X_0$. Indeed, for $x \in X_0$ we have:
    \[\widetilde g \widetilde h(x) = \widetilde g (h(x) - h(0)) = \widetilde g(h(x)) - \widetilde g(h(0)) = gh(x) - gh(0),\]
    where in the second equality we use that $\widetilde g$ is a group homomorphism. Since $g$ and $\widetilde g$ differ by a translation, and translations act trivially on cohomology, we see that $g^*$ and $\widetilde g^*$ act the same way on $\sH^\bullet(X_0; \mathbb{Q})$.
\end{proof}

We can use the modified group action from above to obtain some useful information about the action of $G$ on $X_0$.

Let $X$ be an abelian variety, and let $G$ be a finite group of group automorphisms of $X$. Then $\sH^1(X; \mathbb{Q})$ is a representation of $G$. It turns out that the isotypic decomposition of $\sH^1(X; \mathbb{Q})$ lifts, up to isogeny, to an isotypic decomposition of $X$:
\begin{lemma}{\cite[Proposition 13.6.1]{BirkenhakeLangeCAV}}
    Let $G$ be a finite group acting by group automorphisms on a complex abelian variety $X$. Then for every $\rho \in \Irr_\mathbb{Q}(G)$ there is a $G$-invariant abelian subvariety\index{$X_\rho$} $X_\rho \subseteq X$, with the following properties:
    \begin{enumerate}
        \item The addition map induces a $G$-equivariant isogeny $\prod_{\rho \in \Irr_\mathbb{Q}(G)} X_\rho \to X$.
        \item Under the induced $G$-equivariant isomorphism 
        \begin{equation}
            \label{eq: H^1 of AV isotypic decomposition}
            \sH^1(X;\mathbb{Q}) \cong \bigoplus_{\rho \in \Irr_\mathbb{Q}(G)} \sH^1(X_\rho; \mathbb{Q}),
        \end{equation}
        for every $\rho \in \Irr_\mathbb{Q}(G)$, the $\rho$-isotypic component of $\sH^1(X; \mathbb{Q})$ is mapped to $\sH^1(X_\rho; \mathbb{Q})$.
    \end{enumerate}
\end{lemma}
Since not all parts in this lemma are proven in \cite[Proposition 13.6.1]{BirkenhakeLangeCAV}, we give the proof here.
\begin{proof}
    Consider the decomposition $1 = e_1 + \dots + e_r$ of $1 \in \mathbb{Q}[G]$ into central primitive orthogonal idempotents. Since $G$ acts on $X$, there is an induced map $\mathbb{Q}[G] \to \End_\mathbb{Q}(X) = \End(X) \otimes_\mathbb{Z}\mathbb{Q}$. Every $e_i$ corresponds to a unique $\rho \in \Irr_\mathbb{Q}(G)$. Choose an $m_i \in \mathbb{Z}$ with $m_ie_i \in \mathbb{Z}[G]$. Then $m_ie_iX$ is an abelian subvariety of $X$. By the first part of \cite[Section 13.6]{BirkenhakeLangeCAV}, this construction does not depend on the chosen $m_i$, and \cite[Proposition 13.6.1]{BirkenhakeLangeCAV} shows that the addition map $\prod_\rho X_\rho \to X$ induces an isogeny. Since $G$ acts by group automorphisms, this map is $G$-equivariant. 

    For the second statement, it suffices to show that for every $\rho \in \Irr_\mathbb{Q}(G)$, the representation of $G$ on $\sH^1(X_\rho; \mathbb{Q})$ is $\rho$-isotypic, because the isomorphism \eqref{eq: H^1 of AV isotypic decomposition} is $G$-equivariant.
    
    Below the proof of \cite[Proposition 13.6.1]{BirkenhakeLangeCAV}, it is explained that we can describe $X_\rho$ as follows: for the induced action of $G$ on the tangent space $T_0X$, decompose $T_0X = \bigoplus_{\rho \in \Irr_\mathbb{Q}(G)} U_\rho$ into isotypic components. Let $\Lambda \subseteq T_0 X$ be the lattice with $X \cong T_0X/\Lambda$, then $X_\rho \cong U_\rho / (\Lambda \cap U_\rho)$, and this implies that the representation of $G$ on $\sH^1(X_\rho; \mathbb{Q})$ is $\rho$-isotypic.
\end{proof}

\begin{definition}
    The variety $\prod_\rho X_\rho$ is called the \textbf{isotypic decomposition} of $X$, and the abelian varieties $X_\rho$ are called the \textbf{isotypic components} of $X$.
\end{definition}

If $X$ is equipped with a $G$-invariant holomorphic symplectic form, then all the $X_\rho$ turn out to be symplectic, and in particular, even-dimensional. In particular, this implies that the dimension of every isotypic component of $\sH^1(X; \mathbb{Q})$ is divisible by $4$. The proof of this fact goes in a few steps.

Recall that for an irreducible representation $\rho \in \Irr_\mathbb{C}(G)$ of a finite group $G$, we defined $[\rho] = \{\rho, \rho^*\}$, this set contains either $1$ or $2$ irreducible representations, depending on whether $\rho$ is self-dual or not. For a $G$-representation $V$, we defined the $[\rho]$-isotypic component by 
\begin{equation*}
    V_{[\rho]} := \bigoplus_{\sigma \in [\rho]} V_\sigma = \begin{cases}
    V_{\rho} & \text{ if } \rho \cong \rho^*\\
    V_{\rho} \oplus V_{\rho^*} & \text{ if } \rho \ncong \rho^*.
\end{cases}
\end{equation*}

\begin{lemma}
    \label{lem: symplectic form isotypic decomposition}
    Let $G$ be a finite group, let $V$ be a finite-dimensional representation of $G$ over $\mathbb{C}$ and let $\sigma$ be a $G$-equivariant symplectic form on $V$. Then for every $\rho \in \Irr_\mathbb{C}(G)$, the restriction of $\sigma$ to the $[\rho]$-isotypic component $V_{[\rho]}$ is also symplectic.
\end{lemma}
\begin{proof}
    We have to show that for every $\rho$ the restriction of $\sigma$ to $V_{[\rho]}$ is nondegenerate. Since $\sigma$ is nondegenerate and $V$ is finite-dimensional, it provides an isomorphism $V \to V^*$. This isomorphism is $G$-equivariant, because $\sigma$ is $G$-equivariant. Therefore the isomorphism $V \to V^*$ induced by $\sigma$ is the direct sum of isomorphisms $V_\rho \to (V_{\rho^*})^*$. It follows that $\sigma|_{V_{[\rho]}}$ is nondegenerate for every $\rho \in \Irr_\mathbb{C}(G)$.
\end{proof}

\begin{proposition}
\label{prop: symplectic isotypic AV}
    Let $G$ be a finite group acting freely on an even-dimensional complex abelian variety $X$ and preserving a symplectic form $\sigma$ on $X$. Then every factor $X_\rho$ of the isotypic decomposition $p \colon \prod_{\rho \in \Irr_\mathbb{Q}(G)} X_\rho \to X$ is a symplectic abelian variety.
\end{proposition}
\begin{proof}
    Let $\sigma$ be a $G$-invariant symplectic form on $X$. Since the pullback of a symplectic form under an étale map is again a symplectic form, the pullback $p^*\sigma$ gives a $G$-invariant symplectic form on $\prod_\rho X_\rho$. As a $G$-representation, the tangent space of $\prod_\rho X_\rho$ at the neutral element $0 \in \prod_\rho X_\rho$ splits as $\bigoplus_\rho T_0 X_\rho$. Let $R \subseteq \Irr_\mathbb{C}(G)$ be the set of irreducible representations occuring in $\rho \otimes_\mathbb{Q} \mathbb{C}$. Then $\rho_i \in R$ if and only if $\rho_i^* \in R$, and we have $T_0X_\rho \cong \bigoplus_{\rho_i \in R} (T_0X)_{\rho_i}$ (where $(T_0X)_{\rho_i}$ is the $\rho_i$-isotypic component of $T_0X$). Then Lemma \ref{lem: symplectic form isotypic decomposition} implies that the restriction of $\sigma$ to each $T_0 X_\rho$ is again a symplectic form. It follows from this that $p^*\sigma$ is a sum of symplectic forms on the $X_\rho$, so each $X_\rho$ is symplectic.
\end{proof}

\begin{lemma}
    \label{lem: free automorphism AV}
    Let $X_0$ be a complex abelian variety, and let $g \colon X \to X$ be an automorphism of finite order. If $g$ has no fixed points, then the automorphism $g^*$ of $\sH^1(X_0; \mathbb{Q})$ has the eigenvalue $1$ with multiplicity at least $1$.
\end{lemma}
\begin{proof}
    Let $d = \dim X$. Since $g$ has no fixed points, the Lefschetz trace formula \cite[Page 421]{GriffithsHarris} gives
    \begin{equation}
        \label{eq: Lefschetz trace formula AV}
        \sum_{i=0}^{2d} (-1)^i \Tr(g^* \colon \sH^i(X; \mathbb{Q}) \to \sH^i(X; \mathbb{Q})) = 0.
    \end{equation}
    Let $\lambda_1, \dots, \lambda_{2d}$ be the eigenvalues of the automorphism $g^*$ of $\sH^1(X; \mathbb{Q})$. Since $\sH^\bullet(X; \mathbb{Q}) = \extp^\bullet \sH^1(X; \mathbb{Q})$, we see that the alternating sum of the traces in Equation \eqref{eq: Lefschetz trace formula AV} is equal to
    \[ \prod_{i=1}^{2d} (1 - \lambda_i).\]
    Since this must be zero, at least one of the $\lambda_i$ must be $1$, which finishes the proof.
\end{proof}

If we assume that $X_0$ is symplectic and $g$ is a symplectic automorphism, we can use the isotypic decomposition to show something even stronger:
\begin{corollary}
    \label{cor: symplectic free automorphism AV}
    Let $X_0$ be an even-dimensional complex abelian variety, and let $g \colon X \to X$ be an automorphism of finite order, preserving a holomorphic symplectic form on $X_0$. If $g$ has no fixed points, then the automorphism $g^*$ of $\sH^1(X_0; \mathbb{Q})$ has the eigenvalue $1$ with multiplicity at least $4$.
\end{corollary}
\begin{proof}
    Let $G = \langle g \rangle \subseteq \Aut(X_0)$ be the finite group generated by $g$. Since $g^* \colon \sH^1(X_0; \mathbb{Q}) \to \sH^1(X_0; \mathbb{Q})$ has eigenvalue $1$ with multiplicity at least $1$, there is a $v \in \sH^1(X_0; \mathbb{C})$ which is invariant under $g$, and hence under the action of $G = \langle g \rangle$. It follows that the trivial representation occurs with multiplicity at least $1$ in $\sH^1(X_0; \mathbb{C})$, and therefore also in $\sH^1(X_0; \mathbb{Q})$. 

    Now consider the modified action of $G$ on $X_0$ as in Lemma \ref{lem: tilde g modification}. The modified automorphism $\widetilde g^*$ has the same eigenvalues as $g^*$ on $\sH^1(X_0; \mathbb{Q})$. It follows that the isotypic decomposition $\prod X_\rho \to X_0$ for this modified $G$-action on $X_0$ contains a nontrivial factor $X_{\triv}$ for the trivial representation. Proposition \ref{prop: symplectic isotypic AV} implies that $X_{\triv}$ is symplectic, so we see that $X_{\triv}$ has dimension at least $2$, and this implies that the trivial representation of $\langle g \rangle$ occurs with multiplicity at least $4$ in $\sH^1(X_0; \mathbb{Q})$.
\end{proof}

\section{Orlov's conjecture for symplectic fourfolds}
\label{sec: orlov dimension 4}
In this section, we will prove Conjecture \ref{conj: orlov cohomology} for smooth projective varieties of dimension $4$ admitting a holomorphic symplectic form, even if they do not satisfy condition \ref{condition *}. First, we will use the results of the previous section to understand the different types of symplectic varieties that can exist in dimension $4$, and then we will use LLV algebras and the Verbitsky component to prove the conjecture.

\begin{lemma}
    \label{lem: rep with half of eigenvalues 1}
    Take $n \in \mathbb{N}$, let $G$ be a finite group and let $V$ be a $G$-representation of dimension $2n$ over $\mathbb{Q}$. Write $\rho \colon G \to \GL(V)$ for this representation, and assume that for every $g \in G$, the eigenvalue $1$ occurs with multiplicity at least $n$ in $\rho(g)$. Then the trivial representation of $G$ occurs in the decomposition of $\rho$ into irreducible representations.
\end{lemma}
\begin{proof}
    Let $\chi$ be the character of $\rho$, and let $\mathbf{1}$ be the character of the trivial $G$-representation. For every $g \in G$, let $\lambda_1(g), \dots, \lambda_{2n}(g)$ be the complex eigenvalues of $\rho(g)$. Since $G$ is finite, all the $\lambda_i(g)$ lie on the unit circle. Moreover, for every $g \in G$, at least half of the $\lambda_i(g)$ are equal to $1$. The multiplicity of the trivial representation in $V$ is equal to $\langle \mathbf{1}, \chi\rangle$. Then:
    \begin{equation}
        \label{eq: character inner product}
        \langle \chi, \mathbf{1}\rangle = \frac{1}{\# G} \sum_{g \in G} \Tr(\rho(g)) = \frac{1}{\# G} \sum_{g \in G} \sum_{i=1}^{2n} \lambda_i(g).
    \end{equation}
    The sum of the right hand side is a sum of $2n \cdot \# G$ roots of unity. For the identity element $e \in G$, we know that $\lambda_i(e) = 1$ for all $i$. It follows that more than half of the $\lambda_i(g)$ are equal to $1$, and this implies that the sum in Equation \eqref{eq: character inner product} cannot be zero. 
\end{proof}

\begin{proposition}
    \label{prop: symplectic fourfold types}
    Let $X = (\prod_i X_i)/G$ be as in Setup \ref{stp: symplectic var}, and assume that $X$ has dimension $4$ and that $X$ has a holomorphic symplectic form. Then $X$ is of one of the following forms:
    \begin{enumerate}
        \item We have $X = X_1$ is a $4$-dimensional hyperkähler variety.
        \item We have $X = X_1 \times X_2$ is a product of two K3 surfaces.
        \item We have $X = (X_0 \times X_1)/G$, where $X_0$ is an abelian surface, $X_1$ is a K3 surface and the action of $G$ on $\sH^1(X_0; \mathbb{Q})$ is trivial.
        \item We have $X = X_0/G$, where $X_0$ is a $4$-dimensional abelian variety, and $X_0$ has the isotypic decomposition $X_{\triv} \times X_{\rho} \to X$, where $\rho \in \Irr_\mathbb{Q}(G)$ is nontrivial and $X_{\triv}$ and $X_\rho$ both have dimension $2$.
        \item $X$ is a $4$-dimensional abelian variety (in particular, we may take $X = X_0$).
    \end{enumerate}
\end{proposition}
\begin{proof}
    The product $\prod_i X_i$ has one of the following forms:
    \begin{enumerate}
        \item A $4$-dimensional hyperkähler variety.
        \item A product of two K3 surfaces.
        \item A product of a K3 surface and an abelian surface.
        \item A $4$-dimensional abelian variety.
    \end{enumerate}
    If we are in one of the first two cases, then Proposition \ref{prop: symplectic action product HK not free} implies that $\prod X_i$ does not admit any free symplectic group action, so we must have $X = \prod X_i$ in these two cases.

    Now suppose that we are in the third case, so we have $X = (X_0 \times X_1)/G$, where $X_0$ is an abelian surface, $X_1$ is a K3 surface and $G \acts X_0 \times X_1$ is free. By Lemma~\ref{lem: symplectic HK autom fixed points}, the action of $G$ on $X_1$ is not free, so the action of $G$ on $X_0$ must be free. Since $\dim \sH^1(X_0; \mathbb{Q}) = 4$, Corollary \ref{cor: symplectic free automorphism AV} implies that every $g \in G$ acts trivially on $\sH^1(X_0; \mathbb{Q})$.

    Finally, assume that $\prod_i X_i = X_0$ is an abelian fourfold. Consider the isotypic decomposition $\prod_\rho X_\rho$ of the modified $G$-action as in Lemma~\ref{lem: tilde g modification}. In Proposition~\ref{prop: symplectic isotypic AV} we saw that each $X_\rho$ is symplectic, and Corollary~\ref{cor: symplectic free automorphism AV} and Lemma~\ref{lem: rep with half of eigenvalues 1} imply that the trivial representation has multiplicity at least $4$ in the representation of $G$ in $\sH^1(X_0; \mathbb{Q})$. It follows that $\dim X_{\triv} \geq 2$.

    If $\dim X_{\triv} = 2$, then the isotypic decomposition of $X_0$ has the form $X_{\triv} \times X_\rho$, where $X_\rho$ also has dimension $2$, so we are in the fourth case.

    Otherwise, $X_{\triv}$ has dimension $4$, so $G$ acts trivially on $\sH^1(X_0; \mathbb{Q})$. Then \cite[Proposition 1.2.1]{BirkenhakeLangeCAV} implies that the elements of $G$ act on $X_0$ by translations, so the quotient $X = X_{\triv}/G$ is again an abelian variety.
\end{proof}

Let $V$ be a polarizable $\mathbb{Q}$-Hodge structure of weight $2$. Define $\widetilde V = \mathbb{Q}[0] \oplus V \oplus \mathbb{Q}[4]$, where $V$ sits in degree $2$. Define a Hodge structure on $\widetilde V$ where the two copies of $\mathbb{Q}$ are algebraic. Note that $\widetilde V$ looks like the Hodge structure on the even part of the cohomology of a smooth projective surface.
\begin{lemma}
    \label{lem: twodimensional even hodge structure}
    Let $V$ and $W$ be two polarizable $\mathbb{Q}$-Hodge structures of weight $2$. Suppose that there is a linear isomorphism $f \colon \widetilde V \to \widetilde W$ which preserves the Hodge grading (but not necessarily the normal grading). Then there also is an isomorphism $\widetilde V \to \widetilde W$ which preserves the Hodge bigrading.
\end{lemma}
\begin{proof}
    This follows directly from the fact that the category of polarizable Hodge structures is semisimple.
\end{proof}

Similarly, if $V$ is a polarizable $\mathbb{Q}$-Hodge structure of weight $1$, we can consider $V \oplus V(-1)$, where $V(-1)$ is the Tate twist of $V$, and is a Hodge structure of weight $3$. We give $V \oplus V(-1)$ the grading where $V$ lives in degree $1$ and $V(-1)$ in degree $3$. Note that $V \oplus V(-1)$ looks like the odd part of the cohomology of a smooth projective surface.
\begin{lemma}
    \label{lem: twodimensional odd hodge structure}
    Let $V$ and $W$ be two polarizable $\mathbb{Q}$-Hodge structures of weight $1$, and let $f \colon V \oplus V(-1) \to W \oplus W(-1)$ be a linear isomorphism which preserves the Hodge grading (but not necessarily the normal grading). Then there also is an isomorphism $V \oplus V(-1) \to W \oplus W(-1)$ which preserves the Hodge bigrading.
\end{lemma}
\begin{proof}
    This follows directly from the fact that the category of polarizable Hodge structures is semisimple.
\end{proof}

We can now prove Orlov's conjecture for symplectic fourfolds:
\begin{theorem}
    \label{thm: main theorem dimension 4}
    Let $X$ and $Y$ be two smooth projective varieties of dimension $4$ over $\mathbb{C}$, both equipped with a holomorphic symplectic form, and assume that there is a triangulated equivalence $\Phi \colon D^b(X) \simeq D^b(Y)$. Then there is an isomorphism $\sH^\bullet(X; \mathbb{Q}) \cong \sH^\bullet(Y; \mathbb{Q})$ preserving the Hodge bigrading.
\end{theorem}
\begin{proof}
    By Theorem \ref{thm: symplectic variety is quotient}, we can write $X = (\prod_i X_i)/G$ and $Y = (\prod_j Y_j)/H$ as quotients of a product of abelian and hyperkähler varieties. Moreover, these products are of one of the five possible types as given in Proposition \ref{prop: symplectic fourfold types}. 

    First assume that $X$ is of the first or the fifth type (so $X$ is either a hyperkähler or an abelian variety). Then $Y$ is also a hyperkähler or abelian variety respectively by \cite[Theorem 0.4]{HuybrechtsNieperWisskirchen}, and then we already know that the conjecture is true (either by \cite{TaelmanDerivedEquivalences} or \cite{OrlovAbelian} as seen in the introduction).

    If $X$ is of the second type, then $\llv(X; \mathbb{Q}) \cong \mathfrak{so}(\widetilde \sH(X_1; \mathbb{Q})) \times \mathfrak{so}(\widetilde \sH(X_2; \mathbb{Q}))$. Since $\dim \sH^2(X_i; \mathbb{Q}) = 22$ for $i = 1, 2$, the variety $X$ satisfies condition \ref{condition *}. Hence the same holds for $Y$, so the conjecture is true by Theorem \ref{thm: main theorem}.

    Lastly, assume that $X$ is of the third or fourth type. Since these are the only two types in the classification where $\dim \sH^1(X; \mathbb{Q}) = 4$, the variety $Y$ must also be of the third or fourth type by \cite[Corollary B]{PopaSchnell}. In both these cases, the Lie algebra $\llv(X; \mathbb{Q})$ is a product of two simple factors by Theorem \ref{thm: summary LLV computation}. However, in both these cases, the LLV algebra has a simple factor of type $D_4$, so we cannot apply Theorem~\ref{thm: main theorem}. 

    If $X$ is of the form $(X_0 \times X_1)/G$ where $G$ acts trivially on $\sH^1(X_0; \mathbb{Q})$, and therefore trivially on all of $\sH^\bullet(X_0; \mathbb{Q})$, then the cohomology of $X$ factorizes as
    \[ \sH^\bullet(X; \mathbb{Q}) \cong \sH^\bullet(X_0; \mathbb{Q}) \otimes \sH^\bullet(X_1; \mathbb{Q})^G.\]
    Let $\mathfrak{g}_1 = \llv(X_0; \mathbb{Q})$ and $\mathfrak{g}_2 = \llv_G(X_1; \mathbb{Q})$, then $\llv(X; \mathbb{Q}) \cong \mathfrak{g}_1 \times \mathfrak{g}_2$ by Theorem~\ref{thm: summary LLV computation}. Moreover, note that $\sH^\bullet(X_0; \mathbb{Q}) = \sH^{\ev}(X_0; \mathbb{Q}) \oplus \sH^{\odd}(X_0; \mathbb{Q})$ is a direct sum of two absolutely irreducible $\mathfrak{g}_1$-representations, and $\sH^\bullet(X_1; \mathbb{Q})^G = \SH_G(X_1; \mathbb{Q})$ is an absolutely irreducible $\mathfrak{g}_2$-representation.
    
    If $X$ is of the form $X_0/G$, where $X_0$ has the isotypic decomposition $X_{\triv} \times X_\rho \to X$, where $X_{\triv}$ and $X_\rho$ both have dimension $2$, then we also have a decomposition
    \[\sH^\bullet(X; \mathbb{Q}) \cong \sH^\bullet(X_{\triv}; \mathbb{Q}) \otimes \sH^\bullet(X_\rho; \mathbb{Q})^G.\]
    If we now let $\mathfrak{g}_1 = \llv(X_{\triv}; \mathbb{Q})$ and $\mathfrak{g}_2 = \llv_G(X_\rho; \mathbb{Q})$, then $\llv(X; \mathbb{Q}) \cong \mathfrak{g}_1 \times \mathfrak{g}_2$. As above, we see that $\sH^\bullet(X_0; \mathbb{Q}) = \sH^{\ev}(X_0; \mathbb{Q}) \oplus \sH^{\odd}(X_0; \mathbb{Q})$ is a direct sum of two absolutely irreducible $\mathfrak{g}_1$-representations, and $\sH^\bullet(X_\rho; \mathbb{Q})^G = \SH_G(X_\rho; \mathbb{Q})$ is an absolutely irreducible $\mathfrak{g}_2$-representation.

    So, now suppose that $X$ and $Y$ are two symplectic fourfolds of the third or fourth type, and $\Phi \colon D^b(X) \simeq D^b(Y)$ is an equivalence. The induced isomorphism $\Phi^{\sH} \colon \sH^\bullet(X; \mathbb{Q}) \to \sH^\bullet(Y; \mathbb{Q})$ is a direct sum of isomorphisms 
    \[\Phi^{\sH, \ev} \colon \sH^{\ev}(X; \mathbb{Q}) \to \sH^{\ev}(Y; \mathbb{Q}) \text{ and } \Phi^{\sH, \odd} \colon \sH^{\odd}(X; \mathbb{Q}) \to \sH^{\odd}(Y; \mathbb{Q}).\]
    We saw above that the LLV algebras of $X$ and $Y$ both are a product of two simple factors, say $\llv(X; \mathbb{Q}) \cong \mathfrak{g}_1 \times \mathfrak{g}_2$ and $\llv(Y; \mathbb{Q}) \cong \mathfrak{h}_1 \times \mathfrak{h}_2$. Furthermore, $\sH^{\ev}(X; \mathbb{Q})$ and $\sH^{\odd}(X; \mathbb{Q})$ both are a tensor product of absolutely irreducible representations of $\mathfrak{g}_1$ and $\mathfrak{g}_2$, and similarly $\sH^{\ev}(Y; \mathbb{Q})$ and $\sH^{\odd}(Y; \mathbb{Q})$ both are a tensor product of absolutely irreducible representations of $\mathfrak{h}_1$ and $\mathfrak{h}_2$. For simplicity, write $\sH^{\ev}(X; \mathbb{Q}) = V_1 \otimes V_2$ and $\sH^{\ev}(Y; \mathbb{Q}) = W_1 \otimes W_2$ for these decompositions, and similarly $\sH^{\odd}(X; \mathbb{Q}) = V_1'  \otimes V_2'$ and $\sH^{\odd}(Y; \mathbb{Q}) = W_1' \otimes W_2'$. 
    
    Then Theorem \ref{thm: split Lie alg reps} implies that $\Phi^{\sH, \ev}$ and $\Phi^{\sH, \odd}$ are tensor products of isomorphisms between these smaller factors. These smaller factors $V_i, V_i', W_i, W_i'$ all either have the form $\mathbb{Q}[0] \oplus U \oplus \mathbb{Q}[4]$ for a pure weight $2$ Hodge structure $U$, or the form $U \oplus U(-1)$ for a pure weight $1$ Hodge structure $U$. Then Lemmas~\ref{lem: twodimensional even hodge structure} and \ref{lem: twodimensional odd hodge structure} imply that we have isomorphisms $V_i \cong W_i$ and $V_i' \cong W'_i$ preserving the Hodge bigrading for $i = 1,2$. This implies that we also have an isomorphism $\sH^\bullet(X; \mathbb{Q}) \cong \sH^\bullet(Y; \mathbb{Q})$ which preserves the Hodge bigrading.
\end{proof}

\begin{remark}
    In particular, the above theorem implies the derived invariance of the Hodge numbers for symplectic fourfolds. We should note here that this derived invariance of the Hodge numbers also follows from a result by Abuaf \cite{AbuafHomologicalUnits}. Abuaf \cite[Theorem 1.3(4)]{AbuafHomologicalUnits} shows that for fourfolds, the Hodge numbers $h^{2,0}$ are derived invariants. One can then use the symmetries of the Hodge diamond of a holomorphic symplectic variety, together with the derived invariance of $h^{1,0}$ \cite[Corollary B]{PopaSchnell}, to show that all Hodge numbers are derived invariants.
\end{remark}

\begin{remark}
    The proof also deals with the possibility of a derived equivalence $D^b(X) \simeq D^b(Y)$ where $X$ is of the third type in the classification of Proposition~\ref{prop: symplectic fourfold types} and $Y$ is of the fourth type. I do not know of an example of such an equivalence, but it would be interesting to know whether such a derived equivalence exists. At least on the level of LLV algebras it seems possible. Indeed, take $X = (S \times A)/G$ as in the introduction to this thesis (see just below Theorem~\ref{thm: main theorem}). To construct $Y$, we consider an abelian surface $A$, a $2$-torsion point $a \in A[2]$, and $H = \langle g \rangle$ the cyclic group of order $4$, acting on $A \times A$ by letting $g$ map $(x,y) \in A \times A$ to $(y, x + a)$. Then $Y := (A \times A)/H$ is of the fourth type in Proposition~\ref{prop: symplectic fourfold types}, and the LLV algebras of $X$ and $Y$ both are a product of two simple factors of type $D_4$.
\end{remark}

In the above proof, for the case where $X = X_1 \times X_2$ is a product of two K3 surfaces, we also could have used the same argument as in the third and fourth case. In fact, for this case, we do not have to restrict to a product of two K3 surfaces; with the same argument, one can also prove:
\begin{proposition}
    Let $X_1, \dots, X_n$ and $Y_1, \dots, Y_n$ be K3 surfaces, and let $X = \prod_i X_i$ and $Y = \prod_i Y_i$. If there is a triangulated equivalence $\Phi \colon D^b(X) \simeq D^b(Y)$, then there is an isomorphism $\sH^\bullet(X; \mathbb{Q}) \cong \sH^\bullet(Y; \mathbb{Q})$ which preserves the Hodge bigrading. \qed
\end{proposition}
Note that this also follows directly from Theorem~\ref{thm: main theorem}, since condition~\ref{condition *} is satisfied for $X$ and $Y$. 

Let $X = \prod_i X_i$ and $Y = \prod_j Y_j$ be products of K3 surfaces. Since the cohomology of $X$ and $Y$ is equal to the Verbitsky components, Theorem~\ref{thm: split Lie alg reps} implies that there is a permutation $\sigma \in S_n$ such that the isomorphism $\Phi^{\sH} \colon \sH^\bullet(X; \mathbb{Q}) \to \sH^\bullet(Y; \mathbb{Q})$ is the tensor product of isomorphisms $\Phi^{\sH}_i \colon \sH^\bullet(X_i; \mathbb{Q}) \to \sH^\bullet(Y_{\sigma(i)}; \mathbb{Q})$. It seems to me an interesting question whether the maps $\Phi^{\sH}_i$ can also be defined between the cohomology with $\mathbb{Z}$-coefficients, as is the case if $n = 1$. If this would indeed be possible, and if one could also prove that the obtained isomorphism $\widetilde \sH(X_i; \mathbb{Z}) \cong \widetilde \sH(Y_{\sigma(i)}; \mathbb{Z})$ preserves the lattice structure (i.e. is an isomorphism of the Mukai lattices), then by the derived Torelli theorem \cite{OrlovK3}, one obtains equivalences $D^b(X_i) \simeq D^b(Y_{\sigma(i)})$, thereby partially answering a question by Huybrechts and Nieper-Wisskirchen \cite[Question~0.2]{HuybrechtsNieperWisskirchen}.

\section{Summary}
We now use the results of this chapter to state Theorem \ref{thm: main theorem} in all forms proven so far:
\begin{theorem}
    \label{thm: main theorem generalized}
    Let $X$ and $Y$ be two smooth projective varieties over $\mathbb{C}$, both admitting a holomorphic symplectic form, and let $\Phi \colon D^b(X) \to D^b(Y)$ be a triangulated equivalence.

    Assume that at least one of the following is true:
    \begin{itemize}
        \item The dimension of $X$ is $4$.
        \item The following three properties hold:
        \begin{itemize}
            \item All simple factors of $\llv(X; \mathbb{Q})$ and $\llv(Y; \mathbb{Q})$ of type $D_4$ are of hyperkähler origin.
            \item All simple factors of $\llv(X; \mathbb{Q})$ and $\llv(Y; \mathbb{Q})$ of type $C_2$ are either all of hyperkähler origin or all of abelian origin.
            \item All simple factors of $\llv(X; \mathbb{Q})$ and $\llv(Y; \mathbb{Q})$ of type $A_3$ are either all of hyperkähler origin or all of abelian origin.
        \end{itemize}
    \end{itemize}
    Then there exists an isomorphism $\sH^\bullet(X; \mathbb{Q}) \cong \sH^\bullet(Y; \mathbb{Q})$ which preserves the Hodge bigrading. 
\end{theorem}
\begin{proof}
    If $X$ has dimension $4$, this is Theorem~\ref{thm: main theorem dimension 4}. Otherwise, using the same ideas as in the proofs of Theorems~\ref{thm: * restricted to HK} and~\ref{thm: * restricted to abelian}, we see that the result of Proposition~\ref{prop: matching of LLV factors} still holds, and the result follows similarly to the proof of Theorem~\ref{thm: * restricted to abelian}.
\end{proof}

%% file: chapters/chapter9.tex
An Enriques surface is a surface which is not simply connected, and whose simply connected covering is a $K3$ surface. The canonical bundle $\omega_S$ of an Enriques surface $S$ is not trivial, but $\omega_S^{\otimes 2}$ is trivial. It is known that two Enriques surfaces that are derived equivalent are actually isomorphic \cite[Proposition~6.1]{BridgelandMaciocia2}. There are several ways in which one can generalize the notion of an Enriques surface to higher dimensions, see for instance \cite{OguisoSchroer} and \cite{BossiereNieperWisskirchenSarti}. In this chapter, we use the definition from \cite{OguisoSchroer}:
\begin{definition}
    An \textbf{Enriques variety} is a smooth connected projective variety $X$ over $\mathbb{C}$, which is not simply connected, and whose simply connected cover $\widetilde X \to X$ is a hyperkähler variety.
\end{definition}
If $X$ is an Enriques variety, then its fundamental group is cyclic, say of order $n$, and the canonical bundle of $X$ is torsion of order $n$, see \cite[Lemma 2.3]{OguisoSchroer}. The hyperkähler covering $\widetilde X \to X$ is then the canonical cover of $X$.

In this chapter, we prove Conjecture \ref{conj: orlov cohomology} for Enriques varieties:

\begin{theorem}
    \label{thm: orlov conj for enriques}
    Let $X$ and $Y$ be two Enriques varieties. Let $\Phi \colon D^b(X) \to D^b(Y)$ be an equivalence of triangulated categories. Then for every $i$ there is an isomorphism $\sH^i(X; \mathbb{Q}) \cong \sH^i(Y; \mathbb{Q})$ that preserves the Hodge structure.
\end{theorem}
Let $n$ be the order of the canonical bundle $\omega_X$, and let $G = \mathbb{Z}/n\mathbb{Z}$. Then $\omega_Y$ also has order $n$ (see \cite[Proposition~4.1]{Huybrechts}), and we can write $X \cong \widetilde X / G$ and $Y \cong \widetilde Y/G$, where $\widetilde X$ and $\widetilde Y$ are the canonical covers of $X$ and $Y$ (in particular, $\widetilde X$ and $\widetilde Y$ are hyperkähler varieties). To prove the theorem, we will use a result by Bridgeland-Maciocia \cite{BridgelandMaciocia} to lift the equivalence $\Phi$ to an equivalence $\widetilde \Phi \colon D^b(\widetilde X) \to D^b(\widetilde Y)$. By \cite[Theorem~D]{TaelmanDerivedEquivalences}, we obtain an isomorphism $\sH^\bullet(\widetilde X; \mathbb{Q}) \cong \sH^\bullet(\widetilde Y; \mathbb{Q})$ which preserves the Hodge bigrading. We will show that there in fact exists such an isomorphism which is also equivariant with respect to the $G$-actions. By taking $G$-invariants, we will obtain the desired isomorphism $\sH^\bullet(X; \mathbb{Q}) \cong \sH^\bullet(Y; \mathbb{Q})$.

Since we will use some details of the proof of \cite[Theorem~D]{TaelmanDerivedEquivalences}, we now summarize the proof strategy used there. The reader can find more details in \cite{TaelmanDerivedEquivalences}. A derived equivalence $\Phi \colon D^b(\widetilde X) \to D^b(\widetilde Y)$ induces an isomorphism of LLV algebras $\Phi^{\llv} \colon \llv(\widetilde X; \mathbb{Q}) \to \llv(\widetilde Y; \mathbb{Q})$ and an isomorphism of Verbitsky components $\Phi^{\SH} \colon \SH(\widetilde X; \mathbb{Q}) \to \SH(\widetilde Y; \mathbb{Q})$ which is equivariant with respect to $\Phi^{\llv}$. There is an isomorphism $\llv(\widetilde X; \mathbb{Q}) \cong \mathfrak{so}(\widetilde \sH(\widetilde X; \mathbb{Q}))$ \cite{LooijengaLunts, VerbitskyPhDThesis}, and an injective morphism of $\llv(\widetilde X; \mathbb{Q})$-representations $\Psi_X \colon \SH(\widetilde X; \mathbb{Q}) \to \Sym^d \widetilde \sH(\widetilde X; \mathbb{Q})$, where $d$ is defined by $2d = \dim \widetilde X$ \cite[Proposition~3.5]{TaelmanDerivedEquivalences} (and similarly we have a $\Psi_Y$ for $\widetilde Y$).

Taelman then uses $\Psi_X$ and $\Psi_Y$ to show that there is a similitude $\varphi \colon \widetilde \sH(\widetilde X; \mathbb{Q}) \to \widetilde \sH(\widetilde Y; \mathbb{Q})$ which induces $\Phi^{\llv}$ and $\Phi^{\SH}$. He then uses Witt cancellation to construct an automorphism $\psi$ of $\widetilde \sH(\widetilde Y; \mathbb{Q})$ with the property that $\psi \circ \varphi$ respects the Hodge bigrading. To finish the proof, he uses $\psi$ to construct an automorphism $\Psi^{\sH}$ of $\sH^\bullet(\widetilde Y; \mathbb{Q})$ such that $\Psi^{\sH} \circ \Phi^{\sH}$ respects the Hodge bigrading.

The actions of $G$ on $\widetilde X$ and $\widetilde Y$ induce actions on the cohomology, the LLV algebra, and so forth. In this chapter, we will show that the above constructions of $\Psi_X, \varphi$ and $\psi$ are $G$-equivariant, and use this to show that the induced automorphism $\Psi^{\sH}$ is $G$-equivariant.

\section{Equivariant Witt cancellation}
In this section, we prove an equivariant version of the Witt cancellation by adapting the standard proof as given in \cite[Section 3]{WittCancellation}. This proof only works when cancelling a subspace on which the group acts trivially, which suffices for our application.

\begin{definition}
    Let $G$ be a finite group, let $V$ be a $G$-representation over a field $k$ of characteristic unequal to $2$, and let $q$ be a quadratic form on $V$. We say that $(V,q)$ is a $G$-\textbf{equivariant quadratic space} if for all $g \in G$ and $v \in V$ we have $q(gv) = q(v)$.
\end{definition}

 \begin{theorem}
    \label{thm: weak equiv Witt cancellation}
     Let $(V,q)$ be a $G$-equivariant nondegenerate quadratic space over a field of characteristic different from $2$, and let $\{e_1, \dots, e_n\}$ and $\{f_1, \dots, f_n\}$ be two $q$-orthogonal bases of $V$. Let $m \leq n$ and assume that for all $1 \leq i \leq m$ and $g \in G$ we have $ge_i = e_i, gf_i = f_i$ and $q(e_i) = q(f_i)$. Then there exists a $G$-equivariant isometry between the $G$-equivariant quadratic spaces $\langle e_{m+1}, \dots, e_n \rangle$ and $\langle f_{m+1}, \dots, f_n \rangle$. 
 \end{theorem}

\begin{proof}
    By working inductively, we may assume that $m = 1$. 
    Without the $G$-action, this is \cite[Theorem 3.1]{WittCancellation}. We are finished once we show that their constructed isometry $\tau_{e_1 - f_1}$ is $G$-equivariant. If we write $b$ for the bilinear form associated to $q$, the reflection map $\tau$ is defined by
    $$\tau_u(z) = z - \frac{2 b(u,z)}{q(u)}u,$$
    where $u, z \in V$. We are done if we show that $\tau_u$ is $G$-invariant whenever $gu = u$ for all $g \in G$. This is a direct computation:
    \begin{align*}
        \tau_u(gz) &= gz - \frac{2b(u, gz)}{q(u)}u \\
        &= gz - \frac{2b(g^{-1}u, z)}{q(u)}u\\
        &= g \cdot \left(z - \frac{2b(u, z)}{q(u)}u \right) = g\tau_u(z),
    \end{align*}    
    where we use the fact that $q$ is $G$-equivariant in the second equality, and the third equality follows from the fact that the vector $u$ is $G$-invariant.
\end{proof}

\section{Equivariant constructions}
We will now show that all relevant constructions in the proof of \cite[Theorem~D]{TaelmanDerivedEquivalences} are $G$-equivariant. Examples of this are the element $\varphi$ constructed in \cite[Proposition~4.9]{TaelmanDerivedEquivalences}, and the projection $\GSpin(V, q) \to \SO(V, q)$ for a $G$-equivariant quadratic space $V$. We start with some linear algebra.

\begin{lemma}
    \label{lem: V to Sym^d V}
    Let $V$ be a vector space over a field $k$ of characteristic~$0$ and let $d \geq 1$ be an integer. Denote by $f \colon V \to \Sym^d(V)$ the map $v \mapsto v^d$. Then for $v, w \in V$, we have $f(v) = f(w)$ if and only if there is a $\lambda \in \mu_d(k)$ with $w = \lambda v$.
\end{lemma}
\begin{proof}
    It is clear that $f(v) = f(\lambda v)$ for $\lambda \in \mu_d(k)$. Conversely, suppose that $f(v) = f(w)$. Given a basis $\langle e_1, \dots, e_n \rangle$ of $V$, we know that $\Sym^d(V)$ has a basis given by the set of unordered products 
    $$\{e_{i_1}\cdots e_{i_n} : i_1, \dots, i_n \in \{1, \dots, n\}\}.$$
    If $w$ is not a scalar multiple of $v$, then $V$ has a basis with the first two elements $v$ and $w$, from which we can conclude that $f(v) \neq f(w)$. Otherwise, there is a $\lambda \in k^*$ with $v = \lambda w$. Then $f(\lambda v) = \lambda^d f(v)$, so the assumption that $f(v) = f(w)$ gives $\lambda^d = 1$, which finishes the proof. 
\end{proof}

Recall that the Mukai lattice of a hyperkähler variety $\widetilde X$ is given by $\widetilde \sH(\widetilde X; \mathbb{Q}) = \mathbb{Q}\alpha \oplus \sH^2(\widetilde X; \mathbb{Q}) \oplus \mathbb{Q}\beta$, and that $\llv(\widetilde X; \mathbb{Q}) \cong \mathfrak{so}(\widetilde \sH(\widetilde X; \mathbb{Q}))$. If a group $G$ acts on $\widetilde X$, then $G$ acts on $\widetilde \sH(\widetilde X; \mathbb{Q})$ by taking the induced action of $G$ on $\sH^2(\widetilde X; \mathbb{Q})$ and letting $G$ act trivially on $\alpha$ and $\beta$. Let $d$ be the integer such that $2d = \dim \widetilde X$. By \cite[Proposition~3.5]{TaelmanDerivedEquivalences}, there is a unique injective $\mathfrak{so}(\widetilde \sH(\widetilde X; \mathbb{Q}))$-equivariant homomorphism $\Psi \colon \SH(\widetilde X; \mathbb{Q}) \to \Sym^d \widetilde \sH(\widetilde X; \mathbb{Q})$ which sends $1$ to $\alpha^d/d!$.

\begin{lemma}
    Let $G$ be a finite group acting on a hyperkähler variety $\widetilde X$. Then the map $\Psi$ from above is $G$-equivariant.
\end{lemma}
\begin{proof}
    Since $\SH(\widetilde X; \mathbb{Q})$ is generated by all elements in degree $2$, it suffices to show the following two things for all $g \in G$:
    \begin{enumerate}
        \item $g \Psi(1) = \Psi(g \cdot 1)$.
        \item If $\lambda \in \sH^2(\widetilde X; \mathbb{Q})$ and $\omega \in \SH(\widetilde X; \mathbb{Q})$ with $\Psi(g \omega) = g\Psi(\omega)$, then also $\Psi(g (\lambda \omega)) = g\Psi(\lambda \omega)$.
    \end{enumerate}
    Since $g \cdot 1 = 1$ and $g \alpha = \alpha$, the first statement is clear.

    For the second statement, recall that $\lambda$ gives an element $e_\lambda \in \mathfrak{so}(\widetilde \sH(\widetilde X; \mathbb{Q}))$ with the property that $\lambda \omega = e_\lambda(\omega)$. By also using that $\Psi$ is $\mathfrak{so}(\widetilde \sH(\widetilde X; \mathbb{Q}))$-equivariant and that $g (\lambda \omega) = (g\lambda)(g\omega)$, we see:
    \[ \Psi(g(\lambda\omega)) = \Psi(e_{g\lambda} (g\omega)) = e_{g\lambda}(\Psi(g\omega)).\]
    By then using that $\Psi(g\omega) = g\Psi(\omega)$ and that $e_{g\lambda}(x) = g e_\lambda(g^{-1}x)$ for all $x \in \Sym^d \widetilde \sH(\widetilde X; \mathbb{Q})$, we see that $e_{g\lambda}(\Psi(g\omega)) = ge_\lambda(\Psi(\omega))$, and this is equal to $g\Psi(\lambda \omega)$, which finishes the proof.
\end{proof}

In the following, we will denote for a hyperkähler variety $\widetilde X$ the morphism $\SH(\widetilde X; \mathbb{Q}) \to \Sym^d\widetilde \sH(\widetilde X; \mathbb{Q})$ by $\Psi_X$, and we define $\Psi_Y$ similarly. Let $\Phi \colon D^b(\widetilde X) \to D^b(\widetilde Y)$ be an equivalence of triangulated categories, then there is an induced isomorphism of Verbitsky components $\Phi^{\SH} \colon \SH(\widetilde X; \mathbb{Q}) \to \SH(\widetilde Y; \mathbb{Q})$. By \cite[Proposition~4.9]{TaelmanDerivedEquivalences}, there exist a Hodge similitude $\varphi \colon \widetilde \sH(\widetilde X; \mathbb{Q}) \to \widetilde \sH(\widetilde Y; \mathbb{Q})$ and a scalar $\lambda \in \mathbb{Q}^\times$ such that the diagram 
\begin{equation}
\label{eq: SH psi Sym commutative square}
    \begin{tikzcd}
        \SH(\widetilde X; \mathbb{Q}) \arrow[d, "\Psi_X"] \arrow[rr, "\widetilde \Phi^{\SH}"] && \SH(\widetilde Y; \mathbb{Q}) \arrow[d, "\Psi_Y"]\\
        \Sym^d \widetilde \sH(\widetilde X; \mathbb{Q}) \arrow[rr, "\lambda \Sym^d(\varphi)"] && \Sym^d \widetilde \sH(\widetilde Y; \mathbb{Q})
    \end{tikzcd}
\end{equation}
commutes.

\begin{lemma}
    \label{lem: phi is equivariant}
    Let $\widetilde X$ and $\widetilde Y$ be hyperkähler varieties of complex dimension $2d$, equipped with an action of the cyclic group $G = \mathbb{Z}/n\mathbb{Z}$, and let $\widetilde \Phi \colon D^b(\widetilde X) \to D^b(\widetilde Y)$ be a $G$-equivariant equivalence of triangulated categories. Assume that at least one of $d$ and $n$ is odd, and assume that the induced isomorphism $\widetilde \Phi^{\SH}$ of Verbitsky components is $G$-equivariant. Then the induced map $\varphi \colon \widetilde \sH(\widetilde X; \mathbb{Q}) \to \widetilde \sH(\widetilde Y; \mathbb{Q})$ which makes diagram~\eqref{eq: SH psi Sym commutative square} commute is also $G$-equivariant.
\end{lemma}
\begin{proof}
    Write $\alpha_X \in \widetilde \sH(\widetilde X; \mathbb{Q})$ for the generator of the Mukai lattice in degree $0$, and take $g \in G$. We will first show that $g \varphi(\alpha_X) = \varphi(\alpha_X)$. We know from \cite[Proposition 3.5]{TaelmanDerivedEquivalences} that $\Psi_X(1) = \alpha_X^d/d!$. Since $\alpha$ is $G$-invariant and $\Psi_X$ is $G$-equivariant, we know that $g \Psi_X(1) = \Psi_X(1)$. By using the commutativity of diagram \eqref{eq: SH psi Sym commutative square}, and using that the maps $\Psi_X, \widetilde \Phi^{\SH}$ and $\Psi_Y$ are $G$-equivariant, it follows that $ g\cdot (\lambda \varphi(\alpha_X)^d) = \lambda \varphi(\alpha_X)^d$. If we now use that $\lambda \neq 0$ and that $g \varphi(\alpha_X)^d = (g \varphi(\alpha_X))^d$, we see from Lemma \ref{lem: V to Sym^d V} that there is a $\xi \in \mu_d(\mathbb{Q})$ with $g\varphi(\alpha_X) = \xi\varphi(\alpha_X)$. Additionally, we know that $g^n = 1$, so $\xi^n = 1$. Since $\mu_d(\mathbb{Q})$ is either $\{1\}$ or $\{\pm 1\}$, depending on whether $d$ is odd or even, and $n$ and $d$ can not both be even, we see that $\xi = 1$, so $g\varphi(\alpha_X) = \varphi(\alpha_X)$. 

    Now take an $x \in \sH^2(\widetilde X; \mathbb{Q})$, we will show that $\varphi(gx) = g\varphi(x)$. Consider the operator $e_x \in \llv(\widetilde X; \mathbb{Q})$, which acts as multiplication by $x$ on $\sH^\bullet(\widetilde X; \mathbb{Q})$. Under the isomorphism $\llv(\widetilde X; \mathbb{Q}) \cong \mathfrak{so}(\widetilde \sH(\widetilde X; \mathbb{Q}))$, the operator $e_x$ acts on $\widetilde \sH(\widetilde X; \mathbb{Q})$, and it sends $\alpha_X$ to $x$, see \cite[Section 3.1]{TaelmanDerivedEquivalences}.

    Since $x = e_x(1)$ in $\SH(\widetilde X; \mathbb{Q})$ and $\Psi_X$ is $\llv(\widetilde X; \mathbb{Q})$-equivariant, it follows that 
    \begin{equation} 
        \label{eq: Psi_X(x)}
        \Psi_X(x) = e_x(\alpha_X^d/d!) = x\alpha_X^{d-1}/(d-1)!.
    \end{equation}
    Now take $gx \in \SH(\widetilde X; \mathbb{Q})$ in diagram \eqref{eq: SH psi Sym commutative square}. Since $\widetilde \Phi^{\SH}$ and $\Psi_Y$ are $G$-equivariant, we know that $\Psi_Y(\widetilde \Phi^{\SH}(gx)) = g\Psi_Y(\widetilde \Phi^{\SH}(x))$. Since the diagram is commutative, this implies
    \[\lambda \Sym^d(\varphi)(\Psi_X(gx)) = g \cdot (\lambda \Sym^d(\varphi)(\Psi_X(x))).\]
    If we now use equation \eqref{eq: Psi_X(x)} and the fact that $g\alpha_X = \alpha_X$, we see that $\varphi(gx) = g\varphi(x)$.

    It remains to show that $\varphi(g\beta) = g\varphi(\beta)$ for all $g \in G$. For this, we use the same strategy as above. Choose an $x \in \sH^2(\widetilde X; \mathbb{Q})$ with $b(x, x) \neq 0$. By \cite[Section~3.1]{TaelmanDerivedEquivalences}, we have $e_x(e_x(\alpha_X)) = b(x,x)\beta_X$. It follows that
    \[\Psi_X(x^2) = e_x(e_x(\Psi_X(1))) = b(x,x)\beta_X\alpha_X^{d-1}/(d-1)! + x^2\alpha_X^{d-2}/(d-2)!.\]
    By now applying this formula to $(gx)^2$, we see just as above that $\varphi(g\beta_X) = g\varphi(\beta_X)$, and this finishes the proof.
\end{proof}

Schröer and Oguiso proved that the parity constraints in Lemma~\ref{lem: phi is equivariant} are always satisfied. The proof uses the behaviour of the Euler characteristic $\chi(\mathcal{O}_X)$ under quotients.
\begin{lemma}[{\cite[Proposition 2.4]{OguisoSchroer}}]
    \label{lem: oguiso schroer parity}
    Let $X$ be a hyperkähler variety of dimension $2d$ and assume that the finite group $G = \mathbb{Z}/n\mathbb{Z}$ acts freely on $X$. Then at least one of $d$ and $n$ is odd.\qed
\end{lemma}

Let $G$ be a finite group and let $(V, q)$ be a $G$-equivariant quadratic space. There is an induced action of $G$ on $\SO(V, q)$, where $g \in G$ sends $f \in \SO(V, q)$ to the map which sends $v \in V$ to
\[(gf)(v) := gf(g^{-1}v).\]
An element $\varphi \in \SO(V, q)$ is a $G$-equivariant map $\varphi \colon V \to V$ if and only if $g\varphi = \varphi$.

\begin{proposition}
    \label{prop: GSpin to SO is G-equivariant}
    Let $G$ be a finite group and let $(V, q)$ be a $G$-equivariant quadratic space over $\mathbb{Q}$. Then there is an induced action of $G$ on $\GSpin(V, q)$ such that the map $\pi \colon \GSpin(V, q) \to \SO(V, q)$ is $G$-equivariant.
\end{proposition}
\begin{proof}
    We will prove this using Clifford algebras. Let $C(V, q)$ be the Clifford algebra of $V$; it is defined as $C(V, q) = T(V) / I(V, q)$, where
    \[T(V) = \bigoplus_{n \geq 0} V^{\otimes n}\]
    and $I(V, q) \subseteq T(V)$ is the ideal generated by the elements $v \otimes v - q(v)$ (see \cite[Chapter~20]{Fulton-Harris} for more details about the Clifford algebra). The group $G$ acts on $T(V)$ in the natural way, and the fact that $q$ is $G$-equivariant implies that $I(V, q)$ is $G$-invariant. It follows that there is an induced action of $G$ on $C(V, q)$.

    Recall that $C(V, q)$ is equipped with an involution $x \mapsto x^*$ defined by linearly extending
    \[(v_1 \otimes \dots \otimes v_r)^* := (-1)^rv_r \otimes \dots \otimes v_1\]
    for $v_1, \dots, v_r \in V$. Then $\GSpin(V, q)$ is defined by
    \[\GSpin(V, q) = \{x \in C(V, q)^{\ev} : xVx^* \subseteq V\}.\]
    It follows directly from the definition that $(gx)^* = gx^*$ for all $g \in G$ and $x \in C(V, q)$. From this, we see that if $x \in \GSpin(V, q)$ and $v \in V$, then 
    \[ (gx)v(gx)^* = (gx)(gg^{-1}v)(gx)^* = g(x(g^{-1}v)x^*),\]
    The fact that $x \in \GSpin(V, q)$ then implies that also $x(g^{-1}v)x^* \in V$, and hence also $g(x(g^{-1}v)x^*) \in V$. We also see from this formula that the map $\GSpin(V, q) \to \SO(V, q)$ is indeed $G$-equivariant.
\end{proof}

The action of the LLV algebra $\llv(\widetilde Y; \mathbb{Q})$ on $\sH^\bullet(\widetilde Y; \mathbb{Q})$ integrates to a representation $\rho_0 \colon \uSpin(\widetilde \sH^\bullet(\widetilde Y; \mathbb{Q})) \to \uGL(\sH^\bullet(\widetilde Y; \mathbb{Q}))$. This representation can be extended to a representation
\[\rho \colon \uGSpin(\widetilde \sH(\widetilde Y; \mathbb{Q})) \to \uGL(\sH^\bullet(\widetilde Y; \mathbb{Q})), \]
see \cite[Lemma 5.1]{TaelmanDerivedEquivalences}. The $G$-action on $\sH^\bullet(\widetilde Y; \mathbb{Q})$ induces an action of $G$ on $\GL(\sH^\bullet(\widetilde Y; \mathbb{Q}))$ by conjugation.

\begin{lemma}
    \label{lem: Spin rep on H(tilde X) equivariant}
    The representation $\rho_0 \colon \uSpin(\widetilde \sH(\widetilde Y; \mathbb{Q})) \to \uGL(\sH^\bullet(\widetilde Y; \mathbb{Q}))$ is $G$-equivariant.
\end{lemma}
\begin{proof}
    We want to show that $g\circ \rho_0 = \rho_0 \circ g$ for all $g \in G$. Since $\uSpin(\widetilde \sH(\widetilde Y; \mathbb{Q}))$ is connected, Lemma~\ref{lem: morphism on algebraic groups determined by Lie alg} implies that it suffices to show that the induced map $d\rho_0 \colon \mathfrak{so}(\widetilde \sH(\widetilde Y; \mathbb{Q})) \to \mathfrak{gl}(\sH^\bullet(\widetilde Y; \mathbb{Q}))$ is $G$-equivariant.

    First, we will prove that $d\rho_0$ is $G$-equivariant in degree $2$. Recall from \cite[Section~3.1]{TaelmanDerivedEquivalences} that every degree $2$ element of $\mathfrak{so}(\widetilde \sH(\widetilde Y; \mathbb{Q}))$ is of the form $e_x$ for some $x \in \sH^2(\widetilde Y; \mathbb{Q})$, where the map $e_x \colon \widetilde \sH(\widetilde Y; \mathbb{Q}) \to \widetilde \sH(\widetilde Y; \mathbb{Q})$ sends $\alpha$ to $x$, sends $y \in \sH^2(\widetilde Y; \mathbb{Q})$ to $b(x,y)\beta$ (with $b$ the Beauville-Bogomolov-Fujiki form) and sends $\beta$ to $0$. We then see that $ge_x = e_{gx}$ for all $g \in G$. Furthermore, $d\rho_0$ sends $e_x$ to the endomorphism $e'_x$ of $\sH^\bullet(\widetilde Y; \mathbb{Q}))$, defined by $e'_x(\omega) = x\cdot \omega$ for $\omega \in \sH^\bullet(\widetilde Y; \mathbb{Q})$. Then also $ge'_x = e'_{gx}$ for all $g \in G$. Therefore, we have $d\rho_0(ge_x) = d\rho_0(e_{gx}) = e'_{gx} = ge'_x$, so $d\rho_0$ is indeed $G$-equivariant on the elements in degree $2$. 

    There is a Zariski-dense subset of the degree $-2$ elements of $\mathfrak{so}(\widetilde \sH(\widetilde Y; \mathbb{Q}))$, all of which are of the form $f_x$ for some $x \in \sH^2(\widetilde Y; \mathbb{Q})$ with the Hard Lefschetz property, where $f_x$ is the unique element of $\mathfrak{so}(\widetilde \sH(\widetilde Y; \mathbb{Q}))$ with the property that $\langle e_x, h, f_x\rangle$ is an $\mathfrak{sl}_2$-triple. Moreover, for every $x \in \sH^2(\widetilde Y; \mathbb{Q})$ with the hard Lefschetz property, there is also a unique $f'_x \in \mathfrak{gl}(\sH^\bullet(\widetilde Y; \mathbb{Q}))$ with the property that $\langle e'_x, h, f'_x \rangle$ is an $\mathfrak{sl}_2$-triple (this follows from the proof of \cite[Proposition 2.1]{TaelmanDerivedEquivalences}).

    If we now take a $g \in G$ and an $x \in \sH^2(\widetilde Y; \mathbb{Q})$ with the hard Lefschetz property, then we see, since $e'_{gx} = ge'_x$, that both $\langle e'_{gx}, h, f'_{gx}\rangle$ and $\langle e'_{gx}, h, gf'_x\rangle$ are $\mathfrak{sl}_2$-triples. The Jacobson-Morozov theorem \cite[VIII.11.2, Proposition~2]{BourbakiLie7and8} then implies that $f'_{gx} = gf'_x$. Since $d\rho_0(f_x) = f'_x$, we now conclude just as above that $d\rho_0(gf_x) = gd\rho_0(f_x)$. Since the equality $d\rho_0(gf_x) = gd\rho_0(f_x)$ holds for a Zariski-dense collection of $f_x \in \mathfrak{so}(\widetilde \sH(\widetilde Y; \mathbb{Q}))_{-2}$, it must in fact hold for all $f \in \mathfrak{so}(\widetilde \sH(\widetilde Y; \mathbb{Q}))_{-2}$.

    Since $\mathfrak{so}(\widetilde \sH(\widetilde Y; \mathbb{Q}))$ is generated by its elements in degrees $2$ and $-2$, we conclude that $d\rho_0$ is $G$-equivariant on all of $\mathfrak{so}(\widetilde \sH(\widetilde Y; \mathbb{Q}))$, which finishes the proof.
\end{proof}

In \cite[Lemma~5.1]{TaelmanDerivedEquivalences}, Taelman constructs $\rho$ from $\rho_0$ by using the first short exact sequence from Lemma~\ref{lem: HK Spin short exact sequences}, and writing down a representation of $\mathbb{G}_{m, \mathbb{Q}}$ which commutes with the image of $\rho_0$.

\begin{lemma}
    \label{lem: rep GSpin on cohomology G-equivariant}
    The representation $\rho \colon \uGSpin(\widetilde \sH(\widetilde Y; \mathbb{Q})) \to \uGL(\sH^\bullet(\widetilde Y; \mathbb{Q}))$ is $G$-equivariant.
\end{lemma}
\begin{proof}
    The action of $G$ on $\sH^\bullet(\widetilde Y; \mathbb{Q})$ respects the grading, so the action written down by Taelman in \cite[Lemma~5.1]{TaelmanDerivedEquivalences} is $G$-equivariant. Lemma~\ref{lem: Spin rep on H(tilde X) equivariant} says that the representation $\rho_0 \colon \uSpin(\widetilde \sH(\widetilde Y; \mathbb{Q})) \to \uGL(\sH^\bullet(\widetilde Y; \mathbb{Q}))$ is also $G$-equivariant, so it follows that $\rho$ is $G$-equivariant.
\end{proof}

\section{Proof of Theorem \ref{thm: orlov conj for enriques}}
The proof of Theorem \ref{thm: orlov conj for enriques} follows the proof of \cite[Theorem D]{TaelmanDerivedEquivalences}, where we now check that all relevant constructions are $G$-equivariant.
\begin{proof}[{Proof of Theorem \ref{thm: orlov conj for enriques}}]
    Let $n$ be the order of $\omega_X$, and let $G = \mathbb{Z}/n\mathbb{Z}$. By \cite[Proposition~4.1]{Huybrechts}, the canonical bundle of $Y$ also has order $n$. Let $\widetilde Y$ be the canonical cover of $Y$. A theorem of Bridgeland and Maciocia \cite[Theorem~4.5]{BridgelandMaciocia} implies that there is an automorphism $\mu \colon G \to G$ and an equivalence
    \[\widetilde \Phi \colon D^b(\widetilde X) \to D^b(\widetilde Y),\]
    with the property that $g^* \circ \widetilde \Phi = \widetilde \Phi \circ \mu(g)^*$ for all $g \in G$. By defining a new action of $G$ on $\widetilde X$ by precomposing the map $G \to \Aut(\widetilde X)$ with $\mu^{-1}$, we may assume $\mu$ to be the identity. Then the induced map on cohomology $\widetilde \Phi^{\sH}$ is $G$-equivariant due to the functioriality of the construction of $\widetilde \Phi^{\sH}$ in terms of $\widetilde \Phi$.
    
    There are actions of $G$ on $\llv(\widetilde X; \mathbb{Q})$ and $\llv(\widetilde Y; \mathbb{Q})$ induced by conjugation, and the isomorphism $\widetilde \Phi^{\llv} \colon \llv(\widetilde X; \mathbb{Q}) \to \llv(\widetilde Y; \mathbb{Q})$ is also $G$-equivariant.

    By Lemmas \ref{lem: phi is equivariant} and \ref{lem: oguiso schroer parity}, there exist a scalar $\lambda \in \mathbb{Q}^\times$ and a $G$-equivariant Hodge similitude $\varphi \colon \widetilde \sH(\widetilde X; \mathbb{Q}) \to \widetilde \sH(\widetilde Y; \mathbb{Q})$ such that diagram \eqref{eq: SH psi Sym commutative square} commutes. We now claim that there is a $G$-equivariant Hodge isometry $\psi \in \SO(\widetilde \sH(\widetilde Y; \mathbb{Q}))$ with the property that $\psi\varphi \colon \sH^\bullet(\widetilde X; \mathbb{Q}) \to \sH^\bullet(\widetilde Y; \mathbb{Q})$ preserves the Hodge bigrading. 

    To see this, we use the same strategy as in \cite[Theorem 5.3]{TaelmanDerivedEquivalences}, where we now have to use our equivariant version of Witt cancellation. For completeness, we spell out some details here.

    Write $\widetilde \sH(\widetilde Y; \mathbb{Q}) = N \oplus T$, where $N$ is the algebraic part and $T$ the transcendental part. The subspaces $N$ and $T$ are orthogonal for the Mukai pairing, since the Mukai pairing respects the Hodge structure. Recall that $\widetilde \sH(\widetilde Y; \mathbb{Q}) = \mathbb{Q}\alpha_Y \oplus \sH^2(\widetilde Y; \mathbb{Q}) \oplus \mathbb{Q}\beta_Y$, and we have $\alpha_Y, \beta_Y \in N$. Since $\varphi \colon \widetilde \sH(\widetilde X; \mathbb{Q}) \to \widetilde \sH(\widetilde Y; \mathbb{Q})$ is an equivariant Hodge similitude, it maps the $G$-invariant hyperbolic plane $\mathbb{Q}\alpha_X \oplus \mathbb{Q}\beta_X$ to a $G$-invariant hyperbolic plane $V \subseteq N$. Then equivariant Witt cancellation (Theorem \ref{thm: weak equiv Witt cancellation}) implies that there are a $G$-equivariant isometry $\psi_N$ of $N$ and scalars $s,t \in \mathbb{Q}^\times$ such that $\psi_N \circ \varphi$ sends $\alpha_X$ to $s\alpha_Y$ and $\beta_X$ to $t\beta_Y$.
    
    If $\det \psi_N = -1$, we may compose $\psi_N$ with the isometry of $N$ which sends $\alpha_Y + \beta_Y$ to $\alpha_Y + \beta_Y$ and $\alpha_Y - \beta_Y$ to $-(\alpha_Y - \beta_Y)$, and is the identity on the rest of $N$. This map then also has determinant $-1$, so without loss of generality we may assume that $\psi_N \in \SO(N)$. By extending this by the identity of $T$, we obtain a $G$-equivariant Hodge isometry $\psi \in \SO(\widetilde \sH(Y; \mathbb{Q}))$ with the property that $\psi \circ \varphi \colon \widetilde \sH(\widetilde X; \mathbb{Q}) \to \widetilde \sH(\widetilde Y; \mathbb{Q})$ is a Hodge similitude which also preserves the grading (recall that $\alpha_X, \alpha_Y$ have degree $0$, and $\beta_X, \beta_Y$ have degree $4$).

    Proposition \ref{prop: GSpin to SO is G-equivariant} implies that the map $\GSpin(\widetilde \sH(\widetilde Y; \mathbb{Q})) \to \SO(\widetilde \sH(\widetilde Y; \mathbb{Q}))$ is $G$-equivariant. Since this map is surjective (this follows from Hilbert 90 and the second short exact sequence in Lemma \ref{lem: HK Spin short exact sequences}), we obtain a $G$-equivariant lift $\psi' \in \GSpin(\widetilde \sH(\widetilde Y; \mathbb{Q}))$. By \cite[Lemma 5.1]{TaelmanDerivedEquivalences}, the group $\GSpin(\widetilde \sH(\widetilde Y; \mathbb{Q}))$ acts on $\sH^\bullet(Y; \mathbb{Q})$, and this representation is $G$-equivariant by Lemma \ref{lem: rep GSpin on cohomology G-equivariant}. If we denote the induced action of $\psi'$ on $\sH^\bullet(\widetilde Y; \mathbb{Q})$ by $\Psi^{\sH}$, it follows that $\Psi^{\sH}$ is $G$-equivariant. Therefore, the composition
    \[ \Psi^{\sH} \circ \widetilde \Phi^{\sH} \colon \sH^\bullet(\widetilde X; \mathbb{Q}) \to \sH^\bullet(\widetilde Y; \mathbb{Q})\]
    is also $G$-equivariant. On the level of LLV algebras, it sends $h_{\widetilde X}$ to $h_{\widetilde Y}$ and $h'_{\widetilde X}$ to $h'_{\widetilde Y}$, so it preserves the Hodge bigrading. By now taking $G$-invariants, we obtain an isomorphism 
    \[\sH^\bullet(X; \mathbb{Q}) \to \sH^\bullet(Y; \mathbb{Q})\]
    which preserves the Hodge bigrading, and this finishes the proof.
\end{proof}

\section{Epilogue} Now, at the end of this chapter\footnote{In fact, the end of this thesis, except for the summary, the references and the index of notation.}, one glaring question remains: suppose $X$ and $Y$ are smooth projective varieties with torsion canonical bundle, and whose canonical covers $\widetilde X \to X$ and $\widetilde Y \to Y$ are varieties admitting a holomorphic symplectic form (possibly satisfying condition \ref{condition *}). If there is an equivalence $\Phi \colon D^b(\widetilde X) \simeq D^b(\widetilde Y)$, can we prove that there is an isomorphism $\sH^\bullet(X; \mathbb{Q}) \cong \sH^\bullet(Y; \mathbb{Q})$ which preserves the Hodge bigrading? I don't know the answer to this question\footnote{Given the techniques developed so far, I guess that the answer is yes; at least under suitable hypotheses such as $\widetilde X$ satisfying condition \ref{condition *}}, and leave it as a direction for further research.

%% file: chapters/summary.tex
\section*{Varieties admitting a holomorphic symplectic form: LLV algebras and derived equivalences}

When doing geometry in two or three dimensions, we can often use our intuition, and can study geometric objects by drawing or visualizing them. As soon as we start to do geometry in higher dimensions, we can no longer easily visualize the objects that we study, and it becomes harder to study their properties.

Modern mathematics has found solutions for this. Many ways have been found to get a handle on higher-dimensional geometric objects. This is often done by using simpler mathematical objects, which reflect some geometric information of the geometric object we wish to study. We call such easier objects \textit{invariants}. A well-known example is the cohomology, which roughly counts the number of holes of a shape.

This dissertation is in the area of algebraic geometry, where the geometric objects which we study are so-called \textit{algebraic varieties}. Algebraic varieties are given as the solution sets to polynomial equations, and have a rich structure. This rich structure leads to many interesting invariants of algebraic varieties. These invariants carry information about the specific algebraic structure of the variety. For example, the cohomology of an algebraic variety over $\mathbb{C}$ carries an extra structure, the \textit{Hodge structure}, which contains a lot of information about the variety.

Another important invariant of an algebraic variety is the \textit{derived category}. This is a category which is built in a few steps from the category of coherent sheaves on the variety.

It is possible that two different varieties have the same derived category. We can then wonder how similar these two varieties are. For example, two varieties with the same derived category have the same dimension. It is still unknown in what sense the derived category determines the cohomology and Hodge structure. An important conjecture in this direction is the following conjecture by Orlov:

\begin{conjecture*}[Orlov]
    Let $X$ and $Y$ be two smooth projective varieties over $\mathbb{C}$, and assume that there is an equivalence $D^b(X) \simeq D^b(Y)$ of triangulated categories. Then there is an isomorphism $\sH^\bullet(X; \mathbb{Q}) \cong \sH^\bullet(Y; \mathbb{Q})$ which preserves the grading and Hodge structure.
\end{conjecture*}

In this dissertation, I study a specific kind of varieties, namely smooth projective varieties over $\mathbb{C}$ which admit a holomorphic symplectic form. I prove Orlov's conjecture for many such varieties. 

This result builds upon a result of Taelman \cite{TaelmanDerivedEquivalences}. Using the hard Lefschetz theorem, one can associate to a variety a Lie algebra (the LLV algebra, named after Looijenga, Lunts and Verbitsky \cite{LooijengaLunts, VerbitskyPhDThesis}), and there is a representation of this Lie algebra on the cohomology of the variety. Taelman has proved that for varieties admitting a holomorphic symplectic form this Lie algebra is an invariant of the derived category.

An important part of this thesis is the computation of the LLV algebra of a general variety admitting a holomorphic symplectic form. Using this computation, I prove Orlov's conjecture for many varieties admitting a holomorphic symplectic form. Furthermore, I prove this conjecture for all $4$-dimensional varieties admitting a holomorphic symplectic form.

%% file: chapters/samenvatting.tex
\section*{Variëteiten met een holomorfe symplectische vorm: LLV algebras en afgeleide equivalenties}

Wanneer we meetkunde doen in twee of drie dimensies, kunnen we vaak gebruik maken van onze intuïtie. We kunnen ons iets voorstellen bij de objecten die we bestuderen, en kunnen er tekeningen of modellen van maken. Zodra we meetkunde willen doen in hogere dimensies, verdwijnt de mogelijkheid om deze objecten eenvoudig voor te stellen en wordt het lastiger om hun eigenschappen te bestuderen. 

De moderne wiskunde heeft hier echter oplossingen voor gevonden. Er zijn veel manieren gevonden om grip te krijgen op hoger dimensionale meetkundige objecten. Dit gaat meestal met behulp van simpelere objecten die informatie bevatten over de meetkundige objecten die we bestuderen, en dit noemen we \textit{invarianten}. Een bekend voorbeeld is de cohomologie, die grof gezien het aantal gaten van een meetkundig object telt. 

Dit proefschrift speelt zich af in de algebraïsche meetkunde, waar de bestudeerde meetkundige objecten zogeheten \textit{algebraïsche variëteiten} zijn. Algebraïsche variëteiten zijn gegeven als nulpuntsverzamelingen van polynomiale vergelijkingen, en bezitten een rijke structuur. Dit komt onder andere tot uiting in de vele interessante invarianten die men aan algebraïsche variëteiten kan toekennen. Deze invarianten dragen informatie over de specifieke algebraïsche structuur van de variëteit. De cohomologie van een gladde variëteit over $\mathbb{C}$ bijvoorbeeld draagt nog een extra structuur, de zogenaamde \textit{Hodge structuur}, die veel informatie over de variëteit bevat. 

Een andere belangrijke invariant van een algebraïsche variëteit is de \textit{afgeleide categorie}. Dit is een categorie die in een aantal stappen gebouwd wordt uit de categorie van coherente schoven op de variëteit.

Het is mogelijk dat twee verschillende variëteiten dezelfde afgeleide categorie hebben. We kunnen ons dan vragen in hoeverre deze variëteiten op elkaar lijken. Twee variëteiten met dezelfde afgeleide categorie hebben bijvoorbeeld dezelfde dimensie. Het is nog onbekend in welke mate de afgeleide categorie de cohomologie en Hodge structuur precies bepaalt. Een belangrijk vermoeden in deze richting is het vermoeden van Orlov:
\begin{vermoeden*}[Orlov]
    Zij $X$ en $Y$ twee gladde projectieve variëteiten over $\mathbb{C}$, en neem aan dat er een equivalentie $D^b(X) \simeq D^b(Y)$ van getrianguleerde categorieën bestaat. Dan is er een isomorfisme $\sH^\bullet(X; \mathbb{Q}) \cong \sH^\bullet(Y; \mathbb{Q})$ dat de gradering en Hodge structuur bewaart.
\end{vermoeden*}

In dit proefschrift kijk ik naar een specifiek soort variëteiten, namelijk gladde projectieve variëteiten over $\mathbb{C}$ met een holomorfe symplectische vorm. Voor een groot deel van deze variëteiten bewijs ik dat Orlovs vermoeden inderdaad waar is. 

Dit resultaat bouwt voort op een resultaat van Taelman \cite{TaelmanDerivedEquivalences}. Met behulp van de harde Lefschetz-stelling kan men aan een variëteit een Lie algebra toekennen (de LLV algebra, vernoemd naar Looijenga, Lunts en Verbitsky \cite{LooijengaLunts, VerbitskyPhDThesis}), en deze Lie algebra heeft een representatie op de cohomologie. Taelman heeft bewezen dat voor variëteiten met een holomorfe symplectische vorm deze Lie algebra een invariant van de afgeleide categorie is. 

Een belangrijk onderdeel van dit proefschrift is de berekening van de LLV algebra van een algemene variëteit met een holomorfe symplectische vorm. Deze berekening gebruik ik vervolgens om Orlovs vermoeden te bewijzen voor een groot aantal variëteiten met een holomorfe symplectische vorm. Daarnaast laat ik zien dat Orlovs vermoeden waar is voor alle $4$-dimensionale variëteiten met een holomorfe symplectische vorm.